\documentclass[11pt]{article}
\usepackage{graphicx,psfrag}
\usepackage{amssymb,amsfonts,amsthm}
\usepackage{amsmath,amsthm,amssymb}
\usepackage{amsfonts}

\newcommand{\la}{\langle}
\newcommand{\ra}{\rangle}
\newcommand{\rrr}{{\mathfrak r}}
\newcommand{\vk}{{van Kampen}\ }

\newtheorem{theorem}{Theorem}[section]
\newtheorem{lemma}[theorem]{Lemma}

\newtheorem{cor}[theorem]{Corollary}
\newtheorem{proposition}[theorem]{Proposition}
\theoremstyle{definition}
\newtheorem{definition}[theorem]{Definition}
\newtheorem{definitions}[theorem]{Definitions}
\newtheorem{remarks}[theorem]{Remarks}
\newtheorem{remark}[theorem]{Remark}
\newtheorem{prob}[theorem]{Problem}
\newcommand{\dg}{\mathfrak{d}}
\newcommand{\cgs}{\mathrm{Cayley} (G,S)}
\newcommand{\cgsh}{\mathrm{Cayley} (G,\, S\cup \mathcal{H})}
\newcommand{\dsh}{{\mathrm{dist}}_{S\cup \mathcal{H}}}
\newcommand{\ds}{{\mathrm{dist}}_S}
\newcommand{\luv}{\leq_{uv}}
\newcommand{\La}[1]{\Lambda_{#1}}
\newcommand{\cg}{{\mathcal{G}}}
\newcommand{\su}[1]{\Sat^\mu(#1)}
\newcommand{\dxh }{\dist_{S\cup H}}
\newcommand{\dx }{\dist_S}
\newcommand{\hh}{{\mathcal H}}
\newcommand{\G }{\cgsh}

\newcommand{\calg}{\Gamma}
\newcommand{\St}{\mathrm{St}}
\newcommand{\dn}{\Delta^{(n)}}
\newcommand{\Con}{{\mathrm{Con}}}
\newcommand{\lm}{{\lim}}
\newcommand{\dist}{{\mathrm{dist}}}
\newcommand{\oo}{{\mathcal O}}
\newcommand{\hdist}{\mathrm{hdist}}
\newcommand{\Cay}{\mathrm{Cayley}}

\newcommand{\pp}{{\mathcal P}}

\newcommand{\vv}{{\mathcal V}}
\newcommand{\Sat}{{\mathrm{Sat}}}
\newcommand{\nn}{{\mathcal N}}
\newcommand{\aaa}{{\mathcal A}}
\newcommand{\diam}{{\mathrm{diam}}}

\def\mh{{\mathcal{H}}} 
\def\calp{\mathcal{P}}   
\def\calv{\mathcal{V}}   
\def\calb{\mathcal{B}}   
\def\calc{\mathcal{C}}   

\newcommand {\N}{\mathbb{N}} 
\newcommand {\Z}{\mathbb{Z}}            
\newcommand {\R}{\mathbb{R}} 
\newcommand {\free}{\mathbb{F}} 
\newcommand {\q}{\mathfrak q} 
\newcommand {\g}{\mathfrak g} 
\newcommand {\pgot}{\mathfrak p}
\newcommand {\cf}{\mathfrak c}
\newcommand {\me}{\medskip}
\newcommand {\sm}{\smallskip}
\newcommand {\iv}{^{-1}}
\newcommand{\lio}[1]{\lm^\omega(#1)}
\newcommand{\co}[1]{\Con^\omega(#1)}
\newcommand{\proj}{{\hbox{proj}}}
\newcommand{\ttt}{{\cal T}}

\newcommand {\fn}{\footnote}
\newcommand {\Notat}{\noindent {\it{Notation}}:} 
\begin{document}

\title{Tree-graded spaces and asymptotic cones of groups}
\author{Cornelia Dru\c{t}u$^{a}$ and Mark Sapir$^{b,}$\thanks{The first author
is grateful to the CNRS of France for granting her the
d\'el\'egation CNRS status in the 2003-2004 academic year. The
research of the second author was supported in part by the NSF
grant DMS 0072307 and by the US-Israeli BSF grant 1999298. We
worked on the present paper during the first author's visits to
Vanderbilt University, and the second author's visits to
University of Lille-1 and to the IHES. The authors are grateful to
the administration and faculty of these institutions for their
support and hospitality.}}
\date{}
\maketitle

\begin{center}
\large{\textit{with an Appendix by Denis Osin$^b$ and Mark Sapir}}
\end{center}

\begin{center}

$^a$UFR de Mathématiques, Université de Lille-1 59655 Villeneuve
d'Ascq, France.

$^b$Department of Mathematics, Vanderbilt University, Nashville, TN
37240, U.S.A.
\end{center}

\me

\begin{abstract} We introduce a concept of tree-graded metric
space and we use it to show quasi-isometry invariance of certain
classes of relatively hyperbolic groups, to obtain a
characterization of relatively hyperbolic groups in terms of their
asymptotic cones, to find geometric properties of Cayley graphs of
relatively hyperbolic groups, and to construct the first example of
finitely generated group with a continuum of non-$\pi_1$-equivalent
asymptotic cones. Note that by a result of Kramer, Shelah, Tent and
Thomas, continuum is the maximal possible number of different
asymptotic cones of a finitely generated group, provided that the
Continuum Hypothesis is true.
\end{abstract}

\tableofcontents

\section{Introduction}

An asymptotic cone of a metric space is, roughly speaking, what one
sees when one looks at the space from infinitely far away. More
precisely, any asymptotic cone of a metric space $(X, \dist)$
corresponds to an ultrafilter $\omega$, a sequence of observation
points $e=(e_n)_{n\in \N}$ from $X$ and a sequence of scaling
constants $d=(d_n)_{n\in \N}$ diverging to $\infty$. The cone
$\Con(X;e,d)$ corresponding to $e$ and $d$ is the $\omega$-limit of
the sequence of spaces with basepoints $(X,\dist/d_n , e_n)$ (see
Section \ref{ULAC} for precise definitions).

In particular, if $X$ is the Cayley graph of a group $G$ with a
word metric then the asymptotic cones of $X$ are called asymptotic
cones of $G$.

The concept of asymptotic cone was essentially used by Gromov in
\cite{Gr1} and then formally introduced by van den Dries and
Wilkie \cite{VDW}.

Asymptotic cones have been used to characterize important classes
of groups:

\begin{itemize}
\item A finitely generated group is virtually Abelian if and only if its
asymptotic cones are isometric to the Euclidean space $\R^n$
(\cite{Gr1}, \cite{Pa}).

\item A finitely generated group is virtually nilpotent if and only if its asymptotic cones are locally compact (\cite{Gr1}, \cite{VDW},
\cite{Dr3}).

\item A finitely generated group is hyperbolic if and only if its asymptotic
cones are $\R$-trees (\cite{Gr2}).
\end{itemize}

In \cite{EP} it is shown moreover that asymptotic cones of
non-elementary hyperbolic groups are all isometric to the complete
homogeneous $\R$-tree of valence continuum. The asymptotic cones of
elementary groups are isometric to either a line $\R$ (if the group
is infinite) or to a point. Thus every hyperbolic group has only one
asymptotic cone up to isometry.

Asymptotic cones of quasi-isometric spaces are bi-Lipschitz
equivalent. In particular the topology of an asymptotic cone of a finitely
generated group does not depend on the choice of the generating
set. This was used in \cite{KaL1} and \cite{KaL2} to prove rigidity results for fundamental groups of
Haken manifolds, in \cite{KlL} to prove rigidity for cocompact lattices in
higher rank semisimple groups, and in
\cite{Dr2} to provide an alternative proof of the rigidity  for non-cocompact
lattices in higher rank semisimple groups. For a survey of results on
quasi-isometry invariants and their relations to asymptotic cones
see \cite{Dr3}.

The power of asymptotic cones stems from the fact that they capture
both geometric and logical properties of the group, since a large
subgroup of the ultrapower $G^\omega$ of the group $G$ acts
transitively by isometries on the asymptotic cone
$\Con^\omega(G;e,d)$. Logical aspects of asymptotic cones are
studied and used in the recent papers by Kramer, Shelah, Tent and
Thomas \cite{KSTT}, \cite{KT}.

One of the main properties of asymptotic cones of a metric space
$X$ is that geometry of finite configurations of points in the
asymptotic cone reflects the ``coarse" geometry of similar finite
configurations in $X$. This is the spirit of Gromov-Delzant's
approximation statement \cite{Del} and of the applications  of
$\R$-trees to Rips-Sela theory of equations in hyperbolic groups
and homomorphisms of hyperbolic groups \cite{RipsSela}. This was
also used in Dru\c tu's proof of hyperbolicity of groups with
sub-quadratic isoperimetric inequality \cite{Dr1}.

By a result of Gromov \cite{Gr2} if all asymptotic cones of a
finitely presented group are simply connected then the group has
polynomial isoperimetric function and linear isodiametric
function. Papasoglu proved in \cite{Pp} that groups having
quadratic isoperimetric functions have simply connected asymptotic
cones. In general, asymptotic cones of groups are not necessarily
simply connected \cite{Trofimov}. In fact, if a group $G$ is not
finitely presented then its asymptotic cones cannot all be simply
connected \cite{Gr2,Dr3}. A higher-dimensional version of this
result is obtained by Riley \cite{Ri}. According to the result of
Gromov cited above, examples of finitely presented groups with
non-simply connected asymptotic cones can be found in \cite{Bri}
and \cite{SBR}.

Although asymptotic cones can be completely described in some
cases, the general perception is nevertheless that asymptotic
cones are usually large and ``undescribable''. This might be the
reason of uncharacteristically ``mild'' questions by Gromov
\cite{Gr2}:

\begin{prob} \label{pr1} Which groups can appear as subgroups in
fundamental groups of asymptotic cones of
finitely generated groups?
\end{prob}

\begin{prob} \label{pr2} Is it true that the fundamental group of an
asymptotic cone of a group is either trivial or uncountable?
\end{prob}

In \cite{Gr2}, Gromov also asked the following question.

\begin{prob}\label{pr3} How many non-isometric asymptotic cones can a
finitely generated group have?
\end{prob}

A solution of Problem \ref{pr1} was given by Erschler and Osin
\cite{EO}. They proved that every metric space satisfying some
weak properties can be $\pi_1$- and isometrically embedded into
the asymptotic cone of a finitely generated group. This implies
that every countable group is a subgroup of the fundamental group
of an asymptotic cone of a finitely generated group.

Notice that since asymptotic cones tend to have fundamental groups
of order continuum, this result does not give information about the
structure of the whole fundamental group of an asymptotic cone, or
about how large the class of different asymptotic cones is: there
exists a group of cardinality continuum (for example, the group of
all permutations of a countable set) that contains all countable
groups as subgroups. One of the goals of this paper is to get more
precise information about fundamental groups of asymptotic cones,
and about the whole set of different asymptotic cones of a finitely
generated group.

Problem \ref{pr3} turned out to be related to the Continuum
Hypothesis (i.e. the famous question of whether there exists a set
of cardinality strictly between $\aleph_0$ and $2^{\aleph_0}$).
Namely, in \cite{KSTT}, it is proved that if the Continuum
Hypothesis is not true then any uniform lattice in
$\mathrm{SL}_n(\R)$ has $2^{2^{\aleph_0}}$ non-isometric
asymptotic cones, but if the Continuum Hypothesis is true then any
uniform lattice in $\mathrm{SL}_n(\R)$ has exactly one asymptotic
cone up to isometry, moreover the maximal theoretically possible
number of non-isometric asymptotic cones of a finitely generated
group is continuum. Recall that the Continuum Hypothesis is
independent of the usual axioms of set theory (ZFC).

It is known, however, that even if the Continuum Hypothesis is
true, there exist groups with more than one non-homeomorphic
asymptotic cones \cite{TV}. Nevertheless, it was not known whether
there exists a group with the maximal theoretically possible
number of non-isometric asymptotic cones (continuum).

In \cite{Gr1p}, Gromov introduced a useful generalization of
hyperbolic groups, namely the relatively hyperbolic groups\footnote{These
groups are also called {\em strongly relatively hyperbolic} in
order to distinguish them from weakly relatively hyperbolic groups
in the sense of Farb.}. This class includes:

\begin{enumerate}
\item[(1)] geometrically finite Kleinian groups; these groups
are hyperbolic relative to their cusp subgroups;

\item[(2)] fundamental groups of hyperbolic manifolds of finite volume
 (that is, non-uniform lattices in rank one semisimple
groups with trivial center);  these are hyperbolic relative to
their cusp subgroups;

\item[(3)] hyperbolic groups; these are hyperbolic relative to the
trivial subgroup or more generally to collections of quasi-convex
subgroups satisfying some extra
  conditions;

\item[(4)] free products of groups; these are hyperbolic
  relative to their factors;

  \item[(5)] fundamental groups of non-geometric Haken manifolds with at
least one hyperbolic component; these are hyperbolic relative to
the fundamental groups of the maximal graph-manifold components
and to the fundamental groups of the tori and Klein bottles not
contained in graph-manifold components \cite{Bow4*};

\item[(6)] $\omega$-residually free groups (limit groups in
another terminology); these are hyperbolic relative to the
collection of maximal Abelian non-cyclic subgroups \cite{Dah}.
\end{enumerate}

There exist several characterizations of relatively hyperbolic
groups which are in a sense parallel to the well known
characterizations of hyperbolic groups (see \cite{Bow1*},
\cite{Fa}, \cite{Osin}, \cite{Dah2}, \cite{Ya} and references
therein). But there was no characterization in terms of asymptotic
cones. Also, it was not known whether being relatively hyperbolic
with respect to any kind of subgroups is a quasi-isometry
invariant, except for hyperbolic groups when quasi-isometry
invariance is true.

The following theorems are the main results of the paper (we
formulate these results not in the most general form).

The first theorem gives more information about the possible structure of
fundamental groups of asymptotic cones.

\begin{theorem}[Theorem \ref{thcount} and
Corollary \ref{cth1}] \label{th1}
\begin{itemize}
\item[(1)]
For every countable group $C$, the free product of continuously
many copies of $C$ is the fundamental group of an asymptotic cone
of a 2-generated group.

\item[(2)] There exists a $2$-generated group $\Gamma$ such
that for every finitely presented group $G$, the free product of
continuously many copies of $G$ is the fundamental group of an
asymptotic cone of $\Gamma$.
\end{itemize}
\end{theorem}

The second theorem  answers the question about the number of
asymptotic cones of a finitely generated group.

\begin{theorem}[Theorem \ref{theorem676}] \label{th2}
Regardless of whether the Continuum Hypothesis is true or not,
there exists a finitely generated group $G$ with continuously many
pairwise non-$\pi_1$-equivalent asymptotic cones.
\end{theorem}

The third theorem shows that large classes of relatively hyperbolic
groups are closed under quasi-isometry. We call a finitely generated
group $H$ {\em unconstricted} if one of its asymptotic cones has no
global cut-points.

\begin{theorem}[Corollary \ref{relhipqig}] \label{th3} Let
$G$ be a finitely generated group that is hyperbolic relative to
unconstricted subgroups $H_1, ...,H_m$.

Let $G'$ be a group that is quasi-isometric to $G$. Then $G'$ is
hyperbolic relative to subgroups $H_1',...,H_n'$ each of which is
quasi-isometric to one of $H_1,...,H_m$.
\end{theorem}

The number $m$ of the finite collection of ``parabolic'' subgroups
$\{ H_i \}_{i\in I}$ in Theorem \ref{th3} is not a quasi-isometry
invariant. This can be seen for instance for the fundamental groups
of a finite volume hyperbolic manifold and of a finite covering of
it.

There are previous results showing that some special classes of
relatively hyperbolic groups are closed under quasi-isometry: the
class of fundamental groups of non-geometric Haken manifolds with at
least one hyperbolic component (\cite{KaL1}, \cite{KaL2}) and the
class of non-uniform lattices of isometries of a rank one symmetric
space \cite{Sch}. The class of free products of groups with finite
amalgamated subgroups is closed under quasi-isometry by Stallings'
Ends Theorem (see \cite{PW} for more general results about graphs
graphs of groups with finite edge groups).

The main ingredient in the proof of Theorem \ref{th3} is the
following result, interesting by itself.

\begin{theorem}[Corollary \ref{cutp1}] \label{th3p}
Let $G$ be a finitely generated group that is hyperbolic relative to
subgroups $H_1, ...,H_m$, and let $G'$ be a unconstricted group.
Then the image of $G'$ under any $(L,C)$-quasi-isometry $G'\to G$ is
in an $M$-tubular neighborhood of a coset $gH_i$, $g\in G,
i=1,...,m$, where $M$ depends on $L, C, G$ and $S$ only.
\end{theorem}

Note that the hypothesis of Theorem \ref{th3p} that the group $G$ is
unconstricted clearly cannot be removed. For example, a relatively
hyperbolic group itself is not in a bounded neighborhood of a coset
of any of its ``parabolic" subgroups $H_i$ provided $H_i$ are proper
subgroups.

Theorem \ref{th3p} does not apply in this case because relatively
hyperbolic groups are usually {\em constricted} i.e. they have
global cut-points in every asymptotic cone (see Theorem \ref{th4}
below).

A result similar to Theorem \ref{th3p} is obtained in \cite[$\S
3$]{PW} for $G$ a fundamental group of a graph of groups with finite
edge groups and $S$ a one-ended group. We should note here that
unconstricted groups are 1-ended by Stallings' Ends Theorem. The
converse statement is not true because the asymptotic cones of any
hyperbolic group are $\R$-trees.

Theorem \ref{th3p} in particular gives information about which
unconstricted subgroups can appear as undistorted subgroups in a
relatively hyperbolic group (see Remark \ref{rlast}, (1)). The
following theorem clarifies even more the question of the structure
of undistorted subgroups in relatively hyperbolic groups.

\begin{theorem}[Theorem \ref{thundis}]\label{th3q} Let
$G=\la S\ra$ be a finitely generated group that is hyperbolic
relative to subgroups $H_1,...,H_n$. Let $G'$ be an undistorted
finitely generated subgroup of $G$. Then $G'$ is relatively
hyperbolic with respect to subgroups $H_1',...,H_m'$, where each
$H_i'$ is one of the intersections $G'\cap gH_jg\iv$, $g\in G$,
$j\in \{ 1,2,\dots ,n \}$.
\end{theorem}

We also obtain information about the automorphism group of a
relatively hyperbolic group.

\begin{theorem}[Corollary \ref{829}]\label{th3r}
Let $G$ be a finitely generated group that is relatively hyperbolic
with respect to a unconstricted subgroup $H$.  Let $\mathrm{Fix}(H)$
be the subgroup of the automorphism group of $G$ consisting of the
automorphisms that fix $H$ as a set. Then:
\begin{itemize}
\item[(1)] $\mathrm{Inn}(G)\mathrm{Fix}(H)=\mathrm{Aut}(G)$.

\item[(2)] $\mathrm{Inn}(G)\cap \mathrm{Fix}(H)
= \mathrm{Inn}_{H}(G)$, where $\mathrm{Inn}_{H}(G)$ is by
definition $\{i_h\in \mathrm{Inn}(G) \mid h\in
  H\}\, .$
\item[(3)] There exists a natural homomorphism from
$\mathrm{Out}(G)$ to $\mathrm{Out}(H)$ given by $\phi\mapsto
i_{g_\phi}\phi|_{H}$, where $g_\phi$ is an element of $G$ such
that $i_{g_\phi}\phi\in \mathrm{Fix}(H)$, and $\psi|_{H}$ denotes
the restriction of an automorphism $\psi \in \mathrm{Fix}(H)$ to
$H$.
\end{itemize}
\end{theorem}

We call a finitely generated group {\em wide} if none of its
asymptotic cones has a global cut-point. Wide groups are certainly
unconstricted (the converse statement is not true).

{\it{Here are examples of wide groups:}}
\begin{itemize}

\item Non-virtually cyclic groups satisfying a law (see Corollary \ref{lawcp}).
Recall that a {\it{law}} is a word $w$ in $n$ letters $x_1,\dots ,
x_n$ and {\it{a group satisfying the law}} $w$ is a group $G$ such
that $w=1$ in $G$ whenever $x_1,\dots , x_n$ are replaced by an
arbitrary set of $n$ elements in $G$. For instance Abelian groups
are groups with the law $w=x_1x_2x_1^{-1}x_2^{-1}$. More
generally, solvable groups are groups with a law, and so are
Burnside groups. Also, uniformly amenable groups are groups
satisfying a law (see Corollary \ref{amen}).

While for nilpotent groups the results of Theorems \ref{th3} and
\ref{th3p} are not surprising and were already known in some
particular cases of relatively hyperbolic groups \cite{Sch}, for
solvable non-nilpotent groups and for Burnside groups the situation
is different. For instance the group Sol has asymptotic cones
composed of continuously many Hawaiian earrings \cite{Bu}, so it is
{\it{a priori}} not clear why such a group should have a rigid
behavior with respect to quasi-isometric embeddings into relatively
hyperbolic groups. Burnside groups display a similar picture.

In the case of non-virtually cyclic groups with a law, the
constant $M$ in Theorem \ref{th3p}
 depends only on the law and not on the group $S$ (Corollary \ref{cornice}).

\item Non-virtually cyclic groups with elements of infinite order in
the center (see Theorem \ref{infceter}); the constant $M$ in
Theorem \ref{th3p} is the same for the whole class of such groups
(Theorem \ref{uzg} and Corollary \ref{cuzg}).

\item Groups of isometries acting properly discontinuously and with
compact quotients on products of symmetric spaces and Euclidean
buildings, of rank at least two. The asymptotic cones of such
groups are Euclidean buildings of rank at least two \cite{KlL}.
Most likely the same is true for such groups of isometries so that
the quotients have finite volume, but the proof of this statement
is not straightforward.

\end{itemize}

The main tool in this paper are tree-graded spaces.

\begin{definition}
Let $\free$ be a complete geodesic metric space and let $\pp$ be a
collection of closed geodesic subsets (called {\it{pieces}}).
Suppose that the following two properties are satisfied:

\begin{enumerate}

\item[($T_1$)] Every two different pieces have at most one common
point.

\item[($T_2$)] Every simple geodesic triangle (a simple loop
composed of three geodesics) in $\free$ is contained in one piece.
\end{enumerate}

Then we say that the space $\free$ is {\em tree-graded with respect
to }$\pp$.
\end{definition}


The main interest in the notion of tree-graded space resides in
the following characterization of relatively hyperbolic groups of
which the converse part is proven in Section \ref{secequiv} and
the direct part in the Appendix written by D. Osin and M. Sapir.

\begin{theorem}[Theorem \ref{dir}] \label{th4} A finitely generated group $G$ is
relatively hyperbolic with respect to finitely generated subgroups
$H_1,...,H_n$ if and only if every asymptotic cone
$\Con^\omega(G;e,d)$ is tree-graded with respect to $\omega$-limits
of sequences of cosets of the subgroups $H_i$.
\end{theorem}

Section \ref{tgs} contains many general properties of tree-graded
spaces.

In particular, by Lemma \ref{cutting} any complete homogeneous
geodesic metric space with global cut-points is tree-graded with
respect to a certain uniquely defined collection of pieces which are either singletons or without
cut-points.

We prove in Proposition \ref{equiv}, that the property ($T_2$) in
the definition of tree-graded spaces can be replaced by the
assumption that $\pp$ covers $\free$ and the following property
which can be viewed as a extreme version of the bounded coset
penetration property:

\begin{quotation}
($T_2'$)\quad For every topological arc $\cf:[0,d]\to \free$ and
$t\in [0,d]$, let $\cf[t-a,t+b]$ be a maximal sub-arc of $\cf$
containing $\cf (t)$ and contained in one piece. Then every other
topological arc with the same endpoints as $\cf$ must contain the
points $\cf (t-a)$ and $\cf (t+b)$.
\end{quotation}

Moreover, when ($T_2$) is replaced by ($T_2'$) the condition that
the pieces  are geodesic is no longer needed. Thus, if we do not
ask that the whole space be geodesic either, tree-graded spaces
can be considered in a purely topological setting.

Notice that there are similarities in the study of asymptotic
cones of groups and that of boundaries of groups. Boundaries of
groups do not necessarily have a natural metric, and rarely are
geodesic spaces, but they have a natural topology and they are
also, in many interesting cases, homogeneous spaces with respect
to actions by homeomorphisms. Thus, if the boundary of a group is
homogeneous and has a global cut-point then most likely it is
tree-graded (in the topological sense) with respect to pieces that
do not have cut-points. Such a study of boundaries of groups with
global cut-points appeared, for example, in the work of Bowditch
\cite{Bow2*} on the Bestvina-Mess conjecture. Bowditch developed a
general theory appropriate for the study of topological
homogeneous spaces with global cut-points that is related to the
study of tree-graded spaces that we do in this paper. Results
related to Bowditch's work in this general setting can be found in
\cite{AN}.

As a byproduct of the arguments in Sections \ref{amc} and
\ref{secequiv}, we obtain many facts about the geometry of Cayley
graphs of relatively hyperbolic groups. Recall that given a
finitely generated group $G=\la S\ra$ and a finite collection
$H_1,...,H_n$ of subgroups of it, one can consider the standard
Cayley graph $\Cay(G,S)$ and the modified Cayley graph $\Cay(G,
S\cup\hh)$, where $\hh=\bigsqcup_{i=1}^n (H_i\setminus \{e\})$. The
standard definition of relative hyperbolicity of a group $G$ with
respect to subgroups $H_1,...,H_n$ is given in terms of the
modified Cayley graph $\Cay(G, S\cup\hh)$. Theorem \ref{th4} and
the results of Section \ref{amc} allow us to define the relative
hyperbolicity of $G$ with respect to $H_1,...,H_n$ in terms of
$\Cay(G,S)$ only. This is an important ingredient in our rigidity
results.

An important part in studying tree-graded spaces is played by {\em
saturations} of geodesics. If $G$ is relatively
hyperbolic with respect to $H_1,...,H_n$, $\g$ is a geodesic in
$\Cay(G,S)$ and $M$ is a positive number, then the $M$-\textit{saturation
of} $\g$ is the union of $\g$ and all left cosets of $H_i$ whose $M$-tubular
neighborhoods intersect $\g$. We show that in the study
of relatively hyperbolic groups, saturations play the same role as
the geodesics in the study of hyperbolic groups.

More precisely, we use Bowditch's characterization of hyperbolic
graphs \cite{Bow3*}, and show that tubular neighborhoods of
saturations of geodesics can play the role of ``lines" in that
characterization. In particular, we show that for every geodesic
triangle $[A,B,C]$ in $\Cay(G,S)$ the $M$-tubular neighborhoods of
the saturations of its sides (for some $M$ depending on $G$ and
$S$) have a common point which is at a bounded distance from the
sides of the triangle or a common left coset which is at a bounded
distance from the sides.


We also obtain the following analog for relatively hyperbolic groups
of the Morse Lemma for hyperbolic spaces. Recall that the Morse
lemma states that every quasi-geodesic in a hyperbolic space is at a
bounded distance from a geodesic joining its endpoints. In the
relative hyperbolic version of the lemma we also use the notion of
{\em lift} $\tilde{\pgot}$ of a geodesic $\pgot$ in $\cgsh$. Recall
that the meaning of it is that one replaces each edge in $\pgot$
labelled by an element in $\hh$ by a geodesic in $\cgs$ (see also
Definition \ref{lift}).

We again do not write the statements in the whole generality.

\me

\noindent{\it{Notations}}: Throughout the whole paper,
$\nn_\delta(A)$ denotes the $\delta$-tubular neighborhood of a
subset $A$ in a metric space $X$, that is $\{x\in X\mid
\dist(x,A)<\delta\}$. We denote by $\overline{\nn}_\delta(A)$ its
closure, that is $\{x\mid \dist(x,A)\leq \delta\}$. In the
particular case when $A=\{ x \}$ we also use the notations
$B(x,\delta)$ and $\overline{B(x,\delta)}$ for the tubular
neighborhood and its closure.

\me

\begin{theorem}[Morse property for relatively hyperbolic
groups]\label{morse} Let $G=\la S\ra$ be a group that is hyperbolic
relative to the collection of subgroups $H_1,\dots ,H_m$. Then there
exists a constant $M$ depending only on the generating set $S$ such
that the following holds. Let $\g$ be a geodesic in $\cgs $, let
$\q$ be an $(L,C)$-quasi-geodesic in $\cgs$, and let $\pgot$ be an
$(L,C)$-quasi-geodesic in $\Cay(G,S\cup\hh)$. Suppose that $\g, \q$,
and $\pgot$  have the same endpoints. Then for some $\tau$ depending
only on $L, C, S$:
\begin{itemize}
\item[(1)] $\q$ is
contained in the $\tau$-tubular neighborhood of the $M$-saturation
of $\g$.

\item[(2)] Let
$gH_i$ and $g'H_j$ be two left cosets contained in the
$M$-saturation of $\g$. Let $\q'$ be a sub-quasi-geodesic of $\q$
with endpoints $a\in\nn_\kappa(gH_i)$ and $b\in \nn_\kappa(g'H_j)$
which intersects $\nn_\kappa(gH_i)$ and $\nn_\kappa(g'H_j)$ in sets
of bounded (in terms of $\kappa$) diameter. Then $a$ and $b$ belong
to the $\delta$-tubular neighborhood of $\g$, where $\delta$ depends
only on $L,C, \kappa$.

\item[(3)] In the Cayley graph
$\Cay(G,  S\cup\hh)$, $\q$ is at Hausdorff distance at most $\tau$
from $\pgot$.

\item[(4)] In $\cgs$, $\q$ is
contained in the $\tau$-tubular neighborhood of the
$\tau$-saturation of any lift $\tilde{\pgot}$ of $\pgot$. In its
turn, $\tilde{\pgot}$ is contained in the $\tau$-tubular
neighborhood of the $\tau$-saturation of $\q$.
\end{itemize}
\end{theorem}

The proof of this theorem and more facts about the geometry of
relatively hyperbolic groups are contained in Lemmas \ref{lqqcun},
\ref{limrsat}, \ref{gat}, Proposition \ref{udist} and Proposition
\ref{prudist}.

Theorem \ref{th4} and statements about tree-graded spaces from
Section \ref{tgs} imply that for relatively hyperbolic groups,
Problem \ref{pr2} has a positive answer.

\begin{cor} The fundamental group of an asymptotic cone
of a relatively hyperbolic group $G$ is either trivial or of order
continuum.
\end{cor}

\proof Suppose that the fundamental group of an asymptotic cone of
the group $G$ is non-trivial.  By Theorem \ref{th4}, the
asymptotic cone of $G$ is tree-graded with respect to a set of
pieces that are isometric copies of asymptotic cones of the
parabolic subgroups $H_i$ with the induced metric. The induced
metric on each $H_i$ is equivalent to the natural word metric by
quasi-convexity (see Lemma \ref{qqc}). Moreover, in that set,
every piece appears together with continuously many copies.

The argument in the first part of the proof of Proposition
\ref{pi1} shows that at least one of the pieces has non-trivial
fundamental group $\Gamma$.

The argument in the second part of the proof of Proposition
\ref{pi1} implies that the fundamental group of the asymptotic
cone of $G$ contains the free product of continuously many copies
of $\Gamma$.\endproof

The following statement is another straightforward consequence of
Theorem \ref{th4}.

\begin{cor}\label{cor666} If a group $G$ is hyperbolic relative to
$\{H_1,\dots ,H_m\}$, and each $H_i$ is hyperbolic relative to a
collection of subgroups $\{ H_i^{1}, \dots , H_i^{n_i} \}$ then $G$
is hyperbolic relative to $\{ H_i^{j} \mid i\in \{1,\dots ,m \},\,
j\in \{1,\dots ,n_i \}\}$. \end{cor}

See Problem \ref{prob2.5} below for a discussion of Corollary
\ref{cor666}.

Note that in the alternative geometric definition of relatively
hyperbolic groups given in Theorem \ref{th4} we do not need the
hypothesis that $H_i$ are finitely generated. This follows from the
quasi-convexity of the groups $H_i$ seen as sets in $\cgs$ (Lemma
\ref{qqc}). Moreover, this geometric definition makes sense when $G$
is replaced by a geodesic metric space $X$ and the collection of
cosets of the subgroups $H_i$ is replaced by a collection $\aaa$ of
subsets of $X$. A similar generalization can be considered for
Farb's definition of relative hyperbolicity (including the BCP
condition). Thus, both definitions allow to speak of geodesic spaces
hyperbolic relative to families of subsets. Such spaces, completely
unrelated to groups, do appear naturally. For instance the
complements of unions of disjoint open horoballs in rank one
symmetric spaces are hyperbolic with respect to the boundary
horospheres. Also, the free product of two metric spaces with
basepoints $(X,x_0)$ and $(Y,y_0)$, as defined in \cite[$\S 1$]{PW},
is hyperbolic with respect to all the isometric copies of $X$ and
$Y$. It might be interesting for instance to study actions of groups
on such spaces, hyperbolic with respect to collections of subsets.
To some extent, this is already done in the proof of our Theorem
\ref{relhipqi}, where a particular case of action of a group by
quasi-isometries on an asymptotically tree-graded (=relatively
hyperbolic) space is studied.

Bowditch's characterization of hyperbolic graphs can be easily
generalized to arbitrary geodesic metric spaces. So one can expect
that an analog of Theorem \ref{th4} is true  for arbitrary
geodesic metric spaces.

\subsection{Open problems}\label{last}

\begin{prob} \label{prob0} Is it possible to drop the condition
that $H_i$ are unconstricted from the formulation of Theorem
\ref{th3}?
\end{prob}

An obvious candidate to a counterexample would be, for instance, the
pair of groups $G=A*A*A*A$, where $A=\Z^2$, and $G'=(A*A*A*A)\rtimes
\Z/4\Z$, where $\Z/4\Z$ permutes the factors. The group $G$ is
relatively hyperbolic with respect to $A*A*1*1$ and $1*1*A*A$. It is
easy to check that the group $G'$ is not relatively hyperbolic with
respect to any isomorphic copy of $A*A$. Unfortunately this example
does not work. Indeed, $G'$ is quasi-isometric to $A*A$ by
\cite{PW}, so $G'$ is relatively hyperbolic with respect to a
subgroup that is quasi-isometric to $A*A$, namely itself. Moreover,
it is most likely that $G'$ is hyperbolic relative to a proper
subgroup isomorphic to $A*\Z$ which is also quasi-isometric to $A*A$
by \cite{PW}.

\begin{prob}\label{prob1} Corollary \ref{setgr} shows the following.
Let $G$ be a group, asymptotically tree-graded as a metric space
with respect to a family of subspaces $\aaa$ satisfying the
following conditions:
\begin{itemize}
  \item[(1)] $\mathcal{A}$ is uniformly unconstricted (see Definition \ref{uwandp} for the notion of collection of metric
spaces uniformly unconstricted);

\item[(2)] there exists a constant $c$ such that every point in every $A\in \mathcal{A}$
is at distance at most $c$ from a bi-infinite geodesic in $A$;

  \item[(3)] For a fixed $x_0\in G$ and every $R>0$ the ball
  $B(x_0,R)$ intersects only finitely many $A\in \aaa$.
\end{itemize}
Then the group $G$ is relatively hyperbolic with respect to
subgroups $H_1, ..., H_m$ such that every $H_i$ is quasi-isometric
to some $A\in \aaa$.

Can one remove some of the conditions (1), (2), (3) from this
statement?
\end{prob}

\begin{prob} \label{wideunconstricted} Is every unconstricted group
wide?
\end{prob}

\begin{prob}\label{prob2} Is every constricted group $G$ relatively hyperbolic with
respect to a collection of proper subgroups $\{H_1,\dots , H_m\}$?
Here are some more specific questions. Consider the canonical
representation of every asymptotic cone as a tree-graded space (with
respect to maximal path-connected subsets that are either singletons or without global cut-points, as in
Lemma \ref{cutting}). Is there a family of subsets $\aaa$ of $G$
such that each piece in each asymptotic cone of $G$ is an ultralimit
of a sequence of sets from $\aaa$? Can one take $\aaa$ to be the set
of all left cosets of a (finite) collection of subgroups $\{
H_1,\dots ,H_m\}$?
\end{prob}

Note that a positive answer to Problem \ref{prob2} gives a
positive answer to Problem \ref{prob0}, as being constricted is a
quasi-isometry invariant. Also, it would follow that the rigidity
result Theorem \ref{th3p} holds as soon as $G'$ is not relatively
hyperbolic.

Here is a related question.

\begin{prob} \label{prob2.25} Is every non-virtually cyclic group
without free non-abelian subgroups wide (unconstricted)? Is there a
non-virtually cyclic constricted group with all proper subgroups
cyclic?
\end{prob}

It is easy to notice that in all examples of groups with different
asymptotic cones $\Con^\omega(G;e,d)$, one of the cones
corresponds to a very fast growing sequence $d=(d_n)$.
Equivalently, we can assume that $d_n=n$ but $\omega$ contains
some fast growing sequence of natural numbers $A=\{a_1,a_2,...\}$.
What if we avoid such ultrafilters?  For example, let
$\mathcal{P}$ be the set of all complements of finite sets and of
all complements of sequences $A=\{a_1,a_2,\dots ,a_n,...\}$
($a_1<a_2<...<a_n <...$) which grow faster than linear i.e. $\lim
\frac{a_n}{n}=\infty$. It is easy to see that $\mathcal{P}$ is a
filter. Let $\omega$ be an ultrafilter containing $\mathcal{B}$.
Then no set in $\omega$ grows faster than linear. Let us call
ultrafilters with that property {\em slow}. An asymptotic cone
$\Con^\omega(G,(n))$ corresponding to a slow ultrafilter also will
be called {\em slow}.

\begin{prob} Are there finitely generated groups $G$ with two
bi-Lipschitz non-equivalent slow asymptotic cones? Is it true that
if a slow asymptotic cone of $G$  has (resp. has no) global
cut-points then the group is constricted (resp. wide)? Is it true
that if a slow asymptotic cone of $G$ has global cut-points then $G$
contains non-abelian free subgroups?
\end{prob}

See Section \ref{law} for further discussion of free subgroups of
wide (unconstricted) groups.

The next problem is motivated by Corollary \ref{cor666} above.

\begin{prob}\label{prob2.5} By Corollary \ref{cor666},
one can consider a ``descending process'', finding smaller and
smaller subgroups of a (finitely generated) group $G$ with respect
to which $G$ is relatively hyperbolic. Does this process always
stop? Does every group $G$ contain a finite collection of
unconstricted subgroups with respect to which $G$ is relatively
hyperbolic?
\end{prob}

\begin{prob}\label{prob3} A group $G=\la S\ra$ is weakly hyperbolic relative to
subgroups $H_1,...,H_n$ if the Cayley graph $\Cay(G,S\cup\hh)$ is
hyperbolic. It would be interesting to investigate the behavior of
weak relatively hyperbolic groups up to quasi-isometry. In
particular, it would be interesting to find out if an analog of
Theorem \ref{th3} holds. The arguments used in this paper for the
(strong) relative hyperbolicity no longer work. This can be seen
on the example of $\Z^n$. That group is weakly hyperbolic relative
to $\Z^{n-1}$. But a quasi-isometry $\q :\Z^n \to \Z^n$ can
transform left cosets of $\Z^{n-1}$ into polyhedral or even more
complicated surfaces (see \cite[Introduction]{KlL} for examples).
Nevertheless it is not a real counter-example to a theorem similar
to Theorem \ref{th3} for weak hyperbolic  groups, as every
group quasi-isometric to $\Z^n$ is virtually $\Z^n$.
\end{prob}

\subsection{Plan of the paper}

In Section \ref{tgs}, we establish some basic properties of
tree-graded spaces. In particular, we show that tree-graded
spaces behave ``nicely'' with respect to homeomorphisms.

In Section \ref{ULAC}, we establish general properties of asymptotic
cones and their ultralimits. We show that the ultralimit of a
sequence of asymptotic cones of a metric space $X$ is an asymptotic
cone of $X$ itself.

In Section \ref{amc}, we give an ``internal'' characterization of
{\em asymptotically tree-graded metric spaces}, i.e. pairs of a
metric space $X$ and a collection of subsets $\cal A$, such that
every asymptotic cone $\Con^\omega(X;e,d)$ is tree-graded with
respect to $\omega$-limits of sequences of sets from $\aaa$.

In Section \ref{qiv}, we show that being asymptotically tree-graded
with respect to a family of subsets is a quasi-isometry invariant.
This implies Theorem \ref{th3}.

In Section \ref{ascp}, we show that asymptotic cones of a
non-virtually cyclic group do not have cut-points provided the group
either has an infinite cyclic central subgroup, or satisfies a law.

In Section \ref{exAD}, we modify a construction from the paper
\cite{EO} to prove, in particular, Theorems \ref{th1} and
\ref{th2}.

In Section \ref{secequiv} and in the Appendix (written by D. Osin
and M. Sapir), we prove the characterization of relatively
hyperbolic groups in terms of their asymptotic cones given in
Theorem \ref{th4}. Theorem \ref{th3q} about undistorted subgroups
of relatively hyperbolic groups is also proved in Section
\ref{secequiv}.

\section{Tree-graded spaces}
\label{tgs}

\subsection{Properties of tree-graded spaces}

Let us recall the definition of tree-graded spaces. We say that a
subset $A$ of a geodesic metric space $X$ is a {\em geodesic subset}
if every two points in $A$ can be connected by a geodesic contained
in $A$.

\begin{definition}[tree-graded spaces]\label{tree}
Let $\free$ be a complete geodesic metric space and let $\pp$ be a
collection of closed geodesic subsets (called {\it{pieces}}).
Suppose that the following two properties are satisfied:

\begin{enumerate}

\item[($T_1$)] Every two different pieces have at most one common
point.

\item[($T_2$)] Every simple geodesic triangle (a simple loop
composed of three geodesics) in $\free$ is contained in one piece.
\end{enumerate}

Then we say that the space $\free$ is {\em tree-graded with respect
to }$\pp$.
\end{definition}

\begin{remark}[degenerate triangles] \label{singletons}
We assume that a point is a geodesic triangle composed of geodesics
of length 0. Thus ($T_2$) implies that the pieces cover $\free$.
\end{remark}

The next several lemmas establish some useful properties of
tree-graded spaces. Until Proposition \ref{equiv}, $\free$ is a
tree-graded space with respect to $\pp$.

\begin{lemma}\label{lin}
If all pieces in $\pp$ are $\R$-trees then $\free$ is an
$\R$-tree.
\end{lemma}

\proof It is an immediate consequence of $(T_2)$. \endproof

\begin{lemma}\label{feet} Let $M$ be a piece and $x$ a point outside $M$. If $y$ and
$z$ are points in $M$ such that there exist geodesics $[x,y]$ and
$[x,z]$, joining them to $x$ which intersect $M$ only in $y$ and
$z$, respectively, then $y=z$.
\end{lemma}

\proof Suppose that $y\neq z$. Join $y$ and $z$ by a geodesic
$[y,z]$ in $M$. Let $x'$ be the farthest from $x$ intersection point
of the geodesics $[x,y]$ and $[x,z]$. The triangle $x'yz$ is simple
because by the assumption $[x,y]\cup[x,z]$ intersects with $[y,z]$
only in $y$ and $z$. Therefore that triangle is contained in one
piece $M'$ by ($T_2$). Since $M\cap M'$ contains $[y,z]$, $M=M'$ by
($T_1$), so $x'\in M$, a contradiction since $x'$ belongs both to
$[x,y]$ and to $[x,z]$ but cannot coincide with both $y$ and $z$ at
the same time.
\endproof

\begin{lemma}\label{quadrangle}
Every simple quadrangle (i.e. a simple loop composed of four
geodesics) in $\free$ is contained in one piece.
\end{lemma}

\proof Let $A_1,A_2,A_3$ and $A_4$ be the vertices of the
quadrangle. Suppose that each vertex is not on a geodesic joining
its neighbors, otherwise we have a geodesic triangle and the
statement is trivial. Let $\g$ be a geodesic joining $A_1$ and
$A_3$. Let $P$ be its last intersection point with $[A_1,A_2] \cup
[A_1,A_4]$. Suppose that $P\in [A_1,A_2]$ (the other case is
symmetric). Let $Q$ be the first intersection point of $\g$ with
$[A_2,A_3] \cup [A_3,A_4]$. Replace the arc of $\g$ between $A_1$
and $P$ with the arc of $[A_1,A_2]$ between these two points, and
the arc of $g$ between $Q$ and $A_3$ with the corresponding arc of
$[A_2,A_3] \cup [A_3,A_4]$. Then $\g$ thus modified cuts the
quadrangle into two simple triangles having in common the geodesic
$[P,Q]$. Both triangles are in the same piece by ($T_2$), and so is
the quadrangle.
\endproof

\begin{lemma}\label{cv&proj}
\begin{itemize}
\item[(1)] Each piece is a convex subset of $\free$.

\item[(2)] For every point $x\in \free$ and every piece $M\in
\pp$, there exists a unique point $y\in M$ such that
$\dist(x,M)=\dist(x,y)$. Moreover, every geodesic joining $x$ with
a point of $M$ contains $y$.
\end{itemize}
\end{lemma}

\proof (1) Suppose that there exists a geodesic $\g$ joining two
points of $M$ and not contained in $M$. Let $z$ be a point in
$\g\setminus M$. Then $z$ is on a sub-arc $\g'$ of $\g$
intersecting $M$ only in its endpoints, $a,b$. Lemma \ref{feet}
implies $a=b=z\in M$, a contradiction.

\me

(2) Let $y_n\in M$ be such that $\lm_{n\to \infty}\dist(x,y_n) =
\dist(x,M)$. Since $M$ is closed, we may suppose that every geodesic
$[x,y_n]$ intersects $M$ only in $y_n$. It follows by Lemma
\ref{feet} that $y_1=y_2=\ldots = y$.

Let $z\in M$ and let $\g$ be a geodesic joining $z$ with $x$. Let $z'$
be the last point on $\g$ contained in $M$. Then $z'=y$, by Lemma
\ref{feet}.\endproof

\begin{definition}
We call the point $y$ in part (2) of Lemma \ref{cv&proj}
\textit{the projection of $x$ onto the piece $M$}.
\end{definition}

\begin{lemma}\label{smball} Let $M$ be a piece and $x$ a point outside it with $\dist(x,M)=\delta $, and let $y$  be the projection of $x$ onto
$M$. Then the projection of every  point $z\in \overline{B}(x, \delta )$ onto
$M$ is equal to $y$.
\end{lemma}

\proof Notice that by part (2) of Lemma \ref{cv&proj}
$\overline{B}(x,\delta)\cap M=\{y\}$. Suppose that the projection $z'$ of
$z\in B(x,\delta)$ onto $M$ is different from $y$. Then $z\ne y$, hence $z$ does
not belong to $M$.

Consider a geodesic quadrangle with vertices $x,z,z'$ and $y$. By the definition of projection, the interiors of
$[z,z']\cup[x,y]$ and $[y,z']$ do not intersect.

If there is a common point $p$ of $[x,y]$ and $[z,z']$ then we get a contradiction with Lemma \ref{feet}, so $[x,y]$ and
$[z,z']$ are disjoint. In particular $[z,z']\cup [z',y] \cup [y,x]$ is a topological arc. Since $z\in \overline{B}(x,\delta)\setminus \{ y \}$, the side $[x,z]$ of this quadrangle does not intersect $M$. By part
(1) of Lemma \ref{cv&proj} it follows that $[x,z]$ does not intersect $[y,z']$.

We can replace if necessary $z$ with the
last intersection point of $[z,x]$ with $[z,z']$ and $x$ with the last intersection point of the geodesics $[x,y]$ and
$[x,z]$. We get a simple geodesic
quadrangle $xzz'y$ in which the side $[x,z]$ possibly reduces to a point. By Lemma \ref{quadrangle}, it belongs to one
piece. Since it has $[y,z']$ in common with $M$, that piece is $M$
by ($T_1$). But this contradicts the fact that $[x,z]\cap M=\emptyset$.
\endproof

\begin{cor}\label{projpath} Every continuous path in $\free$
which intersects a piece $M$ in at most one point, projects onto
$M$ in a unique point.
\end{cor}

\proof If the path does not intersect the piece, it suffices to
cover it with balls of radius less than the distance from the path
to the piece and use Lemma \ref{smball}.

If the path intersects $M$ in a point $x$, we may suppose that $x$
is one of its ends and that the interior of the path does not pass
through $x$. Let $z$ be another point on the path and let $y$ be
its projection onto $M$. By the previous argument every point $t$
on the path, $t\neq x$, has the same projection $y$ onto $M$. Let
$\lm_{n\to \infty} t_n=x,\, t_n \neq x $. Then $\lm_{n\to
\infty}\dist(t_n, M)=\lm_{n\to \infty}\dist(t_n,y)=0$. Therefore
$x=y$.
\endproof

\begin{cor} \label{strconv}
\begin{itemize}
  \item[(1)] Every topological arc in $\free$ joining two points in a piece is contained in the piece.
  \item[(2)] Every
non-empty intersection between a topological arc in $\free$ and a
piece is a point or a sub-arc.
\end{itemize}
\end{cor}

\proof (1) If there exists a topological arc $\pgot $ in $\free$
joining two points of a piece $M$ and not contained in $M$, then a
point $z$ in $\pgot \setminus M$ is on a sub-arc $\pgot'$ of $\pgot$
intersecting $M$ only in its endpoints, $a,b$. Corollary
\ref{projpath} implies that both $a$ and $b$ are projections of
$z$ into $M$, contradiction.

(2) immediately follows from (1).\endproof

\begin{cor}\label{projA} Let $A$ be a connected
subset (possibly a   point) in $\free $ which intersects a piece
$M$ in at most one point.
\begin{itemize}
\item[(1)] The subset $A$ projects onto $M$ in a unique point $x$.
\item[(2)] Every path joining a point in $A$ with a point in $M$
contains $x$.
\end{itemize}
\end{cor}

\me

\noindent\textit{Notation:} Let $x\in \free$. We denote by $T_x$ the
set of points $y\in \free$ which can be joined to $x$ by a
topological arc intersecting every piece in at most one point.

\begin{lemma}\label{xy} Let $x\in \free $ and $y\in T_x\, ,\, y\neq x$.   Then every
topological arc with endpoints $x,y$ intersects each piece in at most one point. In particular the arc is contained in $T_x$.
\end{lemma}

\proof Suppose, by contradiction, that there exists a topological
arc $\pgot$ in $\free$ connecting $x,y$ and intersecting a piece $M$ in more than one point. By Corollary \ref{strconv}, $M\cap \pgot$ is a topological arc with
endpoints $a\ne b$. By definition, there also exists an arc $\q$ connecting $x$
and $y$ and touching each piece in at most one point.

Now consider the two paths connecting $x$ and
$M$. The first path $\pgot'$ is a part of $\pgot$ connecting $x$
and $a$. The second path $\q'$ is the composition of the path $\q$
and a portion of $\pgot\iv$ connecting $y$ and $b$. By Corollary \ref{projA}, the path $\q'$ must pass through the
point $a$. Since the portion $[y,b]$ of $\pgot\iv$ does not
contain $a$, the path $\q$ must contain $a$. But then there exists
a part $\q''$ of $\q'$ connecting $a$ and $b$ and intersecting $M$ in
exactly two points. This contradicts part (1) of Corollary
\ref{projA}, as a point in $\q''\setminus \{ a,b \}$ would project onto $M$ in both $a$ and $b$. \endproof

\begin{lemma}\label{txy}
Let $x\in \free $ and $y\in T_x$. Then $T_x=T_y$.
\end{lemma}

\proof It suffices to prove $T_y\subset T_x$. Let $z\in T_y$. By
Lemma \ref{xy}, any geodesics connecting $y$ with $x$ or $z$
intersects every piece in at most one point. Let $t$ be the
farthest from $y$ intersection point between two geodesics
$\pgot=[y,x]$ and $\q=[y,z]$. Then
$\rrr=[x,t]\cup [t,z]$ is a topological arc. The arc $\rrr$
intersects every piece in at most one point. Indeed, if $\rrr$
intersects a piece $M$ in two points $a, b$ then it intersects it in a subarc by Corollary \ref{strconv}, so at least one of the two segments $[x,t],[t,z]$ intersects $M$ in an arc, contradiction. Thus $z\in
T_x$.
\endproof

\begin{lemma}\label{txtree}
Let $x\in \free $.
\begin{itemize}
\item[(1)] Every topological arc joining two distinct points in
$T_x$ is contained in $T_x$.

\item[(2)] The subset $T_x$ is a real tree.
\end{itemize}
\end{lemma}

\proof (1) is an immediate consequence of the two previous lemmas.

(2) First we prove that for every $y,z\in T_x$ there exists a
unique geodesic joining $y$ and $z$, also contained in $T_x$.
Since $\free$ is a geodesic space, there exists a geodesic in
$\free$ joining $x$ and $y$. By the first part of the lemma, this
geodesic is contained in $T_x$. Suppose there are two distinct geodesics $\g,\g'$ in $T_x$ joining
$y$ and $z$. A point on $\g$ which is not on $\g'$ is contained in
a simple bigon composed of a sub-arc of $\g$ and a sub-arc of
$\g'$. This bigon, by $(T_2)$, is contained in a piece. This
contradicts Lemma \ref{xy}.

Now consider a geodesic triangle $yzt$ in $T_x$. Deleting, if
necessary, a common sub-arc we can suppose that $[y,z]\cap
[y,t]=\{ y\}$. If $y\not\in [z,t]$ then let $z'$ be the nearest to
$y$ point of $[y,z]\cap [z,t]$ and let $t'$ be the nearest to $y$
point of $[y,t]\cap [z,t]$. The triangle $yz't'$ is simple,
therefore it is contained in one piece by $(T_2)$. This again
contradicts Lemma \ref{xy}. Thus $y\in [z,t]$.
\endproof

\me

\noindent{\em Convention:} We assume that a 1-point metric space
has a cut-point.

\me
\begin{lemma}\label{cut}
Let $A$ be a path connected subset of $\free$ without a cut-point.
Then $A$ is contained in a piece. In particular every simple loop
is contained in a piece.
\end{lemma}

\proof By our convention, $A$ contains at least two points. Fix a
point $x\in A$. The set $A$ cannot be contained in the real tree
$T_x$, because otherwise it would have a cut-point. Therefore, a
topological arc joining in $A$ the point $x$ and some $y\in A$
intersects a piece $M$ in a sub-arc $\pgot$. Suppose that
$A\not\subset M$. Let $z\in A\setminus M$ and let $z'$ be the
projection of $z$ onto $M$. Corollary \ref{projA} implies that
every continuous path joining $z$ to any point $\alpha $ of $\pgot
$
contains $z'$. In particular $z'\in A$, and $z$ and $\alpha$ are
in two distinct connected components of $\free \setminus \{ z'
\}$. Thus, $z'$ is a cut-point of $A$, a contradiction.\endproof

\begin{proposition}\label{homeom}
Let $\free$ and $\free'$ be two tree-graded spaces with respect to
the sets of pieces $\pp$ and $\pp'$, respectively. Let $\Psi\colon
\free\to \free'$ be a homeomorphism.  Suppose that all pieces in
$\pp$ and $\pp'$  do not have cut-points. Then $\Psi$ sends any
piece from $\pp$ onto a piece from $\pp'$, and $\Psi (T_x)=T_{\Psi
(x)}$ for every $x\in \free$.
\end{proposition}

\proof Indeed, for every piece $M$ in $\free$, $\Psi(M)$ is a path
connected subset of $\free'$ without cut-points. Therefore
$\Psi(M)$ is inside a piece $M'$ of $\free'$ by Lemma \ref{cut}.
Applying the same argument to $\Psi\iv$, we have that $\Psi\iv(M')$ is contained in a piece $M''$. Then $M\subseteq \Psi\iv(M')\subseteq M''$, hence $M=M''$ and $\Psi(M)=M'$.
\endproof

\begin{proposition}\label{equiv} Condition ($T_2$) in the
definition of tree-graded spaces can be replaced the assumption that
pieces cover $\free$ plus any one of the following conditions:
\begin{enumerate}
\item[{\rm ($T_2'$)}]
 For every topological arc $\cf:[0,d]\to \free$ and $t\in
[0,d]$, let $\cf[t-a,t+b]$ be a maximal sub-arc of $\cf$ containing
$\cf (t)$ and contained in one piece. Then every other topological
arc with the same endpoints as $\cf$ must contain the points
$\cf (t-a)$ and $\cf (t+b)$.
\begin{figure}[!ht]
\centering
\unitlength .85mm 
\linethickness{0.4pt}
\ifx\plotpoint\undefined\newsavebox{\plotpoint}\fi 
\begin{picture}(128.04,45.33)(0,0)
\put(93.47,23.86){\line(0,1){.97}}
\put(93.44,24.83){\line(0,1){.968}}
\put(93.38,25.8){\line(0,1){.964}}
\multiput(93.27,26.77)(-.0306,.19171){5}{\line(0,1){.19171}}
\multiput(93.11,27.73)(-.0327,.15844){6}{\line(0,1){.15844}}
\multiput(92.92,28.68)(-.02987,.1176){8}{\line(0,1){.1176}}
\multiput(92.68,29.62)(-.03125,.10323){9}{\line(0,1){.10323}}
\multiput(92.4,30.55)(-.0323,.09154){10}{\line(0,1){.09154}}
\multiput(92.08,31.46)(-.03309,.0818){11}{\line(0,1){.0818}}
\multiput(91.71,32.36)(-.03369,.07354){12}{\line(0,1){.07354}}
\multiput(91.31,33.24)(-.031701,.061663){14}{\line(0,1){.061663}}
\multiput(90.86,34.11)(-.032159,.056156){15}{\line(0,1){.056156}}
\multiput(90.38,34.95)(-.032499,.051229){16}{\line(0,1){.051229}}
\multiput(89.86,35.77)(-.032736,.046783){17}{\line(0,1){.046783}}
\multiput(89.3,36.56)(-.032883,.042741){18}{\line(0,1){.042741}}
\multiput(88.71,37.33)(-.032951,.039042){19}{\line(0,1){.039042}}
\multiput(88.09,38.08)(-.032948,.035637){20}{\line(0,1){.035637}}
\multiput(87.43,38.79)(-.032882,.032486){21}{\line(-1,0){.032882}}
\multiput(86.74,39.47)(-.036033,.032515){20}{\line(-1,0){.036033}}
\multiput(86.02,40.12)(-.039438,.032476){19}{\line(-1,0){.039438}}
\multiput(85.27,40.74)(-.043136,.032364){18}{\line(-1,0){.043136}}
\multiput(84.49,41.32)(-.047176,.032167){17}{\line(-1,0){.047176}}
\multiput(83.69,41.87)(-.051618,.031877){16}{\line(-1,0){.051618}}
\multiput(82.86,42.38)(-.060579,.033726){14}{\line(-1,0){.060579}}
\multiput(82.01,42.85)(-.066814,.033334){13}{\line(-1,0){.066814}}
\multiput(81.15,43.28)(-.07394,.0328){12}{\line(-1,0){.07394}}
\multiput(80.26,43.68)(-.0822,.0321){11}{\line(-1,0){.0822}}
\multiput(79.35,44.03)(-.09192,.03119){10}{\line(-1,0){.09192}}
\multiput(78.44,44.34)(-.1036,.03){9}{\line(-1,0){.1036}}
\multiput(77.5,44.61)(-.1348,.03251){7}{\line(-1,0){.1348}}
\multiput(76.56,44.84)(-.15883,.03078){6}{\line(-1,0){.15883}}
\multiput(75.61,45.02)(-.19207,.02828){5}{\line(-1,0){.19207}}
\put(74.65,45.17){\line(-1,0){.966}}
\put(73.68,45.26){\line(-1,0){.969}}
\put(72.71,45.32){\line(-1,0){.971}}
\put(71.74,45.33){\line(-1,0){.97}}
\put(70.77,45.29){\line(-1,0){.968}}
\multiput(69.8,45.22)(-.2408,-.0303){4}{\line(-1,0){.2408}}
\multiput(68.84,45.09)(-.19133,-.03292){5}{\line(-1,0){.19133}}
\multiput(67.88,44.93)(-.13546,-.02967){7}{\line(-1,0){.13546}}
\multiput(66.93,44.72)(-.11723,-.03129){8}{\line(-1,0){.11723}}
\multiput(66,44.47)(-.10284,-.0325){9}{\line(-1,0){.10284}}
\multiput(65.07,44.18)(-.09114,-.0334){10}{\line(-1,0){.09114}}
\multiput(64.16,43.85)(-.07461,-.03124){12}{\line(-1,0){.07461}}
\multiput(63.26,43.47)(-.067501,-.031922){13}{\line(-1,0){.067501}}
\multiput(62.39,43.06)(-.061275,-.032445){14}{\line(-1,0){.061275}}
\multiput(61.53,42.6)(-.055762,-.032836){15}{\line(-1,0){.055762}}
\multiput(60.69,42.11)(-.050832,-.033116){16}{\line(-1,0){.050832}}
\multiput(59.88,41.58)(-.046384,-.033299){17}{\line(-1,0){.046384}}
\multiput(59.09,41.01)(-.04234,-.033398){18}{\line(-1,0){.04234}}
\multiput(58.33,40.41)(-.03864,-.033421){19}{\line(-1,0){.03864}}
\multiput(57.59,39.78)(-.035235,-.033377){20}{\line(-1,0){.035235}}
\multiput(56.89,39.11)(-.03369,-.034936){20}{\line(0,-1){.034936}}
\multiput(56.22,38.41)(-.032077,-.036423){20}{\line(0,-1){.036423}}
\multiput(55.57,37.68)(-.031997,-.039828){19}{\line(0,-1){.039828}}
\multiput(54.97,36.93)(-.033712,-.046085){17}{\line(0,-1){.046085}}
\multiput(54.39,36.14)(-.033569,-.050534){16}{\line(0,-1){.050534}}
\multiput(53.86,35.33)(-.033333,-.055467){15}{\line(0,-1){.055467}}
\multiput(53.36,34.5)(-.032991,-.060983){14}{\line(0,-1){.060983}}
\multiput(52.89,33.65)(-.032523,-.067213){13}{\line(0,-1){.067213}}
\multiput(52.47,32.77)(-.03191,-.07433){12}{\line(0,-1){.07433}}
\multiput(52.09,31.88)(-.0311,-.08258){11}{\line(0,-1){.08258}}
\multiput(51.75,30.97)(-.03341,-.10255){9}{\line(0,-1){.10255}}
\multiput(51.45,30.05)(-.03234,-.11695){8}{\line(0,-1){.11695}}
\multiput(51.19,29.12)(-.03088,-.13519){7}{\line(0,-1){.13519}}
\multiput(50.97,28.17)(-.02886,-.15919){6}{\line(0,-1){.15919}}
\multiput(50.8,27.21)(-.0324,-.2405){4}{\line(0,-1){.2405}}
\put(50.67,26.25){\line(0,-1){.967}}
\put(50.58,25.28){\line(0,-1){3.877}}
\multiput(50.68,21.41)(.0332,-.2404){4}{\line(0,-1){.2404}}
\multiput(50.81,20.45)(.02936,-.15909){6}{\line(0,-1){.15909}}
\multiput(50.98,19.49)(.0313,-.13509){7}{\line(0,-1){.13509}}
\multiput(51.2,18.55)(.03271,-.11684){8}{\line(0,-1){.11684}}
\multiput(51.47,17.61)(.03374,-.10244){9}{\line(0,-1){.10244}}
\multiput(51.77,16.69)(.03137,-.08248){11}{\line(0,-1){.08248}}
\multiput(52.11,15.78)(.03214,-.07423){12}{\line(0,-1){.07423}}
\multiput(52.5,14.89)(.032736,-.067109){13}{\line(0,-1){.067109}}
\multiput(52.93,14.02)(.033184,-.060878){14}{\line(0,-1){.060878}}
\multiput(53.39,13.17)(.033509,-.055361){15}{\line(0,-1){.055361}}
\multiput(53.89,12.34)(.033729,-.050428){16}{\line(0,-1){.050428}}
\multiput(54.43,11.53)(.031977,-.043423){18}{\line(0,-1){.043423}}
\multiput(55.01,10.75)(.032123,-.039726){19}{\line(0,-1){.039726}}
\multiput(55.62,9.99)(.032192,-.036322){20}{\line(0,-1){.036322}}
\multiput(56.26,9.27)(.032191,-.033171){21}{\line(0,-1){.033171}}
\multiput(56.94,8.57)(.035341,-.033265){20}{\line(1,0){.035341}}
\multiput(57.65,7.91)(.038746,-.033298){19}{\line(1,0){.038746}}
\multiput(58.38,7.27)(.042446,-.033263){18}{\line(1,0){.042446}}
\multiput(59.15,6.67)(.046489,-.033152){17}{\line(1,0){.046489}}
\multiput(59.94,6.11)(.050937,-.032955){16}{\line(1,0){.050937}}
\multiput(60.75,5.58)(.055866,-.03266){15}{\line(1,0){.055866}}
\multiput(61.59,5.09)(.061377,-.032251){14}{\line(1,0){.061377}}
\multiput(62.45,4.64)(.067601,-.031708){13}{\line(1,0){.067601}}
\multiput(63.33,4.23)(.07471,-.031){12}{\line(1,0){.07471}}
\multiput(64.22,3.86)(.09125,-.03311){10}{\line(1,0){.09125}}
\multiput(65.14,3.53)(.10294,-.03217){9}{\line(1,0){.10294}}
\multiput(66.06,3.24)(.11733,-.03092){8}{\line(1,0){.11733}}
\multiput(67,2.99)(.13555,-.02924){7}{\line(1,0){.13555}}
\multiput(67.95,2.78)(.19143,-.03231){5}{\line(1,0){.19143}}
\put(68.91,2.62){\line(1,0){.963}}
\put(69.87,2.5){\line(1,0){.968}}
\put(70.84,2.43){\line(1,0){.97}}
\put(71.81,2.4){\line(1,0){.971}}
\put(72.78,2.41){\line(1,0){.969}}
\put(73.75,2.47){\line(1,0){.965}}
\multiput(74.71,2.57)(.19198,.02889){5}{\line(1,0){.19198}}
\multiput(75.67,2.72)(.15873,.03128){6}{\line(1,0){.15873}}
\multiput(76.63,2.9)(.1347,.03294){7}{\line(1,0){.1347}}
\multiput(77.57,3.13)(.1035,.03033){9}{\line(1,0){.1035}}
\multiput(78.5,3.41)(.09182,.03148){10}{\line(1,0){.09182}}
\multiput(79.42,3.72)(.0821,.03236){11}{\line(1,0){.0821}}
\multiput(80.32,4.08)(.07384,.03304){12}{\line(1,0){.07384}}
\multiput(81.21,4.47)(.066708,.033545){13}{\line(1,0){.066708}}
\multiput(82.07,4.91)(.056441,.031657){15}{\line(1,0){.056441}}
\multiput(82.92,5.38)(.051517,.03204){16}{\line(1,0){.051517}}
\multiput(83.75,5.9)(.047074,.032317){17}{\line(1,0){.047074}}
\multiput(84.55,6.45)(.043033,.0325){18}{\line(1,0){.043033}}
\multiput(85.32,7.03)(.039335,.032601){19}{\line(1,0){.039335}}
\multiput(86.07,7.65)(.03593,.032629){20}{\line(1,0){.03593}}
\multiput(86.79,8.3)(.032779,.03259){21}{\line(1,0){.032779}}
\multiput(87.47,8.99)(.032835,.035741){20}{\line(0,1){.035741}}
\multiput(88.13,9.7)(.032827,.039146){19}{\line(0,1){.039146}}
\multiput(88.75,10.45)(.032748,.042845){18}{\line(0,1){.042845}}
\multiput(89.34,11.22)(.032587,.046887){17}{\line(0,1){.046887}}
\multiput(89.9,12.02)(.032336,.051332){16}{\line(0,1){.051332}}
\multiput(90.42,12.84)(.031982,.056257){15}{\line(0,1){.056257}}
\multiput(90.9,13.68)(.031506,.061763){14}{\line(0,1){.061763}}
\multiput(91.34,14.54)(.03346,.07364){12}{\line(0,1){.07364}}
\multiput(91.74,15.43)(.03283,.08191){11}{\line(0,1){.08191}}
\multiput(92.1,16.33)(.03201,.09164){10}{\line(0,1){.09164}}
\multiput(92.42,17.25)(.03092,.10333){9}{\line(0,1){.10333}}
\multiput(92.7,18.18)(.03371,.13451){7}{\line(0,1){.13451}}
\multiput(92.93,19.12)(.0322,.15854){6}{\line(0,1){.15854}}
\multiput(93.13,20.07)(.02999,.19181){5}{\line(0,1){.19181}}
\put(93.28,21.03){\line(0,1){.965}}
\put(93.38,21.99){\line(0,1){1.871}}
\put(57,21.86){\circle*{1.8}} \put(84.75,22.11){\circle*{.71}}
\put(85,22.36){\circle*{1}}
\qbezier(14.5,19.11)(33.88,15.11)(56.75,22.11)
\qbezier(57,22.86)(70.75,28.36)(84.5,22.86)
\qbezier(84.75,22.86)(110.88,15.74)(127.5,23.11)
\put(14.5,19.11){\circle*{1.58}}
\put(127.25,23.11){\circle*{1.58}} \thicklines
\qbezier(14.25,19.36)(25.75,22.86)(35.25,18.36)
\qbezier(57,21.86)(72.88,43.24)(85.25,22.11)
\put(22.75,15.36){\makebox(0,0)[cc]{$\cf$}}
\put(50.61,22.76){\makebox(0,0)[cc]{$\cf(t-a)$}}
\put(79.81,20.55){\makebox(0,0)[cc]{$\cf(t+b)$}}
\put(84.93,22.6){\circle*{1.69}}
\qbezier(35.32,18.08)(47.35,11.98)(56.66,21.86)
\qbezier(85.04,22.7)(103.43,-5.99)(127.29,22.91)
\end{picture}
\caption{Property ($T_2'$).}
\label{fig1}
\end{figure}
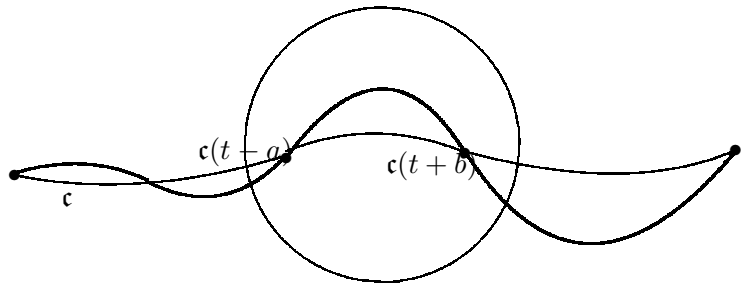
\item[{\rm ($T_2''$)}] Every
simple loop in $\free$ is contained in one piece.
\end{enumerate}
\end{proposition}

\proof Obviously ($T_1$) and ($T_2'$) imply ($T_2$). Therefore it
is enough to establish the implications $(T_1)\&
(T_2'')\Rightarrow (T_2')$ and $(T_1)\&(T_2)\Rightarrow (T_2'')$.
The second of these implications is given by Lemma \ref{cut}.

Suppose that ($T_1$) and ($T_2''$) hold for some space $\free$
with respect to some set of pieces $\pp$.

Let $\cf:[0,d]\to \free$ be a topological arc, $t\in [0,d]$, and
$a, b$ as in ($T_2'$). If $\cf':[0,d']\to \free$ is another
topological arc with the same endpoints as $c$, then
$K=\cf^{-1}(\cf'[0,d'])$ is a compact set containing $0$ and $d$.
Suppose that, say, $t-a \not\in K$. Let $\alpha$ be the supremum of
$K\cap [0, t-a]$ and $\beta$ be the infimum of $K\cap [t-a, d]$.
Then  $\alpha < t-a < \beta $. Since $\alpha, \beta\in K$, there
exist $\alpha', \beta'\in [0,d']$ such that
$\cf'(\alpha')=\cf(\alpha), \cf'(\beta')=\cf(\beta)$. The
restriction of $\cf$ to $[\alpha, \beta]$ and
the restriction of $\cf'$ to
$[\alpha',\beta']$ form a simple loop which is contained in one
piece by ($T_2''$). In particular $\cf([\alpha, \beta])$ is
contained in one piece. Since $[t-a,t+b]$ is the maximal interval containing
$t$ such that the restriction of $\cf$ to that interval is
contained in one piece, it follows that $b+a\ne 0$. Therefore the intersection of
the intervals $[\alpha,\beta]$ and $[t-a,t+b]$ has a non-empty
interior. Hence the pieces containing $\cf([\alpha,\beta])$ and
$\cf([t-a,t+b])$ must coincide by property ($T_1$). But this
contradicts the maximality of the interval $[t-a,t+b]$.
\endproof

\begin{remark}\label{geodsub} If a collection of subsets $\calp$ of a
geodesic metric space $X$ satisfy $(T_1)$ and $(T_2'')$, and each
set in $\calp$ is path connected then each set in $\calp$ is a
geodesic subspace. Thus if one replaces property $(T_2)$ by the
stronger property $(T_2'')$ in Definition \ref{tree} then one can
weaken the condition on $\calp$.
\end{remark}

\proof Let $M\in \calp$, let $x,y$ be two points in $M$ and let
$\mathfrak r$ be a topological arc joining $x$ and $y$ in $M$.
Suppose that a geodesic $\mathfrak g$ connecting $x$ and $y$ in
$X$ is not contained in $M$. Let $z\in {\mathfrak g}\setminus M$.
There exists a simple non-trivial bigon with one side a sub-arc in
$\mathfrak r$ and the other a sub-arc in $\mathfrak g$ containing
$z$. Property $(T_2'')$ implies that this bigon is contained in a
piece, and property $(T_1)$ implies that this piece is $M$. Hence
$z$ is in $M$, a contradiction.\endproof


\begin{lemma}\label{ctree}
For every $x\in \free$, $T_x$ is a closed subset of $\free$.
\end{lemma}

\proof Let $(y_n)$ be a sequence in $T_x$ converging to a point
$y$. Suppose that the geodesic $[x,y]$ intersects a piece $M$ in a
maximal non-trivial sub-arc $[\alpha , \beta]$. We can assume that the geodesic
$[y_n,y]$ intersects $[y_n,x]$ only in $y_n$. Otherwise we can
replace $y_n$ with the farthest from it intersection point between
these two geodesics. By property ($T_2'$) the arc $[x, y_n]\cup
[y_n, y]$ must contain $[\alpha, \beta]$. Since $y_n\in T_x$, it
follows by Lemma \ref{xy} that $[\alpha , \beta]\subset [y_n,y]$
and so $\dist(y_n,y)\geq \dist(\alpha , \beta) >0 $. This
contradicts $\dist(y_n,y)\to 0$. We conclude that $[x,y]$
intersects every piece in at most one point and that $y\in
T_x$.\endproof


\begin{lemma}\label{projc}
The projection of $\free$ onto any of the pieces is a metric
retraction.
\end{lemma}

 \proof Let $M$ be a piece, $x,y$ two
 points in $\free$ and $[x,y]$ a geodesic
 joining them. If $[x,y]\cap M = \emptyset$ then $[x,y]$ projects onto one point
 $z$, by Corollary \ref{projpath}, and $d(x,y)\geq
 d(z,z)=0$.

If $[x,y]\cap M = [\alpha ,\beta ]$ then $\alpha$ is the projection of $x$ onto $M$ and $\beta $ is the projection of $y$ onto $M$, by Corollary \ref{projpath}. Obviously $d(x,y)\geq
 d(\alpha ,\beta)$. \endproof

\begin{lemma}\label{lcount} Let $\pgot\colon [0, l]\to \free$
be a path in a tree-graded space $\free$. Let $U_\pgot$ be the
union of open subintervals $(a,b)\subset [0,l]$ such that the
restriction of $\pgot$ onto $(a,b)$ belongs to one piece (we
include the trees $T_x$ into the set of pieces). Then $U_\pgot$ is
an open and dense subset of $[0,l]$.
\end{lemma}

\proof Suppose that $U_\pgot$ is not dense. Then there exists a
non-trivial interval $(c,d)$ in the complement $[0,l]\setminus
U_\pgot$. Suppose that the restriction $\pgot'$ of $\pgot$ on
$(c,d)$ intersects a piece $P$ in two points
$y=\pgot(t_1),z=\pgot(t_2)$. We can assume that $y$ is not in the
image of $(t_1,t_2]$ under $\pgot$. Since $y\not\in U_\pgot$ there
is a non-empty interval $(t_1,t_3)$ such that the restriction of
$\pgot$ onto that interval does not intersect $P$. Let $t>t_1$ be
the smallest number in $(t_1,t_2]$ such that $z'=\pgot(t)$ is in
$P$. Then $z'\ne y$. Applying Corollary \ref{projA} to the
restriction of $\pgot$ onto $[t_1,t]$, we get a contradiction.
This means that $\pgot'$ intersects every piece in at most one
point. Therefore $\pgot'$ is contained in a tree $T_x$ for some
$x$, a contradiction.
\endproof

\begin{proposition}\label{pi1} Let $\free$ be a tree-graded space with
the set of pieces $\pp$. If the pieces in $\pp$ are locally
uniformly contractible then $\pi_1(\free)$ is the free product of
$\pi_1(M)$, $M\in \pp$.
\end{proposition}

\proof We include all trees $T_x$ into $\pp$. Fix a base point $x$
in $\free$ and for every piece $M_i\in \pp$ let $y_i$ be the
projection of $x$ onto $M_i$, and let $\g_i$ be a geodesic
connecting $x$ and $y_i$. We identify $\pi_1(M_i)$ with the
subgroup $G_i=\g_i\pi_1(M_i,y_i)\g_i^{-1}$ of $\pi_1(\free, x)$.
Consider an arbitrary loop $\pgot\colon [0,l]\to\free$ in $\free$
based at $x$. Let $\pgot'$ be the image of $\pgot$. Let
$\pp_\pgot$ be the set of pieces from $\pp$ which are intersected
by $\pgot'$ in more than one point. By Lemma \ref{lcount} the set
$\pp_\pgot$ is countable.

Let $M\in \pp_\pgot$. The projection $\pgot_M$ of $\pgot'$ onto
$M$ is a loop containing the intersection $\pgot'\cap M$. Let us
prove that $\pgot_M=\pgot'\cap M$. If there exists a point $z\in
\pgot_M\setminus \pgot'$ then $z$ is a projection of some point
$y\in \pgot'\setminus M$ onto $M$. By Corollary \ref{projA}, a
subpath of $\pgot$ joining $y$ with a point in $\pgot'\cap P$ must
contain $z$, a contradiction.

Therefore $\pgot'$ is a union of at most countably many loops
$\pgot_i$, $i\in \N$, contained in pieces from $\pp_\pgot$. By
uniform local contractibility of the pieces, all but finitely many
loops $\pgot_i$ are contractible inside the corresponding pieces.
Consequently, in the fundamental group $\pi_1(\free)$, $\pgot$ is
a product of finitely many loops from $G_i$. Hence $\pi_1(\free,
x)$ is generated by the subgroups $G_i$.

It remains to prove that for every finite sequence of loops
$\pgot_i\in G_i$, $i=1,...,k$, if $M_i\ne M_{j}$ for $i\neq j$,
and if the loops $\pgot_i$ are not null-homotopic in $M_i$, then
the loop $\pgot_1\pgot_2...\pgot_n$ is not null-homotopic in
$\free$. Suppose that $\pgot$ is null-homotopic, and that $\gamma:
t\to \pgot(t)$ is the homotopy, $\pgot(0)=\pgot$, $\pgot(1)$ is a
point. Let $\pi_i$ be the projection of $\free$ onto $M_i$. Lemma
\ref{projc} implies that $\pi_i \circ \gamma : t\to \pgot_i'(t)$
is a homotopy which continuously deforms $\pgot_i'$ in $M_i$ into
a point. Hence each of the loops $\pgot_i$ is null-homotopic, a
contradiction.
\endproof

\subsection{Modifying the set of pieces}

\begin{lemma}[gluing pieces together]\label{retr}
Let $\free$ be a space which is tree-graded with respect to $\pp\{ M_k \mid k\in K \}$.
\begin{itemize}
\item[(1)] Let $Y=\bigcup_{k\in F}M_k$ be a finite connected union
of pieces. Then $\free$ is tree-graded with respect to $\pp'=\{
M_k \mid k\in K\setminus F \} \cup \{Y \} $.

\item[(2)] Let $\cf$ be a topological arc in $\free$ (possibly a
point) and let $Y(\cf)$ be a set of the form $\cf\cup \bigcup_{j\in J}M_j$, where $J$ is a subset
of $K$ such that every $M_j$ with $j\in J$ has a non-empty
intersection with $\cf$, and $J$ contains all $i\in K$ such that
$M_i\cap \cf$ is a non-trivial arc.

Then $\free$ is tree-graded with respect to $\pp'=\{ M_k \mid k\in
K\setminus J \} \cup \{Y(\cf) \} $.

\item[(3)] Let $\{\cf_i \, ;\, i\in F\}$ be a finite collection of
topological arcs in $\free$ and let $Y(\cf_i)=\cf_i\cup
\bigcup_{j\in J_i}M_j$ be sets defined as in (2). If
$Y=\bigcup_{i\in F}Y(\cf_i)$ is connected then $\free$ is
tree-graded with respect to $\pp'=\{ M_k \mid k\in K\setminus
\bigcup_{i\in F}J_i \} \cup \{Y\} $.
\end{itemize}
\end{lemma}

\begin{remark}
In particular all properties on projections on pieces obtained
till now hold for sets $Y$ defined as in (1)-(3). We shall call
sets of the form $Y(\cf)$ {\em sets of type $Y$}.
\end{remark}

\proof \textbf{(1)} We first prove that $Y$ is convex. Every $y,
y'\in Y$ can be joined by a topological arc $\cf:[0,d]\to Y $. By
Corollary \ref{strconv}, we may write $\cf[0,d] = \bigcup_{k\in
F'}\left[ \cf[0,d]\cap M_k \right]$, where $F'\subset F$ and
$\cf[0,d]\cap M_k$ is a point or an arc. Property $(T_1)$ implies
that every two such arcs have at most one point in common.
Therefore there exists a finite sequence $t_0=0 <t_1<t_2< \dots <
t_{n-1}< t_n=d $ such that $\cf[t_i,t_{i+1}]=\cf[0,d]\cap
M_{k(i)}\, ,\, k(i) \in F',$ for every $i\in \{ 0,1,\dots n-1 \}$.
Property $(T_2')$ implies that every geodesic between $y$ and $y'$
must contain $\cf(t_1)\, ,\, \cf(t_2)\, ,\dots \cf(t_{n-1})$.
Hence every such geodesic is of the form $[y, \cf(t_1)]\cup
[\cf(t_1), \cf(t_2)]\cup \dots \cup [\cf(t_{n-1}), y]$, so by
Corollary \ref{strconv} it is contained in $Y$.

For every $k\in K\setminus F$, $M_k\cap Y$, if non-empty, is a
convex set composed of finitely many points. Hence it is a point. This and the previous discussion imply that $\free$ is
tree-graded with respect to $\pp'$.

\me

 \textbf{(2)} In order to prove that $Y$ is convex, let $\g$ be a
geodesic joining two points $x,y\in Y$. We show that $\g$
is inside $Y$.

\me

\noindent\textbf{Case I.} Suppose that $x,y\in \cf$. Consider a
point $z=\g(t)$ in $\g$. Take the maximal interval $[t-a,t+b]$
such that $\g([t-a, t+b])$ is contained in one piece $M$. If
$a+b\ne 0$ then by property $(T_2')$ the path $\cf$ must pass
through $\g(t-a)$ and $\g(t+b)$. By part (1) of Corollary
\ref{strconv} the (non-trivial) subarc of $\cf$ joining
$\g(t-a)$ and $\g(t+b)$ is contained in $M$. Then $M$ is one of
the pieces contained in $Y$. Therefore $z\in Y$. If $a+b=0$ then
again by $(T_2')$ the curve $\cf$ must pass through $z$, so $z\in
Y$. We conclude that in
both cases $z\in Y$.

\me

\noindent\textbf{Case II.}  Suppose that $x\in \cf$ and $y\in
M\setminus \cf$, where $M$ is a piece in $Y$. By the definition of
$Y$, $M$ has a non-trivial intersection with $\cf$. If $x\in M$,
we can use the convexity of $M$ (Corollary \ref{strconv}). So
suppose that $x\not\in M$.

Let $\alpha$ be the projection of $x$ onto $M$. By Corollary
\ref{projA}, part (2), $\alpha \in \cf$. Then the sub-arc $\cf'$
of $\cf$ with endpoints $x$ and $\alpha$ forms together with the
geodesic $[\alpha ,y]\subseteq M$ a topological arc. Property
$(T_2')$ implies that $\alpha \in \g$. Corollary \ref{strconv},
part (1), implies that the portion of $\g$ between $\alpha$ and
$y$ is contained in $Y$. For the remaining part of $\g$ we apply
the result in Case I of the proof (since both endpoints of that
part of $\g$ belong to $\cf$).

\me

\noindent\textbf{Case III.}  Suppose that $x\in M_1\setminus \cf$
and that $y\in M_2\setminus \cf$. Let $\alpha$ be the projection
of $x$ onto $M_2$. As before, we obtain that $\alpha \in \cf$,
$\alpha \in \g$ and that the portion of $\g$ between $\alpha$ and
$y$ is contained in $M_2$, hence in $Y$. For the remaining part of
$\g$ we apply the result of Case II.

\me

\textbf{(3)} We argue by induction on the size $k$ of the set $F$.
The statement is true for $k=1$ by part (2) of this Proposition.
Suppose it is true for some $k\ge 1$. Let us prove it for $k+1$.
We have two cases.

\me

\noindent{\textbf{Case I.}} Suppose that there exist $i,j\in F,
i\neq j$, such that the intersection $\cf_i\cap Y(\cf_j)$ is not
empty. According to part (2) of the Proposition and Corollary
\ref{strconv}, part (2), the intersection is a sub-arc, and
$\free$ is tree-graded with respect to $\pp_j'=\{ M_k \mid k\in
K\setminus J_j \} \cup \{Y(\cf_j) \} $. Let $Y'(\cf_i)=Y(\cf_i)\cup Y(\cf_j)$. Then $Y'(\cf_i)$ is a set defined as in part (2) of the Proposition but with $\pp$ replaced
by $\pp_j'$. Thus we can
write $Y=Y'(\cf_i) \cup \bigcup_{s\in F\setminus \{i,j \}}
Y(\cf_s)$ and use the induction hypothesis.

\me

\noindent{\textbf{Case II.}} For every $i,j\in F, i\neq j,$ we
have $\cf_i\cap Y(\cf_j)= \emptyset$.

Then there are no pieces that appear in both $Y(\cf_i)$ and $Y(\cf_j)$ for $i\ne j\in F$. Hence by $(T_1)$, for every $k\in
J_i$, $l\in J_j$, $M_k\cap M_l$ consists of at most one point. By
part (2) of the Proposition and Corollary \ref{projA} that
point must be equal to the projection of $\cf_i$ onto $Y(\cf_j)$.
Therefore $Y(\cf_i)\cap Y(\cf_j)$
is either empty or one point. This implies that $\free$ is
tree-graded with respect to $\pp''=\{ M_k \mid k\in K\setminus
\bigcup_{i\in F}J_i \} \cup \{Y(\cf_i)\mid i\in F \} $. It remains
to apply part (1) of the Proposition. \endproof

\me


\begin{definition}\label{dbouq}
Let $(M_1,x_1),\, (M_2,x_2),\dots,\, (M_k,x_k)$ be finitely many
metric spaces with fixed basepoints. The \textit{bouquet} of these
spaces, denoted by $\bigvee_{i=1}^k (M_i,x_i)$, is the metric
space obtained from the disjoint union of all $M_i$ by identifying
all the points $x_i$. We call the point $x$ thus obtained
\textit{the cut-point of the bouquet}. The metric on
$\bigvee_{i=1}^k (M_i,x_i)$ is induced by the metrics on $M_i$ in
the obvious way.
\end{definition}

Clearly each $M_i$ is a closed subset of the bouquet
$\bigvee_{i=1}^k (M_i,x_i)$. It is also clear that the bouquet is
a geodesic metric space if and only if all $M_i$ are geodesic
metric spaces.

\begin{lemma}[cutting pieces by cut-points]\label{bouq}
Let $\free$ be a space which is tree-graded with respect to
$\pp=\{ M_k \mid k\in K \}$. Let $I\subset K$ be such that for
every $i\in I$ the piece $M_i$ is the bouquet of finitely many
subsets of it, $\{ M_i^j \}_{j\in F_i}$, and its cut-point is
$x_i$.

Then $\free$ is tree-graded with respect to the set
$$
\pp'= \{ M_k \mid k\in K\setminus I \} \cup \{ M_i^j \mid j\in
F_i,\, i\in I \}\, .
$$
\end{lemma}

\proof Since $M_i^j\cap M_k \subset M_i\cap M_k$ for $i\in I,\,
k\in K\setminus I$, and $M_i^j\cap M_t^s \subset M_i\cap M_t$ for
$i\neq t,\, i,t\in I$, property $(T_1)$ for $(\free , \pp' )$ is
an immediate consequence of property $(T_1)$ for $(\free , \pp )$.

Let $\Delta$ be a simple geodesic triangle. Property $(T_2)$ for
$(\free , \pp )$ implies that either $\Delta \subset M_k$ for some
$k\in K\setminus I$ or $\Delta \subset M_i$ for some $i\in I$. We
only need to consider the second case. Assume that $\Delta $ has a
point in $M_i^{j_1}$ and a point in $M_i^{j_2}$, with $j_1\neq j_2$.
Then $x_i$ is a cut-point for $\Delta $. This contradicts the fact
that $\Delta $ is a simple loop. We conclude that there exists $j\in
F_i$ such that $M_i^j$ contains~$\Delta$. Thus $\pp'$ satisfies
$(T_2)$. \endproof

Lemma \ref{txy} implies that two trees $T_x$ and $T_y$ are either
disjoint or coincident. Let $\{ T_i \mid i\in I \}$ be the
collection of all the trees $\{ T_x \mid x\in \free \}$.

\begin{remark} \label{treesrem}
The set $\pp'=\pp\cup \{ T_i \mid i\in I \}$ also satisfies
properties $(T_1)$ and $(T_2)$. Therefore all the properties and
arguments done for $\free$ and $\pp$ up to now also hold for
$\free $ and $\pp'$. In this case, $T_x=\{ x \}$ for every $x\in
\free$. The disadvantage of this point of view is that trees $T_x$
always have cut-points.
\end{remark}

\subsection{Geodesics in tree-graded spaces}

\Notat \quad For every path $\pgot$ in a metric space $X$, we denote
the start of $\pgot$ by $\pgot_-$ and the end of $\pgot $ by $\pgot_+$. \me

\me

\begin{lemma}\label{sir}
Let $\g=\g_1\g_2\dots\g_{2m}$ be a curve in a tree-graded space
$\free$ which is a composition of geodesics. Suppose that all
geodesics $\g_{2k}$ with $k\in \{ 1,\dots ,m-1 \}$ are non-trivial
and for every $k\in \{ 1,\dots ,m \}$ the geodesic $g_{2k}$ is
contained in a piece $M_k$ while for every $k\in \{ 0,1,\dots ,m-1
\}$ the geodesic $\g_{2k+1}$ intersects $M_k$ and $M_{k+1}$ only
in its respective endpoints. In addition assume that if
$\g_{2k+1}$ is empty then $M_k\ne M_{k+1}$. Then $\g$ is a
geodesic.
\end{lemma}

\proof Suppose that $\g$ is not simple. By $(T_2'')$, any simple
loop formed by a portion of $\g$ has to be contained in one piece
$M$. On the other hand the loop must contain the whole neighborhood
of one vertex $(\g_i)_+=(\g_{i+1})_-$ in $\g$. Let $k$ be such that
$\{ \g_i, \g_{i+1}  \} = \{ \g_{2k}, \g_{2k\pm 1}  \}$. The
intersection of $M$ and $M_k$ contains a sub-arc of $\g_{2k}$,
whence $M=M_k$. At the same time, $M$ contains a subarc of
$\g_{2k\pm 1}$ or  (if $\g_{2k\pm 1}$ is empty) of $\g_{2k-2}$. In
all cases we immediately get a contradiction.

Therefore $\g$ is simple and has two distinct endpoints $x,y$.
Consider any geodesic $\mathfrak r$ joining $x$ and $y$. By $(T_2')$
$\mathfrak r$ contains all the endpoints of all geodesics $\g_i$.
Therefore the length of $\g$ coincides with the length of $\mathfrak
r$ and $\g$ is itself a geodesic. \endproof

\begin{cor}\label{transproj}
Let $M$ and $M'$ be two distinct pieces in a tree-graded space
$\free$. Suppose that $M'$ projects onto $M$ in $x$ and $M$
projects on $M'$ in $y$. Let $A$ be a set in $\free$ that projects
onto $M'$ in $z \neq y$. Then $A$ projects onto $M$ in $x$ and
$\dist(A,M)\geq \dist(M',M)$.
\end{cor}

\proof Let $a\in A$ and let $[a,z]$, $[z, y]$ and $[y,x]$ be
geodesics. Then $\g_a=[a,z]\cup [z, y]\cup [y, x]$ is a geodesic,
according to Lemma \ref{sir}. It cannot intersect $M$ in a
sub-geodesic, because $[z,y]\cup [y, x]$ intersects $M$ in $x$. Hence
$\g_a\cap M=\{x\}$ and $x$ is the projection of $a$ onto $M$. Also
$\dist(a,x)\geq \dist(y, x)$.
\endproof



\subsection{Cut-points and tree-graded spaces}

\begin{remark}[about singletons] \label{singletons11}
Notice that if $\free$ is tree graded with respect to $\pp$ then we
can always add some or all one-point subsets (singletons) of $\free$
to $\pp$, and $\free$ will be tree-graded with respect to a bigger
set of pieces. To avoid using extra pieces, we shall always assume
that pieces cannot contain other pieces. Property ($T_1$) guarantees
that this only restricts using singletons as pieces.
\end{remark}

Property $(T_2')$ implies that any tree-graded space containing
more than one piece has a global cut-point. Here we shall show
that any geodesic metric space with cut-points has a uniquely
determined collection of pieces with respect to which it is
tree-graded.

In order to do this, we need to define a partial order relation on
the set of collections of subsets of a space. If $\pp$ and $\pp'$
are collections of subsets of $X$ and a space $X$ is tree-graded
with respect to both $\pp$ and $\pp'$, we write $\pp \prec \pp'$ if
for every set $M \in \pp$ there exists $M'\in \pp'$ such that
$M\subset M'$. The relation $\prec$ is a partial order because by
Remark \ref{singletons}, pieces of $\pp$ (resp. $\pp'$) cannot
contain each other.

\begin{lemma}\label{cutting}
Let $X$ be a complete geodesic metric space containing at least
two points and let $\calc$ be a non-empty set of global cut-points
in $X$.
\begin{itemize}
\item[(a)] There exists the largest in the sense of $\prec$  collection $\calp$
of subsets of $X$ such that
\begin{itemize}
\item $X$ is tree-graded with respect to $\pp$;
\item any piece in $\pp$ is either a singleton or a set
with no global cut-point from $\calc$.
\end{itemize}
Moreover the intersection of any two distinct pieces from $\calp$ is
either empty or a point from $\calc$.

\item[(b)] Let $X$ be a homogeneous space with a cut-point. Then every point in
$X$ is a cut-point, so let $\calc=X$.  Let $\pp$ be the set of
pieces defined in part (a).  Then for every $M\in \calp$ every $x\in
M$ is the projection of a point $y\in X\setminus M$ onto $M$.
\end{itemize}
\end{lemma}

\proof (a) Let $\calp$ be the set of all maximal path connected
subsets $M$ with the property that either $|M|=1$ or cut-points of
$M$ do not belong to $\calc$. The existence of maximal subsets with
this property immediately follows from Zorn's lemma.

Any $M\in \calp$ is closed. Indeed, let $\bar M$ be the closure of
$M$ in $X$ and suppose that $\bar M\ne M$. Let $a\in \bar M\setminus
M$. There exists a sequence of points $(a_n)$ in $M$ converging to
$a$. Let $M'$ be the union of $M$ and geodesics $[a,a_n]$,
$n=1,2,...$ (one geodesic for each $n$). By construction, the set
$M'$ is path connected. Let us prove that cut-points of $M'$ do not
belong to $\calc$. This will contradict the maximality of $M$.

Let $c\in \calc\cap M'$,  $x,y\in M'\setminus \{ c\}$. We want to
connect $x$ and $y$ with a path avoiding $c$. If $x,y \in
M\setminus \{ c\}$ then we are done.

Suppose that $x \in M\setminus \{ c\}$ and $y\in [a_n,a]$ for
some $n$. The point $x$ can be connected by
some path $\pgot_k\subseteq M$ avoiding $c$
with $a_k$ for every $k\in \N$.

If $c\not \in [a_n,y]$ then the path $\pgot_n\cup [a_n,y]\subseteq
M'$ avoids $c$ and we are done.

If $c \in [a_n,y]$ then $\dist (c,a)>\dist (y,a)$. In particular
$c$ is not in $[a,a_m]$ for $m$ large enough. Then we join $y$
with $x$ by a path $[y,a]\cup [a,a_m]\cup \pgot_m$ avoiding $c$.

It remains to consider the case when $x \in [a_m,a]$ and $y\in
[a_n,a]$ for some $m, n$. If $c\not \in [a_m,x]$ then we can replace
$x$ with $a_m$ and use the previous argument. Likewise if $c\not
\in [a_n,y]$. If $c\in [a_m,x]\cap [a_n,y]$ then we join $x$ and
$y$ in $X\setminus \{ c\}$ by $[x,a]\cup [a,y]$.

Let $M_1, M_2$ be distinct sets from $\calp$, $c\in \calc$.
Suppose that $M_1\cap M_2$ contains a point $x$ that is different
from $c$. Then any point $z_i\in M_i$, $z_i\ne c$, $i=1,2$, can be
joined with $x$ by a path in $M_i$ avoiding $c$. Hence $z_1$ and
$z_2$ can be joined in $M_1\cup M_2$ by a path avoiding $c$.
Consequently if $M_1\cap M_2$ contains more than one point or contains a
point not from $\calc$, we get a contradiction with the
maximality of $M_i$. Thus $\calp$ satisfies $(T_1)$ and the
intersection of any two sets from $\calp$ is in $\calc$ or empty.

To prove $(T_2'')$ notice that every non-trivial simple loop is
path connected and does not have cut-points, hence it is contained
in some $M$.

The fact that each piece $M\in \calp$ is a geodesic subset follows
from Remark \ref{geodsub}.

Suppose that $X$ is tree-graded with respect to another collection
of pieces $\pp'$ that contains only singletons and pieces without
cut-points from $\calc$. Let $M'\in\pp'$. Then $M'$ is contained in
a maximal path-connected subset which is either a singleton or
without cut-point in $\calc$, that is $M'\subset M$ for some $M\in
\pp$. Thus $\pp'\prec \pp$. Hence $\pp$ is the largest in the sense
of $\prec$ collection of subsets of $X$ satisfying the conditions of
part (a).

(b) Let $M\in \calp$. Since $M\neq X$ it follows that one point
$x_0\in M $ is the projection on $M$ of a point $y_0\in X\setminus
M$. If $M$ is a point this ends the proof. Suppose in the sequel
that $M$ has at least two points. Let $[y_0,x_0]$ be a geodesic
joining $y_0$ and $x_0$ and let $[x_0,z_0]$ be a geodesic in $M$. By
the definition of the projection, $[y_0,x_0]\cap [x_0,z_0]=\{x_0\}$.
Let $x$ be an arbitrary point in $M$. Consider an isometry $g$ such
that $g(x_0)=x$. Let $[y,x]$ and $[x,z]$ be the respective images of
$[y_0,x_0]$ and $[x_0,z_0]$ under $g$. If $g(M)=M$ then $x$ is the
projection of $y$ on $M$. Suppose $g(M)\neq M$. Then $g(M)\cap M =
\{ x \}$, hence $[x,z]\subset g(M)$ intersects $M$ in $x$. Corollary
\ref{projA} implies that $z$ projects on $M$ in $x$.\endproof

\begin{remarks}\label{quest}
(1) In general not every point in $\calc$ is the intersection point
of two distinct pieces. An example is an $\R$-tree without endpoints
$X$, $\calc=X$, in which case $\calp$ is the set of all singleton
subsets of $X$.

(2) Lemma \ref{cutting} implies that every asymptotic cone of a
group which has a cut-point is tree-graded with respect to a
uniquely determined collection of pieces each of which is either a
singleton or a closed geodesic subset without cut-points.
\end{remarks}

\section{Ultralimits and asymptotic cones}\label{ULAC}

\subsection{Preliminaries}

Most of the interesting examples of tree-graded spaces that we
know are asymptotic cones of groups. In this section, we start
with giving the definitions of ultralimit, asymptotic cone and
related objects (most of these definitions are well known). We
show that the collection of asymptotic cones of a space is closed
under ultralimits. We also show that simple geodesic triangles in
ultralimits and asymptotic cones can be approximated by
ultralimits of polygons with certain properties. As a consequence
we show that the family of tree-graded spaces is also
 closed under ultralimits. These results play a central part in
 the theorems obtained in Sections \ref{amc} and \ref{exAD}.

\noindent \textit{Convention}: In the sequel $I$ will denote an
arbitrary countable set.

\begin{definition}[ultrafilter]
A (non-principal\fn{We shall only use non-principal ultrafilters
in this paper, so the word non-principal will be omitted.})
ultrafilter $\omega$ over $I$ is a set of subsets of $I$
satisfying the following conditions:

\begin{enumerate}

\item If $A, B\in \omega$ then $A\cap B\in \omega$;

\item If $A\in \omega$, $A\subseteq B\subseteq I$, then $B\in
\omega$;

\item For every $A\subseteq I$ either $A\in\omega$ or $I\setminus
A\in\omega$;

\item No finite subset of $I$ is in $\omega$.
\end{enumerate}
\end{definition}

Equivalently $\omega$ is a finitely additive  measure on the class
$\mathcal{P}(I)$ of subsets of $I$ such that each subset has
measure either $0$ or $1$ and all finite sets have measure 0. If
some statement $P(n)$ holds for all $n$ from a set $X$ belonging
to an ultrafilter $\omega$, we say that $P(n)$ holds {\em
$\omega$-almost surely}.

\begin{remark}\label{udisj}
By definition $\omega$ has the property that $\omega
(\sqcup_{i=1}^m A_i)=1$  (here $\sqcup$ stands for disjoint union)
implies that there exists $i_0\in \{ 1,2,\dots ,m\}$ such that
$\omega (A_{i_0})=1$ and $\omega (A_{i})=0$ for every $i\neq i_0$.
This can be reformulated as follows: let $P_1(n),\, P_2(n),\,
\dots ,\, P_m(n)$ be properties such that for any $n\in I$ no two
of them can be true simultaneously. If the disjunction of these
properties holds $\omega$-almost surely then there exists $i\in \{
1,2,\dots , m\}$ such that $\omega$-almost surely $P_i(n)$ holds
and all $P_j(n)$ with $j\neq i$ do not hold.
\end{remark}

\begin{definition}[$\omega$-limit] Let $\omega$ be an ultrafilter over $I$. For every
sequence of points $(x_n)_{n\in I}$ in a topological space $X$,
its $\omega$-limit $\lm_\omega x_n$ is a point $x$ in $X$ such
that for every neighborhood $U$ of $x$ the relation $x_n\in U$
holds $\omega$-almost surely.
\end{definition}

\begin{remark} If $\omega$-limit $\lm_\omega x_n$ exists then it
is unique, provided the space $X$ is Hausdorff. Every sequence of
elements in a compact space has an $\omega$-limit \cite{Bou}.
\end{remark}

\begin{definition}[ultraproduct] For every sequence of sets $(X_n)_{n\in I}$ the
{\em ultraproduct} $\Pi X_n/\omega$ corresponding to an
ultrafilter $\omega$ consists of  equivalence classes of sequences
$(x_n)_{n\in I}$, $x_n\in X_n$, where two sequences $(x_n)$ and
$(y_n)$ are identified if $x_n=y_n$ $\omega$-almost surely. The
equivalence class of a sequence $(x_n)$ in $\Pi X_n/\omega$ is
denoted by $(x_n)^\omega$. In particular, if all $X_n$ are equal
to the same $X$, the ultraproduct is called the {\em ultrapower}
of $X$ and is denoted by $X^\omega$.
\end{definition}

Recall that if $G_n$, $n\ge 1$, are groups then $\Pi G_n/\omega$
is again a group with the operation
$(x_n)^\omega(y_n)^\omega=(x_ny_n)^\omega$.

\begin{definition}[$\omega$-limit of metric spaces] Let
$(X_n,\dist_n)$, $n\in I$, be a sequence of metric spaces and let
$\omega$ be an ultrafilter over $I$. Consider the ultraproduct
$\Pi X_n/\omega$ and an {\em observation point} $e=(e_n)^\omega$
in $\Pi X_n/\omega$. For every two points $x=(x_n)^\omega,
y=(y_n)^\omega$ in $\Pi X_n/\omega$ let
$$D(x,y)=\lm_\omega \dist_n(x_n,y_n)\, .$$ The function
$D$ is a pseudo-metric on $\Pi X_n/\omega$ (i.e. it satisfies the
triangle inequality and the property $D(x,x)=0$, but for some
$x\ne y$, the number $D(x,y)$ can be $0$ or $\infty$). Let $\Pi_e
X_n/\omega$ be the subset of $\Pi X_n/\omega$ consisting of
elements which are finite distance from $e$ with respect to $D$.
The {\em $\omega$-limit} $\lio{X_n}_e$ {\em of the metric spaces}
$(X_n,\dist_n)$ {\em relative to the observation point} $e$ is the
metric space obtained from $\Pi_e X_n/\omega$ by identifying all
pairs of points $x,y$ with $D(x,y)=0$. The equivalence class of a
sequence $(x_n)$ in $\lio{X_n}_e$ is denoted by $\lio{x_n}$.
\end{definition}

\begin{remark}[changing the observation point]\label{equality}
It is easy to see that if $e,e'\in \Pi X_n/\omega$ and
$D(e,e')<\infty$ then $\lio{X_n}_e=\lio{X_n}_{e'}$.
\end{remark}

\begin{definition}[asymptotic cone] Let $(X,\dist)$ be a metric
space, $\omega$ be an ultrafilter over a set $I$, $e=(e_n)^\omega$
be an observation point. Consider a sequence of numbers
$d=(d_n)_{n\in I}$ called {\em scaling constants} satisfying
$\lm_\omega d_n=\infty$.

In the ultrapower $X^\omega$ we define the subset $X_e^\omega
=\Pi_e X_n/\omega$, where $(X_n, \dist_n)=(X,\dist /d_n)$. We call
it \textit{ultrapower of $X$ with respect to the observation point
$e$}.

The $\omega$-limit $\lio{X,\frac{\dist}{d_n}}_e$ is called an {\em
asymptotic cone of $X$.} It is denoted by $\co{X;e,d}$  (see
\cite{Gr1}, \cite{Gr2}, \cite{VDW}).
\end{definition}

\begin{definition}
For a sequence $(A_n), n\in I,$ of subsets of $(X,\dist)$ we denote
by $\lio{A_n}$ the subset of $\co{X; e, d}$ that consists of all the
elements $\lio{x_n}$ such that $x_n\in A_n$ $\omega$-almost surely.
Notice that if $\lim_\omega \frac{\dist(e_n,A_n)}{d_n}=\infty $ then
the set $\lio{A_n}$ is empty.
\end{definition}

\begin{remark} \label{rkcomp} It is proved in \cite{VDW} that any asymptotic cone
of a metric space is complete. The same proof gives that
$\lio{A_n}$ is always a closed subset of the asymptotic cone
$\co{X;e,d}$.
\end{remark}

\begin{definition}[quasi-isometries]
A {\em quasi-isometric embedding} of a metric space $(X,\dist_X)$
into a metric space $(Y,\dist_Y)$ is a map $\q\colon X\to Y$ such
that

$$
\frac{1}{L}\dist_X(x,x')-C\leq \dist_Y(\q(x),\q(x'))\leq
L\dist_X(x,x')+C, \hbox{ for all } x, x'\in X.
$$

In particular if  $(X,\dist_X)$ is an interval of the real line
$\R$ then $\q$ is called a {\em quasi-geodesic} or an $(L,C)$-{\em
quasi-geodesic}.

A \textit{quasi-isometry} is a quasi-isometric embedding $\q\colon
X\to Y$ such that there exists a quasi-isometric embedding
$\q'\colon Y\to X$ with the property that $\q\circ \q'$ and
$\q'\circ \q$ are at finite distance from the identity maps.
\end{definition}

\begin{remark}[quasi-injectivity] Although a quasi-isometric
embedding is not necessarily injective, a weaker version of
injectivity holds: If $\q$ is an $(L,C)$-quasi-isometric embedding
then $\dist (x,y)>LC$ implies $\dist (\q (x), \q(y))>0$.
\end{remark}

\begin{definition}[Lipschitz maps] Let $L\ge 1$. A map
$\q\colon (X,\dist_X)\to (Y,\dist_Y)$ is called {\em Lipschitz} if
$$\dist_Y(\q(x),\q(x'))\le L\dist_X(x,x')$$ for every $x,x'\in X$.
The map $\q$ is called {\em bi-Lipschitz} if it also satisfies
$$\dist_Y(\q(x),\q(x'))\ge \frac{1}{L}\dist_X(x,x')\, .$$
\end{definition}

\begin{remark} \label{grr0} Let $(X_n)$ and $(Y_n)$ be sequences of metric
spaces, $e_n\in X_n$, $e_n'\in Y_n$ ($n\in I$). Then it is easy to
see that any sequence $\q_n\colon X_n\to Y_n$ of
$(L_n,C_n)$-quasi-isometries with $\q_n(e_n)=e_n'$, $n\in I$,
induces an $(L,C)$-quasi-isometry $\q\colon
\lio{X_n}_e\to\lio{Y_n}_{e'}$ where $e=(e_n)^\omega$,
$e'=(e_n')^\omega$, and $L=\lm_\omega L_n,\, C=\lm_\omega C_n$
provided $L<\infty,\, C<\infty$. Moreover, the $\omega$-limit of
the images $\q_n(X_n)$ coincides with the image of $\q$.
\end{remark}

\begin{remark}\label{grr05} Let $\q_n\colon [0,\ell_n]\to X$ be a
sequence of $(L,C)$-quasi-geodesics in a geodesic metric space
$(X,\dist)$. Then the $\omega$-limit $\lio{\q_n([0,\ell_n])}$ in
any asymptotic cone $\co{X,e,d}$ is either empty, or a
bi-Lipschitz arc or a bi-Lipschitz ray or a bi-Lipschitz line.
This immediately follows from Remark \ref{grr0}.
\end{remark}

\begin{remark}\label{grr2} Any quasi-isometric embedding $\q$ of $(X,\dist_X)$
into $(Y,\dist_Y)$ induces a bi-Lipschitz embedding of $\co{X;
e,d}$ into $\co{Y; (q(e_n)),d}$ for every $\omega$, $e$ and $d$
\cite{Gr2}.
\end{remark}

\me

Every finitely generated group $G=\la X\ra$ can be considered a
metric space where the distance between two elements $a, b$ is the
length of the shortest group word in $X$ representing $a\iv b$. The asymptotic cones of $G$ corresponding to different observation points are isometric
\cite{Gr2}. Thus when we consider an asymptotic cone of a finitely
generated group, we shall always assume that the observation point
$e$ is $(1)^\omega$.

Let $G_n$, $n\in I$, be the metric space $G$ with metric
$\frac{\dist}{d_n}$ for some sequence of scaling constants
$(d_n)_{n\in I}$. The set
$\Pi_e G_n/\omega$ denoted by $G$ is a subgroup of the
ultrapower $G^\omega$.

\begin{remark}\label{grr3}
Notice \cite{Gr2} that the group $G^\omega_e$ acts on $\co{G;e,d}$
by isometries: $$(g_n)^\omega \lio{x_n}=\lio{g_nx_n}.$$ This
action is transitive, so, in particular, every asymptotic cone of
a group is homogeneous.

More generally if a group $G$ acts by isometries on a metric space
$(X,\dist)$ and there exists a bounded subset $B \subset X$ such
that $X=GB$ then all asymptotic cones of $X$ are homogeneous
metric spaces.
\end{remark}

\begin{definition}[asymptotic properties]\label{apr}
We say that a space \textit{has a certain property asymptotically}
if each of its asymptotic cones has this property. For example, a
space may be asymptotically CAT(0), asymptotically without
cut-point etc.
\end{definition}

\begin{definition}[asymptotically tree-graded spaces]\label{asco}
Let $(X,\dist)$ be a metric space and let $\mathcal{A}=\{ A_i \mid
i\in I \}$ be a collection of subsets of $X$. In every asymptotic
cone $\co{X;e,d}$, we consider the collection of subsets
$$
\mathcal{A}_\omega = \left\{ {\lio{ A_{i_n}}}\mid (i_n)^\omega\in
I^\omega \hbox{ such that the sequence } \left(
\frac{\dist(e_{n},A_{i_n})}{d_{n}}\right) \hbox{ is
bounded}\right\}\, .
$$

We say that $X$ is \textit{asymptotically tree-graded with respect
to} $\aaa$ if every asymptotic cone $\co{X; e,d}$ is tree-graded
with respect to $\aaa_\omega $.
\end{definition}

Corollary \ref{cor451} will show that there is no need to vary the
ultrafilter in Definition \ref{asco}: if a space is tree-graded
with respect to a collection of subsets for one ultrafilter, it is
tree-graded for any other with respect to the same collection of
subsets.

\subsection{Ultralimits of asymptotic cones are asymptotic cones}

\begin{definition}[an ultraproduct of ultrafilters]\label{defprod}
Let $\omega$ be an ultrafilter over $I$ and let $\mu=(\mu_n)_{n\in
I}$ be a sequence of ultrafilters over $I$. We consider each
$\mu_n$ as a measure on the set $\{n\}\times I$ and $\omega$ as a
measure on $I$.

For every subset $A\subseteq I\times I$ we set $\omega\mu(A)$
equal to the $\omega$-measure of the set of all $n\in I$ such that
$\mu_n(A\cap (\{n\}\times I))=1$.

In other words $$\omega\mu(A)=\int \mu_n\left(A\cap(\{n\}\times I)\right)\:
d\omega(n).$$
\end{definition}

Notice that this is a generalization of the standard notion of
product of ultrafilters (see \cite[Definition 3.2 in Chapter
VI]{She}).

\begin{lemma} (cf \cite[Lemma 3.6 in Chapter VI]{She})
$\omega\mu$ is an ultrafilter over $I\times I$.
\end{lemma}

\proof It suffices to prove that $\omega\mu$ is finitely additive
and that it takes the zero value on finite sets.

Let $A$ and $B$ be two disjoint subsets of $I \times I$. Then for
every $n\in I$ the sets $A\cap (\{n\}\times I)$ and
$B\cap(\{n\}\times I)$ are disjoint. Hence (by the additivity of
$\mu_n$) for every $n\in I$
$$\mu_n((A\cup B)\cap (\{n\}\times I))=\mu_n(A\cap
(\{n\}\times I))+\mu_n(B\cap (\{n\}\times I)).$$ Therefore (by the
additivity of $\omega$) $$\omega\mu(A\sqcup
B)=\omega\mu(A)+\omega\mu(B).$$

Let now $A$ be a finite subset of $I\times I$. Then the set of
numbers $n$ for which $\mu_n(A\cap (\{n\}\times I))=1$ is empty.
So $\omega\mu(A)=0$ by definition.\endproof

\begin{lemma}[double ultralimit of sequences]\label{doublels}
Let $\omega, \mu_n$, $n\in I$, be as in Definition \ref{defprod}.
Let $r^{(n)}_k$ be an uniformly bounded double indexed sequence of
real numbers, $k, n\in I$. Then
\begin{equation}\label{ulseq}
\lm_{\omega\mu} r^{(n)}_k =\lm_\omega \lm_{\mu_n} r^{(n)}_k
\end{equation}
(the internal limit is taken with respect to $k$).
\end{lemma}

\proof Let $r=\lm_{\omega\mu} r^{(n)}_k $. It follows that, for
every $\varepsilon >0$,
$$
\omega\mu\left\{ (n,k) \mid r^{(n)}_k \in (r-\varepsilon ,
r+\varepsilon ) \right\} =1 \Leftrightarrow
$$
$$
\omega \left\{n\in I \mid \mu_n  \left\{ k \mid r^{(n)}_k \in
(r-\varepsilon , r+\varepsilon ) \right\} =1 \right\} =1 \, .
$$

It follows that
$$
\omega \left\{ n\in I \mid \lm_{\mu_n} r^{(n)}_k \in
[r-\varepsilon , r+\varepsilon ] \right\} =1\, ,
$$ which implies that
$$
\lm_\omega \lm_{\mu_n} r^{(n)}_k\in [r-\varepsilon , r+\varepsilon
] \, .
$$

Since this is true for every $\varepsilon >0$ we conclude that
$\lm_\omega \lm_{\mu_n} r^{(n)}_k=r$.\endproof

Lemma \ref{doublels} immediately implies:

\begin{proposition}[double ultralimit of metric spaces]\label{doublel}
Let $\omega$ and $\mu$ be as in Definition \ref{defprod}. Let
$\left(X^{(n)}_k, \dist^{(n)}_k\right)$ be a double indexed
sequence of metric spaces, $k,n\in I$, and let $e$ be a double
indexed sequence of points $e^{(n)}_k \in X^{(n)}_k$. We denote by
$e^{(n)}$ the sequence $\left( e^{(n)}_k \right)_{k\in I}$.

The map

\begin{equation}\label{isodl}
{\lm}_{\omega\mu} \left( x^{(n)}_k \right) \mapsto {\lm}_\omega
\left( {\lm}_{\mu_n}\left( x^{(n)}_k\right)\right)\, ,
\end{equation}
 is an isometry from
${\lm}^{\omega\mu} \left( X^{(n)}_k \right)_e$ onto
${\lm}^\omega\left({\lm}^{\mu_n}\left(X^{(n)}_k\right)_{e^{(n)}}\right)_{e'},
$ where $e_n'={\lm}^{\mu_n} (e^{(n)})$

\end{proposition}

\begin{cor}[ultralimits of cones are cones]\label{ulc}
Let $X$ be a metric space. Let $\omega$ and $\mu$ be as above. For
every $n\in I$ let $e^{(n)}=\left(e^{(n)}_k\right)_{k\in I}$ be an
observation point, $d^{(n)}=\left(d^{(n)}_k\right)_{k\in I}$ be a
sequence of scaling constants satisfying $\lm_{\mu_n} d^{(n)}_k
=\infty$ for every $n\in I$. Let
$\mathrm{Con}^{\mu_n}\left(X;e^{(n)},d^{(n)}\right)$ be the
corresponding asymptotic cone of $X$. Then the map
\begin{equation}\label{isoc}
\lm_{\omega\mu} \left( x^{(n)}_k \right) \mapsto \lm_\omega
\left(\lm_{\mu_n }\left( x^{(n)}_k\right) \right),
\end{equation}
 is an isometry from $\mathrm{Con}^{\omega\mu}(X;e,d)$ onto $$\lm^\omega \left(
\mathrm{Con}^{\mu_n}\left(X;e^{(n)},d^{(n)}\right)
\right)_{(\lm^{\mu_n } (e^{(n)}))},$$ where $e=\left(
e^{(n)}_k\right)_{(n,k)\in I \times I}$ and $d=\left(
d^{(n)}_k\right)_{(n,k)\in I \times I}$.
\end{cor}

\proof Let us prove that $\lm_{\omega\mu}d^{(n)}_k =\infty$. Let
$M>0$. For every $n\in I$ we have that $\lm_{\mu_n } d^{(n)}_k
=\infty $, whence $\mu_n \left\{ k\in I \mid d^{(n)}_k >M
\right\}=1$. It follows that $\left\{ n\in I \mid \mu_n \left\{
k\in I \mid d^{(n)}_k >M \right\}=1 \right\}=I$, therefore its
$\omega$-measure is $1$. We conclude that $\omega\mu\left\{
(n,k)\mid d^{(n)}_k >M  \right\}=1$.

It remains to apply Proposition \ref{doublel} to the  sequence of
metric spaces $\left( X, \frac{1}{d^{(n)}_k}\dist \right)$ and
to~$e$.\endproof

\subsection{Another definition of asymptotic cones}\label{adef}

In \cite{Gr2}, \cite{VDW} and some other papers, a more
restrictive definition of asymptotic cones is used. In that
definition, the set $I$ is equal to $\N$ and the scaling constant
$d_n$ must be equal to $n$ for every $n$. We shall call these
asymptotic cones {\em restrictive}.

It is easy to see that every restrictive asymptotic cone is an
asymptotic cone in our sense. The converse statement can well be
false although we do not have any explicit examples.

Also for every ultrafilter $\omega$ over $I$ and every sequence of
scaling constants $d=(d_n)_{n\in I}$, there exists an ultrafilter
$\mu$ over $\N$ such that the asymptotic cone $\Con^\omega(X;e,d)$
contains an isometric copy of the restrictive asymptotic cone
$\Con^\mu(X;e,(n))$. Indeed, let $\phi$ be a map $I\to \N$ such
that $\phi(i)=[d_i]$. Now define the ultrafilter $\mu$ on $\N$ by
$\mu(A)=\omega(\phi\iv(A))$ for every set $A\subseteq \N$. The
embedding $\Con^\mu(X;e,(n))\to \Con^\omega(X;e,d)$ is defined by
$\lim^\mu (x_n)\mapsto \lim^\omega \left( x_{\phi(i)}
\right)_{i\in I}$.

\begin{remark} \label{rem676} In the particular case when the sets $\{i\in I\mid [d_i]=k\}$
are of uniformly bounded (finite) size, this embedding is a
surjective isometry \cite{Ri}.
\end{remark}

The restrictive definition of asymptotic cones is, in our opinion,
less natural because the $\omega$-limit of restrictive asymptotic
cones is not canonically represented as a restrictive asymptotic
cone (see Corollary \ref{ulc}). Conceivably, it may even not be a
restrictive asymptotic cone in general. The next statement shows
that it is a restrictive asymptotic cone in some particular cases.

\begin{proposition} \label{prop676} Let $\nu_n$, $n\in \N$ be a sequence of
ultrafilters over $\N$. Let $(I_n)$ be sequence of pairwise
disjoint subsets of $\N$ such that $\nu_n(I_n)=1$. Let
$C_n=\Con^{\nu_n}(X;e^{(n)},(n))$, $n\in \N$, be a restrictive
asymptotic cone of a metric space $X$. Then the $\omega$-limit of
asymptotic cones $C_n$ is a restrictive asymptotic cone.
\end{proposition}

\proof Let $\mu_n$ be the restriction of $\nu_n$ onto $I_n$, $n\in
\N$. Then $C_n$ is isometric to $\Con^{\mu_n}(X;e^{(n)},d^{(n)})$
where $d^{(n)}$ is the sequence of all numbers from $I_n$ in the
increasing order. By Corollary \ref{ulc},
$\lio{C_n}_{\lm^{\mu_n}(e^{(n)})}$ is the asymptotic cone
$\Con^{\omega\mu}(X;e,d)$ where $e=\left(
e^{(n)}_k\right)_{(n,k)\in \N \times \N}$ and $d=\left(
d^{(n)}_k\right)_{(n,k)\in \N \times \N}$. For every natural
number $a$ the set of pairs $(n,k)$ such that $d^{(n)}_k=a$
contains at most one element because the subsets $I_n\subseteq \N$
are disjoint. It remains to apply Remark \ref{rem676}.
\endproof

\subsection{Simple triangles in ultralimits of metric spaces}

\begin{definition}[$k$-gons] We say that a metric space $P$ is a
geodesic (quasi-geodesic) $k$-gon if it is a union of $k$
geodesics (quasi-geodesics) $\q_1,...,\q_k$ such that
$(\q_i)_+=(\q_{i+1})_-$ for every $i=1,...,k$ (here $k+1$ is
identified with $1$).

For every $i=1,...,k$, we denote the polygonal curve $P\setminus
\left( \q_{i-1}\cup \q_i\right)$ by $\oo_{x_i}(P)$, where $x_i=(\q_{i-1})_+ =(\q_i)_-$. When there is no possibility of confusion we simply
denote it by $\oo_{x_i}$.
\end{definition}

\begin{lemma}\label{Dd}
(1) Let $P_n$, $n\in \N,$ be a sequence of geodesic $k$-gons in
metric spaces $(X_n, \dist_n)$. Let $\omega$ be an ultrafilter
over $\N$, such that $\lio{P_n}=P$, where $P$ is a simple geodesic
$k$-gon in the metric space $\lio{X_n}_e$ with metric $\dist$. Let
${\mathcal{V}}_n$ be the set of vertices of $P_n$ in the clockwise
order. Let $D_n$ be the supremum over all points $x$ contained in
two distinct edges of $P_n$ of the distances $\dist\left(x,
{\mathcal{V}}_n\right).$ Then $\lm_\omega D_n=0$.

(2) Let $P$ be a simple $k$-gon in $(X,\dist)$. For every
$\delta>0$ we define $D_\delta = D_\delta(P)$ to be the
 supremum over all $k$-gons $P_\delta$ in $X$
that are at Hausdorff distance at most $\delta$ from $P$ and over
all points $x$ contained in two distinct edges of $P_\delta$ of
the distances $\dist\left(x, {\mathcal{V}}_\delta\right)$, where
${\mathcal{V}}_\delta$ is the set of vertices of $P_\delta$. Then
$\lm_{\delta \to 0} D_\delta =0$.
\end{lemma}

\proof (1) Since the $\omega$-limit of the diameters of $P_n$ is
the diameter of $P$, it follows that the diameters of $P_n$ are
uniformly bounded $\omega$-almost surely. In particular $D_n$ is
uniformly bounded $\omega$-almost surely, therefore its
$\omega$-limit exists and it is finite. Suppose that $\lm_\omega
D_n=2D>0$. Then $\omega$-almost surely there exists $x_n$
contained in two distinct edges of $P_n$ such that
$\dist_n\left(x_n, {\mathcal{V}}_n\right)> D.$ Without loss of
generality we may suppose that $x_n\in [A_n,B_n]\cap [B_n,C_n]$
for every $n$, where $[A_n,B_n], [B_n,C_n]$ are two consecutive
edges of $P_n$ such that $\lio{[A_n,B_n]}=[A,B]\, ,\,
\lio{[B_n,C_n]}=[B,C]$, where $[A,B], [B,C]$ are two consecutive
edges of $P$. Then $\lio{x_n}\in [A,B]\cap [B,C]$, which by
simplicity of $P$ implies that $\lio{x_n}=B$. On the other hand we
have that $\dist_n \left(x_n, {\mathcal{V}}_n\right)> D$, which
implies that $\dist (\lio{x_n},B)\geq D$. We have obtained a
contradiction.

\medskip

(2) Assume that $\lm_{\delta \to 0} D_\delta =2D>0$. It follows
that there exists a sequence $(P_n)$ of $k$-gons endowed with
metrics such that their Hausdorff distance to $P$ tends to zero
and such that there exists $x_n$ contained in two distinct edges
of $P_n$ and at distance at least $D$ of the vertices of $P_n$.
According to \cite{KaL1}, it follows that $\lio{P_n}=P$ for every
ultrafilter $\omega$. On the other hand $D_n >D$ for all $n\in
\N$. We thus obtain a contradiction of (1).\endproof

\me

\begin{proposition}[limits of simple polygons]\label{approxul}

Consider an ultrafilter $\omega$ over $\N$ and a sequence of
metric spaces, $(X_n, \dist_n)$, $n\in \N$. Let $e\in \Pi
X_n/\omega$ be an observation point. For every simple geodesic
triangle $\Delta$ in $\lio{X_n}_e$, for every sufficiently small
$\varepsilon >0$ there exists $k_0=k_0(\varepsilon )$ and a simple
geodesic triangle $\Delta_\varepsilon$ with the properties:
\begin{itemize}
  \item[(a)] The Hausdorff distance between $\Delta$ and
  $\Delta_\varepsilon$ does not exceed $\varepsilon $;
  \item[(b)] $\Delta_\varepsilon $ contains the midpoints
  of the edges of $\Delta $;
  \item[(c)] The triangle
$\Delta_\varepsilon$ can be written as $\lio{P^{\varepsilon}_n}$,
where each $P^{\varepsilon}_n$ is a geodesic $k$-gon in $X_n$,
$k\leq k_0$, $P^{\varepsilon}_n$ is simple and the lengths of all
edges of $P^{\varepsilon}_n$ are $O(1)$ $\omega$-almost surely.
\end{itemize}
\end{proposition}

\proof Let $A,B,C$ be the vertices of $\Delta$, in the clockwise
order, and let $M_{AB},M_{BC}$ and $M_{AC}$ be the midpoints of
$[A,B],[B,C]$ and $[A,C]$, respectively.

We construct $\Delta_\varepsilon$ in several steps.

\medskip

\noindent {\bf Step I. Constructing not necessarily simple
geodesic triangles $\Delta_\varepsilon$.}

\medskip

 For every small
$\delta
>0$ we divide each of the halves of edges of $\Delta $ determined by a vertex and a midpoint
into finitely many segments of length at least $\delta$ and at
most $2\delta $. Let $\vv$ be the set of endpoints of all these
segments, endowed with the natural cyclic order. We call $\vv$ a
$\delta$-\textit{partition of} $\Delta$. We assume that $\{
A,B,C,M_{AB},M_{BC},M_{AC}\}\subset \vv$. Every $t\in \vv$ can be
written as $t=\lio{t_n}$, hence $\vv=\lio{\vv_n}$, where each
$\vv_n$ is endowed with a cyclic order. Let $P_n$ be a geodesic
$k$-gon with vertices $\vv_n$, where $k=|\vv|$. The limit set
$\Delta_\delta =\lio{P_n}$ is a geodesic triangle with vertices
$A,B,C$ and at Hausdorff distance at most $\delta $ from $\Delta$.

\me

\Notat \quad Let $E,F$ be two points on an edge of
$\Delta_\delta$. We denote the part of the geodesic side of
$\Delta_\delta$ between $E$ and $F$ in $\Delta_\delta$ by
$[E,F]_\delta$. If $E,F$ are two points on an edge of $\Delta$, we
denote the part of the side of $\Delta$ between $E$ and $F$  by
$[E,F]$. This is to avoid confusion between different geodesics
joining two such points.

\medskip

\noindent \textbf{Step II. Making $\Delta_\varepsilon $ simple.}

\medskip

For every $\delta>0$, we consider $D_\delta = D_\delta (\Delta )$
given by Lemma \ref{Dd}. Let
$$\alpha (\Delta)\inf\left\{\dist\left(x, \oo_x(\Delta)\right) \mid x\in \{ A,B,C\}
\right\}\, .$$ By Lemma \ref{Dd} we have $\lm_{\delta \to
0}D_\delta =0$. Therefore, for $\delta$ small enough we have
\begin{equation}\label{delta}
2D_\delta +4\delta < \alpha (\Delta )\mbox{ and }D_\delta +2\delta
\leq \frac{1}{10} \min \left\{ \dist(A,B),\dist(B,C),\dist(C,A)
\right\}\, .
\end{equation}

Fix a $\delta$ satisfying (\ref{delta}), a $\delta$-partition
$\vv$ of $\Delta$, and a corresponding triangle $\Delta_\delta
=\lio{P_n}$.

Let $A_1$ and $A_2$ be the nearest to $A$ points of $\vv\setminus
\nn_{D_\delta+\delta}(A)$ on the edges $[A,B]$ and $[A,C]$,
respectively. For an appropriate choice of $\Delta_\delta$, we may
suppose that $\dist(A,A_1)=\dist(A,A_2)$. We note that
$\dist(A,A_1) \in [D_\delta +\delta , D_\delta +2 \delta]$.
Similarly we take $B_1 \in [B,C]\cap \vv \, ,\, B_2\in [B,A]\cap
\vv$ and $C_1\in [C,A]\cap \vv \, ,\,  C_2\in [C,B]\cap \vv $ with
$\dist(B,B_1)=\dist(B,B_2)\in [D_\delta +\delta , D_\delta +2
\delta]$ and $\dist(C,C_1)=\dist(C,C_2)\in [D_\delta +\delta ,
D_\delta +2 \delta]$.

Suppose that $[A_1, B_2]_\delta$ and $[B_1,C_2]_\delta $ have a
point $E$ in common. The definition of $D_\delta$ implies that
$E\in \nn_{D_\delta} (\{ A,B,C\} )$. On the other hand $E\in [A_1,
B_2]_\delta$ implies $E\not\in \nn_{D_\delta} (\{ A,B\} )$ and
$E\in [B_1, C_2]_\delta$ implies $E\not\in \nn_{D_\delta} (\{
B,C\} )$, a contradiction.

We conclude, by repeating the previous argument, that the segments
$[A_1, B_2]_\delta$, $[B_1,C_2]_\delta$ and $[C_1, A_2]_\delta$
are pairwise disjoint. Since $\dist (A,A_1) \, ,\, \dist(B,B_2)
\leq D_\delta+2\delta \leq \frac{1}{10}\dist(A,B)$, it follows
that $M_{AB}$ is contained in $[A_1, B_2]_\delta $. Likewise,
$M_{BC}$ and $M_{AC}$ are contained in $[B_1,C_2]_\delta$ and
$[C_1, A_2]_\delta$, respectively.

Let $d_A$ be the supremum of $\dist(E,A)$ for all $E$ satisfying
two conditions: $E\in [A_1,A]_\delta$ and
$\dist(A_2,E)+\dist(E,A)=\dist(A_2,A)$. Since these two conditions
define a closed set, it follows that there exists $A'\in
[A_1,A]_\delta $ such that
$\dist(A_2,A')+\dist(A',A)=\dist(A_2,A)$ and $\dist (A,A')=d_A$.
Obviously $A'\not\in \{ A_1, A_2\}$. In other words, $A'$ is the
farthest from $A$ point in $[A_1,A]_\delta$ which is contained in
a geodesic joining $A_2$ and $A$. Hence $A'$ has the property that
every geodesic joining it with $A_2$ intersects $[A_1,A']_\delta$
only in $A'$. Similarly we find points $B'\in [B_1,B]_\delta$ and
$C'\in [C_1,C]_\delta $.

Recall that $\Delta_\delta =\lio{P_n}$. Let $P_n^A$ be a sequence
of polygonal lines in $P_n$ with endpoints $A_n',B^2_n$, having as
limit $[A',B_2]_\delta$. Likewise let $P_n^B$ and $P_n^C$ be
sequences of polygonal lines in $P_n$, with endpoints $B_n',C^2_n$
and $C_n',A^2_n$, having as limits $[B',C_2]_\delta$ and
$[C',A_2]_\delta$, respectively. We consider the new sequence of
polygons $P_n'=P_n^A \cup [B^2_n,B_n']\cup P_n^B \cup
[C^2_n,C_n']\cup P_n^C \cup [A^2_n,A_n']$. The limit set
$\lio{P_n'}$ is $[A',B_2]_\delta \cup \g_{B_2B'}\cup
[B',C_2]_\delta \cup \g_{C_2C'}\cup [C',A_2]_\delta \cup
\g_{A_2A'}$ where $\g_{B_2B'}=\lio{[B^2_n,B_n']}$ is a geodesic
and likewise for $\g_{C_2C'}, \g_{A_2A'}$.

We have $\dist(C',A)=\dist(C', A_2)+\dist(A_2,A)=\dist(C',
A_2)+\dist(A_2,A')+\dist(A',A)$. It follows that by joining the
pairs of points $(C', A_2)$, $(A_2,A')$ and $(A',A)$ by geodesics
we obtain a geodesic from $C'$ to $A$. In particular
$[C',A_2]_\delta \cup \g_{A_2A'}$ is a geodesic. Likewise,
$[A',B_2]_\delta \cup \g_{B_2B'}$ and $[B',C_2]_\delta \cup
\g_{C_2C'}$ are geodesics. Therefore $\lio{P_n'}$ is a geodesic
triangle $\Delta_\delta'$ with vertices $A',B',C'$. By
construction the Hausdorff distance between $\Delta_\delta'$ and
$\Delta_\delta$ is at most $D_\delta +2\delta $, hence the
Hausdorff distance between $\Delta_\delta'$ and $\Delta$ is at
most $D_\delta +3\delta $.

Suppose that two edges of $\Delta_\delta'$ have a common point
$E$. Suppose the two edges are $[A',B_2]_\delta \cup \g_{B_2B'}$
and $[B',C_2]_\delta \cup \g_{C_2C'}$. If $E\in [A',A_1]_\delta $
then $\dist (A,E)\leq D_\delta +2\delta $. On the other hand $E\in
[B',C_2]_\delta \cup \g_{C_2C'}$ implies $E\in \nn_{D_\delta
+2\delta }([B,C])$. Hence $\dist (A,[B,C])\leq 2D_\delta +4\delta
< \alpha (\Delta )$, a contradiction.

If $E\in \g_{C_2C'}$ then $\dist (C,E)\leq D_\delta +2\delta $
which together with $E\in [A',B_2]_\delta \cup \g_{B_2B'} \subset
\nn_{D_\delta +2\delta } ([A,B])$ implies $\dist (C,[A,B])\leq
2D_\delta +4\delta <\alpha (\Delta )$, a contradiction.

If $E\in [A_1,B_2]_\delta $ then $E\not \in [B_1 ,C_2]_\delta$.
Also since $\dist (B,E)\geq \dist (B,B_2)=\dist (B,B_1)$ it
follows that $E\not \in [B' ,B_1]_\delta$, a contradiction.

If $E\in \g_{B_2B'}$ then an argument similar to the previous
gives $E\not \in [B_1,C_2]_\delta $. We conclude that $E\in [B'
,B_1]_\delta$. By the choice of $B'$ we have $E=B'$.

We conclude that $\Delta_\delta'$ is a simple geodesic triangle,
containing the midpoints of the edges of $\Delta$, at Hausdorff
distance at most $D_\delta +3\delta$ from $\Delta $, and
$\Delta_\delta'=\lm^\omega (P_n')$, where $P_n'$ is a geodesic
$m$-gon, with $m\leq k+3$.

\medskip

\noindent \textbf{Step III. Making polygons simple.}

\medskip

Let $D_n$ be the supremum over all points $x$ contained in two
distinct edges of $P_n'$ of the distances from $x$ to the vertices
of $P_n'$. Applying Lemma \ref{Dd}, (1), to $(P_n')$ and to
$\Delta_\delta'=\lm^\omega (P_n')$ we obtain that $D_n$ tends to
zero as $n\to \infty $. Let $v_n$ be a vertex of $P_n'$. We
consider the farthest point $v_n'$ in the ball $B(v_n, 2 D_n)$
contained in both edges of endpoint the vertex $v_n$. Cut the
bigon of vertices $v_n,v_n'$ from the polygon, and repeat this
operation for every vertex $v_n$ of $P_n'$. As a result, we obtain
a new polygon $P_n''$ which is simple and at Hausdorff distance at
most $2 D_n$ from $P_n'$. It follows that $\lm^\omega
(P_n'')=\lm^\omega (P_n')=\Delta_\delta'$.\endproof

\begin{theorem}[being tree-graded is closed under ultralimits]\label{ultg}
For every $n\in \N$ let $\free_n$ be a complete geodesic metric
space which is tree-graded with respect to a collection $\pp_n$ of
closed geodesic subsets of $\free_n$. Let $\omega$ be an
ultrafilter over $\N$ and let $e\in \Pi \free_n/\omega$ be an
observation point. The ultralimit $\lio{\free_n}_e$ is tree-graded
with respect to the collection of limit sets
$$
\pp_\omega = \left\{ \lio{M_n}\mid M_n \in \pp_n\, ,\,  \dist
(e_n,M_n) \mbox{ bounded uniformly in }n\right\}\, .
$$
\end{theorem}

\proof {\bf Property $(T_1)$.} Let $\lio{M_n}\, ,\, \lio{M_n'}\in
\pp_\omega$ be such that there exist two distinct points $x_\omega
, y_\omega $ in $\lio{M_n}\cap \lio{M_n'}$. It follows that
$x_\omega = \lio{x_n}=\lio{x_n'}$ and $y_\omega \lio{y_n}=\lio{y_n'}$, where $x_n,y_n\in M_n,\, x_n',y_n'\in
M_n'$, $\dist(x_n,x_n')=o(1),\, \dist(y_n,y_n')=o(1)$, while
$\dist (x_n,y_n)=O(1),\, \dist (x_n',y_n')=O(1)$.

By contradiction suppose that $M_n \neq M_n'\; \omega$-almost
surely. Then property $(T_2)$ of the space $\free_n$ and Corollary
\ref{projA} imply that $M_n$ projects into $M_n'$ in a unique
point $z_n$ and that $z_n\in [x_n,x_n'] \cap [y_n,y_n']$. It
follows that $\dist (x_n,z_n)=o(1)$ and $\dist (y_n,z_n)=o(1)$,
therefore that $\dist (x_n,y_n)=o(1)$. This contradiction implies
that $M_n = M_n'\; \omega$-almost surely, so
$\lio{M_n}=\lio{M_n'}$.

\medskip

{\bf Property $(T_2)$.} Let $\Delta $ be a simple geodesic
triangle in $\lio{\free_n}_e$. Consider an arbitrary sufficiently
small $\varepsilon>0$ and apply Proposition \ref{approxul}. We
obtain a simple geodesic triangle $\Delta_\varepsilon$ satisfying
properties (a), (b), (c) in the conclusion of the Proposition. In
particular $\Delta_\varepsilon = \lio{P_n^\varepsilon}$, where
$P_n^\varepsilon $ is a simple geodesic polygon in $\free_n$.
Property $(T_2'')$ applied to $\free_n$ implies that
$P_n^\varepsilon$ is contained in one piece $M_n$. Consequently
$\Delta_\varepsilon \subset \lio{M_n}$. Property (b) of
$\Delta_\varepsilon$ implies that $\lio{M_n}$ contains the three
distinct middle points of the edges of $\Delta $. This and
property $(T_1)$ already proven imply that all triangles
$\Delta_\varepsilon$ are contained in the same $\lio{M_n}$.
Property (a) and the fact that $\lio{M_n}$ is closed imply that
$\Delta \subset \lio{M_n}$.\endproof

\begin{definition}
Let $P$ be a polygon with quasi-geodesic edges and with set of
vertices $\mathcal{V}$. Points in $P\setminus \mathcal{V}$ are
called \textit{interior points of }$P$. Let $p\in P$. The
\textit{inscribed radius in }$p$ with respect to $P$ is either the
distance from $p$ to the set $\mathcal{O}_p$, if $p$ is a vertex,
or the distance from $p$ to the set $P\setminus \q$ if $p$ is
contained in the interior of the edge $\q$.
\end{definition}

\medskip

\unitlength .3mm 
\linethickness{0.4pt}
\ifx\plotpoint\undefined\newsavebox{\plotpoint}\fi 


\begin{figure}[!ht]
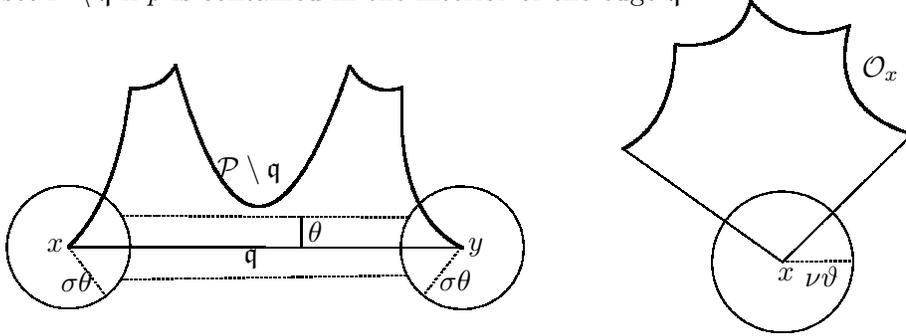

\centering
\caption{Properties ($F_1$) and ($F_2$).}
\label{fig2}
\end{figure}

\begin{definition}[fat polygons]
Let $\vartheta>0$, $\sigma \geq 1$ and $\nu \geq 4\sigma$. We call
a $k$-gon $P$ with quasi-geodesic edges \textit{$(\vartheta ,
\sigma, \nu )$-fat} if the following properties hold:
\begin{enumerate}
  \item[$(F_1)$](\textbf{large comparison
angles, large inscribed radii in interior points}) for every edge
$\q $ with endpoints $\{ x,y \}$ we have $$\dist(\q \setminus
\nn_{\sigma \vartheta}(\{x,y\}), P \setminus \q)\geq \vartheta;$$

  \item[$(F_2)$]
(\textbf{large edges, large inscribed radii in vertices}) for
every vertex $x$ we have $$\dist(x, \oo_x)\geq \nu \vartheta.$$
\end{enumerate}
\end{definition}

\begin{remarks}
1. For almost all applications, we can assume that $\sigma$ in that
definition is equal to 2, so the ``fatness" really depends on two
parameters, $\vartheta$ ad $\nu$. We need $\sigma$ to make fatness
preserved under quasi-isometry (see Theorem \ref{qi}).

2. Property ($F_1$) implies that in each of the vertices $x,y$
certain comparison angles are at least $\frac{1}{\sigma }$ and that
in the interior points of $P$ outside $\nn_{\sigma \vartheta }
(\calv )$ the inscribed radius is at least $\vartheta$.

3. Property ($F_2$) ensures that for every edge $\q$ the set $\q
\setminus \nn_{\sigma \vartheta}(\{x,y\})$ has diameter at least
$2\sigma\vartheta $, in particular it is never empty. It also
ensures that the inscribed radius in every vertex is at least $\nu
\vartheta$.

\end{remarks}

\begin{proposition}[triangles in an asymptotic cone are $\omega$-limits of fat polygons]\label{approx}
For every simple geodesic triangle $\Delta$ in $\co{X; e,d}$, for
every sufficiently small $\varepsilon >0$ there exists
$k_0=k_0(\varepsilon )$ and a simple geodesic triangle
$\Delta_\varepsilon$ with the properties:
\begin{itemize}
  \item[(a)] The Hausdorff distance between $\Delta$ and
  $\Delta_\varepsilon$ does not exceed $\varepsilon$;
  \item[(b)] $\Delta_\varepsilon $ contains the midpoints
  of the edges of $\Delta $;
  \item[(c)] For every $\vartheta>0$ and $\nu \geq 8$, the triangle
$\Delta_\varepsilon$ can be written as $\lio{P^{\varepsilon}_n}$,
where each $P^{\varepsilon}_n$ is a geodesic $k$-gon in $X$,
$k\leq k_0$, and $P^{\varepsilon}_n$ is $(\vartheta , 2, \nu
)$-fat $\omega$-almost surely.
\end{itemize}
\end{proposition}

\proof Proposition \ref{approxul} applied to $\left( X\, ,\,
\frac{1}{d_n}\dist \right)$, $\omega ,\, e$ and $\Delta $ implies
that for every $\varepsilon >0 $ there exists
$k_0=k_0(\varepsilon)$ and $\Delta_\varepsilon$ satisfying (a) and
(b) and such that $\Delta_\varepsilon =\lio{P_n}$, where $P_n$ are
simple geodesic $k$-gons in $X$, $3\leq k\leq k_0$, such that the
lengths of all edges in $P_n$ are $O(d_n)$ $\omega$-almost surely.
Remark \ref{udisj} implies that there exists $m\in \{ 3,\dots ,
k_0 \}$ such that $P_n$ have $m$ edges $\omega$-almost surely. Let
$\vartheta>0$ and $\nu \geq 8$. We modify the sequence of polygons
$(P_n)$ so that their limit set stays the same while the polygons
become $(\vartheta , 2, \nu )$-fat.

\medskip

Let $\vv_n=\{ v_1^n,v_2^n,\dots ,v_{m}^n \}$ be the set of
vertices of $P_n$ in the clockwise order. We denote the limit set
$\lio{\vv_n}$ by $\vv$, and we endow it with the clockwise order
on $\Delta_\varepsilon$. There exists $\varrho
>0$ such that for every $v\in
\vv$, the distance between $v$ and $\oo_v (\Delta_\varepsilon)$ is
at least $2\varrho $, where $\oo_v (\Delta_\varepsilon)$ is taken
in $\Delta_\varepsilon$ considered as a polygon with vertices
$\vv$. It follows that $\omega$-almost surely for every $i\in \{
1,2,\dots , m\}$ we have $\dist\left(v_i^n , \oo_{v_i^n
}(P_n)\right)\geq \varrho d_n$. In particular, $\omega$-almost
surely all the edges of $P_n$ have length at least $\varrho d_n$.

\medskip

\noindent \textit{Convention}: In what follows we use the notation
$[v_{i}^n,v_{i+1}^n]$ for a generic edge of $P_n$, where $i+1$ is
taken modulo $m$.

\medskip

Let $\epsilon_n$ be the supremum of distances $\dist\left( x,
\vv_n \right)$ for all $x \in [v_{i}^n,v_{i+1}^n] \cap \nn_{\nu
\vartheta }\left( [v_{j}^n,v_{j+1}^n] \right), i\neq j,\, i,j \in
\{ 1,2,\dots , m\}$. Suppose that $\lm_\omega
\frac{\epsilon_n}{d_n}=2\eta
>0$. Then there exist $x_n \in [v_{i}^n,v_{i+1}^n] \cap
\nn_{\nu \vartheta }\left( [v_{j}^n,v_{j+1}^n] \right),\, i\neq
j,\, i,j \in \{ 1,2,\dots , m\} $, with $\dist \left( x_n, \vv_n
\right)\geq \eta d_n$ $\omega$-almost surely. Taking the
$\omega$-limit, we get a contradiction with the fact that
$\Delta_\varepsilon$ is simple. Therefore $\lm_\omega
\frac{\epsilon_n}{d_n}=0$.

\medskip

\Notat \quad We denote by ${\mathfrak N}$ the set of all $n\in \N$
such that for every $i\in \{ 1,2,\dots , m\}$ we have $\dist
\left(v_i^n , \oo_{v_i^n }\right)\geq \varrho d_n$ and such that
$\varrho d_n \geq 2\epsilon_n +2 +(2\nu +1) \vartheta $. Obviously
${\mathfrak N}\in\omega$.

\medskip

 Let $[v_{i-1}^n,v_i^n]$ and
$[v_i^n,v_{i+1}^n]$ be two consecutive edges of $P_n$. Let
$\bar{v}_i^n$ be the farthest point of $v_i^n$ in
$[v_{i-1}^n,v_i^n] \cap  \nn_{\epsilon_n+1}(v_i^n)$  contained in
the $\nu \vartheta$-tubular neighborhood of a different edge
$\pgot$ of $P_n$. The edge $\pgot$ has to be at a distance at most
$\epsilon_n +1 +\nu \vartheta $ from $v_i^n$. It follows that for
every $n\in \mathfrak N$ the edge $\pgot$ must be
$[v_i^n,v_{i+1}^n]$. Therefore $\bar{v}_i^n$ is the farthest from
$v_i^n$ point in $[v_{i-1}^n,v_i^n] $ contained in
$\nn_{\nu\vartheta}([v_i^n,v_{i+1}^n])$. Let $\tilde{v}_i^n$ be
the farthest from $v_i^n$ point $t_n\in [v_i^n,v_{i+1}^n]$ such
that $\dist (\bar{v}_i^n \, ,\, t_n)\leq \nu\vartheta$. It follows
that $\dist (\bar{v}_i^n, \tilde{v}_i^n)= \nu\vartheta$. We modify
$P_n$ by replacing $[\bar{v}_i^n,v_i^n]\cup [v_i^n,\tilde{v}_i^n]$
with a geodesic $[\bar{v}_i^n,\tilde{v}_i^n]$. We repeat the
argument for each of the vertices of $P_n$, and in the end we
obtain a sequence of polygons $P'_n$ with at most $2m$ edges each.
As the Hausdorff distance between $P'_n$ and $P_n$ is at most
$\epsilon_n+1+\nu \vartheta $, $\lio{P'_n}=\lio{P_n}$.

 Let us show that for $n\in \mathfrak N$, $P'_n$ is $(\vartheta , 2, \nu
)$-fat.

\medskip

\noindent \textbf{Verification of property $(F_1)$ for $n\in
\mathfrak N$.}

\medskip

There are two types of edges in $P_n'$, the edges of the form
$[\tilde{v}_{i}^n,\bar{v}_{i+1}^n]$, which we shall call
\textit{restricted edges}, and the edges of the form
$[\bar{v}_{i}^n,\tilde{v}_{i}^n]$, which we shall call
\textit{added edges}. We denote by $RE_n$ the union of the
restricted edges of $P_n'$ and by $AE_n$ the union of the added
edges of $P_n'$.

Let $[\tilde{v}_{i}^n,\bar{v}_{i+1}^n]$ be a restricted edge. We
first show that for $n\in {\mathfrak N}$,

      $$
\dist \left(  [\tilde{v}_{i}^n,\bar{v}_{i+1}^n] \setminus
\nn_{2\vartheta}(\{\tilde{v}_{i}^n,\bar{v}_{i+1}^n \}), RE_n
\setminus [\tilde{v}_{i}^n,\bar{v}_{i+1}^n]\right)\geq \vartheta.
      $$

Suppose there exists $y$ in $ [\tilde{v}_{i}^n,\bar{v}_{i+1}^n]
\setminus \nn_{2\vartheta}(\{ \tilde{v}_{i}^n,\bar{v}_{i+1}^n \})$
contained in $\nn_{\vartheta}([\tilde{v}_j^n,\bar{v}_{j+1}^n])$
which is inside $\nn_{\vartheta }([v_j^n,v_{j+1}^n])$, with $j\neq
i$. Then $y\in \nn_{\epsilon_n +1 }(\{ v_i^n , v_{i+1}^n \})$. The
choice of $\bar{v}_{i+1}^n$ implies that $y\in
\nn_{\epsilon_n+1}(v_i^n) $. Therefore $\dist(v_i^n,
[v_j^n,v_{j+1}^n] )\leq \epsilon_n +1 + \vartheta$. The previous
inequality implies that $j=i-1$ for $n\in \mathfrak N$.

Hence there exists $t\in [\tilde{v}_{i-1}^n,\bar{v}_{i}^n]$ such
that $\dist (t,y)< \vartheta $. By the definition of
$\bar{v}_{i}^n$ we have $t=\bar{v}_{i}^n$. This contradicts the
choice of $\tilde{v}_{i}^n$.

       Now let us show that for $n\in \mathfrak N$,
      $$
      \dist \left(  [\tilde{v}_{i}^n,\bar{v}_{i+1}^n] \setminus \nn_{2\vartheta
      }(\{ \tilde{v}_{i}^n,\bar{v}_{i+1}^n \}) ,AE_n \right)\geq  \vartheta.
      $$

      Suppose there exists $z$ in $ [\tilde{v}_{i}^n,\bar{v}_{i+1}^n]
\setminus \nn_{2\vartheta}(\{ \tilde{v}_{i}^n,\bar{v}_{i+1}^n \})$
contained in $\nn_{\vartheta }([\bar{v}_j^n,\tilde{v}_{j}^n])$. It
follows that $z$ belongs to the $(\epsilon_n+\nu
\vartheta+1)$-neighborhood of $v_j^n$ and that $\dist (v_j^n,
[v_{i}^n,v_{i+1}^n]) \leq \epsilon_n+\nu \vartheta+1$. For $n\in
{\mathfrak N}$ this implies that $j\in \{ i,i+1\}$. Suppose $j=i$
(the other case is similar). Let $t\in
[\bar{v}_i^n,\tilde{v}_{i}^n]$ with $\dist (t,z)\leq \vartheta$.
Then $\dist (\tilde{v}_{i}^n,t)\geq \dist(\tilde{v}_{i}^n,z)-\dist
(t,z) \geq 2\vartheta -\vartheta \geq \dist (t,z)$. It follows
that $\dist (\bar{v}_{i}^n,z)\leq \dist (\bar{v}_{i}^n,t)+\dist
(t,z)\leq \dist (\bar{v}_{i}^n,t)+\dist (\tilde{v}_{i}^n,t) \dist (\bar{v}_{i}^n,\tilde{v}_{i}^n)= \nu \vartheta $. This
contradicts the choice of $\tilde{v}_{i}^n$.

\medskip

Now consider an added edge $[\bar{v}_i^n,\tilde{v}_{i}^n]\subset
B(v_i^n, \epsilon_n+1+\nu \vartheta )$. Let $n\in \mathfrak N$. If
there exists $u\in [\bar{v}_i^n,\tilde{v}_{i}^n] \setminus
\nn_{2\vartheta }(\{ \bar{v}_i^n,\tilde{v}_{i}^n \})$ contained in
$\nn_\vartheta ([\bar{v}_j^n,\tilde{v}_{j}^n])$ with $j\neq i$
then $u\in \nn_{\epsilon_n+1+(\nu +1) \vartheta}(v_j^n)$. It
follows that $\dist(v_i^n, v_j^n)\leq \dist(v_i^n, u)+\dist (u,
v_j^n)\leq 2\epsilon_n +2+(2\nu +1)\vartheta $. This contradicts
the fact that $n\in \mathfrak N$.

If there exists $s\in [\bar{v}_i^n,\tilde{v}_{i}^n] \setminus
\nn_{2\vartheta }(\{ \bar{v}_i^n,\tilde{v}_{i}^n \})$ contained in
the $\vartheta$-tubular neighborhood of
$[\tilde{v}_{j}^n,\bar{v}_{j+1}^n]$ then $v_i^n\in \nn_{(\nu +1
)\vartheta +\epsilon_n +1} ([ v_{j}^n,v_{j+1}^n ])$, which
together with the hypothesis $n\in \mathfrak N$ implies that $j\in
\{ i-1, i \}$. The fact that $\dist(s, \tilde{v}_{i}^n)\geq
2\vartheta $ together with the choice of $\tilde{v}_{i}^n$ implies
that $\dist (s, [\tilde{v}_{i}^n,\bar{v}_{i+1}^n] )\geq 2\vartheta
$.  The fact that $\dist (s, \bar{v}_{i}^n)\geq 2\vartheta $
together with the choice of $\bar{v}_{i}^n$ implies that $\dist
(s, [\tilde{v}_{i-1}^n,\bar{v}_{i}^n] )\geq 2\vartheta $.
Therefore $j\not\in \{ i-1, i \}$, a contradiction.

\medskip

\noindent \textbf{Verification of property $(F_2)$ for $n\in
\mathfrak N$.}

\medskip

Let $\bar v=\bar{v}_i^n$ be a vertex of $P_n'$ and let $v=v_i^n$.
We have that $\oo_{\bar v}(P_n')=(RE_n\setminus
[\tilde{v}_{i-1}^n,v]) \cup (AE_n \setminus [\bar
v,\tilde{v}_{i}^n])$. The set $RE_n\setminus
[\tilde{v}_{i-1}^n,\bar v]$ is composed of
$[\tilde{v}_{i}^n,\bar{v}_{i+1}^n]$ and of a part $RE_n'$
contained in $\oo_{v}(P_n)$. By construction we have $\dist (\bar
v, [\tilde{v}_{i}^n,\bar{v}_{i+1}^n] )\geq \nu \vartheta $. On the
other hand $\dist(\bar v, RE_n' )\geq \dist (v, RE_n' )-\dist
(\bar{v}, v)\geq \dist (v, \oo_{v}(P_n) ) - \epsilon_n-1 \geq
\varrho d_n - \epsilon_n-1 $, which is larger that $\nu \vartheta$
for $n\in \mathfrak N$.

Since $AE_n \setminus [\bar{v},\tilde{v}_{i}^n] \subset
\nn_{\epsilon_n+1+\nu \vartheta }(\vv_n\setminus \{ v\})$ it
follows that $$\dist (\bar{v}, AE_n \setminus
[\bar{v},\tilde{v}_{i}^n] )\geq \dist (v, \vv_n\setminus \{ v\} )
- \epsilon_n -1-(\epsilon_n+1+\nu \vartheta )\geq \varrho d_n
-(2\epsilon_n+2+\nu \vartheta)\geq \nu \vartheta $$ for $n\in
\mathfrak N$.

Now let $\tilde v=\tilde{v}_i^n$ be a vertex of $P_n'$. We have
that $\oo_{\tilde{v}}(P_n')=(RE_n\setminus
[\tilde{v},\bar{v}_{i+1}^n]) \cup (AE_n \setminus
[\bar{v},\tilde{v}])$. As before, we show that $\dist (\tilde{v},
AE_n \setminus [\bar{v},\tilde{v}])\geq \nu \vartheta $ for $n\in
\mathfrak N$.

The set $RE_n\setminus [\tilde{v},\bar{v}_{i+1}^n]$ is composed of
$[\tilde{v}_{i-1}^n,\bar{v}]$ and of $RE_n'$. As above, $\dist
(\tilde{v}, RE_n')\geq \nu \vartheta $ for $n\in \mathfrak N$. The
distance $\dist (\tilde{v}, [\tilde{v}_{i-1}^n,\bar{v}] )$ is at
least $\nu \vartheta$ by the choice of $\bar{v}$.

We conclude that for $n\in \mathfrak N$ the polygon $P_n'$ is
$(\vartheta , 2, \nu )$-fat.\endproof

\section{A characterization of asymptotically tree-graded
spaces}\label{amc}

In this section, we find metric conditions for a metric space to
be asymptotically tree-graded with respect to a family of subsets.

\begin{theorem}[a characterization of asymptotically tree-graded spaces]\label{tgi}
Let $(X,\dist)$ be a geodesic metric space and let $\aaa=\{ A_i
\mid i\in I \}$ be a collection of subsets of $X$. The metric
space $X$ is asymptotically tree-graded with respect to $\aaa$ if
and only if the following properties are satisfied:
\begin{itemize}
  \item[$(\alpha_1)$] For every $\delta >0$ the diameters
of the intersections $\nn_{\delta}(A_i)\cap \nn_{\delta}(A_j)$ are
uniformly bounded for all $i\ne j$.

\item[$(\alpha_2)$] For every $\theta$ from $\left[ 0, \frac{1}{2}
\right)$ there exists a number $M>0$ such that for every geodesic
$\q$ of length $\ell$ and every $A\in \aaa$ with
$\q(0),\q(\ell)\in \nn_{\theta\ell}(A)$ we have $\q([0, \ell
])\cap \nn_{M}(A)\neq \emptyset$.

\item[$(\alpha_3)$] For every $k\ge 2$ there exist $\vartheta>0$,
$\nu \geq 8 $ and $\chi
>0$ such that every $k$-gon $P$ in $X$ with
geodesic edges which is $(\vartheta , 2, \nu )$-fat satisfies $P
\subset \nn_{\chi}(A)$ for some $A\in \aaa$.
\end{itemize}
\end{theorem}

\begin{figure}[!ht]
\centering
\unitlength .6 mm 
\linethickness{0.4pt}
\ifx\plotpoint\undefined\newsavebox{\plotpoint}\fi 
\begin{picture}(168,91)(0,0)
\put(97.88,24.25){\oval(79.75,25.5)[]}
\multiput(59.93,33.93)(-.031868,.062637){14}{\line(0,1){.062637}}
\multiput(59.04,35.68)(-.031868,.062637){14}{\line(0,1){.062637}}
\multiput(58.15,37.44)(-.031868,.062637){14}{\line(0,1){.062637}}
\multiput(57.25,39.19)(-.031868,.062637){14}{\line(0,1){.062637}}
\multiput(56.36,40.95)(-.031868,.062637){14}{\line(0,1){.062637}}
\multiput(55.47,42.7)(-.031868,.062637){14}{\line(0,1){.062637}}
\multiput(54.58,44.45)(-.031868,.062637){14}{\line(0,1){.062637}}
\multiput(53.68,46.21)(-.031868,.062637){14}{\line(0,1){.062637}}
\multiput(52.79,47.96)(-.031868,.062637){14}{\line(0,1){.062637}}
\multiput(51.9,49.71)(-.031868,.062637){14}{\line(0,1){.062637}}
\multiput(51.01,51.47)(-.031868,.062637){14}{\line(0,1){.062637}}
\multiput(50.11,53.22)(-.031868,.062637){14}{\line(0,1){.062637}}
\multiput(49.22,54.98)(-.031868,.062637){14}{\line(0,1){.062637}}
\multiput(48.33,56.73)(-.031868,.062637){14}{\line(0,1){.062637}}
\multiput(47.44,58.48)(-.031868,.062637){14}{\line(0,1){.062637}}
\multiput(46.55,60.24)(-.031868,.062637){14}{\line(0,1){.062637}}
\multiput(45.65,61.99)(-.031868,.062637){14}{\line(0,1){.062637}}
\multiput(44.76,63.75)(-.031868,.062637){14}{\line(0,1){.062637}}
\multiput(43.87,65.5)(-.031868,.062637){14}{\line(0,1){.062637}}
\multiput(42.98,67.25)(-.031868,.062637){14}{\line(0,1){.062637}}
\multiput(42.08,69.01)(-.031868,.062637){14}{\line(0,1){.062637}}
\multiput(41.19,70.76)(-.031868,.062637){14}{\line(0,1){.062637}}
\multiput(40.3,72.51)(-.031868,.062637){14}{\line(0,1){.062637}}
\multiput(39.41,74.27)(-.031868,.062637){14}{\line(0,1){.062637}}
\multiput(38.51,76.02)(-.031868,.062637){14}{\line(0,1){.062637}}
\multiput(37.62,77.78)(-.031868,.062637){14}{\line(0,1){.062637}}
\multiput(36.73,79.53)(-.031868,.062637){14}{\line(0,1){.062637}}
\multiput(35.84,81.28)(-.031868,.062637){14}{\line(0,1){.062637}}
\multiput(34.95,83.04)(-.031868,.062637){14}{\line(0,1){.062637}}
\multiput(34.05,84.79)(-.031868,.062637){14}{\line(0,1){.062637}}
\multiput(33.16,86.55)(-.031868,.062637){14}{\line(0,1){.062637}}
\multiput(32.27,88.3)(-.031868,.062637){14}{\line(0,1){.062637}}
\multiput(31.38,90.05)(-.031868,.062637){14}{\line(0,1){.062637}}
\multiput(134.93,33.93)(.032836,.055721){15}{\line(0,1){.055721}}
\multiput(135.91,35.6)(.032836,.055721){15}{\line(0,1){.055721}}
\multiput(136.9,37.27)(.032836,.055721){15}{\line(0,1){.055721}}
\multiput(137.88,38.94)(.032836,.055721){15}{\line(0,1){.055721}}
\multiput(138.87,40.62)(.032836,.055721){15}{\line(0,1){.055721}}
\multiput(139.86,42.29)(.032836,.055721){15}{\line(0,1){.055721}}
\multiput(140.84,43.96)(.032836,.055721){15}{\line(0,1){.055721}}
\multiput(141.83,45.63)(.032836,.055721){15}{\line(0,1){.055721}}
\multiput(142.81,47.3)(.032836,.055721){15}{\line(0,1){.055721}}
\multiput(143.8,48.97)(.032836,.055721){15}{\line(0,1){.055721}}
\multiput(144.78,50.65)(.032836,.055721){15}{\line(0,1){.055721}}
\multiput(145.77,52.32)(.032836,.055721){15}{\line(0,1){.055721}}
\multiput(146.75,53.99)(.032836,.055721){15}{\line(0,1){.055721}}
\multiput(147.74,55.66)(.032836,.055721){15}{\line(0,1){.055721}}
\multiput(148.72,57.33)(.032836,.055721){15}{\line(0,1){.055721}}
\multiput(149.71,59)(.032836,.055721){15}{\line(0,1){.055721}}
\multiput(150.69,60.68)(.032836,.055721){15}{\line(0,1){.055721}}
\multiput(151.68,62.35)(.032836,.055721){15}{\line(0,1){.055721}}
\multiput(152.66,64.02)(.032836,.055721){15}{\line(0,1){.055721}}
\multiput(153.65,65.69)(.032836,.055721){15}{\line(0,1){.055721}}
\multiput(154.63,67.36)(.032836,.055721){15}{\line(0,1){.055721}}
\multiput(155.62,69.03)(.032836,.055721){15}{\line(0,1){.055721}}
\multiput(156.6,70.71)(.032836,.055721){15}{\line(0,1){.055721}}
\multiput(157.59,72.38)(.032836,.055721){15}{\line(0,1){.055721}}
\multiput(158.57,74.05)(.032836,.055721){15}{\line(0,1){.055721}}
\multiput(159.56,75.72)(.032836,.055721){15}{\line(0,1){.055721}}
\multiput(160.54,77.39)(.032836,.055721){15}{\line(0,1){.055721}}
\multiput(161.53,79.06)(.032836,.055721){15}{\line(0,1){.055721}}
\multiput(162.51,80.74)(.032836,.055721){15}{\line(0,1){.055721}}
\multiput(163.5,82.41)(.032836,.055721){15}{\line(0,1){.055721}}
\multiput(164.48,84.08)(.032836,.055721){15}{\line(0,1){.055721}}
\multiput(165.47,85.75)(.032836,.055721){15}{\line(0,1){.055721}}
\multiput(166.45,87.42)(.032836,.055721){15}{\line(0,1){.055721}}
\multiput(167.44,89.09)(.032836,.055721){15}{\line(0,1){.055721}}
\qbezier(31,90)(93,-4)(167,90)
\put(136.93,23.93){\line(1,0){.909}}
\put(138.75,23.93){\line(1,0){.909}}
\put(140.57,23.93){\line(1,0){.909}}
\put(142.38,23.93){\line(1,0){.909}}
\put(144.2,23.93){\line(1,0){.909}}
\put(146.02,23.93){\line(1,0){.909}}
\put(94,23){\makebox(0,0)[cc]{$A$}}
\put(66,60){\makebox(0,0)[cc]{$\q$}}
\put(35,60){\makebox(0,0)[cc]{$\theta l$}}
\put(157,56){\makebox(0,0)[cc]{$\theta l$}}
\put(142,18){\makebox(0,0)[cc]{$M$}}
\put(60.5,15){\oval(23,20)[lb]}
\put(58.93,4.93){\line(1,0){.988}}
\put(60.91,4.95){\line(1,0){.988}}
\put(62.88,4.98){\line(1,0){.988}} \put(64.86,5){\line(1,0){.988}}
\put(66.83,5.03){\line(1,0){.988}}
\put(68.81,5.05){\line(1,0){.988}}
\put(70.78,5.08){\line(1,0){.988}}
\put(72.76,5.1){\line(1,0){.988}}
\put(74.73,5.12){\line(1,0){.988}}
\put(76.71,5.15){\line(1,0){.988}}
\put(78.69,5.17){\line(1,0){.988}}
\put(80.66,5.2){\line(1,0){.988}}
\put(82.64,5.22){\line(1,0){.988}}
\put(84.61,5.25){\line(1,0){.988}}
\put(86.59,5.27){\line(1,0){.988}}
\put(88.56,5.3){\line(1,0){.988}}
\put(90.54,5.32){\line(1,0){.988}}
\put(92.52,5.34){\line(1,0){.988}}
\put(94.49,5.37){\line(1,0){.988}}
\put(96.47,5.39){\line(1,0){.988}}
\put(98.44,5.42){\line(1,0){.988}}
\put(100.42,5.44){\line(1,0){.988}}
\put(102.39,5.47){\line(1,0){.988}}
\put(104.37,5.49){\line(1,0){.988}}
\put(106.34,5.52){\line(1,0){.988}}
\put(108.32,5.54){\line(1,0){.988}}
\put(110.3,5.56){\line(1,0){.988}}
\put(112.27,5.59){\line(1,0){.988}}
\put(114.25,5.61){\line(1,0){.988}}
\put(116.22,5.64){\line(1,0){.988}}
\put(118.2,5.66){\line(1,0){.988}}
\put(120.17,5.69){\line(1,0){.988}}
\put(122.15,5.71){\line(1,0){.988}}
\put(124.12,5.73){\line(1,0){.988}}
\put(126.1,5.76){\line(1,0){.988}}
\put(128.08,5.78){\line(1,0){.988}}
\put(130.05,5.81){\line(1,0){.988}}
\put(132.03,5.83){\line(1,0){.988}}
\put(134,5.86){\line(1,0){.988}}
\put(135.98,5.88){\line(1,0){.988}}
\put(137.95,5.91){\line(1,0){.988}}
\put(136.5,16){\oval(21,20)[rb]}
\put(56.93,48.93){\line(1,0){.987}}
\put(58.9,48.93){\line(1,0){.987}}
\put(60.88,48.93){\line(1,0){.987}}
\put(62.85,48.93){\line(1,0){.987}}
\put(64.83,48.93){\line(1,0){.987}}
\put(66.8,48.93){\line(1,0){.987}}
\put(68.78,48.93){\line(1,0){.987}}
\put(70.75,48.93){\line(1,0){.987}}
\put(72.72,48.93){\line(1,0){.987}}
\put(74.7,48.93){\line(1,0){.987}}
\put(76.67,48.93){\line(1,0){.987}}
\put(78.65,48.93){\line(1,0){.987}}
\put(80.62,48.93){\line(1,0){.987}}
\put(82.6,48.93){\line(1,0){.987}}
\put(84.57,48.93){\line(1,0){.987}}
\put(86.55,48.93){\line(1,0){.987}}
\put(88.52,48.93){\line(1,0){.987}}
\put(90.49,48.93){\line(1,0){.987}}
\put(92.47,48.93){\line(1,0){.987}}
\put(94.44,48.93){\line(1,0){.987}}
\put(96.42,48.93){\line(1,0){.987}}
\put(98.39,48.93){\line(1,0){.987}}
\put(100.37,48.93){\line(1,0){.987}}
\put(102.34,48.93){\line(1,0){.987}}
\put(104.31,48.93){\line(1,0){.987}}
\put(106.29,48.93){\line(1,0){.987}}
\put(108.26,48.93){\line(1,0){.987}}
\put(110.24,48.93){\line(1,0){.987}}
\put(112.21,48.93){\line(1,0){.987}}
\put(114.19,48.93){\line(1,0){.987}}
\put(116.16,48.93){\line(1,0){.987}}
\put(118.13,48.93){\line(1,0){.987}}
\put(120.11,48.93){\line(1,0){.987}}
\put(122.08,48.93){\line(1,0){.987}}
\put(124.06,48.93){\line(1,0){.987}}
\put(126.03,48.93){\line(1,0){.987}}
\put(128.01,48.93){\line(1,0){.987}}
\put(129.98,48.93){\line(1,0){.987}}
\put(131.96,48.93){\line(1,0){.987}}
\put(136,40.5){\oval(22,17)[rt]}
\put(146.93,14.93){\line(0,1){.962}}
\put(146.93,16.85){\line(0,1){.962}}
\put(146.93,18.78){\line(0,1){.962}}
\put(146.93,20.7){\line(0,1){.962}}
\put(146.93,22.62){\line(0,1){.962}}
\put(146.93,24.55){\line(0,1){.962}}
\put(146.93,26.47){\line(0,1){.962}}
\put(146.93,28.39){\line(0,1){.962}}
\put(146.93,30.31){\line(0,1){.962}}
\put(146.93,32.24){\line(0,1){.962}}
\put(146.93,34.16){\line(0,1){.962}}
\put(146.93,36.08){\line(0,1){.962}}
\put(146.93,38.01){\line(0,1){.962}}
\put(48.93,12.93){\line(0,1){.963}}
\put(48.93,14.86){\line(0,1){.963}}
\put(48.93,16.78){\line(0,1){.963}}
\put(48.93,18.71){\line(0,1){.963}}
\put(48.93,20.63){\line(0,1){.963}}
\put(48.93,22.56){\line(0,1){.963}}
\put(48.93,24.49){\line(0,1){.963}}
\put(48.93,26.41){\line(0,1){.963}}
\put(48.93,28.34){\line(0,1){.963}}
\put(48.93,30.26){\line(0,1){.963}}
\put(48.93,32.19){\line(0,1){.963}}
\put(48.93,34.11){\line(0,1){.963}}
\put(48.93,36.04){\line(0,1){.963}}
\put(48.93,37.97){\line(0,1){.963}}
\put(58.5,41.5){\oval(19,15)[lt]}
\end{picture}
\caption{Property ($\alpha_2$).} \label{fig3}
\end{figure}
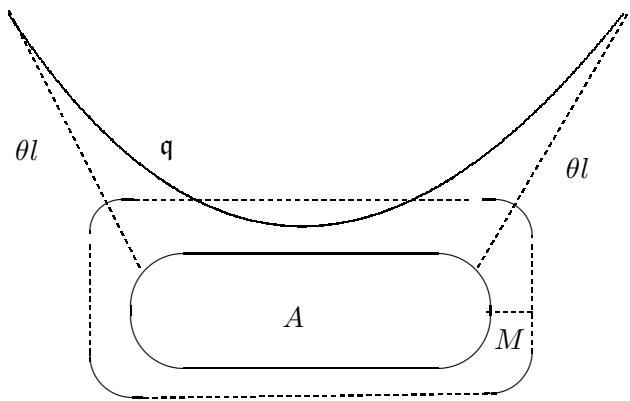

\begin{remarks}\label{tgi3}

(1) If the space $X$ is asymptotically uniquely geodesic (for
instance asymptotically $CAT(0)$) then in $(\alpha_3)$ it is
enough to consider $k=3$ (only triangles).

(2) From the proof of Theorem \ref{tgi}, it will be clear that
conditions $(\alpha_2)$, $(\alpha_3)$ can be replaced by the
following stronger conditions:

\begin{itemize}

\item[$(\alpha_2')$] For every $L\geq 1$, $C\ge 0$, and $\theta
\in \left[ 0, \frac{1}{2} \right)$ there exists $M>0$ such that
for every $(L,C)$-quasi-geodesic $\q$ defined on $[0,\ell ]$ and
every $A\in \aaa$ such that
  $\q(0),\q(\ell )\in \nn_{\theta\ell/L}
  (A)$ we have $\q([0, \ell ])\cap \nn_{M}
  (A)\neq \emptyset$;

\item[$(\alpha_3')$] For every $L\geq 1$, $C\geq 0$ and $k\geq 2$,
 and for every $\sigma \geq 1$
and $\nu \geq 4\sigma $, there exist $\vartheta_0>0$ such that for
every $\vartheta \geq \vartheta_0$ every $k$-gon $P$ with
$(L,C)$-quasi-geodesic edges which is $(\vartheta\,  ,\, \sigma \,
,\, \nu)$-fat is contained in $\nn_{\chi}(A)$ for some $A\in
\aaa$, where $\chi = \sigma L^2 \vartheta +c$ with $c$ a constant
independent of $\vartheta $.
\end{itemize}

(3) Also from the proof of Theorem \ref{tgi}, it will be clear
that for every $\epsilon\le\frac12$ the condition $(\alpha_2)$ can
be replaced by the following weaker condition :

\begin{itemize}
\item[$(\alpha_2^\epsilon)$] For every $\theta$ from $\left[ 0,
\epsilon \right)$ there exists a number $M>0$ such that for every
geodesic $\q$ of length $\ell$ and every $A\in \aaa$ with
$\q(0),\q(\ell)\in \nn_{\theta\ell}(A)$ we have $\q([0, \ell
])\cap \nn_{M}(A)\neq \emptyset$.
\end{itemize}
(Notice that condition $(\alpha_2)$ is the same as the condition
$(\alpha_2^{\frac12})$.)

\me

(4) If $\aaa=\{ A_i\mid i\in I\}$ satisfies conditions
$(\alpha_1)$, $(\alpha_2)$, $(\alpha_3)$, then the family
$\nn_c(\aaa)=\{\nn_c(A_i)\mid i\in I\}$ also satisfies these
conditions, for every $c>0$.
\end{remarks}

{\em Proof of Theorem \ref{tgi}.} First we show that conditions
$(\alpha_1)$, $(\alpha_2^\epsilon)$ (for an arbitrary
$\epsilon\le\frac12$) and $(\alpha_3)$ imply that $X$ is
asymptotically tree-graded with respect to $\aaa$.


\begin{lemma}[$(\alpha_1)$ and $(\alpha_2^\epsilon)$ imply uniform quasi-convexity]\label{12qc}  Let $(X,d)$ be a geodesic metric space and let $\aaa=\{ A_i \mid
i\in I \}$ be a collection of subsets of $X$ satisfying properties
$(\alpha_1)$ and $(\alpha_2^\epsilon)$ for some $\epsilon$. Let
$M_0=M_0(\theta)$ be the number from property
$(\alpha_2^\epsilon)$ corresponding to $\theta = \frac{2}{3}\epsilon$.

There exists $t>0$ such that for every $A\in\aaa$, $M\geq M_0$ and
$x,y\in \nn_M(A)$, every geodesic joining $x$ and $y$ in $X$ is
contained in $\nn_{tM}(A)$.
\end{lemma}

\proof Suppose, by contradiction, that for every $n\in \N $ there
exist $M_n\geq M_0$, $x_n,y_n \in \nn_{M_n}(A_n)$ and a geodesic
$[x_n,y_n]$ not contained in $\nn_{n M_n}(A_n)$. For every $n\ge
1$ let $D_n$ be the infimum of the distances between points
$x,y\in \nn_{M_n}(A)$ for some $A\in\aaa$ such that $[x,y]\not
\subset\nn_{nM_n}(A)$ for some geodesic $[x,y]$.

We note that $D_n\geq 2(n-1)M_n\geq 2(n-1)M_0$, hence $\lm_{n\to
\infty}D_n=\infty$. For every $n\ge 1$, choose $x_n,y_n \in
\nn_{M_n}(A_n)$ such that $\dist(x_n,y_n)=D_n+1$. Also choose
$a_n,b_n\in [x_n,y_n]$ such that
$\dist(x_n,a_n)=\dist(y_n,b_n)=\frac{\theta(D_n+1)}{2}$. Then
$\dist (a_n, A_n)\leq \dist(a_n,x_n)+\dist(x_n,A_n)\leq
\frac{\theta(D_n+1)}{2} + M_n \leq \frac{\theta(D_n+1)}{2} +
\frac{D_n+1}{2(n-1)}$. Likewise $\dist (b_n, A_n)\leq
\frac{\theta(D_n+1)}{2} + \frac{D_n+1}{2(n-1)}$. On the other hand
$\dist(a_n,b_n)\geq \dist(x_n,y_n)-\dist(x_n,a_n)-
\dist(y_n,b_n)\geq (1-\theta)(D_n+1)$. For $n$ large enough we
have $\frac{\theta}{2} + \frac{1}{2(n-1)}\leq \frac{2}{3}\theta$.
We apply $(\alpha_2^\epsilon)$ with $\theta =\frac{2}{3}\epsilon$
to $[a_n,b_n]$ and we deduce that there exists $z_n \in [a_n,b_n]
\cap \nn_{M_0}(A_n)$. We have that either
$[x_n,z_n]\not\subset\nn_{nM_n}(A_n)$ or $[z_n,y_n]\not \subset
\nn_{nM_n}(A_n)$, while $\dist(x_n,z_n),\dist(z_n,y_n)\leq
(1-\frac{\theta}{2})(D_n+1)< D_n$ for $n$ large enough. This
contradicts the choice of $D_n$. \endproof

\begin{lemma}\label{connected}
 Let $(X,d)$ be a geodesic metric space and let $\aaa=\{ A_i \mid
i\in I \}$ be a collection of subsets of $X$ satisfying properties
$(\alpha_1)$ and $(\alpha_2^\epsilon)$ for some $\epsilon$. Then
in every asymptotic cone $\co{X;e,d}$, every set $\lio{A_n}$ is
connected and a geodesic subspace.
\end{lemma}

\proof Indeed, consider any two points $x=\lio{x_n}, y=\lio{y_n}$
in $\lio{A_n}$, and geodesics $\q_n$ connecting $x_n,y_n$ in $X$.
Then by Lemma \ref{12qc}, $\q_n$ is inside $\nn_M(A_n)$ for some
fixed $M$. Therefore the geodesic $\lio{\q_n}$ is inside
$\lio{\nn_M(A_n)}=\lio{A_n}$.
\endproof

\begin{lemma}\label{at1} Let $(X,d)$ be a geodesic metric space and let $\aaa=\{ A_i \mid
i\in I \}$ be a collection of subsets of $X$ satisfying properties
$(\alpha_1)$ and $(\alpha_2^\epsilon)$. Then in every asymptotic
cone $Con^\omega(X;e,d)$ the collection of subsets $\aaa_\omega$
satisfies $(T_1)$.
\end{lemma}

\proof Suppose that, in an asymptotic cone $\co{X; e,d}$ of $X$,
the intersection $\lio{A_{i_n}} \cap \lio{A_{j_n}}$ contains two
distinct points $\lio{x_n}, \lio{y_n}$ but $A_{i_n}\neq A_{j_n}$
$\omega$-almost surely. For every $n\ge 1$ consider a geodesic
$[x_n,y_n]$. Its length $\ell_n$ is $O(d_n)$ while $\delta_n$
defined as the maximum of the distances $\dist(x_n, A_{i_n}),
\dist(x_n,A_{j_n}), \dist(y_n,A_{i_n}),\dist(y_n,A_{j_n})$, is
$o(d_n)$. According to Lemma \ref{12qc}, $[x_n,y_n]$ is contained
in $\nn_{t\delta_n}(A_{i_n}) \cap \nn_{t\delta_n} (A_{j_n})$ for
some $t>0$.

Consider $a_n, b_n \in [x_n,y_n]$ at distance $6t\delta_n$ from
$x_n$ and $y_n$, respectively. Property $(\alpha_2^\epsilon)$ can
be applied twice, to $[x_n,a_n]\subset [x_n,y_n]$ and $A_{i_n}$
(resp. $A_{j_n}$) for $n$ large enough. It implies that there
exist $z_n \in [x_n,a_n] \cap \nn_{M_0}(A_{i_n})$ and $z_n' \in
[x_n,a_n] \cap \nn_{M_0}(A_{j_n})$ (where $M_0$ is the same as in
Lemma \ref{12qc}). A similar argument for $[b_n,y_n]\subset
[x_n,y_n]$ and $A_{i_n}$ (resp. $A_{j_n}$) implies that there
exist $u_n \in [b_n,y_n] \cap \nn_{M_0}(A_{i_n})$ and $u_n' \in
[b_n,y_n] \cap \nn_{M_0}(A_{j_n})$. Hence $[a_n,b_n]\subset
[z_n,u_n]\subset \nn_{t M_0}(A_{i_n})$ and $[a_n,b_n]\subset
[z_n',u_n']\subset \nn_{t M_0}(A_{j_n})$. It follows that
$[a_n,b_n]\subset \nn_{t M_0}(A_{i_n}) \cap \nn_{tM_0}(A_{j_n})$,
while $\dist(a_n,b_n)=O(d_n)$. This contradicts property
$(\alpha_1)$.\endproof

\begin{lemma}\label{a2T2} {\bf (asymptotic $(T_1)$ and $(\alpha_3)$ implies asymptotic $(T_2)$)}
Let $(X,\dist)$ be a geodesic metric space and let $\aaa=\{ A_i
\mid i\in I \}$ be a collection of subsets of $X$. Suppose that
property $(\alpha_3)$ holds. Then every simple geodesic triangle
in any asymptotic cone $\co{X;e,d}$ is contained in one of the
sets from $\aaa_\omega$.
\end{lemma}

\proof Let $\Delta$ be a simple geodesic triangle in $\co{X;
e,d}$. Let $\varepsilon_m=\frac{1}{2^m}$ be fixed, for every large
enough integer $m$. By Proposition \ref{approx}, we can find $k_0$
and a simple triangle $\Delta_m=\Delta_{\varepsilon_m}
  =\lio{P^m_n}$ satisfying properties (a),(b), and (c) for
   $\vartheta $ and $\nu \geq 8$ given by $(\alpha_3)$ for $k_0(\varepsilon_m )$. It follows that $\omega$-almost surely, $P^m_n$ are contained in
$\nn_{\chi} (A_n)$ for some $A_n\in \aaa$. We conclude that
$\Delta_m \subset A_\omega =\lm^\omega (A_n)$. By property (b) all
triangles $\Delta_m$ have at least 3 distinct points in common
(e.g. the midpoints of the edges of $\Delta$). This and property
$(T_1)$ of the collection $\aaa_\omega$ imply that the set
$A_\omega$ is independent of $m$. Since $\Delta$ is a Hausdorff
limit of $\Delta_m$ and $A_\omega$ is closed (see Remark
\ref{rkcomp}), we deduce that $\Delta \subset A_\omega$.\endproof

Lemmas \ref{connected}, \ref{at1} and \ref{a2T2} show that
$(\alpha_1)$, $(\alpha_2^\epsilon)$, $(\alpha_3)$ imply that the
space $X$ is asymptotically tree-graded. Now we prove the
(stronger version of the) converse statement.

\begin{lemma}[asymptotic
$(T_1)$ implies $(\alpha_1)$] \label{T1a1} Let $(X,\dist)$ be a
geodesic metric space asymptotically satisfying $(T_1)$ with
respect to $\aaa$. Then $X$ satisfies $(\alpha_1)$ with respect to
$\aaa$.
\end{lemma}

\proof By contradiction, suppose $X$ asymptotically satisfies
$(T_1)$ but for some $\delta
>0$ there exists a sequence of pairs of points $x_n , y_n$ in $\nn_{\delta }(A_{i_n})\cap
\nn_{\delta}(A_{j_n}),$ where $A_{i_n}$ and $A_{j_n}$ are distinct
sets in $\aaa$, with $\lm_{n\to\infty}\dist(x_n,y_n)=\infty$. Set
the observation point $e$ to be $(x_n)^\omega$, and let
$d_n=\dist(x_n,y_n)$ for every $n\ge 1$. Then $M_1=\lio{A_{i_n}}$
and $M_2=\lio{A_{j_n}}$ are not empty, so these are distinct
pieces in $\co{X; e,d}$. The limits $x=\lio{x_n}$ and
$y=\lio{y_n}$ are distinct points in $\co{X;e,d}$ that belong to
both $M_1$ and $M_2$. This contradicts $(T_1)$. \endproof

\begin{definition}[almost closest points] Let $x\in X$,
$A, B\subseteq X$. A point $y\in
A$ is called an {\em almost closest to $x$ point in $A$} if
$\dist(x,y)\le \dist(x,A)+1$. Points $a\in A$, $b\in B$ are called
{\em almost closest representatives of $A$ and $B$} if
$\dist(a,b)\le \dist(A,B)+1$.
\end{definition}

\begin{definition}[almost projection]
Let $x$ be a point in $X$ and $A\subset X$. The \textit{almost
projection of $x$ on $A$} is the set of almost closest to $x$
points in $A$. For every subset $B$ of $X$ we define the {\em
almost projection} $\proj_A(B)$ of $B$ onto $A$ as $\bigcup_{b\in
B}\proj_A(b)$.
\end{definition}

\begin{remark}
If all $A\in \aaa$ were closed sets and the space $X$ was proper
(i.e. all balls in $X$ compact) then we could use closest points
and usual projections instead of almost closest points and almost
projections.
\end{remark}

\begin{lemma}\label{eproj}
If the space $X$ is asymptotically tree-graded with respect to
$\aaa$ then for every $x\in X$, $A\in \aaa$, with $\dist(x,A)=2d$
$$\diam(\proj_A(\nn_d(x))=o(d).$$
\end{lemma}

\proof Suppose there exists $\varepsilon>0$ and $x_n\in X, A_n\in
\aaa$ with $\dist(x_n,A_n)=2d_n, \lm_{n\to \infty} d_n=\infty$,
and the projection $\proj_{A_n}(\nn_{d_n}(x_n))$ is of diameter at
least $\varepsilon d_n$. Let $e=(x_n)$ and $d=(d_n)$. In the
asymptotic cone $\co{X; e,d}$, we have the point $x=\lio{x_n}$ at
distance 2 of $A=\lio{A_n}$, two points $y,z\in \nn_1(x)$, and two
points $y', z'$ in $A$ such that $y', z'$ are the respective
projections of $y, z$ onto $A$, but $\dist(y',z')\ge \varepsilon$.
This contradicts Lemma \ref{smball}.
\endproof

\begin{lemma}[asymptotically tree-graded
implies $(\alpha_2')$]\label{T1T2a2}  Let $(X,\dist)$ be a
geodesic metric space which is asymptotically tree-graded with
respect to $\aaa$. Then $X$ satisfies $(\alpha_2')$.
\end{lemma}

\proof Fix $L\ge 1, C\ge 0$. By contradiction, suppose that for
some fixed $\theta\in \left[ 0, \frac{1}{2} \right)$ there exists
a sequence of $(L,C)$-quasi-geodesics $\q_n\colon [0,\ell_n ]\to
X$ and a sequence of sets $A_n\in \aaa$, such that $\q_n(0),
\q_n(\ell_n)\in \nn_{\theta \ell_n/L}(A_n)$ and $\dist\left(
\q_n([0,\ell_n]), A_n \right)=2D_n, \, \lm_{n\to\infty}D_n
=\infty$. Since $\dist\left( \q_n([0,\ell_n]), A_n \right) \leq L
\ell_n +\frac{\theta \ell_n}{L}$ this implies
$\lm_{n\to\infty}\ell_n=\infty $.

Let $t_0=0 <t_1<\cdots <t_{m-1}<t_m=\ell_n$ be such that
$\frac{D_n-C}{2L}\leq \dist(t_i,t_{i+1})\leq \frac{D_n -C}{L}$ for
all $i\in \{ 0,1,\dots ,m-1\}$. We have $m \leq
\frac{3L\ell_n}{D_n}$ for large enough $n$. Let $y_i$ be an almost
projection of $\q_n(t_i)$ onto $A_n$. According to Lemma
\ref{eproj}, $\dist(y_i,y_{i+1})=o(D_n)$. Consequently
$\dist(\q_n(0), \q_n(\ell_n))\leq
\dist(\q_n(0),y_0)+\Sigma_{i=0}^{m-1}\dist(y_i,y_{i+1})+\dist(y_m,
\q_n(\ell_n))\leq \frac{2\theta \ell_n}{L}+m\cdot o(D_n)\leq
\frac{2\theta \ell_n}{L}+3Lo(1)\ell_n$. On the other hand
$\dist(\q_n(0), \q_n(\ell_n))\geq \frac{\ell_n}{L}-C$. This is a
contradiction with $\theta<\frac12$.\endproof

It remains to prove that being asymptotically tree-graded implies
$(\alpha_3)$.

\begin{definition}[almost geodesics] If an $(L,C)$-quasi-geodesic $\q$ is $L$-Lipschitz then $\q$ will be called an {\em
$(L,C)$-almost geodesic}.
\end{definition}

\begin{remark} \label{buragos}Every
$(L,C)$-quasi-geodesic in a geodesic metric space is at bounded
(in terms of $L, C$) distance from an $(L+C,C)$-almost geodesic
with the same end points \cite[Proposition 8.3.4]{Buragos}.
\end{remark}

\begin{lemma}[$\aaa$ is uniformly quasi-convex with respect to quasi-geodesics]
\label{qqc}

Let $X$ be a geodesic metric space which is asymptotically
tree-graded with respect to a collection of subsets $\aaa$. For
every $L\geq 1$ and $C\geq 0$, there exists $t\geq 1$ such that
for every $d\geq 1$ and for every $A\in \aaa$, every
$(L,C)$-quasi-geodesic joining two points in $\nn_d(A)$ is
contained in $\nn_{td}(A)$.
\end{lemma}

\proof Suppose by contradiction that there exists a sequence $\q_n
:[0,\ell_n]\to X $ of $(L,C)$-quasi-geodesics with endpoints
$x_n,y_n \in \nn_{d_n}(A_n)$ such that there exists $z_n \in
\q_n([0,\ell_n]) $ with $k_n=\dist(z_n,A_n)\geq nd_n\ge n$. By
Remark \ref{buragos}, we can assume that each $\q_n$ is an
$(L+C,C)$-almost geodesic. This allows us to choose $z_n \in
\q_n([0,\ell_n])$ so that $\dist(z_n,A_n)$ is maximal. In $\co{X;
(z_n),(k_n)}$, the limit set $\q=\lio{\q_n}$ is either a
topological arc with endpoints in $\lio{A_n}$ and not contained in
$\lio{A_n}$, or a bi-Lipschitz ray with origin in $\lio{A_n}$ or a
bi-Lipschitz line (Remark \ref{grr05}). Notice also that $\q$ is
contained in $\nn_1(\lio{A_n})$. In all three cases we obtain a
contradiction with  Corollary \ref{projpath}.\endproof

Let $(X,\dist)$ be a geodesic space that is asymptotically
tree-graded with respect to the collection of subsets $\aaa$.

\me

\noindent \textit{Notation}: For every $L\geq 1,C\geq 0$, we
denote by $M(L,C)$ the constant given by $(\alpha_2')$ for
$\theta=\frac13$. We also denote by $\dist$ the distance function
in any of the asymptotic cones of $X$.

\me

\noindent \textit{Conventions}: To simplify the notations and
statements, in the sequel we shall not mention the constants
$L\geq 1$ and $C\geq 0$ for each quasi-geodesic anymore. We assume
that all constants provided by the following lemmas in the section
depend on $L$ and $C$.

\begin{lemma}\label{Aprime} Let $\q_n\colon [0,\ell_n]\to X$, $n\ge 1$, be
a sequence of $(L,C)$-quasi-geodesics in $X$ and let $A_n$, $n\ge
1$, be a sequence of sets in $\aaa$. Suppose that
$\dist(\q_n(0),A_n)=o(\ell_n), \dist(\q_n(\ell_n),A_n)=o(\ell_n)$
$\omega$-almost surely. Then there exists $t^1_n\in [0,
\frac13\ell_n]$, $t_n^2\in[\frac23\ell_n,\ell_n]$ such that
$\q_n(t_n^i)\in \nn_M(A_n)$, $i=1, 2$, where $M=M(L,C)$,
$\omega$-almost surely.
\end{lemma}

\proof By Lemma \ref{qqc}, the quasi-geodesic $\q_n$ is inside
$\nn_{t_n}(A_n)$ for $t_n=o(\ell_n)$. It remains to apply
$(\alpha_2')$ to the quasi-geodesics $\q_n([0,\frac13 \ell_n])$
and $\q_n([\frac23\ell_n,\ell_n])$.
\endproof

\begin{lemma}[linear divergence]\label{lemaA}
For every $\varepsilon >0$ and every $M\geq M(L,C)$ there exists
$t_\varepsilon
>0$ such that if $A\in \aaa$, $\q$ is a quasi-geodesic
with origin $a\in \nn_M (A)$, such that $\q\cap \nn_M (A) = \{ a\}
$ and $t\geq t_\varepsilon $ then
$$
\dist(\q(t),A)> (1-\varepsilon )\dist(\q(t),a)\, .
$$
\end{lemma}

\proof We suppose that for some $\varepsilon >0$ there exists a
sequence $A_n\in \aaa$, a sequence $\q_n$ of quasi-geodesics with
origin $a_n\in \nn_M (A_n)$ such that $\q_n\cap \nn_M (A_n) = \{
a_n\} $, and a sequence of numbers $t_n\to \infty $ with the
property
$$
\dist(\q_n(t_n),A_n)\leq (1-\varepsilon )\dist(\q_n(t_n),a_n)\, .
$$

In $\co{X; (a_n), (t_n)}$, we obtain the points $a=\lio{a_n} \in
\lio{A_n}$ and $b=\lio{\q_n(t_n)}$, joined by the bi-Lipschitz arc
$\q([0,1])=\lio{\q_n([0,t_n])}$, such that
$$\dist(b, \lio{A_n})\leq (1-\varepsilon)\dist(b,a).$$ It follows
that the projection of $b$ on $\lio{A_n}$ is a point $c\neq a$.
Corollary \ref{projA} implies that $\q([0,1])$ contains $c$ and
Corollary \ref{strconv} implies that a sub-arc $\q([0,2\beta])$ of
$\q([0,1])$ is contained in $\lio{A_n}$. We apply Lemma
\ref{Aprime} to the sub-quasi-geodesic $\q_n([0,\beta t_n])$ and
obtain that this sub-quasi-geodesics intersects $\nn_M(A_n)$ in a
point different from $a_n$, a contradiction.\endproof

\begin{lemma}\label{lemaB}
For every $\varepsilon>0$, $\delta>0$ and $M\geq M(L,C)$ there
exists $D>0$ such that for every $A\in \aaa$ and every two
quasi-geodesics $\q_i\colon [0,\ell_i]\to X$, $i=1,2$, that
connect $a\in \nn_M (A)$ with two points $b_1$ and $b_2$
respectively, if the diameter of $\q_1\cap \nn_M (A)$ does not
exceed $\delta$, $b_2\in \nn_M (A)$, and $\dist(a,b_2)\geq D$ then
$$
\dist(b_1,b_2)\geq \frac{1}{L+\varepsilon }(\ell_1+\ell_2 ).
$$
\end{lemma}

\proof Suppose there exist sequences $\q^{(n)}_i
:[0,\ell^{(n)}_i]\to X$, $i=1,2$, $n\ge 1$, of pairs of
quasi-geodesics joining $a^{(n)}\in \nn_M(A_n)$ to $b^{(n)}_i$
such that $\q^{(n)}_1\cap \nn_M(A_n)$ has diameter at most
$\delta$, $b^{(n)}_2\in \nn_M(A_n)$,
$\lm_{n\to\infty}\dist(a^{(n)},b^{(n)}_2)=\infty$, but
\begin{equation}\label{ra}
\dist(b^{(n)}_1,b^{(n)}_2)\leq \frac{1}{L+\varepsilon
}(\ell^{(n)}_1+\ell^{(n)}_2).
\end{equation}

Denote $\dist(a^{(n)},b^{(n)}_1)$ by $f_n$ and
$\dist(a^{(n)},b^{(n)}_2)$ by $d_n$. Since $\ell^{(n)}_1\leq L
(f_n+C)$, $\ell^{(n)}_2\le L(d_n+C)$, for every large enough $n$
the inequality (\ref{ra}) implies that
\begin{equation}\label{comp}
\dist(b^{(n)}_1,b^{(n)}_2)\leq (1-\gamma)(f_n+d_n).
\end{equation}
for some $\gamma>0$.

\noindent \textbf{Case I.}\quad Suppose that
$\lio{\frac{f_n}{d_n}} < \infty $. In the asymptotic cone $\co{X;
(a_n), (d_n)}$,  the two points $\lio{b^{(n)}_i}$, $i=1,2$, are
joined by the Lipschitz arc $\lio{\q^{(n)}_1}\cup
\lio{\q^{(n)}_2}$ (it is Lipschitz as any union of two Lipschitz
arcs). Lemma \ref{lemaA} implies that
$$\lio{\q^{(n)}_1}\cap \lio{\q^{(n)}_2}
={\lio{a^{(n)}}}$$ (here we use the fact that the diameters of the
intersections $\q^{(n)}_1$ with $\nn_M(A_n)$ are uniformly
bounded, so we can cut a comparatively little piece of each
$\q^{(n)}_1$ to make it satisfy the conditions of Lemma
\ref{lemaA}).

Thus the points $\lio{b^{(n)}_i}$ are joined by the simple arc
$\lio{\q^{(n)}_1} \cup \lio{\q^{(n)}_2}.$ This and property
$(T_2')$ imply that every geodesic joining $\lio{b^{(n)}_1}$ and
$\lio{b^{(n)}_2}$ contains $\lio{a^{(n)}}$. Therefore
$$\dist(\lio{b^{(n)}_1}, \lio{b^{(n)}_2} )=\dist(\lio{b^{(n)}_1},
\lio{a^{(n)}})+\dist(\lio{a^{(n)}}, \lio{b^{(n)}_2}).$$ This
contradicts the inequality (\ref{comp}).

\me

\noindent \textbf{Case II.}\quad Suppose that $\lm_\omega
\frac{f_n}{d_n} =\infty $. In the asymptotic cone $\co{X;
(a^{(n)}),(f_n)}$, we denote $a=\lio{a^{(n)}}=\lio{b^{(n)}_2}\in
\lio{A_n}$ and $b=\lio{b^{(n)}_1}$. Then inequality (\ref{comp})
implies that $\dist(a, b)\leq (1-\gamma)\dist(a, b)$, a
contradiction.\endproof

\begin{lemma}\label{lemaC}
For every $M\geq M(L,C)$, $\varepsilon >0$ and $\delta>0$ there
exists $D'>0$ such that for every $A\in \aaa$, and every two
quasi-geodesics $\q_i\colon [0,\ell_i]\to X$, $i=1,2$, joining $a$
in $\nn_M (A)$ with $b_i$, if the diameter of $\q_1\cap \nn_M (A)$
does not exceed $\delta$, $b_2\in \nn_M (A)$, $\dist(a,b_2)\geq
D'$, then the union $\q_1\sqcup \q_2$ of these two quasi-geodesics
is an $(L+\varepsilon,K )$-quasi-geodesic, where $K = 2D'$.
\end{lemma}

\begin{figure}[!ht]
\centering
\unitlength .7mm 
\linethickness{0.4pt}
\ifx\plotpoint\undefined\newsavebox{\plotpoint}\fi 
\begin{picture}(147,101)(0,0)
\put(97.88,24.25){\oval(79.75,25.5)[]}
\put(136.93,23.93){\line(1,0){.909}}
\put(138.75,23.93){\line(1,0){.909}}
\put(140.57,23.93){\line(1,0){.909}}
\put(142.38,23.93){\line(1,0){.909}}
\put(144.2,23.93){\line(1,0){.909}}
\put(146.02,23.93){\line(1,0){.909}}
\put(94,23){\makebox(0,0)[cc]{$A$}}
\put(142,18){\makebox(0,0)[cc]{$M$}}
\put(60.5,15){\oval(23,20)[lb]}
\put(58.93,4.93){\line(1,0){.988}}
\put(60.91,4.95){\line(1,0){.988}}
\put(62.88,4.98){\line(1,0){.988}} \put(64.86,5){\line(1,0){.988}}
\put(66.83,5.03){\line(1,0){.988}}
\put(68.81,5.05){\line(1,0){.988}}
\put(70.78,5.08){\line(1,0){.988}}
\put(72.76,5.1){\line(1,0){.988}}
\put(74.73,5.12){\line(1,0){.988}}
\put(76.71,5.15){\line(1,0){.988}}
\put(78.69,5.17){\line(1,0){.988}}
\put(80.66,5.2){\line(1,0){.988}}
\put(82.64,5.22){\line(1,0){.988}}
\put(84.61,5.25){\line(1,0){.988}}
\put(86.59,5.27){\line(1,0){.988}}
\put(88.56,5.3){\line(1,0){.988}}
\put(90.54,5.32){\line(1,0){.988}}
\put(92.52,5.34){\line(1,0){.988}}
\put(94.49,5.37){\line(1,0){.988}}
\put(96.47,5.39){\line(1,0){.988}}
\put(98.44,5.42){\line(1,0){.988}}
\put(100.42,5.44){\line(1,0){.988}}
\put(102.39,5.47){\line(1,0){.988}}
\put(104.37,5.49){\line(1,0){.988}}
\put(106.34,5.52){\line(1,0){.988}}
\put(108.32,5.54){\line(1,0){.988}}
\put(110.3,5.56){\line(1,0){.988}}
\put(112.27,5.59){\line(1,0){.988}}
\put(114.25,5.61){\line(1,0){.988}}
\put(116.22,5.64){\line(1,0){.988}}
\put(118.2,5.66){\line(1,0){.988}}
\put(120.17,5.69){\line(1,0){.988}}
\put(122.15,5.71){\line(1,0){.988}}
\put(124.12,5.73){\line(1,0){.988}}
\put(126.1,5.76){\line(1,0){.988}}
\put(128.08,5.78){\line(1,0){.988}}
\put(130.05,5.81){\line(1,0){.988}}
\put(132.03,5.83){\line(1,0){.988}}
\put(134,5.86){\line(1,0){.988}}
\put(135.98,5.88){\line(1,0){.988}}
\put(137.95,5.91){\line(1,0){.988}}
\put(136.5,16){\oval(21,20)[rb]}
\put(56.93,48.93){\line(1,0){.987}}
\put(58.9,48.93){\line(1,0){.987}}
\put(60.88,48.93){\line(1,0){.987}}
\put(62.85,48.93){\line(1,0){.987}}
\put(64.83,48.93){\line(1,0){.987}}
\put(66.8,48.93){\line(1,0){.987}}
\put(68.78,48.93){\line(1,0){.987}}
\put(70.75,48.93){\line(1,0){.987}}
\put(72.72,48.93){\line(1,0){.987}}
\put(74.7,48.93){\line(1,0){.987}}
\put(76.67,48.93){\line(1,0){.987}}
\put(78.65,48.93){\line(1,0){.987}}
\put(80.62,48.93){\line(1,0){.987}}
\put(82.6,48.93){\line(1,0){.987}}
\put(84.57,48.93){\line(1,0){.987}}
\put(86.55,48.93){\line(1,0){.987}}
\put(88.52,48.93){\line(1,0){.987}}
\put(90.49,48.93){\line(1,0){.987}}
\put(92.47,48.93){\line(1,0){.987}}
\put(94.44,48.93){\line(1,0){.987}}
\put(96.42,48.93){\line(1,0){.987}}
\put(98.39,48.93){\line(1,0){.987}}
\put(100.37,48.93){\line(1,0){.987}}
\put(102.34,48.93){\line(1,0){.987}}
\put(104.31,48.93){\line(1,0){.987}}
\put(106.29,48.93){\line(1,0){.987}}
\put(108.26,48.93){\line(1,0){.987}}
\put(110.24,48.93){\line(1,0){.987}}
\put(112.21,48.93){\line(1,0){.987}}
\put(114.19,48.93){\line(1,0){.987}}
\put(116.16,48.93){\line(1,0){.987}}
\put(118.13,48.93){\line(1,0){.987}}
\put(120.11,48.93){\line(1,0){.987}}
\put(122.08,48.93){\line(1,0){.987}}
\put(124.06,48.93){\line(1,0){.987}}
\put(126.03,48.93){\line(1,0){.987}}
\put(128.01,48.93){\line(1,0){.987}}
\put(129.98,48.93){\line(1,0){.987}}
\put(131.96,48.93){\line(1,0){.987}}
\put(136,40.5){\oval(22,17)[rt]}
\put(146.93,14.93){\line(0,1){.962}}
\put(146.93,16.85){\line(0,1){.962}}
\put(146.93,18.78){\line(0,1){.962}}
\put(146.93,20.7){\line(0,1){.962}}
\put(146.93,22.62){\line(0,1){.962}}
\put(146.93,24.55){\line(0,1){.962}}
\put(146.93,26.47){\line(0,1){.962}}
\put(146.93,28.39){\line(0,1){.962}}
\put(146.93,30.31){\line(0,1){.962}}
\put(146.93,32.24){\line(0,1){.962}}
\put(146.93,34.16){\line(0,1){.962}}
\put(146.93,36.08){\line(0,1){.962}}
\put(146.93,38.01){\line(0,1){.962}}
\put(48.93,12.93){\line(0,1){.963}}
\put(48.93,14.86){\line(0,1){.963}}
\put(48.93,16.78){\line(0,1){.963}}
\put(48.93,18.71){\line(0,1){.963}}
\put(48.93,20.63){\line(0,1){.963}}
\put(48.93,22.56){\line(0,1){.963}}
\put(48.93,24.49){\line(0,1){.963}}
\put(48.93,26.41){\line(0,1){.963}}
\put(48.93,28.34){\line(0,1){.963}}
\put(48.93,30.26){\line(0,1){.963}}
\put(48.93,32.19){\line(0,1){.963}}
\put(48.93,34.11){\line(0,1){.963}}
\put(48.93,36.04){\line(0,1){.963}}
\put(48.93,37.97){\line(0,1){.963}}
\put(58.5,41.5){\oval(19,15)[lt]}
\put(59,44){\makebox(0,0)[cc]{$a$}}
\put(135,42){\makebox(0,0)[cc]{$b_2$}}
\put(46,100){\makebox(0,0)[cc]{$b_1$}}
\qbezier(49,101)(53.5,99.5)(52,92) \qbezier(52,92)(50,84.5)(52,83)
\qbezier(52,83)(57.5,80)(53,73) \qbezier(53,73)(48.5,65.5)(54,60)
\qbezier(54,60)(61,53.5)(56,53) \qbezier(56,53)(50.5,53)(55,45)
\put(57,85){\makebox(0,0)[cc]{$\q_1$}}
\qbezier(55,45)(66,54)(73,47) \qbezier(73,47)(85,36)(89,41)
\qbezier(89,41)(92,48)(99,43) \qbezier(99,43)(109,36.5)(115,52)
\qbezier(115,52)(118,60.5)(125,47)
\qbezier(125,47)(128.5,41.5)(130,42)
\put(124,56){\makebox(0,0)[cc]{$\q_2$}}
\end{picture}
\caption{Lemma \ref{lemaC}} \label{fig4}
\end{figure}

\proof Let $\q=\q_1\sqcup \q_2\colon [0, \ell_1+\ell_2]\to X$. For
every $[t_1,t_2]\subset [0,\ell_1 +\ell_2]$ we have
$$
\dist(\q(t_1), \q(t_2))\leq L(t_2-t_1 )+2C\,
$$ by the triangular inequality. This implies
$\dist(\q(t_1),\q(t_2))\le (L+\varepsilon)(t_2-t_1)+K$, for $K
\geq 2C$. We need to prove that for some well chosen $K$ we have
\begin{equation}
\frac{1}{L+\varepsilon}(t_2-t_1)-K\le\dist(\q(t_1),\q(t_2))\, .\label{ineq1}\end{equation}

We consider the constant $D$ given by Lemma \ref{lemaB} and set
$D' = 2L^2(D +C)+C$ and $K=2D'$. The hypothesis $\dist(a,b_2)\geq
D'$ implies that $\ell_2 \geq 2L (D+C)$.

Let $[t_1,t_2]\subset [0,\ell_1 +\ell_2]$. If $t_2-t_1$ is smaller
than $2L(D+C)$ then (\ref{ineq1}) obviously holds. Suppose that
$t_2-t_1\geq 2L(D+C)$. If $[t_1,t_2] \cap [\ell_1 ,\ell_1 +\ell_2
]$ is an interval of length at least $L(D+C)$ then the distance
between $\mathfrak{q}(\ell_1)$ and $\q(t_2)$ is bigger than $D$.
Lemma \ref{lemaB} implies (\ref{ineq1}).

The same inequality is true if $(t_1,t_2)$ does not contain
$\ell_1$. Suppose that $[t_1,t_2] \cap [\ell_1 ,\ell_1 +\ell_2 ]$
is a nontrivial interval of length at most $L(D+C)$. Then
 $$\dist(\q(t_1), \q(t_2)) \geq \dist(\q(t_1),
 \q(\ell_1))-\dist(\q(t_2),
 \q(\ell_1))\geq \frac{1}{L}(\ell_1 -t_1)-D'
 \geq \frac{1}{L}(t_2 -t_1)-2D'$$ and (\ref{ineq1}) holds.
\endproof

\begin{definition}[saturations]\label{satt}
For every $(L,C)$-quasi-geodesic $\q$ in $X$ we define the
\textit{saturation} $\Sat(\q)$ as the union of $\q$ and all
$A\in\aaa$ with $\nn_M(A)\cap \q\ne\emptyset$.
\end{definition}

\begin{lemma}\label{lsat}
Let $\q_n$ be a sequence of $(L,C)$-quasi-geodesics in $X$. In
every asymptotic cone $\co{X;e,d}$ if the limit $\lio{\Sat\left(
\q_n\right)}$ is not empty then it is either a piece $\lio{A_n}$
from $\aaa_\omega$, or the union of ${\mathfrak p}=\lio{\q_n}$
and a collection of pieces from $\aaa_\omega$ such that each piece
intersects $\lio{\q_n}$ in at least one point and all pieces from
$\aaa_\omega$ that intersect $\lio{\q_n}$ in a non-trivial sub-arc
are in the collection (recall that by Corollary \ref{strconv} if a
piece in a tree-graded space intersects an arc in more than two
points then it intersects the arc by a sub-arc).
\end{lemma}

\proof {\textbf{Case I.}}\quad Suppose that $\lm_\omega
\frac{\dist(e_n,\q_n)}{d_n} < \infty$. Let $u_n\in \q_n$ be an
almost closest point to $e_n$ in $\q_n$.

Suppose that a piece $A =\lio{A_n}$ intersects $\q=\lio{\q_n}$ in
an arc $\q([t_1,t_2])$, $t_1<t_2$. This arc is a limit of
sub-quasi-geodesics $\q_n'$ of $\q_n$ defined on intervals of
length $(t_2-t_1)d_n$. The ends of $\q_n'$ are at distance
$o(d_n)$ from $A_n$ $\omega$-almost surely. Lemma \ref{Aprime}
implies that $\omega$-almost surely $A_n\subseteq \Sat(\q_n)$
since $\diam(\nn_M(A_n)\cap \q_n)=O(d_n)$.

Suppose $A$ is such that $A_n\subseteq \Sat(\q_n)$ and $\lm_\omega
\frac{\dist(e_n,A_n)}{d_n}<\infty$. Let $a_n$ be an almost nearest
point to $u_n$ in $\q_n\cap \nn_{M}(A_n)$. Lemma \ref{qqc} implies
that the sub-arc $\q_n'$ of $\q_n$ with endpoints $u_n$ and $a_n$
is contained $\omega$-almost surely in $\nn_{t_n}(A_n)$ for some
number $t_n=O(d_n)$. If
$\lm_\omega\frac{\dist(u_n,a_n)}{d_n}=\infty$ then by applying
Lemma \ref{Aprime} we obtain ($\omega$-almost surely) a point in
$\q_n\cap \nn_{M}(A_n)$ nearer to $u_n$ than $a_n$ by a distance
$O(d_n)$, a contradiction. Hence $\lm_\omega
\dist(u_n,a_n)/d_n<\infty$. Then $a=\lio{a_n}$ exists and is an
intersection point of $A$ with $\q$.

\me

\noindent {\textbf{Case II.}}\quad Suppose that $\lm_\omega
\frac{\dist(e_n,\q_n)}{d_n}= \infty$. Let $A_n\subset \Sat(\q_n)$
be such that $\lm_\omega \frac{\dist(e_n,A_n)}{d_n}< \infty$. We
have $A=\lio{A_n} \subseteq \lio{\Sat(\q_n)}$. Suppose there
exists $B =\lio{B_n}\subset \lio{\Sat(\q_n)}$ with $B\neq A$
whence $B_n\neq A_n$ $\omega$-almost surely.

For every $n\ge 1$, let $y_n$ be an almost closest to $e_n$ point
in $A_n$. Also pick $b_n=\q_n(t_n)\in\nn_M(B_n)$. If
$\dist(t_n,\q_n\iv(\nn_M(A_n)))=0$ then we set $s_n=t_n$.
Otherwise let $s_n$ be the almost closest to $t_n$ number in
$\q_n\iv(\nn_M(A_n))$. We assume that $s_n\le t_n$ otherwise we
can reverse the orientation of $\q_n$. Then the diameter of the
intersection of $\q_n([s_n,t_n])$ with $\nn_M(A_n)$ is bounded in
terms of $L, C$. By Lemma \ref{lemaC}, $\rrr_n=[y_n,\q_n(s_n)]$
$\cup \q_n([s_n,t_n])$ is an $(L+\varepsilon,K)$-quasi-geodesic
where $[y_n,\q_n(s_n)]$ is any geodesic connecting $y_n$ and
$\q_n(s_n)$ in $X$.

Notice that $\dist(y_n,B_n)\le O(d_n)$, $\q_n(t_n)\in B_n$. Then
by Lemma \ref{qqc}, $\rrr_n\subseteq \nn_{O(d_n)}(B_n)$
$\omega$-almost surely. Applying Lemma \ref{Aprime} we find $y_n',
a_n'$ in $[y_n,\q_n(s_n)]$ with $\dist(y_n',a_n')=O(d_n)$ which
belong to both $\nn_M(A_n)$ and $\nn_M(B_n)$. This contradicts
property $(\alpha_1)$.

Thus we can conclude that there is no sequence $B_n\subset
\Sat(\q_n)$ with $B_n\neq A_n$ $\omega$-almost surely, such that
$\lm_\omega \frac{\dist(e_n,B_n)}{d_n}< \infty$. Hence in this
case $\lio{\Sat(\q_n)}=A$.\endproof

\begin{lemma}\label{diam}
For every $d>0$,every $(L,C)$-quasi-geodesic $\q$ and every
$A\in\aaa$, $\nn_M(A)\cap\q=\emptyset$, the diameter of
$\nn_d(A)\cap \nn_d(\Sat(\q))$ is bounded in terms of $d, L, C$.
\end{lemma}

\proof Suppose that for some $d>0$ and some $(L,C)$ there exist
sequences of $(L,C)$-quasi-geodesics $\q_n$, of sets $A_n\in\aaa$,
$A_n \not\subset \Sat(\q_n)$, and of points $x_n, y_n\in
\nn_d(A_n)\cap \nn_d(\Sat(\q_n))$ such that the sequence
$\dist(x_n,y_n)=p_n$ is unbounded. Consider the corresponding
asymptotic cone $\co{X;(x_n),(p_n)}$.  The limit sets $\lio{A_n}$
and $\lio{\Sat(\q_n)}$ contain points $x=\lio{x_n}$ and
$y=\lio{y_n}$ in common, $\dist(x,y)=1$. By Lemma \ref{lsat},
either $\lio{\Sat(\q_n)}$ is $\lio{A_n'}$ with $A_n'\in\aaa$,
$A_n'\ne A_n$ $\omega$-almost surely, or $\lio{\Sat(\q_n)}$ is
equal to $Y(\q)$ where $\q$ is the arc $\lio{\q_n}$, and
$\lio{A_n}\not\subset \lio{\Sat(\q_n)}$. In the first case we get
a contradiction with property ($T_1$) for $\aaa$. In the second
case we get a contradiction with Lemma \ref{retr}, part (2).
\endproof

\begin{lemma}[uniform variant of Lemma \ref{eproj} for saturations]\label{leprojun}
For every $x\in X$ and every $(L,C)$-quasi-geodesic $\q$ in $X$
with $\dist(x,\Sat(\q))=2d$,
$$\diam(\proj_{\Sat(\q)}(\nn_d(x))=o(d).$$

\end{lemma}

\proof By contradiction, suppose that there exists a sequence of
quasi-geodesics $\q_n$ and points $x_n$ with $\lm_\omega
\dist(x_n,\Sat(\q_n))=2d_n$ such that $\lm_\omega d_n=\infty$, and
the  almost projection of $\nn_{d_n}(x_n)$ on $\Sat(\q_n)$ has
diameter at least $td_n$ for some fixed $t$. In the asymptotic
cone $\co{X, (x_n), (d_n)}$ we have, according to Lemma
\ref{lsat}, that $\lio{\Sat(\q_n)}$ is either one piece or a set
of type $Y$. We apply Lemma \ref{retr}, part (2), and get a
contradiction.\endproof

\begin{lemma}[uniform property $(\alpha_2')$ for saturations]\label{ptgi2u}
For every $\lambda \geq 1$, $\kappa \geq 0$ and $\theta \in \left[
0, \frac{1}{2} \right) $ there exists $R$ such that for every
$(\lambda ,\kappa )$-quasi-geodesic $\cf:[0,\ell ]\to X $ joining
two points in $\nn_{\theta \ell /L} \left(\Sat(\q) \right)$, where
$\q$ is a quasi-geodesic, we have $\cf([0,\ell ])\cap \nn_{R}
(\Sat(\q))\neq \emptyset$ (in particular, the constant $R$ does
not depend on $\q$).
\end{lemma}

\proof One can simply repeat the argument of Lemma \ref{T1T2a2}
but use Lemma \ref{leprojun} instead of Lemma
\ref{eproj}.\endproof

\begin{lemma}[uniform quasi-convexity of saturations]\label{lqqcun}
For every $\lambda \geq 1$, $\kappa \geq 0$, there exists $\tau $
such that for every $R\ge1$, for every quasi-geodesic $\q$, the
saturation $\Sat(\q)$ has the property that every $(\lambda , \kappa
)$-quasi-geodesic $\cf$ joining two points in its $R$-tubular
neighborhood is entirely contained in its $\tau R$-tubular
neighborhood.
\end{lemma}

\proof By Remark \ref{buragos}, it is enough to prove the
statement for $(\lambda,\kappa)$-almost geodesics $\cf$. Suppose
there exists a sequence of quasi-geodesics $\q_n$, a sequence of
numbers $R_n\ge 1$, a sequence ${\mathfrak c}_n$ of $(\lambda ,
\kappa)$-almost geodesics joining the points $x_n, y_n$ in the
$R_n$-tubular neighborhood of $\Sat(\q_n)$ such that ${\mathfrak
c}_n$ is not contained in the $nR_n$-tubular neighborhood of
$\Sat(\q_n)$.

Let $z_n \in {\mathfrak c}_n$ be such that $d_n=\dist(z_n,
\Sat(\q_n))$ is maximal. By Lemma \ref{lsat}, in the asymptotic
cone $\co{X;(z_n),(d_n)}$, we have that $S=\lio{\Sat(\q_n)}$ is
either one piece or a set $Y(\q)$ of type $Y$. On the other hand
by Remark \ref{grr05} $\lio{{\mathfrak c}_n}$ is either a
topological arc with endpoints in $S$ and not contained in it, or
a bi-Lipschitz ray with origin in $S$ or a bi-Lipschitz line. In
addition, $\lio{c_n}$ is contained in $\nn_1(S)$. In all three
cases Lemma \ref{retr}, part (2), and Corollary \ref{projpath}
give a contradiction.\endproof

\me

\begin{lemma}[saturations of  polygonal lines]\label{limrsat}
Let $X$ be a geodesic metric space. Then the following is true for
every $k\ge 1$.

\begin{itemize}
\item[(1)] For every $n\ge 1$, let $\bigcup_{i=1}^k \q_i^{(n)}$ be
a polygonal line composed of $(L,C)$-quasi-geodesics $\q_i^{(n)}$.
Then in every asymptotic cone the limit set $\lio{\bigcup_{i=1}^k
\Sat(\q^{(n)}_i)}=\bigcup_{i=1}^k \lio{\Sat(\q^{(n)}_i)}$ is
either a piece or a connected union of sets of type $Y$ (as in
Lemma \ref{retr}, part (3)).

\item[(2)] The results in Lemmas \ref{leprojun}, \ref{ptgi2u}, \ref{lqqcun} are true
if we replace $\Sat(\q)$ with $\bigcup_{i=1}^k
\Sat\left(\q_i\right)$, where $\bigcup_{i=1}^k \q_i$ is a
polygonal line composed of $(L,C)$-quasi-geodesics.

\item[(3)] For every $\delta >0$, for every polygonal line
$\bigcup_{i=1}^k \q_i$ composed of $(L,C)$-quasi-geodesics, and
every $A\in \aaa$ such that $A\not\subset\bigcup_{i=1}^k \Sat
\left( \q_i \right)$, the intersection $\nn_\delta(A)\cap
\nn_\delta \left( \bigcup_{i=1}^k \Sat \left( \q_i \right)
\right)$ has a uniformly bounded diameter in terms of $A$,
$\q_1,\ldots,\q_k$.
\end{itemize}
\end{lemma}

\proof We prove simultaneously (1), (2) and (3) by induction on
$k$. For $k=1$ all three statements are true. Suppose they are
true for $i\leq k$. We prove them for $k+1$. We note that (1)
implies (2) in the same way as Lemma \ref{lsat} implies the cited
Lemmas, and the implication (1) $\Rightarrow$ (3) follows from
Lemma \ref{retr}, part (3) (the argument is essentially the same
as in Lemma \ref{diam}). Thus it is enough to prove part (1).

Let $\co{X; e,d}$ be an asymptotic cone. We suppose that
$$\lm_\omega \frac{\dist\left(e_n , \bigcup_{i=1}^{k+1}
\Sat\left({q}^{(n)}_i\right) \right)}{d_n} < \infty $$ (otherwise
the $\omega$-limit is empty). There are two possible situations.

\noindent {\bf Case I.} Suppose that there exists an integer $i$
between $2$ and $k$ such that $$\lm_\omega\frac{\dist\left(e_n,
\Sat\left(\q^{(n)}_i\right) \right)}{d_n} < \infty.$$ By the
inductive hypothesis $\lio{\bigcup_{j=1}^i
\Sat\left(\q^{(n)}_j\right)}$ is a set of type $Y$, and so is the
set $$\lio{\bigcup_{j=i}^{k+1}\Sat\left(\q^{(n)}_j\right)}.$$
These two sets have a common non-empty subset $\lio{
\Sat\left(\q^{(n)}_i\right)}$. Since a connected union of two sets
of type $Y$ is again a set of type $Y$, statement (1) follows.

\me

\noindent \textbf{Case II.} Suppose that for every $i$ between $2$
and $k$, we have $$\lm_\omega\frac{\dist\left(e_n,
\Sat\left(\q^{(n)}_i\right) \right)}{d_n} = \infty.$$ If the same
is true either for $i=1$ or for $i=k+1$ one can apply Lemma
\ref{lsat}. Thus suppose that for $i=1,k+1$, we have
$$\lm_\omega\frac{\dist\left(e_n, \Sat\left(\q^{(n)}_i\right)
\right)}{d_n} < \infty.$$

By Lemma \ref{lsat}, for $i=1, k+1$, for the limit set
$\lio{\Sat(\q^{(n)}_i)}$ one of the following two possibilities
occurs:

\me

(\textbf{A}$_i$) it is equal to $\lio{A_n}$, where $A_n\in \aaa,
A_n \subseteq\Sat(\q^{(n)}_i)$;

\me

(\textbf{B}$_i$) it is equal to $Y\left(\q_i\right)$ as in Lemma
\ref{retr}, part (2), where $\q_i=\lio{\q^{(n)}_i}$.

\me

It remains to show that the union $\lio{\Sat(\q^{(n)}_1)}\cup
\lio{\Sat(\q^{(n)}_{k+1})}$ is connected.

Suppose that we are in the situation (\textbf{B}$_1$). Let $u_n\in
\q^{(n)}_1$ be an almost nearest point from $e_n$. Then
$\dist(u_n,e_n)=O(d_n)$. Let $v_n\in \bigcup_{j=2}^{k+1}
\Sat\left(\q^{(n)}_j\right)$ be an almost nearest point to $e_n$.
By our assumption, $\omega$-almost surely $v_n \in
\Sat\left(\q_{k+1}^{(n)}\right)$ and $\dist(v_n,e_n)=O(d_n)$.
Hence $\dist(u_n,v_n)=O(d_n)$. Let $R_k$ be the constant given by
the variant of Lemma \ref{ptgi2u} for polygonal lines composed of
$k$ $(L,C)$-quasi-geodesics with $(\lambda,\kappa)=(L,C)$,
$\theta=\frac13$ (that $R_k$ exists by the induction hypothesis).
Let $a_n$ be an almost nearest point from $u_n$ in $\q^{(n)}_1
\cap \nn_{R_k}\left(\bigcup_{j=2}^{k+1}
\Sat\left(\q^{(n)}_j\right)\right)$. Let ${\mathfrak p}^{(n)} $ be
the sub-quasi-geodesic of $\q^{(n)}_1$ with endpoints $u_n$ and
$a_n$. According to the part (2) of the proposition (which by the
induction assumption is true for $k$), ${\mathfrak p}^{(n)}
\subset \nn_{td_n }\left( \bigcup_{j=2}^{k+1}
\Sat(\q^{(n)}_j)\right)$ for some $t$ independent on $n$. If
$\dist(u_n,a_n)\gg d_n$ then according to Lemma \ref{ptgi2u} there
exists another point on ${\mathfrak p}^{(n)} \cap \nn_{R_k}\left(
\bigcup_{j=2}^{k+1} \Sat\left(\q^{(n)}_j\right)\right)$ whose
distance from $u_n$ is smaller than $\dist(a_n,u_n)$ by $O(d_n)$,
a contradiction. Therefore $\dist(u_n,a_n)\le O(d_n)$ and the
limit point $\lio{a_n}$ is a common point of $\q_1$ and
$\lio{\bigcup_{i=2}^{k+1}\Sat\left(\q^{(n)}_i\right)}=\lio{\Sat\left(\q^{(n)}_{k+1}\right)}$.

The same argument works if we are in the situation
(\textbf{B}$_{k+1})$. Therefore we suppose that we are in the
situations (\textbf{A}$_{1}$) and (\textbf{A}$_{k+1}$). We have
that $\lio{\Sat\left(\q^{(n)}_i\right)}, i=1,k+1,$ is equal to
$\lio{A^{(n)}_i}$, where $A^{(n)}_i\in \aaa, A^{(n)}_i \subseteq
\Sat\left(\q^{(n)}_i\right)$. Suppose that $A^{(n)}_1 \neq
A^{(n)}_{k+1}$ $\omega$-almost surely.  Let $v^{(n)}_i\in
\Sat\left(\q^{(n)}_i\right)$ be an almost nearest point from
$e_n$. By hypothesis $v^{(n)}_i \in A^{(n)}_i$.

The two assumptions: $$\lm_\omega \frac{\dist (e_n, \Sat (\q
_i^{(n)}))}{d_n}=\infty,$$ $i\in \{ 2,\dots ,k\}$, and
$$\lio{\Sat (\q_{k+1}^{(n)})}=\lio{
A_{k+1}^{(n)}}$$ imply that $A_1^{(n)}\not \subset
\bigcup_{i=2}^{k+1}\Sat(\q_i^{(n)})$ $\omega$-almost surely.

Suppose that $[0,\ell_1^{(n)}]$ is the domain of $\q_1^{(n)}$. The
following two cases may occur.

\noindent {\bf Case I.} If the distance from $\ell_1^{(n)}$ to the
pre-image $(\q_1^{(n)})^{-1}(A_1^{(n)})$ is at most $LC+1$ then we
denote $\q_1^{(n)}(\ell_1^{(n)})$ by $a_n$. We have that $\dist
(a_n, \q_1^{(n)}\cap A_1^{(n)})\leq L^2C+L+C$, which implies by
Lemma \ref{qqc} that a geodesic $\pgot_n=[v_1^{(n)},a_n]$ is
contained in the $t(L^2C+L+C)$-tubular neighborhood of
$A_1^{(n)}.$

\noindent {\bf Case II.} If the distance from $\ell_1^{(n)}$ to
$(\q_1^{(n)})^{-1}(A_1^{(n)})$ is larger than $LC+1$, then we
consider $t_n\in [0,\ell_1^{(n)}]$ at distance $LC+1$ of
$(\q_1^{(n)})^{-1}(A_1^{(n)})$ such that all points in
$[t_n,\ell_1^{(n)}]$ are at distance at least $LC+1$ of
$(\q_1^{(n)})^{-1}(A_1^{(n)})$. We denote by $a_n$ the point
$\q_1^{(n)}(t_n)$. According to Lemma \ref{qqc} we have that a
geodesic $[v_1^{(n)},a_n]$ is contained in the
$t(L^2C+L+C)$-tubular neighborhood of $A_1^{(n)}.$

By our assumption, $\lm_\omega\frac{\dist
(v_1^{(n)},a_n)}{d_n}=\infty$. Lemma \ref{lemaC} implies that
$[v_1^{(n)},a_n]$ and the restriction of $\q_1^{(n)}$ to
$[t_n,\ell_1^{(n)}]$ form an $(L+\varepsilon , K)$-quasi-geodesic
$\omega$-almost surely. We denote it by $\pgot_n$.

\medskip

Both in Case I and in Case II we have obtained an $(L+\varepsilon
, K)$-quasi-geodesic $\pgot_n$ with one of the endpoints
$v_1^{(n)}$ and the other one contained in $\q_2^{(n)}$. The
distance from $v_1^{(n)}$ to $\bigcup_{i=2}^{k+1}\Sat
(\q_i^{(n)})$ does not exceed $\dist (v_1^{(n)},v_{k+1}^{(n)})$,
hence it is at most $O(d_n)$. It follows that $\pgot_n \subset
\nn_{O(d_n)}\left( \bigcup_{i=2}^{k+1}\Sat (\q_i^{(n)}) \right)$.
In particular $[v_1^{(n)},a_n]$ is contained in the same tubular
neighborhood. Since the length $\lambda_n$ of $[v^{(n)}_1,a_n]$
satisfies $\lm_\omega \frac{\lambda_n}{d_n}=\infty $, by applying
Lemmas \ref{ptgi2u} and \ref{lqqcun} we obtain that a sub-segment
$[\alpha_n,\beta_n]$ of $[v_n^1,a_n]$ of length
$\frac{\lambda_n}{2}$ is contained in $\nn_{\tau R}\left(
\bigcup_{i=2}^{k+1} \Sat\left(\q^i_n\right) \right)$, where $R$ is
an universal constant. On the other hand we have
$[\alpha_n,\beta_n]\subset \nn_{t(L^2C+L+C)} (A_1^{(n)})$. This
contradicts the inductive hypothesis (3). We conclude that if we
are in situation (\textbf{A}$_1$) then $\lm_\omega
\frac{\dist\left( e_n\, ,\,
\Sat\left(\q^{k+1}_n\right)
\right)}{d_n}=\infty $.\endproof

\me

\begin{cor}
Let $\Delta $ be a quasi-geodesic triangle. Then every edge
$\mathfrak a $ of $\Delta $ is contained in an $M$-tubular
neighborhood of $Sat ({\mathfrak b}) \cup Sat ({\mathfrak c})$,
where ${\mathfrak b}$ and ${\mathfrak c}$ are the two other edges
of $\Delta $ and $M$ is an universal constant.
\end{cor}

\begin{lemma}\label{gat}
For every $R>0, k\in \N $ and $ \delta >0$ there exists $\varkappa
>0$ such that if $\bigcup_{i=1}^k \q_i$ is a polygonal line composed of
quasi-geodesics and $A,B\in \aaa, A\cup B\subset \bigcup_{i=1}^k
\Sat(\q_i), A\neq B$, the following holds. Let $a\in \nn_R(A)$ and
$b\in \nn_R(B)$ be two points that can be joined by a
quasi-geodesic ${\mathfrak p}$ such that ${\mathfrak p} \cap
\nn_R(A)$ and ${\mathfrak p} \cap \nn_R(B)$ has diameter at most
$\delta$. Then $\{ a,b \}\subset \nn_{\varkappa} \left(
\bigcup_{i=1}^k \q_i \right) $.
\end{lemma}

\proof Suppose $\q_i$ is defined on the interval $[0,\ell_i]$. Let
$\mathfrak{r} :[0,\sum_{i=1}^k \ell_i] \to X$ be the map defined
by $\mathfrak{r}(\sum_{i=1}^{j-1} \ell_i+t)=\q_j(t)$, for all
$t\in [0,\ell_j]$ and all $j\in \{ 2, \dots , k \}$. It satisfies
\begin{equation}\label{eq454}
\dist (\mathfrak{r}(t), \mathfrak{r}(s))\leq L|t-s|+kC \, .
\end{equation}

Let $x$ be a point in $\mathfrak{r}\cap \nn_M(B)$ and $t_x\in
[0,\sum_{i=1}^k \ell_i]$ such that $\mathfrak{r}(t_x)=x$. We have
two cases.

\medskip

\noindent (a) If the distance from $t_x$ to the pre-image
${\mathfrak r}\iv(\nn_M(A))$ does not exceed $LC+1$ then $x\in
\nn_{M+L^2C+L+kC}(A)$ by (\ref{eq454}). By Lemma \ref{lemaC}, if
$\dist (a,x)$ is larger than $D'$ then the union of $\pgot$ and a
geodesic $[a,x]$ form an $(L+\varepsilon , K)$-quasi-geodesic,
with endpoints in $\nn_{R+M}(B)$. It follows that this
quasi-geodesic and in particular $[a,x]$ are contained in
$\nn_{t(M+R)}(B)$. On the other hand $[a,x]$ is contained in
$\nn_{t(M+R+L^2C+L+kC)}(A)$. If $\dist (a,x)$ is larger than the
diameter given by $(\alpha_1)$ for $\delta = t(M+R+L^2C+L+kC)$
then we obtain a contradiction with $(\alpha_1)$.

\medskip

\noindent (a) Suppose that the distance from $t_x$ to ${\mathfrak
r}\iv(\nn_M (A))$ is larger than $LC+1$. Consider $s_0$ at
distance $LC+1$ from $\mathfrak{r}^{-1}(\nn_M (A))$ such that
every $s$ between $s_0$ and $t_x$ is at distance at least $LC+1$
from $\mathfrak{r}^{-1}(\nn_M (A))$. It follows that
$\mathfrak{r}([s_0,t_x])$ or $\mathfrak{r}([t_x,s_0])$ is
disjoint of $\nn_M(A)$. Let $y=\mathfrak{r}(s_0)$. The restriction
$\mathfrak{r}'$ of $\mathfrak{r}$ to $[s_0, t_x]$ or $[t_x,s_0]$
can be written as $\bigcup_{j=1}^m \q_j'$, where $m\le k$ and each
$\q_j'$ coincides with one of the $\q_i$'s or a restriction of
it. We note that $A\not\subset \Sat (\mathfrak{r}')$.

If the distance from $a$ to $y$ is larger than the constant $D'$
given by Lemma \ref{lemaC} then $\pgot$ and a geodesic $[a,y]$
form an $(L+\varepsilon , K)$-quasi-geodesic. Lemma \ref{limrsat},
part (2), implies that this quasi-geodesic, and in particular
$[a,y]$, is contained in the $\tau R$-tubular neighborhood of
$\Sat (\mathfrak{r}')$. On the other hand, $[a,y]$ is contained in
the $t (R+ M+L^2C+L+kC)$-tubular neighborhood of $A$. For $\dist
(a,y)$ larger than the diameter given by Lemma \ref{limrsat}, (3),
for $\delta = \max \left ( t (R+ M+L^2C+L+kC), \tau R\right)$ we
obtain a contradiction.
\endproof

\begin{lemma}\label{pi3g}
Suppose that a metric space  $X$ is asymptotically tree-graded with
respect to $\aaa$. Then $X$ satisfies $(\alpha_3')$.
\end{lemma}

\proof Let $k\geq 2$, $\sigma \geq 1$ and $\nu \geq 4\sigma$. Fix
a sufficiently large number $\vartheta$ (it will be clear later in
the proof how large $\vartheta$ should be). Let $P$ be a $k$-gon
with quasi-geodesic edges that is $(\vartheta \, ,\,  \sigma \,
,\, \nu ) $-fat. Changing if necessary the polygon by a finite
Hausdorff distance, we may suppose that its edges are
$(L+C,C)$-almost geodesics.

Let $\q:[0,\ell ]\to X$ be an edge with endpoints $\q (0)=x,\q
(\ell )=y$. We denote $\q_1,\q_2,\dots ,\q_{k-1}$ the other edges
in the clockwise order. By Lemma \ref{limrsat}, part (2),
$$
\q\subset \nn_{\tau R}\left( \bigcup_{i=1}^{k-1}
\Sat\left(\q_i\right) \right).$$ We take $\vartheta > \tau R$.
Then for every point $z\in \q\setminus {\nn}_{\sigma \vartheta }(
\{ x,y \})$ there exists $A\subset \Sat\left(\q_i\right), i\in \{
1,2,\dots ,k-1 \}$ such that $z\in \nn_{\tau R}(A)$. If such a
point $z$ is contained in $\nn_{\tau R}(A)\cap \nn_{\tau R}(B),\,
A\neq B,$ then Lemma \ref{gat} implies that $z\in \nn_{\varkappa}
(\bigcup_{i=1}^{k-1} \q_i)$, where $\varkappa$ depends on $\tau R$
and $k$. If we choose $\vartheta
>\varkappa$ then this gives a contradiction.

Let $t_\q$ be the supremum of the numbers $t\in [0,\ell ]$
contained in $\q^{-1}\left( \nn_{\sigma \vartheta } (x)\right)$.
Let $s_\q$ be the infimum of the numbers in $[t_\q, \ell ]$
contained in $\q^{-1}\left( \nn_{\sigma \vartheta }(y)\right)$.
Let $a_\q =\q (t_\q )$ and $b_\q = \q (s_\q )$. We note that
$\dist (a_\q, x)=\sigma \vartheta $ and $\dist (b_\q, y)=\sigma
\vartheta $. According to the argument in the paragraph above, $\q
([t_\q , s_\q ])$ is covered by the family of open sets $\nn_{\tau
R}(A)$, with $A\subset \Sat\left(\q_i\right), i\in \{ 1,2,\dots
,k-1 \}$, and the traces of these sets on $\q ([t_\q , s_\q ])$
are pairwise disjoint. The connectedness of $\q ([t_\q , s_\q ])$
implies that there exists $A$ as above such that $\q ([t_\q , s_\q
])\subset \nn_{\tau R}(A)$.

Thus, for every edge $\q$ a sub-arc $\q' :[t_\q, s_\q ]\to X$ with
endpoints $a_\q,b_\q$ is contained in $\nn_{\tau R}(A)$ for some
$A\subset \Sat\left(\q_i\right), i\in \{ 1,2,\dots ,k-1 \}$ ($A$
may depend on $\q$). We note that $t_\q$ and $\ell -s_\q$ are less
than $\sigma \vartheta L+C$, hence $\q |_{[0,t_\q ]}\in
\nn_{\sigma\vartheta L^2+LC+C}(a_\q)$ and $\q |_{[s_\q,\ell]}\in
\nn_{\sigma \vartheta L^2+LC+C}(b_\q)$.

Suppose that we have two consecutive edges $\q_1,\q_2$ with
endpoints $x,y$ and $y,z$ respectively, such that $\q_1' \subset
\nn_{\tau R } (A)$ and $\q_2' \subset \nn_{\tau R }(B)$, $A\neq
B$. We denote $\q_3,\q_4,\dots ,\q_{k}$ the other edges in the
clockwise order. We have $\q_i':[t_{\q_i}, s_{\q_i}]\to X$ with
endpoints $a_{\q_i}, b_{\q_i}$. Suppose $b_{\q_1}=\q_1' \cap
\nn_{\sigma \vartheta }(y)$ and $a_{\q_2}=\q_2' \cap \nn_{\sigma
\vartheta }(y)$.

Let $\bar\q_1$ be the restriction of $\q_1'$ to $[t_{\q_1},
t_{\q_1}+3L\tau R ]$ and $\tilde\q_1=[x, a_{\q_1}]\cup \bar\q_1$.
We note that since $\dist (a_{\q_1} , b_{\q_1})\geq \dist
(x,y)-2\sigma\vartheta \geq \nu\vartheta -2\sigma\vartheta \geq
2\sigma\vartheta$, we have $s_{\q_1} -t_{\q_1}\geq
\frac{2\sigma\vartheta }{L}-C$, so for $\vartheta $ large enough
we have $s_{\q_1} -t_{\q_1} \geq 10 L\tau R $ and the restriction
$\bar\q_1$ makes sense.

Likewise we construct $\tilde\q_2 = \bar\q_2 \cup [b_{\q_2}, z]$,
where $\bar\q_2$ is the restriction of $\q_2'$ to the last
sub-interval of length $3L\tau R $.

Let $[a,b]$ be a geodesic joining the points $a=a_{\q_2}$ and
$b=b_{\q_1}$. It has length at most $2\sigma \vartheta$. Let
$[a',b']\subset [a,b]$ be a sub-geodesic which intersects
$\nn_{\tau R} (A)$ in $a'$ and $\nn_{\tau R} (B)$ in $b'$
(eventually reduced to a point). Notice that $A\subseteq
\Sat(\tilde\q_1)$, $B\subseteq \Sat(\tilde\q_2)$. Lemma \ref{gat}
applied to the polygonal line $\tilde\q_2\cup \bigcup_{i=3}^k \q_i
\cup \tilde\q_1$ and to the points $a',b'$ implies that $\{
a',b'\} \subset \nn_{\varkappa} \left( \tilde\q_2\cup
\bigcup_{i=3}^k \q_i \cup \tilde\q_1\right)$, where $\varkappa $
depends on $\tau R$. Since $\dist( y, \{ a',b' \} )$ is at most
$2\sigma \vartheta$, it follows that $y\in \nn_{\varkappa +2\sigma
\vartheta } \left( \tilde\q_2\cup \bigcup_{i=3}^k \q_i \cup
\tilde\q_1\right)\subset \nn_{\varkappa + 3\sigma \vartheta
+3L^2\tau R +C} \left(\bigcup_{i=3}^k \q_i \right)$. On the other
hand property ($F_2$) implies that $\dist (y, \bigcup_{i=3}^k
\q_i)\geq \nu \vartheta \geq 4\sigma \vartheta $. For $\vartheta$
large enough this gives a contradiction.

We conclude that there exists $A\in \aaa$ such that
$\bigcup_{i=1}^k \q_i' \subset {\nn}_{\tau R }(A)$. Hence ${P}$ is
inside the $(\tau R +\sigma \vartheta L^2+LC+C)$-tubular
neighborhood of $A$. \endproof

The following corollary immediately follows from the proof of
Theorem \ref{tgi}

\begin{cor}[there is no need to vary the ultrafilter in Definition \ref{asco}]
\label{cor451} Let $X$ be a metric space, let $\aaa$ be a
collection of subsets in $X$. Let $\omega$ be any ultrafilter over
$\N$. Suppose that every asymptotic cone $\co{X;e,d}$ is
tree-graded with respect to the collection of sets $\lio{A_n}$,
$A_n\in \aaa$. Then $X$ is asymptotically tree-graded with respect
to $\aaa$.
\end{cor}

\section{Quasi-isometric behavior}\label{qiv}

One of the main reasons we are interested in the property of being
asymptotically tree-graded is the rigid behavior of this property
with respect to quasi-isometry.

\subsection{Asymptotically tree-graded spaces}
\begin{theorem}[being asymptotically tree-graded is a geometric property.] \label{qi}
Let $X$ be a metric space and let $\aaa$ be a collection of
subsets of $X$. Let $\q$ be a quasi-isometry $X\to X'$. Then:
\begin{itemize}
    \item[(1)] If $X$
satisfies properties $(\alpha_1)$ and $(\alpha_2')$ with respect
to $\aaa$ then $X'$ satisfies properties $(\alpha_1)$ and
$(\alpha_2^\epsilon)$, for a sufficiently small $\epsilon$, with
respect to $\q(\aaa)=\{\q(A)\mid A\in\aaa\}$.
    \item[(2)] If $X$ satisfies $(\alpha_3')$ with respect to $\aaa$
then $X'$ satisfies $(\alpha_3)$ with respect to $\q ( \aaa )$.
    \item[(3)] If $X$ is asymptotically tree-graded with
respect to $\aaa$ then $X'$ is asymptotically tree-graded with
respect to $\q(\aaa)$.
\end{itemize}
\end{theorem}

\proof (1) follows from Theorem \ref{tgi} and Remark \ref{tgi3}.

\me

(2) Assume that $\q$ is an $(L,C)$-quasi-isometry and that $\bar{
\q} :X'\to X$ is an $(L,C)$-quasi-isometry so that $\bar{\q} \circ
\q$ and $\q \circ \bar{\q}$ are at distance at most $C$ from the
respective identity maps.

Fix an arbitrary integer $k\geq 2$. Let $\sigma = 2L^2+1$ and $\nu =
4\sigma $. Property $(\alpha_3')$ in $X$ implies that for the
constants $L,C$ of the quasi-isometries, for the given $k$, $\sigma
$ and $\nu $ there exists $\vartheta _0$ such that for every
$\vartheta \geq \vartheta_0$ a $k$-gon with $(L,C)$-quasi-geodesic
edges in $X$ which is $(\vartheta\, ,\, \sigma \, ,\, \nu)$-fat is
contained in $\nn_\chi (A)$, where $A\in \aaa$ and $\chi =\chi
(L,C,k,\sigma ,\nu , \vartheta)$.

Let $\vartheta_1=\max (\vartheta _0 , 2L^2C+C )$ and let $\theta = L (\vartheta _1 + C)$. Let $P$ be a geodesic $k$-gon in $X'$ which
is $(\theta , 2, \nu
 )$-fat. Then $\bar{\q }(P)$ is a $k$-gon in $X$ with $(L,C)$-quasi-geodesic edges which is
$(\vartheta_1\, ,\, \sigma \, ,\, \nu)$-fat. Consequently,
$\bar{\q }(P) \subset \nn_\chi (A)$, where $A\in \aaa$ and $\chi
=\chi (L,C,k,\sigma ,\nu , \vartheta_1)$. It follows that
$P\subset \nn_C (\q \circ \bar{\q} (P))\subset \nn_{L\chi +2C }(\q
(A)) $.

\me

(3) The statement follows from (1) and (2). It also follows
immediately from the definition of asymptotically tree graded
spaces. Indeed, it is easy to see that $\omega$-limits of sequences
of subsets commute with quasi-isometries. Since quasi-isometric
spaces have bi-Lipschitz equivalent asymptotic cones (Remark
\ref{grr2}) it remains to note that a metric space that is
bi-Lipschitz equivalent to a space that is tree-graded with respect
to $\pp$, is itself tree-graded with respect to the images of the
sets in $\pp$ under the bi-Lipschitz map.\endproof

\begin{definition}\label{wandp}
Let $B$ be a geodesic metric space. We say that $B$ is \textit{wide}
if every asymptotic cone of $B$ does not have global cut-points.

We say that $B$ is \textit{constricted} if every asymptotic cone of
$B$ has a global cut point.

We say that $B$ is \textit{unconstricted} if there exists an
ultrafilter $\omega$ and a sequence $d=(d_n)$ of scaling constants
satisfying $\lim_\omega d_n=\infty$ such that for every observation
point $e=(e_n)^\omega$ the asymptotic cone $\co{B;e,d}$ has no
cut-points.
\end{definition}

\begin{remarks}
\begin{enumerate}
    \item[(1)] Note that ``unconstricted'' is in general more than the negation of
``constricted'', as the latter only means that there exists one
asymptotic cone without cut-points. The two notions coincide for
finitely generated groups, according to the comment following Remark
\ref{grr2}.
    \item[(2)] Note also that most probably ``wide'' is stronger than
    ``unconstricted'', but we do not have an example of an unconstricted group which is not
    wide (see Problem \ref{wideunconstricted}).

\end{enumerate}
\end{remarks}

Definion \ref{wandp} has the following uniform version.

\begin{definition}\label{uwandp}
Let $\mathcal{B}$ be a family of geodesic metric spaces. We say that
$\mathcal{B}$ is \textit{uniformly wide} if for every sequence $B_n$
of metric spaces in $\mathcal{B}$ with metrics $\dist_n$ and
basepoints $b_n\in B_n$, for every ultrafilter $\omega$ and for
every sequence of scaling constants $(d_n)$ with $\lim_\omega
d_n=\infty$, the ultralimit $\lio{B_n,\frac{1}{d_n}\dist_n}_b$ is
without cut-points.

We say that $\mathcal{B}$ is \textit{uniformly unconstricted} if for
every sequence $B_n$ of metric spaces in $\mathcal{B}$ with metrics
$\dist_n$, there exists an ultrafilter $\omega$ and a sequence of
scaling constants $d=(d_n)$ with $\lim_\omega d_n=\infty$ such that
for every sequence of basepoints $b_n\in B_n$, the ultralimit
$\lio{B_n,\frac{1}{d_n}\dist_n}_b$ is without cut-points.
\end{definition}

\begin{remarks}
\begin{itemize}
  \item[(a)] All metric spaces in a family that is uniformly wide (uniformly
  unconstricted) are wide (unconstricted).
  \item[(b)] If $\mathcal{B}$ is a family of wide metric spaces containing only
finitely many pairwise non-isometric spaces then $\mathcal{B}$ is
uniformly wide.
 \item[(c)] For examples of groups that are wide or unconstricted and of
families of groups that are uniformly wide or unconstricted, see
Section \ref{ascp}.
\end{itemize}
\end{remarks}

\begin{proposition}\label{cutp}
Let metric space $X$ be asymptotically tree-graded with respect to a
collection of subsets $\aaa$. Let $\mathcal{B}$ be a family of
metric spaces which is uniformly unconstricted. Suppose that for
some constant $c$, every point in every space $B\in \mathcal{B}$ is
at distance at most $c$ from an infinite geodesic in $B$. Then for
every $(L,C)$ there exists $M=M(L,C, \mathcal{B})$ such that for
every $B\in \mathcal{B}$ and every $(L,C)$-quasi-isometric embedding
$\q\colon B \to X$ there exists $A\in \aaa$ such that $\q(B) \subset
\nn_M (A)$.
\end{proposition}

\proof We argue by contradiction and assume that there is a sequence
of metric spaces $B_n \in \mathcal{B}$ and a sequence of
$(L,C)$-quasi-isometric embeddings $\q_n:B_n \to X$ such that
$\q_n(B_n) \not\subset \nn_n (A)$ for all $A\in \aaa$. By definition
there  exists an ultrafilter $\omega$ and a sequence $d=(d_n)$ with
$\lim_\omega d_n=\infty$ such that for every sequence of basepoints
$b_n\in B_n$, the ultralimit $\lio{B_n,\frac{1}{d_n}\dist_n}_b$ is
without cut-points. Fix a point $b_n\in B_n$. Let $e=(\q_n(b_n))$.
In $\co{X; e, d}$, the limit set $\lio{\q_n(B_n)}$ is a bi-Lipschitz
image of $\lio{B_n,\frac{1}{d_n}\dist_n}_b$, therefore it is without
cut-points. Lemma \ref{cut} implies that
\begin{equation}\label{inclz}
 \lio{\q_n(B_n)}\subset \lio{A_n}\, , \mbox{ where }A_n\in
\aaa\, .
\end{equation}

Consider a sequence $u_n\in B_n$ such that $\lm_\omega
\frac{\dist_n(b_n,u_n)}{d_n}< \infty $. Each $u_n$ is contained in
$\nn_c(\g_n)$, where $\g_n$ is a bi-infinte geodesic in $B_n$.
Suppose that $\g_n$ is parameterized with respect to the arc-length
in $\left( B_n\, ,\, \frac{1}{d_n}\dist_n \right)$ and so that
$\dist_n (u_n, \g_n (0))< c$. The inclusion in (\ref{inclz}) implies
that for every $t\in \R$, $\lm_\omega \frac{\dist_n(\q_n(\g_n (t))\,
,\, A_n)}{d_n}=0$. Therefore for every $s<t$ we have $\omega$-a.s.
that the image by $\q_n$ of the segment $\g_n([s,t])$ contains a
point in $\nn_{M_0} (A_n)$, where $M_0$ is the constant given by
$(\alpha_2')$, for $L$ and $C$. By taking first $s<t<0$ then
$0<s<t$, we may deduce that there exist $\alpha_n <0<\beta_n$ such
that $\q_n(\g_n(\alpha_n))\, ,\, \q_n(\g_n (\beta_n)) \in \nn_{M_0}
(A_n)$. We conclude that $\q_n(\g_n (0)) \in \nn_{\tau M_0} (A_n)$,
by Lemma \ref{qqc}. Hence $\q_n(u_n)\in \nn_{M} (A_n)$
$\omega$-almost surely, where $M=\tau M_0 +Lc+C$.

Let $x_n \in B_n$ be such that $\q_n(x_n)\in \q_n(B_n) \setminus
\nn_n(A_n)$ and let $[b_n,x_n]$ be a geodesic in $B_n$. The previous
argument implies that $\lm_\omega \frac{\dist_n(b_n,x_n)}{d_n} =
\infty $ and that for every $t$ the point $b_n(t)$ on $[b_n,x_n]$ at
distance $td_n$ of $b_n$ has the image by $\q_n$ contained in
$\nn_{M} (A_n)$ $\omega$-almost surely. Let $y_n$ be the farthest
point from $b_n$ in the closure of $[b_n,x_n]\cap \q_n^{-1}( \nn_{M}
(A_n))$. We have that $\lm_\omega \frac{\dist_n\left(b_n,
y_n\right)}{n}=\infty $. Also, $y_n \in [b_n, x_n]\cap
\overline{\q_n^{-1}( \nn_{M} (A_n))}$ implies that for every
$\varepsilon >0$ the distance from $\q_n (y_n)$ to $A_n$ is at most
$M+L \varepsilon +C$. Hence $\q_n (y_n) \in \nn_{M+C+1} (A_n)$. On
the other hand, $b_n \in \nn_{M} (A_n)$ $\omega$-almost surely.
According to Lemma \ref{qqc}, $\q_n([b_n,y_n])\subset \nn_{\tau
(M+C+1)} (A_n)$.
 In $\co{X; (\q_n(y_n)), d}$, $\q=\lio{{\mathfrak
 q}_n([b_n,y_n])}$ is a bi-Lipschitz ray contained in $A
 =\lio{A_n}$ and in $\lio{\q_n (B_n)}$. Since $\lio{\q_n(B_n)}$ is
the image of a bi-Lipschitz embedding of the ultralimit
$\lio{B_n,\frac{1}{d_n}\dist_n}_y$, it is without cut-points,
therefore it is contained in a piece $A'
 =\lio{A_n'}$. Property $(T_1)$ implies that $A=A'$. In
 particular $\lio{\q_n([y_n,x_n])}\subset A $. The same argument
 as before yields that every sequence $v_n\in B_n$ such that $\lm_\omega
\frac{\dist_n(y_n,v_n)}{d_n}< \infty $ satisfies $\q_n(v_n)\in
\nn_{M} (A_n)$ $\omega$-almost surely. Hence, $\dist\left(
\lio{\q_n(y_n)}, \lio{\q_n(x_n)} \right)=\infty$ and there exists
$v_n \in [y_n,x_n]$ such that
 $\dist\left( \lio{\q_n(y_n)}, \lio{\q_n(v_n)} \right)>0$ and $\q_n (v_n) \in \nn_{M}
 (A_n)$, which contradicts the choice of $y_n$.\endproof

\begin{remark} \label{cutprem} The condition that every point is contained in the $c$-tubular
neighborhood of a bi-infinite geodesic is satisfied for instance if
$B$ is a geodesic complete locally compact homogeneous metric space
of infinite diameter. In particular it is true for Cayley graphs of
infinite finitely generated groups.
\end{remark}

\begin{cor}\label{cutp1}
Let $X$ be asymptotically tree-graded with respect to a collection
of subsets $\aaa$. Let $B$ be an unconstricted metric space. Then
every $(L,C)$-quasi-isometric embedding $\q\colon B\to X$ maps $B$
into an $M$-neighborhood of a piece $A\in\aaa$, where $M$ depends
only on $L$, $C$ and $B$.
\end{cor}

\me

\Notat\quad We shall denote the Hausdorff distance between two
sets $A$, $B$ in a metric space by $\hdist(A,B)$.

\subsection{Asymptotically tree-graded groups}

\begin{definition}
We say that a finitely generated group $G$ is
\textit{asymptotically tree-graded with respect to the family of
subgroups} $\{ H_1,H_2,\dots ,H_k\}$ if the Cayley graph $\Cay(G)$
with respect to some (and hence every) finite set of generators is
asymptotically tree-graded with respect to the collection of left
cosets $\{gH_i \mid g\in G \, ,\, i=1,2,\dots ,k\}.$
\end{definition}

\begin{remark}\label{infind}
If $\{ H_1,H_2,\dots ,H_k\}\neq \{ G\}$ and if every $H_i$ is
infinite then every $H_i$ has infinite index in $G$.
\end{remark}

\proof Indeed, a finite index subgroup is at bounded distance of
the whole group, which would contradict ($\alpha_1$).
\endproof

\begin{proposition} \label{propfg} Let $G=\la S\ra$ be a group that is asymptotically tree-graded
with respect to subgroups $H_1,...,H_n$. Then each of the
subgroups $H_i$ is finitely generated.
\end{proposition}

\proof Take $h\in H_i$ and consider a geodesic $\g$ in $\Cay(G,S)$
connecting $1$ and $h$. By Lemma \ref{qqc} there exists a constant
$M>0$ such that $\g$ is in the $M$-tubular neighborhood of $H_i$.
Let $v_1,...,v_k$ be the consecutive vertices of $\g$. For each
$j=1,...,k$ consider a vertex $w_j$ in $H_i$ at distance $\le M$
from $v_j$. Then the distance between $w_j$ and $w_{j+1}$ is at
most $2M+1$, $j=1,...,k-1$. Hence each element $w_j\iv w_{j+1}$
belongs to $H_i$ and is of length at most $2M+1$. Since $h$ is a
product of these elements, we can conclude that $H_j$ is generated
by all its elements of length at most $2M+1$.
\endproof

\begin{remark}
Corollary \ref{cutp1} gives certain restrictions on the groups
that can be quasi-iso\-metrical\-ly embedded into asymptotically
tree-graded groups. For instance, if $G$ is a group asymptotically
tree-graded with respect to a finite family of free Abelian groups
of rank at most $r$, no free Abelian group of rank at least $r+1$
can be quasi-isometrically embedded into $G$.
\end{remark}

\me

\begin{theorem}\label{relhipqi}
Let $X$ be a space that is asymptotically tree-graded with respect
to a collection of subspaces $\mathcal{A}$. Assume that
\begin{itemize}
  \item[(1)] $\mathcal{A}$ is uniformly unconstricted;
  \item[(2)] for some constant $c$ every
point in every $A\in \mathcal{A}$ is at distance at most $c$ from a
bi-infinite geodesic in $A$;
  \item[(3)] For a fixed $x_0\in X$ and every $R>0$ the ball
  $B(x_0,R)$ intersects finitely many $A\in \aaa$.
\end{itemize}

Let $G$ be a finitely generated group which is quasi-isometric to
$X$. Then there exist subsets $A_1,\dots,A_m\in \aaa$ and
subgroups $H_1,\dots ,H_m$ of $G$ such that
\begin{itemize}
  \item[(I)] every $A\in \aaa$ is quasi-isometric to $A_i$ for
  some $i\in \{ 1,2,\dots , m\}$;
  \item[(II)] $H_i$ is quasi-isometric to $A_i$ for
  every $i\in \{ 1,2,\dots , m\}$;
  \item[(III)] $G$ is asymptotically tree-graded with
respect to the family of subgroups $\{ H_1,H_2,\dots ,H_m\}$.
\end{itemize}
\end{theorem}

\proof First we show (in the next lemma) that there is a natural
quasi-transitive quasi-action of $G$ on $X$ by quasi-isometries.

\noindent \textit{Notation}: Let $g\in G$. We denote by
$\mathbf{g}$ the multiplication on the left by $g$ in $G$.

\begin{lemma}\label{qiact}
Let $\q\colon G\to X$ and $\bar{\q}\colon X \to G$ be
$(L_0,C_0)$-quasi-isometries such that $\q \circ \bar{\q}$ is at
distance $C_0$ from the identity map on $X$ and the same is true
for $\bar{\q} \circ \q$ with respect to the identity map on $G$.
\begin{itemize}
  \item [(1)] For every $g\in G$ the map
  $\q_g=\q \circ {\mathbf g}\circ \bar{\q}\colon X\to X$ is an
$(L,C)$-quasi-isometry, where $L=L_0^2$ and $C=L_0C_0+C_0$.
  \item [(2)] For $g, h \in G$ the map $\q_g \circ \q_h$ is at
distance at most $C$ from the map $\q_{gh }$.
  \item [(3)] For every $g\in G$ the map $\q_g \circ
\q_{g^{-1}}$ is at distance at most $C+C_0$ from the identity.
  \item [(4)] For every $x,y\in X$ there exists $g\in G$ such that
  $\dist (x, \q_{g}(y))\leq C_0$.
\end{itemize}

\end{lemma}

\proof Statement (1) follows from the fact that $\mathbf{g}$ acts as
an isometry on $G$. Statement (2) is a consequence of the fact that
$\bar{\q} \circ \q$ is at distance at most $C_0$ from the identity
map on $G$. For (3) we use (2) and the fact that $\q \circ \bar{\q}$
is at distance at most $C_0$ from the identity map on $X$.

(4) Let $g=\bar{\q}(x)$ and $h = \bar{\q}(y) $. Then
$\q_{hg^{-1}}(x)=\q(h)=\q(\bar{\q}(y))$, which is at distance at
most $C_0$ from $y$.
\endproof

\me

\Notat \quad Let $H$ be a subgroup in $G$ and let $x\in X$. We
denote by $Hx$ the set $\{ \q_h (x) \mid h\in H\}$.

\me

Proposition \ref{cutp}, Remark \ref{cutprem} and hypothesis (1)
imply that there exists $M=M(L,C)$ such that for every $A\in \aaa$
and every $(L,C)$-quasi-isometric embedding $\q :A \to X$ there
exists $A'\in \aaa$ such that $\q (A)\subset \nn_M(A')$.

\begin{lemma}\label{neigh}
\begin{itemize}
  \item[(1)] If $A,A'\in \aaa$ satisfy $A\subset \nn_r(A')$ for some $r>0$ then
$A=A'$.
  \item[(2)] Let $\q :X \to X$ and
  $\bar{\q}$ be $(L,C)$-quasi-isometries such that $\q \circ
  \bar{\q}$ and $\bar{\q} \circ \q $ are at distance at most $K
  $ from the identity map on $X$. If $A,A'\in \aaa$ are such that
  $\q(A)\subset \nn_r(A')$ or $A'\subset \nn_r(\q (A))$ for some $r>0$
   then $\q(A)\subset \nn_M(A')$, $\bar{\q}(A')\subset \nn_M(A)$ and
   $\hdist (\q(A),A')\, ,\, \hdist (\bar{\q}(A'),A) \leq LM+C+K
  $.
\end{itemize}
\end{lemma}

\proof (1) follows from property $(\alpha_1)$ and hypothesis (2)
of Theorem \ref{relhipqi}.

(2) Suppose $A'\subset \nn_r(\q (A))$. By Proposition \ref{cutp},
there exists $\bar{A}$ such that $\q(A)\subset \nn_M(\bar{A})$.
Then $A'\subset \nn_{r+M}(\bar{A})$, which implies that
$A'=\bar{A}$. We may therefore reduce the problem to the case when
$\q(A)\subset \nn_r(A')$.

The set $\bar{\q}(A')$ is contained in $\nn_M (A'')$ for some
$A''\in \aaa$. Also $\bar{\q} \circ \q (A)\subset
\nn_{Lr+C}(\bar{\q}(A'))$, which implies that $A\subset
\nn_{Lr+C+M+K}(A'')$. This and (1) imply that $A=A''$. It follows
that $\bar{\q}(A')\subset \nn_M (A)$, which implies that
$A'\subset \nn_{LM+C+K} (\q (A))$.

Proposition \ref{cutp} implies that there exists $\tilde{A}\in A$
such that $\q(A)\subset \nn_M(\tilde{A})$. Hence $A'\subset
\nn_{(L+1)M+C+K} (\tilde{A})$, so $A'=\tilde{A}$. We conclude that
$\q(A)\subset \nn_M(A')$ and $$\hdist (\q(A),A')\, ,\, \hdist
(\bar{\q}(A'),A) \leq LM+C+K .$$\endproof

\me

\Notat \quad We denote the constant $LM+2C+C_0$ by $D$.

\me

\begin{definition}
For every $r>0$ and every $A\in \mathcal{A}$ we define the
$r$-\textit{stabilizer of }$A$ as
$$
\St_{r}(A)=\{g \in G \mid \hdist(\q_g(A),A)\leq r\}\, .
$$
\end{definition}

\begin{cor}\label{stab}
\begin{itemize}
  \item[(a)] For every $g\in G$ and $A,A'\in \aaa$ such that $\q_g(A)\subset
\nn_r (A')$ or $A'\subset \nn_r (\q_g(A))$, where $r>0$, we have
$\hdist (\q_g(A),A')\leq D$.
  \item[(b)] For every $A\in \mathcal{A}$ and for every $r >D$, $\St_{r
}(A)=\St_{D}(A)$. Consequently $\St_{D}(A)$ is a subgroup of $G$.
  \item[(c)] Let $A,B\in \mathcal{A}$ and $g \in G $ be such that
$\hdist(\q_g(A),B)$ is finite. Then $$\St_{D} (A)=g^{-1} \St_{D}
(B)g.$$
\end{itemize}
\end{cor}

\proof Statement (a) is a reformulation in this particular case of
part 2  of Lemma \ref{neigh}, and (b) is a consequence of (a).

(c) For every $r>0$ there exists $R$ large enough so that we have
$\St_r(B)\subset g\St_{R} (A)g^{-1}$.

Applying the previous result again for $g^{-1} , B,A$, together
with (b), we obtain the desired equality.\endproof

\me

Let $\mathcal{F}=\{ A_1,\dots ,A_k\}$ be the collection of all the
sets in $\aaa$ that intersect $B(x_0, M+C_0)$. We show that this
set satisfies (I). Let $A$ be an arbitrary set in $\aaa$ and let
$a\in A$. There exists $g\in G$ such that $\q_g (a)\in
B(x_0,C_0)$, by Lemma \ref{qiact}, (4). On the other hand, there
exists $A'\in \aaa$ such that $\q_g (A)\subset \nn_M(A')$. It
follows that $A'$ intersects $B(x_0,C_0+M)$, hence it is in
$\mathcal{F}$. Corollary \ref{stab}, (a), implies that $\hdist
(\q_g(A),A')\leq D $, consequently $A$ is quasi-isometric to $A'$.

For every $i\in \{1, 2,\dots ,k \}$ define
$$
I(A_i)=\{ j\in \{ 1,2,\dots ,k \} \mid \hbox{there exists } g \in
G \mbox{ such that } \hdist(\q_g(A_i), A_j)\leq D \}.
$$

For every $j\in I(A_i)$ we fix $g_j \in G$ such that
$\hdist(\q_{g_j} A_i, A_j)\leq D $. Let $\Gamma(A_i)=\{g_j
\}_{j\in I(A_i) }$ and let $K(A_i)=\max_{j\in I(A_i)} \dist
(g_j\bar{\q}(x_0), \bar{\q}(x_0))$.

We define the constant $K=L_0\max_{i\in \{ 1,2,\dots ,k \}}
K(A_i)+(2L_0+1)\delta_0$, where $\delta_0=L_0C_0+2C_0$.

The following argument uses an idea from \cite[$\S 5.1$]{KaL2}.

\begin{lemma}\label{qtranz}
For every $A\in \mathcal{A}$ the $D$-stabilizer of $A$ acts
$K$-transitively on $A$, that is $A$ is contained in the
$K$-tubular neighborhood of every orbit $\St_{D }(A)a$, where
$a\in A$.
\end{lemma}

\proof Let $a$ and $b$ be two arbitrary points in $A$. Lemma
\ref{qiact}, (4), implies that there exist $g,\gamma \in G$ such
that $\q_g (a),\, \q_\gamma (b) \in B(x_0,C_0)$. This implies that
\begin{equation}\label{ingrup}
\dist (\mathbf{g}\circ \bar{\q}(a), \bar{\q}(x_0))\leq \delta_0\,
,\: \dist (\mathbf{\gamma} \circ \bar{\q}(b), \bar{\q}(x_0))\leq
\delta_0\, .
\end{equation}

There exist $i,j\in \{ 1,2,\dots ,k \}$ such that $\hdist (\q_g
(A), A_i),\, \hdist (\q_\gamma (A), A_j)\leq D$. Then $\q_{\gamma
g^{-1}} (A_i)$ is at finite Hausdorff distance from $A_j$, which
implies that $\hdist (\q_{\gamma g^{-1}} (A_i), A_j)\leq D$ and
that $j\in I(A_i)$. Let $g_j$ be such that $\hdist (\q_{g_j}
(A_i), A_j)\leq D$. It follows that $g\gamma^{-1}g_j\in
\St_D(A_i)$. The relation $\hdist (\q_g (A), A_i)\leq D$ and
Corollary \ref{stab}, (c), imply that $\gamma^{-1}g_jg\in
\St_D(A)$. We have that
$$
\dist (\q_{\gamma^{-1}g_jg}(a),b)\leq L_0 \dist
(\gamma^{-1}g_jg\bar{\q}(a), \bar{\q}(b))+C_0+L_0C_0\leq L_0 \dist
(g_jg\bar{\q}(a), \gamma\bar{\q}(b))+\delta_0\, .
$$

This and inequalities (\ref{ingrup}) imply that
$$
\dist (\q_{\gamma^{-1}g_jg}(a),b)\leq L_0 \dist (g_j\bar{\q}(x_0),
\bar{\q}(x_0))+(2L_0+1)\delta_0 \leq K\, .
$$
\endproof

\begin{cor}
 For every $A\in \mathcal{A}$ the normalizer of $\St_{D}(A)$ in $G$
is $\St_{D}(A)$.
\end{cor}

\proof Let $g \in G$ be such that $\St_{D} (A)=g\iv \St_{D} (A)g$.
Let $B\in \mathcal{A}$ be such that $\hdist(\q_g(A),B)\leq D$.
Corollary \ref{stab}, (c), implies that $\St_{D}(A)=\St_{D}(B)=S$.
Let $a\in A$ and $b\in B$. We have $\hdist(Sa, Sb)\leq
L\dist(a,b)+C$ and also $\hdist(A,Sa)\leq K$ and $\hdist(B,Sb)\leq
K$, therefore $\hdist(A,B)\leq 2K+L\dist(a,b)+C$. Lemma
\ref{neigh}, (1), implies that $B=A$ and $g \in
\St_{D}(A)$.\endproof

\begin{lemma}\label{hds}
For every $i\in \{ 1,2,\dots ,m\}$ we have
$$
\hdist(\bar{\q }(A_i), \St_D(A_i))\leq \kappa\: ,
$$ where $\kappa$ is a constant depending on $L_0,C_0, M$ and
$\dist(\q (1),x_0)$.
\end{lemma}

\proof Let $x_i\in A_i \cap B(x_0,M+C_0)$. For every $g\in
\St_D(A_i)$ we have $\dist(\q_g (x_i), A_i)\leq D$, hence $\dist
(\mathbf{g}\circ\bar{\q}(x_i), \bar{\q}(A_i))\leq L_0D+2C_0 $. It
follows that $\dist (g, \bar{\q}(A_i))\leq L_0D+2C_0+\dist
(1,\bar{\q}(x_i))$. Or $\dist (1,\bar{\q}(x_i))\leq L_0 \dist (\q
(1), x_i)+(L_0+1)C_0\leq L_0M+ (2L_0+1)C_0 +L_0 \dist (\q (1),
x_0)$.

Let $\bar{\q}(b)\in \bar{\q}(A_i)$. According to Lemma
\ref{qtranz}, there exists $g\in \St_M(A_i)$ such that
$$\dist(b,\q_g (x_i))\leq K.$$ Hence $\dist(\bar{\q}(b),
\mathbf{g}\circ\bar{\q}(x_i))\leq L_0K + 2C_0$ and
$\dist(\bar{\q}(b),g)\leq L_0K + 2C_0+\dist
(1,\bar{\q}(x_i))$.\endproof

\begin{cor}\label{hda}
Let $A\in \mathcal{A}$. There exists $g \in G$ and $i\in \{
1,2,\dots ,m\}$ such that $$\hdist(\bar{\q }(A), g\St_{D}
(A_i))\leq \kappa +L_0D+2C_0.$$
\end{cor}

\me

We continue the {\em proof of Theorem \ref{relhipqi}}. Consider
the minimal subset $\{B_1,...,B_m\}$ of $\{A_1,...,A_k\}$ such
that for each $A_i$ there exists $B_{j_i}$ and $\gamma_i$ such
that $\hdist(A_i, \q_{\gamma_i}(B_{j_i}))\le D$. Let
$\mathcal{B}=\{ B_1,\dots , B_m\}$. We denote $S_i=\St_{D}(B_i)$,
$i\in \{ 1,2,\dots ,m \}$. Let us show that $G$ is asymptotically
tree-graded with respect to $S_1,...,S_m$.

Indeed, by Theorem \ref{qi}, $\Cay(G)$ is asymptotically
tree-graded with respect to $\{ \bar{\q}(A), A\in \aaa \}$.
Corollary \ref{hda} implies that each $\bar{\q}(A)$ is at
uniformly bounded Hausdorff distance from $g\St_D (A_i)$ for some
$i\in \{ 1,2,\dots ,k \}$ and $g\in G$. Corollary \ref{stab}, (c),
implies that $\St_D (A_i)= \gamma_i S_{j_i} \gamma_i^{-1}$, with
the notations introduced previously. It follows that $\hdist
(g\St_D (A_i), g\gamma_i S_{j_i} )\leq \max_{i\in \{1,\dots , k
\}} \dist (1,\gamma_i^{-1} )$. We conclude that $\bar{\q}(A)$ is
at uniformly bounded Hausdorff distance from $g\gamma_i S_{j_i}$.
Thus $G$ is asymptotically tree-graded with respect to
$S_1,...,S_m$.
\endproof

\me

\begin{cor}\label{relhipqig}
Let $G$ be a group that is asymptotically tree-graded with respect
to the family of subgroups $\{ H_1,H_2,\dots ,H_k\}$, where $H_i$ is
an unconstricted infinite group for every $i\in \{ 1,2,\dots ,k \}$.
Let $G'$ be a finitely generated group which is quasi-isometric to
$G$.  Then $G'$ is asymptotically tree-graded with respect to a
finite collection of subgroups $\{S_1,\dots,S_m\}$ such that each
$S_i$ is quasi-isometric to one of the $H_j$.
\end{cor}

\begin{remarks}
If the groups $H_i$ in Corollary \ref{relhipqig} are contained in
classes of groups stable with respect to quasi-isometries (for
instance the class of virtually nilpotent groups of a fixed
degree, some classes of virtually solvable groups) then $S_i$ are
in the same classes. \me
\end{remarks}

\begin{cor}\label{setgr} If a group is asymptotically tree-graded
with respect to a family of subsets $\aaa$ satisfying conditions
(1), (2), (3) in Theorem \ref{relhipqi}, then it is asymptotically
tree-graded with respect to subgroups $H_1, ..., H_m$ such that
every $H_i$ is quasi-isometric to some $A\in \aaa$.
\end{cor}

\begin{remark}\label{522}
\begin{itemize}
  \item[(a)] If in Theorem \ref{relhipqi} we have that the cardinality of
  $\mathcal{A}$ is at least two then
for every $i\in \{1,2,\dots ,m\}$, $H_i$ has infinite index in
$G$.
  \item[(b)] If in Corollary \ref{relhipqig} we have $\{H_1,\dots ,H_k \}\neq \{
G\}$ then for every $j\in \{1,2,\dots ,m\}$, $S_j$ has infinite
index in $G'$.
\end{itemize}
\end{remark}

\proof (a) Suppose that $\{H_1,\dots , H_k \}= \{ G\}$. According
to the proof of Theorem \ref{relhipqi}, it follows that
$G=\mathrm{St}_D (B)$ for some $B\in \aaa$. Lemma \ref{hds} then
implies that $\hdist(\bar{\q }(B), G)\leq \kappa$, whence $\hdist(
B , X)\leq 3C_0+L_0\kappa$. This contradicts the property
$(\alpha_1)$ satisfied by $\aaa $.

Therefore  $\{H_1,\dots , H_k \}\neq \{ G\}$. Now the statement
follows from Remark \ref{infind}.

Statement (b) follows from (a).\endproof

\section{Cut-points in asymptotic cones of groups}\label{ascp}

Theorem \ref{relhipqi} shows that we need to study unconstricted
groups. In this section we provide two classes of examples of such
groups. We begin with some general remarks. Let $G$ be a finitely
generated group such that an asymptotic cone $\co{G;e,d}$ has a
cut-point, where $e=(1), d=(d_n)$. Lemma \ref{cutting} implies that
$\co{G;e,d}$ is a
 tree-graded space with respect to a set of pieces $\pp$ such that each piece is
 either a point or a geodesic subset
 without cut-point. In particular, if all the pieces are points then the cone is a
 tree. By homogeneity in this case it can be either a
 line or a tree in which every point is a branching point with the same degree.

 The case when one asymptotic cone is a line turns out to be quite
 particular. More precisely, we have the following general
 results.

\begin{proposition}\label{ur}
Let $\mathcal{G}$ be a family of finitely generated non-virtually
cyclic groups. Then for any sequence of groups $G_n\in\mathcal{G}$
endowed with word metrics $\dist_n$, any sequence $(\lambda_n)$ of
positive numbers with $\lim \lambda_n=0$, any $e\in \Pi G_n$ and any
ultrafilter $\omega$, the ultralimit $\lio{G_n,\lambda_n\dist_n}_e$
is neither a point nor a (real) line.
\end{proposition}

\proof We may assume without loss of generality that $e_n=1$ for
every $n$. If an ultralimit $\lio{G_n,\lambda_n \dist_n}_e$ is a
point then $G_n$ are finite $\omega$-almost surely, which is a
contradiction.

Suppose that an ultralimit $\lio{G_n,\lambda_n \dist_n}_e$ is a
line. Since $G_n$ are all infinite, it follows that for any $n\in
\N$, there exists in $\Cay (G_n)$ a geodesic line $\g_n$ through
$1$. Then $\lio{\g_n}=\lio{G_n}$. Suppose by contradiction that
$\omega$-almost surely $G_n \not\subset \nn_{1/\lambda_n} (\g_n)$.
Then $\omega$-almost surely there exists $z_n\in G_n$ at distance
at least $1/\lambda_n$ of $\g_n$. Let ${\mathfrak l}_n$ be a
geodesic joining $z_n$ to $z_n'\in \g_n$ and of length $\dist_n
(z_n, \g_n)$. For every point $t\in {\frak l}_n$ we have $\dist_n
(t, z_n')=\dist_n (t, \g_n)$. By homogeneity we may suppose that
$z_n'=1$.

In the ultralimit $\lio{G_n,\lambda_n\dist_n}_e$, ${\frak l}_\omega
= \lio{{\frak l}_n}$ is either a geodesic segment of length at least
$1$ with one endpoint $\lio{1}$, or a geodesic ray of origin
$\lio{1}$. If ${\frak l}_\omega$ has a point in common with
$\lio{\g_n}$ that is different from $\lio{1}$, then $\omega$-almost
surely there exists $t_n\in {\frak l}_n$ at distance of order
$\frac{1}{\lambda_n}$ of $1$ and at distance $o\left(
\frac{1}{\lambda_n}\right)$ of $\g_n$. This contradicts the equality
$\dist_n (t_n, 1)=\dist_n (t_n, \g_n)$. Hence ${\frak l}_\omega \cap
\lio{\g_n} = \left\{ \lio{1} \right\}$. But in this case
$\lio{G_n}\neq \lio{\g_n}$, contradiction.

It follows that $\omega$-almost surely $G_n \subset
\nn_{1/\lambda_n} (\g_n)$, which implies that $G_n$ is hyperbolic
with boundary of cardinality $2$, consequently virtually cyclic. We
have obtained a contradiction.\endproof

\begin{cor}\label{ptrof}
A finitely generated group with one asymptotic cone a point or a
line is virtually cyclic.
\end{cor}

\subsection{Groups with central infinite cyclic subgroups}

Let $G$ be a finitely generated group containing a central infinite
cyclic subgroup $H=\la a\ra$. We fix a finite set of generators $X$
of $G$ and the corresponding word metric $\dist$ on $G$.

\begin{lemma}\label{zg}
For every asymptotic cone $\co{G;e,d}$ of $G$ and every $\epsilon
>0$, there exists an element $h=(h_n)^\omega$ in $G^\omega_e \cap H^\omega $ which acts
isometrically on $\co{G;e,d}$, such that for every $x\in
\co{G;e,d}$, $\dist (hx,x)=\epsilon$.
\end{lemma}

\proof Let $w$ be a word in $X$ representing $a$ in $G$. It is
obvious that for every $r>0$ there exists $h=a^n\in H$ such that
$|h|$ is in the interval $[r-|w|, r+|w|]$. For every $n\ge 1$ we
consider $h_n\in H$ such that $|h_n|\in [\epsilon d_n-|w|\, ,\,
\epsilon d_n+|w|]$. According to Remark \ref{grr3}, the element
$h=(h_n)^\omega$ in $G^\omega_e$ acts as an isometry on
$\co{G;e,d}$. Moreover, for every $g=\lio{g_n}\in \co{G;e,d}$ we
have that $\dist(hg, g)=\lim_\omega \frac{\dist (h_ng_n\, ,\,
g_n)}{d_n}=\lim_\omega \frac{\dist (g_nh_n\, ,\,
g_n)}{d_n}=\lim_\omega \frac{|h_n|}{d_n}=\epsilon $.\endproof

\begin{lemma}\label{centre}
If an asymptotic cone $C$ of $G$ has a cut-point then $C$ is
isometric to a point or a (real) line.
\end{lemma}

\proof Let $C=\co{G;e,d}$ be an asymptotic cone that has a
cut-point, where $e=(1), d=(d_n)$. Then $C$ is tree-graded with
respect to a collection $\calp$ of pieces that are either points
or geodesic sets without cut-points. Let $h$ in $G^\omega_e \cap
H^\omega$ be as in Lemma \ref{zg} for $\epsilon=1$.

If all sets in $\calp$ are points then $C$ is an $\R$-tree. If
this tree contains a vertex of degree $>2$, then it does not admit
an isometry $h$ such that $\dist(h(x),x)=1$ for every $x$. Thus in
this case $C$ is isometric to $\R$ or to a point.

So we may suppose that $\calp$ contains pieces that are not
points. Let $M$ be such a piece.

\me

\textbf{Case I}. Suppose $h(M)=M$. Let $x$ be an arbitrary point in
$M$. By Lemma \ref{cutting}, part (b), there exists $y\in C\setminus
M$ such that $x$ is the projection of $y$ on $M$. Let $\delta =
\dist(x,y)$. Since $h$ acts as an isometry, it follows that
$y'=h(y)$ projects on $M$ in $x'=h(x)$ and that $\delta =\dist
(x',y')$. We have $\dist (x,x')=\dist (y,y')=1$. On the other hand
Lemma \ref{sir} implies that $[y,x]\cup [x,x']\cup [x',y']$ is a
geodesic. Consequently $\dist (y,y')=1+2\delta$, a contradiction.

\me

\textbf{Case II}. Suppose $h(M)\neq M$. Then $h(M)$ is another piece
of the tree-graded space $C$, by Proposition \ref{homeom}. Let $x$
be the projection of $h(M)$ on $M$ and let $y$ be the projection of
$M$ on $h(M)$. Let $z\in M \setminus \{x\}$ and $z'=h(z)$. By moving
$z$ a little, for instance along the geodesic $[z,x]$, we can ensure
that $z'\ne y$. Every geodesic joining $z$ and $z'$ contains $x$ and
$y$, by Lemma \ref{cv&proj}. Let $t$ be a point in $C\setminus M$
that projects on $M$ in $z$ (it exists by Lemma \ref{cutting}, part
(b)). The projection of $t'=h(t)$ onto $h(M)$ is then $z'$. Lemma
\ref{sir} implies that $[t,z]\cup [z, x] \cup [x, y]\cup [y,z' ]\cup
[z',t']$ is a geodesic, whence $\dist(t,t')=1+2\dist(t,z)$. This
contradicts the fact that $\dist(t,t')=1$.
\endproof

\begin{theorem}\label{infceter}
Let $G$ be a non-virtually cyclic finitely generated group that
has a central infinite cyclic subgroup $H$. Then $G$ is wide.
\end{theorem}

\proof By contradiction suppose that $G$ is not wide. Lemma
\ref{centre} implies that one of the asymptotic cones of $G$ is a
line or a point. Corollary \ref{ptrof} implies that $G$ is
virtually cyclic, a contradiction. \endproof

\begin{cor}
Let $G$ be a non-virtually cyclic group, that is asymptotically
tree-graded with respect to certain proper subgroups. Then every
finitely generated subgroup in the center $Z(G)$ is finite.
\end{cor}

Theorem \ref{infceter} has the following uniform version.

\begin{theorem}\label{uzg}
Let $\cg$ be the family of all finitely generated non-virtually
cyclic groups with a central infinite cyclic subgroup. The family
$\cg$ is uniformly unconstricted.
\end{theorem}

\proof Consider $G_n$ a sequence of groups in $\cg$, $\dist_n$ a
word metric on $G_n$ and $H_n=\la a_n \ra $ a central infinite
cyclic subgroup of $G_n$. Let $d_n\geq n \dist_n (1,a_n)$ for all
$n$. An argument as in the proof of Lemma \ref{zg} implies that
for every sequence of observation points $e$ and for every
$\epsilon
>0$, the ultralimit $\lio{G_n, \dist_n/d_n}_e$ has as isometry $h$
moving every point by $\epsilon$. With an argument analogous to
the one in the proof of Lemma \ref{centre} we deduce that
$\lio{G_n, \dist_n/d_n}_e$ is a line or a point. This contradicts
Proposition \ref{ur}.\endproof

\begin{cor}\label{cuzg}
Let $X$ be a metric space asymptotically tree-graded with respect to
a collection of subsets $\aaa$. For every $(L,C)$ there exists
$M=M(L,C)$ such that for every $(L,C)$-quasi-isometric embedding
$\q\colon G \to X$ of a finitely generated non-virtually cyclic
group $G$ with a central infinite cyclic subgroup, there exists
$A\in \aaa$ such that $\q(G) \subset \nn_M (A)$.
\end{cor}

\subsection{Groups satisfying a law}

\label{law}

\begin{proposition}\label{freegr}
Let space $\free$ be tree-graded with respect to a collection
$\calp$ of proper subsets. Suppose that  $\free$ is not an
$\R$-tree, and let $G$ be a group acting transitively on $\free$.
Then $G$ contains a non-abelian free subgroup.
\end{proposition}

\begin{remark} If $\mathbb F$ is an $\mathbb R$-tree, $G$ may not
contain non-abelian free subgroups even if it acts transitively on
$\mathbb F$. Indeed, let $G$ be the group of upper triangular
$2\times 2$-matrices with determinant $1$ acting by isometries on
the hyperbolic plane ${\mathbb H}^2$. The action is transitive.
Therefore the (solvable) group $G^\omega_e$ acts transitively on an
asymptotic cone of  ${\mathbb H}^2$ which is an ${\mathbb R}$-tree.
\end{remark}

\proof[Proof of Proposition \ref{freegr}]. By Lemma \ref{cutting} we
can assume that every piece in $\calp$ is either a point or does not
have a cut-point. Since $\mathbb F$ is not a tree, we can assume
that $\calp$ contains a non-singleton piece $M$.

\begin{lemma}\label{transl}
Let $a$ and $b$ be two distinct points in $M$. There exists an
isometry $g\in G$ such that the following property holds:
\begin{itemize}
\item $a\neq g(b)$, the projection of $g(M)$ onto $M$
is $a$ and the projection of $M$ onto $g(M)$ is $g(b)$.
\end{itemize}
We shall denote this property of $g$ by $P(a,b,M)$.
\end{lemma}

\proof There are two cases:

(A) There exist two distinct pieces in $\calp$ that intersect.

(B) Any two distinct pieces in $\calp$ are disjoint.

By homogeneity, in case (A), every point is contained in two
distinct pieces. In case (B) let $x,y$ be two distinct points in
$M$. There exists an isometry $g\in G$ such that $g(x)=y$. Since
$g(M)$ intersects $M$ in $y$ it follows that $g(M)=M$. We conclude
that in this case the stabilizer of $M$ in $G$ acts transitively
on $M$.

Suppose we are in case (A). Then we can construct a geodesic
$\g\colon [0,s]\to \free $ such that $s=\Sigma_{i=1}^\infty s_n$
with $0<s_n<\frac{1}{n^2}$ and $\g \left[ \Sigma_{i=0}^n
s_i,\Sigma_{i=0}^{n+1}s_i\right]\subset M_n$ for some pieces
$M_n$, where $M_n\neq M_{n+1}$ for all $n\in \N \cup \{ 0 \}$.
Here $s_0=0$. Such a geodesic exists by Lemma \ref{sir}. We call
such a geodesic \textit{fractal at the arrival point}. By gluing
together two geodesics fractal at their respective arrival points,
$\g \cup \g'$, and making sure that the two respective initial
pieces, $M_0$ and $M_0'$, are distinct, we obtain a geodesic
\textit{fractal at the departure and arrival points} or
\textit{bifractal}. By homogeneity, every point in $\free$ is the
endpoint of a bifractal geodesic.

Let $[a,c]$ be a bifractal geodesic. Corollary \ref{strconv}, (b),
implies that $[a,c]$ can intersect $M$ in $a$ or in a non-trivial
sub-geodesic $[a,c']$. Since $[a,c]$ is fractal at the departure
point the latter case cannot occur. It follows that the
intersection of $[a,c]$ and $M$ is $\{a\}$. There exists an
isometry $g\in G$ such that $g(b)=c$. Since $[a,c]$ is fractal at
the arrival point also, it follows that $[a,c]\cap g(M)=\{c\}$.
For every $x\in g(M)$ we have that $[a,c]\cup [c,x]$ is a
geodesic, by Lemma \ref{sir}. In particular $a$ is the projection
of $g(M)$ on $M$. A symmetric argument gives that $c=g(b)$ is the
projection of $M$ on $g(M)$.

Now suppose that case (B) holds. Lemma \ref{cutting}, part (b),
implies that $a$ is the projection of a point $x\in \free \setminus
M$. Let $g$ be an isometry in $G$ such that $g(b)=x$. If $[a,x]$
intersects $g(M)$ in $x$ then we repeat the previous argument.
Assume $[a,x]\cap g(M)=[x',x]$. By the hypothesis in case (B),
$x'\neq a$. We have $x'=g(b')$ for some $b'\in M$. Since the
stabilizer of $M$ in $G$ acts transitively on $M$, there exists $g'$
in it such that $g'(b)=b'$. We have that $gg'(M)=g(M)$ projects onto
$M$ in $a$ and $M$ projects onto $gg'(M)$ in $x'=gg'(b)$.\endproof

\me

\Notat \quad  For every $t\in M$ let $\Pi_t(M)$ be the set of
  points $x$ in $\free \setminus M $ that project onto $M$ in $t$.

\me

\begin{lemma}\label{prop}
Let $g$ satisfy property $P(a, b,M)$. Then:
\begin{itemize}
  \item[(a)] the isometry $g^{-1}$ satisfies property
  $P(b,a,M)$;
  \item[(b)] for every $t\neq b$ we have $g(\Pi_t(M))\subset
  \Pi_a(M)$.
\end{itemize}
\end{lemma}

\proof (a) We apply the isometry $g^{-1}$ to the situation in
$P(a,b,M)$.

(b) The set $g(\Pi_t(M))$ projects on $g(M)$ in $g(t) \neq g(b)$.
This, property $P(a,b,M)$ and Corollary \ref{transproj} imply that
$g(\Pi_t(M))$ projects onto $M$ in $a$ and that $\dist
(g(\Pi_t(M)), M)\geq \dist (g(M),M)>0$.\endproof

We now finish the {\em proof of Proposition \ref{freegr}}. Let
$a,b,c,d$ be four pairwise distinct elements in $M$. Lemma
\ref{transl} implies that there exist $g\in G$ satisfying
$P(a,b,M)$ and $h$ satisfying $P(c,d,M)$.

\begin{figure}[!ht]
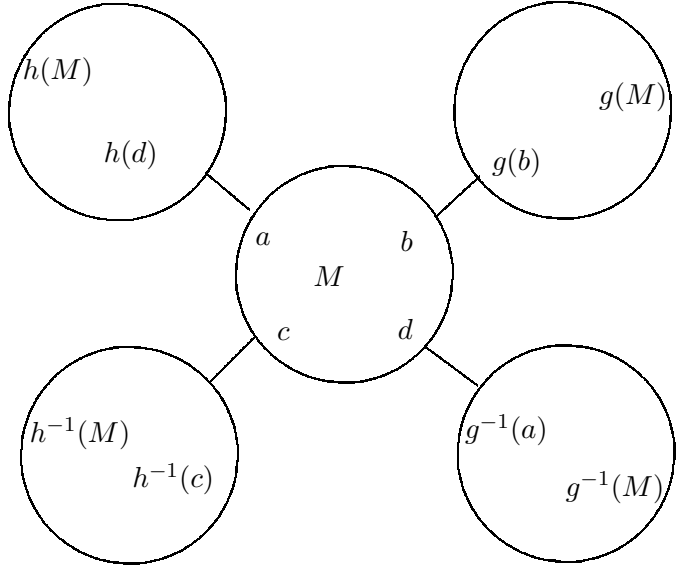

\centering
\unitlength .9mm 
\linethickness{0.4pt}
\ifx\plotpoint\undefined\newsavebox{\plotpoint}\fi 

\caption{Action of the elements $g, g^{-1}, h, h^{-1}$.}
\label{fig5}
\end{figure}

Then $g^{-1}$ is
satisfying $P(b,a,M)$ and $h^{-1}$ is satisfying $P(d,c,M)$ by
Lemma \ref{prop}. In particular $g(M)\subset \Pi_a(M)$,
$g^{-1}(M)\subset \Pi_b(M)$, $h(M)\subset \Pi_c(M)$,
$h^{-1}(M)\subset \Pi_d(M)$.

Since $b\not \in \{a,c,d \}$, Lemma \ref{prop}, part (b), implies
that $g(\Pi_a(M)\cup \Pi_c(M)\cup \Pi_d(M))\subset \Pi_a(M)$. The
isometries $g^{-1}, h, h^{-1}$ satisfy similar properties. The
Tits ping-pong argument allows to conclude that $g$ and $h$
generate a free group.\endproof

\begin{theorem}
Let $\mathcal{G}$ be a family of finitely generated non-virtually
cyclic groups satisfying a law. Then $\mathcal{G}$ is uniformly
wide.
\end{theorem}

\proof Suppose that an ultralimit ${\cal
C}=\lio{G_n,\frac{1}{d_n}\dist_n}_e$ has a cut-point, where
$\lim_\omega d_n=\infty$. Then by Lemma \ref{cutting} and
Proposition \ref{ur}, ${\cal C}$ is a tree graded space, not reduced
to a point nor isometric to $\R$. The group ${\cal G}=\Pi_e (G_n
,\frac{1}{d_n}\dist_n)/ \omega $ acts transitively on ${\cal C}$. If
${\cal C}$ is not an $\R$-tree, Proposition \ref{freegr} implies
that $\cal G$ contains a non-abelian free subgroup, and so it cannot
satisfy a non-trivial law, a contradiction.

Suppose that ${\cal C}$ is an $\mathbb R$-tree (but not a real
line). By \cite[Proposition 3.7, p.111]{Chis}, if ${\cal G}$ does
not fix an end of ${\cal C}$, ${\cal G}$ contains a non-abelian free
subgroup, a contradiction. Therefore we can assume that ${\cal G}$
fixes an end of ${\cal C}$. This means that ${\cal G}$
asymptotically fixes a ray $s(t)$, $t\in [0,\infty)$, starting at
$e$.  We shall now show that this assumption leads to a
contradiction.

Since the action of ${\cal G}$ on $\cal C$ is transitive, the ball
of radius $1$ in ${\cal C}$ around $e$ contains at least $9$
disjoint isometric copies of the ball of radius $1/4$ (of course,
here $9$ can be replaced by any positive integer). This implies that
$\omega$-almost surely for all $n$, the number of elements in the
ball of radius $d_n$ in the Cayley graph of $G_n$ is at least 9
times bigger than the number of elements in the ball of radius
$d_n/4$.

For $x\in \{1, 1.25, 1.5, 1.75\}$ Let $s(x)=(u_n(x))^\omega$, for
some $u_n(x)\in G_n$. Take any $g=(g_n)^\omega \in {\cal G}$ such
that $\dist_n(g_n,1)\le d_n$. Then $\dist(g\cdot 1,e)\le 1$. Note
that the image $g\cdot s$ is a ray which must be asymptotically
equal to $s$. Therefore the intersection $g\cdot s$ and $s$ contains
the subray $s(t), t\in [1,\infty)$. Since $g$ acts asymptotically on
this ray by translation, either $g\cdot s(1)$ or $g^{-1}\cdot s(1)$
belongs to the interval $s(t), t\in [1,2]$ of this subray. Therefore
either $g\cdot s(1)$ or $g^{-1} \cdot s(1)$ is within distance $1/4$
from  $s(x)$ for some $x\in \{1, 1.25, 1.5, 1.75\}$. This implies
that $\omega$-almost surely for any $n$, and any $g_n\in G_n$ with
$\dist_n(g_n,1)\le d_n$, for some $x\in \{1, 1.25, 1.5, 1.75\}$, and
a choice of $\epsilon\in \{1,-1\}$, we have
$$\dist_n(u_n(x)^{-1}g_n^\epsilon u_n(1),1)\le d_n/4.$$ This implies
that the $\omega$-almost surely for every $n$ the ball of radius
$d_n$ in the Cayley graph of $G_n$ contains at most 8 times more
elements than the ball of radius $d_n/4$, a contradiction with the
statement from the previous paragraph. \endproof

\me

\textit{Examples}: Solvable groups of a given degree, Burnside
groups of a fixed exponent and uniformly amenable groups (see
Corollary \ref{amen} below) are examples of groups satisfying a law.

\me

\begin{cor}\label{lawcp}
Let $G$ be a finitely generated non-virtually cyclic group
satisfying a law. Then $G$ is wide.
\end{cor}

\begin{cor} \label{cornice}
Let metric space $X$ be asymptotically tree-graded with respect to a
collection of subsets $\mathcal{A}$. For every non-trivial group law
and every $(L,C)$ there exists a constant $M$ depending on $(L,C)$
and on the law such that the following holds. Any
(L,C)-quasi-isometric embedding of a finitely generated
non-virtually cyclic group satisfying the law into $X$ has the image
in $\nn_M (A)$ for some $A\in \mathcal{A}$.
\end{cor}

The following statement is probably well known but we did not find a
proper reference.

\begin{lemma} \label{utr} Let $\omega$ be any ultrafilter, $G$ any group.
The group $G$ satisfies a law if and only if its ultrapower $\Pi
G/\omega$ does not contain free non-abelian subgroups.
\end{lemma}

\proof Clearly, if $\Pi G/\omega$ contains a free non-abelian
subgroup then $G$ does not satisfy any law. Conversely assume that
$G$ does not satisfy any law. Let us list all words in two
variables: $u_1, u_2,...,$ and form a sequence of words $v_1=u_1,
v_2=[u_1,u_2], v_3=[u_1,u_2,u_3], ...$ (iterated commutators). We
can choose the sequence $u_1, u_2,...$ in such a way that all words
$v_i$ are non-trivial. Since $G$ does not satisfy a law, for every
$i$ there exists a pair $(x_i,y_i)$ in $G$ such that $v_i(x_i,y_i)$
is not 1 in $G$. Let $x=(x_i)^\omega,\, y=(y_i)^\omega$ be elements
in the ultrapower. Suppose that the subgroup $\la x,y\ra$ of $\Pi
G/\omega$ has a relation. That relation is some word $u_i$ in two
variables. Hence $u_i(x_j,y_j)=1$ $\omega$-almost surely. In
particular, since $\omega$ is a non-principal ultrafilter, for some
$j>i$, $u_i(x_j,y_j)=1$. But then $v_j(x_j,y_j)=1$ since $u_i$ is a
factor in the commutator $v_j$, a contradiction.\endproof

Recall that a discrete group $G$ is (F\"olner) \textit{amenable} if
for every finite subset $K$ of $G$ and every $\epsilon \in
 (0,1)$ there exists a finite subset $F\subset G$ satisfying:
$$
|KF| < (1+\epsilon ) |F|.
$$

The group $G$ is \textit{uniformly amenable} if, in addition, one
can bound the size of $F$ in terms of $\epsilon$ and $|K|$, i.e.
there exists a function $\phi : (0,1) \times \N \to \N$ such that
$$
|F|\leq \phi(\epsilon,
    |K|)\, .
$$
For details on the latter notion see \cite{Kel}, \cite{Boz} and
\cite{Wys}. The following result has also been obtained in
\cite[Corollary 5.9]{Kel},
 we give a proof here for the sake of completeness.

\begin{cor}\label{amen}
A uniformly amenable finitely generated group satisfies a law and so
it is wide if it is not virtually cyclic.
\end{cor}

\proof Indeed, by \cite{Wys}, if $G$ is uniformly amenable then any
ultrapower $\Pi G/\omega$ is F\"olner amenable. Hence we can apply
Lemma \ref{utr} if we prove that any subgroup $S$ of an arbitrary
F\"olner amenable group $H$ is F\"olner amenable.

The argument is fairly standard and well known, we present it here
only for the sake of completeness. Take an arbitrary small $\epsilon
>0$. Take $K$ a finite subset in $S$. There exists a subset $F$ in $H$ such
that $|KF| < (1+\epsilon )|F|$. Consider a graph whose vertices are
the elements of the set $F$, and whose edges correspond to the pairs
of points $(f_1,f_2)\in F\times F$ such that $f_2=kf_1$, where $k\in
K$. Let $C$ be a connected component of this graph with set of
vertices $\mathcal{V}_C$. Then $K\mathcal{V}_C$ does not intersect
the sets of vertices of other connected components. Hence there
exists a connected component $C$ such that $|K\mathcal{V}_C|
<(1+\epsilon )|\mathcal{V}_C|$ (otherwise if all these inequalities
have to be reversed, the sum of them gives a contradiction with the
choice of $F$). Without loss of generality, we can assume that
$\mathcal{V}_C$ contains $1$. Otherwise we can shift it to $1$ by
multiplying on the right by $c^{-1}$ for some $c \in \mathcal{V}_C$.
Then $\mathcal{V}_C$ can be identified with a finite subset of $S$.
Therefore $S$ contains a subset $\mathcal{V}_C$ such that
$|K\mathcal{V}_C| < (1+\epsilon )|\mathcal{V}_C|$.\endproof

\begin{remark}\label{dlH}
The amenability defined by the existence of a left invariant mean on
the set of functions uniformly continuous to the left is not
inherited by subgroups in general. If $\mathcal{H}$ is a separable
infinite dimensional Hilbert space and $G=U(\mathcal{H})$ is the
group of unitary operators on $\mathcal{H}$ endowed with the weak
operator topology, then $G$ is amenable in the above sense
\cite{dlH}. On the other hand, if we take $\mathcal{H}=\ell^2
(F_2)$, with $F_2$ the free group of two generators, then $G$
contains $F_2$ \cite[Remark G.3.7]{BHV}.
\end{remark}

\section{Fundamental groups of asymptotic cones}\label{exAD}

In \cite{EO}, A. Erschler and D. Osin constructed (modifying an
idea from \cite{Ols}), for every ``sufficiently good" metric space
$M$, a two-generated group $G$ with the property that $M$
$\pi_1$-embeds isometrically into an asymptotic cone $\co G$. Thus
any countable group is a subgroup of the fundamental group of some
asymptotic cone of a finitely generated group. In this section we
modify, in turn, the construction from \cite{EO} to show that the
fundamental group of an asymptotic cone can be isomorphic to the
uncountable free power of any countable group. Moreover, that
asymptotic cone can be completely described as a tree-graded
space. In particular, if, say, $M$ is compact and locally
contractible then there exists a 2-generated group one of whose
asymptotic cones is tree-graded with respect to pieces isometric
to $M$. We also construct a 2-generated recursively presented
group with the maximal possible (under the continuum hypothesis)
number of non-homeomorphic asymptotic cones.

\subsection{Preliminaries on nets}

Let $(X,\dist)$ be a metric space. We recall some notions and
results from \cite{GLP}.

\begin{definition} A $\delta$-{\it separated set} $A$ in $X$ is a set such that
for every $x_1,x_2\in A$, $\dist(x_1,x_2 )\geq \delta $. A
$\delta$-{\it net} in $X$ is a set $B$ such that $X\in \nn_{\delta
}(B)$.
\end{definition}

\begin{remark}
A maximal $\delta$-separated set in $X$ is a $\delta$-net in $X$.
\end{remark}

\proof Let $N$ be a maximal $\delta$-separated set in $X$. For
every $x\in X\setminus N$, the set $N\cup \{ x\}$ is no longer
$\delta$-separated, by maximality of $N$. Hence there exists $y\in
N$ such that $\dist(x,y)< \delta$.
\endproof

\begin{definition}
We call a maximal $\delta$-separated set in $X$ a
$\delta$-\textit{snet}.
\end{definition}

Note that if $X$ is compact then every snet is finite, hence every
separated set is finite.

\begin{remark}\label{cresc}
Let $(X,\dist)$ be a metric space and let $(M_n)_{n\in \N }$ be an
increasing sequence of subsets of X. Let $(\delta_n)_{n\in \N }$
be a decreasing sequence of positive numbers converging to zero.
There exists an increasing sequence
$$
N_1\subset N_2 \subset  \cdots \subset N_n \subset \cdots \, ,
$$ such that $N_n$ is a $\delta_n$-snet in $(M_n , \dist)$.
\end{remark}

\proof There exists a $\delta_1$-snet in $M_1$, which we denote
$N_1$. It is a $\delta_1$-separated set in $M_2$. Let $N_2$ be a
$\delta_2$-snet in $M_2$ containing $N_1$. Then $N_2$ is a
$\delta_2$-separated set in $M_3$. Inductively we construct an
increasing sequence $(N_n)_{n\in \N}$.\endproof

\medskip

\Notat \quad Let $A$ be a subset in a metric space. We denote by
$\Gamma_\kappa(A)$ the metric graph with set of vertices $A$ and
set of edges
$$
\{ (a_1,a_2 )\mid a_1,a_2\in A,\, 0< \dist(a_1,a_2)\leq \kappa
\}\, ,
$$
such that the edge $(a_1,a_2 )$ is of length $\dist(a_1,a_2)$. We
shall denote the length of every edge $e$ by $|e|$. We endow
$\Gamma_\kappa(A)$ with its length metric.

\medskip

\Notat \quad Let $(X ,{\rm{dist}} )$ be a proper geodesic metric
space, let $O$ be a fixed point in it and let $\zeta\in (0,1)$. We
denote by $B_n=\overline{B(O, n)}$ the closed ball of radius $n$
around $O$. We consider an increasing sequence of subsets in $X$,
$$
\{ O\} \subset N_1\subset N_2 \subset  \cdots \subset N_n \subset
\cdots \, ,
$$ such that $N_n$ is an $\zeta^n$-snet in $B_n$. Let $\calg_n$ be the finite graph
$\calg_{\zeta^{[n/2]} }(N_{n})$, endowed with its length metric
${\rm{dist}}_n$ (here $[\frac n2]$ is the integer part of $\frac
n2$).

\me

We recall that two metric spaces with fixed basepoints
$(X,\dist_X, x)$ and $(Y,\dist_Y ,y)$ are said to be
\textit{isometric} if there exists an isometry $\phi :X\to Y$ such
that $\phi (x)=y$.

\me

\begin{lemma}\label{lim} In the notation as above:
\begin{itemize}
  \item[(1)] for every $n\ge 2$, for every $x,y\in N_n$ we have
  \begin{equation}\label{eq2'}
\dist(x,y) \leq \dist_n (x,y)\leq \left( 1+6\zeta^k \right) \left(
\dist(x,y) + 2\zeta^k \right) + 2\zeta^k \, ,
\end{equation} where $k=[\frac n2]$;
  \item[(2)] for every observation point $e\in \Pi N_n/\omega$,
the spaces $\lio{N_{n}, {\rm{dist}}_n}_e$, $\lio{\calg_{n},
{\rm{dist}}_n}_e$ and $\lio{B_{n}, {\rm{dist}}}_e$ with the
basepoints $\lio{e}$ are isometric.
  \item[(3)] the spaces $\lio{N_{n}, {\rm{dist}}_n}$, $\lio{\calg_{n},
{\rm{dist}}_n}$ with the basepoints $\lio{O}$ and $(X,
{\rm{dist}})$ with the basepoint $O$ are isometric.

\end{itemize}
\end{lemma}

\proof  (1) Let $x,y$ be two fixed points in $N_{n}$. If
$\dist(x,y) \leq \zeta^{[n/2]}$ then by construction $\dist(x,y)
=\dist_n(x,y)$ and both inequalities in (\ref{eq2'}) are true. Let
us suppose that $\dist(x,y) > \zeta^{[n/2]}$.

The distance $\dist_n(x,y)$ in $\calg_n$ is the length of some
path composed of the edges  $e_1e_2...e_s$, where $x=(e_1)_-$ and
$y=(e_s)_+$. It follows that
$$\dist_n(x,y)=\sum_{i=1}^{s} |e_i|\geq \dist(x,y).$$ We conclude that
$$
\dist_n(x,y)\geq \dist(x,y)\, .
$$

We also note that
\begin{equation}\label{add}
\dist_n(x,y)\geq \dist_m(x,y) \hbox{ for every } m\geq n\, ,
\end{equation}
since $N_n\subseteq N_m$.

\medskip

The distance $\dist(x,y)$ is the length of a geodesic
$\cf\colon[0,\dist(x,y)]\to X$. Since $x,y\in N_n\subset
\overline{B(O,n)}$, the image of this geodesic is entirely
contained in $\overline{B(O,2n)}$. Let $t_0=0,t_1,t_2,\dots
,t_m=\dist(x,y)$ be a sequence of numbers in $[0,\dist(x,y) ]$
such that $0< t_{i+1}-t_i \leq \frac{\zeta^{n}}{2}, \hbox{ for
every } i\in \{ 1,2,\dots ,m-1 \}$ and $m\leq \frac{2
\dist(x,y)}{\zeta^{n}}+1$. Since $\dist(x,y)> \zeta^{[n/2]}>
\zeta^n $, we can write $m\leq \frac{3 \dist(x,y)}{\zeta^{n}}$.
Let $x_i=\cf(t_i),\, i\in \{ 0,1,2,\dots ,m \}$. For every $i\in
\{ 0,1,2,\dots ,m \}$ there exists $w_i\in N_{2n}$ such that
$\dist(x_i,w_i)\leq \zeta^{2n}$. We note that $w_0=x, w_m=y$. We
can write
\begin{equation}\label{eq1}
\dist(x,y)=\sum_{i=0}^{m-1} \dist(x_i,x_{i+1})\geq
\sum_{i=0}^{m-1} [ \dist(w_i,w_{i+1}) - 2 \zeta^{2n}]\sum_{i=0}^{m-1} \dist(w_i,w_{i+1} ) -2m\zeta^{2n} \, .
\end{equation}

We have $\dist(w_i,w_{i+1})\leq \dist(x_i,w_i) +
\dist(x_i,x_{i+1})+\dist(x_{i+1},w_{i+1})\leq 2\zeta^{2n} +
\frac{\zeta^n}{2} \leq \zeta^n$ for $n$ large enough. Therefore
$w_i,w_{i+1}$ are connected in $\calg_{2n}$ by an edge of length
$\dist(w_i,w_{i+1} )$. We conclude that
$$\sum_{i=0}^{m-1} \dist(w_i,w_{i+1})=\sum_{i=0}^{m-1}
\dist_{2n}(w_i,w_{i+1} )\geq \dist_{2n}(w_0,w_m )=\dist_{2n}
(x,y).$$ This and (\ref{eq1}) implies that
$$
\dist(x,y) \geq \dist_{2n} (x,y) -6\dist(x,y) \zeta^n \, .
$$

We have obtained that
\begin{equation}\label{eq2}
\frac{1}{1+6\zeta^n}\dist_{2n}(x,y) \leq \dist(x,y) \leq
\dist_n(x,y), \hbox{ for all } x,y\in N_n.
\end{equation}

Let again $x,y$ be two points in $N_{n}$, $k=[n/2]$. There exist
$x',y'\in N_k\subset N_n$ such that $\dist(x, x'),\dist(y, y')\leq
\zeta^k $.
 This implies that $\dist(x, x')=\dist_n (x, x')\leq \zeta^k$ and likewise
 $\dist(y, y')=\dist_n (y, y')\leq \zeta^k $. Hence $\dist_n (x, y)\leq \dist_n (x', y') + 2\zeta^k$.

Inequalities (\ref{add}) and (\ref{eq2}) imply $$\dist_n
(x',y')\leq \dist_{2k}(x',y' )\leq (1+6\zeta^k) \dist(x',y' )\leq
(1+6\zeta^k) (\dist(x,y) + 2\zeta^k ).$$ This gives (\ref{eq2'}).

\me

(2) \quad We have $N_n\subset \calg_{n}\subset \nn_{\zeta^{[n/2]}}
\left( N_{n}\right)$. Therefore $\lm^\omega \left( \calg_{n},
{\rm{dist}}_n \right)_e=\lm^\omega \left( N_{n}, {\rm{dist}}_n
\right)_e$. Thus it is enough to prove that $\lm^\omega \left(
N_{n}, {\rm{dist}}_n \right)_e$ and $\lio{B_{n}, {\rm{dist}}}_e$
with the basepoints $\lio{e}$ are isometric.

We define the map
\begin{equation}\label{eq3}
\Psi\colon \lio{x_n} \mapsto \lio{x_n}
\end{equation}
from $\lm^\omega \left( N_{n}, {\rm{dist}}_n \right)_e$ to
$\lio{B_{n}, {\rm{dist}}}_e$. Inequalities (\ref{eq2'}) imply that
the map $\Psi$ is well defined and that it is an isometric
embedding.

We prove that $\Psi$ is surjective. Let $(y_n)^\omega \in \Pi_e
B_n/\omega$. For every $y_n$ there exists $x_n\in N_n$ such that
$\dist(x_n,y_n) \leq \zeta^n$. Since the sequence
$(\dist(y_n,e_n))$ is bounded, the sequence $(\dist(x_n,e_n))$ is
also bounded by the second inequality in (\ref{eq2'}), and so is
the sequence $(\dist_n (x_n,e_n))$. We have $\lio{x_n}\in
\lm^\omega \left( N_{n}, \dist_n \right)_e$ and $\Psi
(\lio{x_n})=\lio{x_n}$. As $\lm_\omega \dist(x_n,y_n)=0$ we
conclude that $\lio{x_n}=\lio{y_n}$.

\medskip

(3)\quad According to (2) it suffices to prove that $\lio{B_{n},
{\rm{dist}}}_O$ with the basepoint $\lio{O}$ and $X$ with the
basepoint $O$ are isometric. Let $x\in X$. For $n$ large enough,
$x\in \overline{B(O,n)}$. We define the map
\begin{equation}\label{eq4}
\Phi\colon x\mapsto \lio{x}
\end{equation}
from $X$ to $\lio{B_n}_O$.

The map $\Phi$ is clearly an isometric embedding. Let us show that
$\Phi $ is surjective. Let $(x_n)_{n\in \N }$ be such that $x_n\in
B_n$ and such that $\dist(O,x_n)$ is uniformly bounded by a
constant $C$. It follows that $x_n \in \overline{B(O, C)}$ for all
$n\in \N $. Since the space $X$ is proper, $\overline{B(O, C)}$ is
compact and there exists an $\omega$-limit $x$ of $(x_n)$. It
follows that $\lm_\omega \dist(x_n,x)=0$, which implies that
$\lio{x_n}=\lio{x}=\Phi (x)$.\endproof

\me

\Notat \quad We shall denote the point $\lio{O}$ also by $O$. This
should not cause any confusion.

\me

\begin{remark}\label{1.6}
The hypothesis that $X$ is proper is essential for the
surjectivity of $\Phi$ in the proof of part (3) of Lemma
\ref{lim}.
\end{remark}

\begin{definition} For every proper geodesic metric space $(X,\dist)$ with a fixed basepoint $O$,
and every sequence of points $e=(e_n)^\omega$, $e_n\in
B_n=\overline{B(O,n)}$, we shall call the limit $\lio{B_n}_e$ an
{\em ultraball} of $X$ with center $O$ and observation point $e$.
\end{definition}

\begin{remark} \label{uball}
Notice that the ultraballs $\lio{B_n}_e$ and $\lio{B_n}_{e'}$ with
observation points $e=(e_n)^\omega$ and $e'=(e_n')^\omega$, such
that $\dist(e_n,e_n')$ is uniformly bounded $\omega$-almost
surely, are the same spaces with different basepoints (see Remark
\ref{equality}).
\end{remark}

\begin{remark} It is easy to prove, using results from \cite[$\S
I.3$]{BGS} and \cite{KaL1}, that an ultraball of a complete
homogeneous locally compact CAT(0)-space is either the whole space
or a horoball in it (for a definition see \cite{BH}). In
particular the ultraballs of the Euclidean space $\R^n$ are $\R^n$
itself and all its half-spaces.
\end{remark}

We are now going to construct a proper geodesic metric space with
basepoint $(Y_C,\dist , O)$ whose fundamental group is any
prescribed countable group $C$, and such that every ultraball with
center $O$ of $Y_C$ either is isometric to the space $Y_C$ itself
or is simply connected.

Let $C=\la S \mid R \ra $ be a countable group. We assume that
$S=\{s_n \mid n\in \N\}=C$, and that $R$ is just the
multiplication table of $C$, i.e. that all relations in $R$ are
triangular. For every $n\in \N$, consider $X_n$ the part of the
cone $z^2=x^2+y^2$ in $\R^3$ which is above the plane $z=n-1$. The
intersection of this (truncated) cone with the plane $z=n-1$ will
be called its {\em base}. Cut a slit in $X_n$ of length $n\pi$, in the intersection of $X_n$ with the plane $z=2n$. This slit has simple closed curve boundary of length $2n\pi$, same as the length of the base of $X_{n+1}$. The resulting space is denoted by $Y_n$. The vertex of $Y_1$
is denoted by $O$.

Now consider the following construction. We start with the space
$Y_1$, glue in the space $Y_2$ so that the base hole of $Y_2$ is
isometrically identified with the boundary of the slit cut in $Y_1$, glue in
$Y_3$ so that the base hole of $Y_3$ is identified with the boundary of the
slit in $Y_2$, etc. The resulting space with the natural
gluing metric is denoted by $Y$. Now enumerate all relations in
$R=\{r_1,r_2,...\}$. For every $m=1,2,...$, $r_m$ has the form
$x_ix_jx_k\iv$. Choose a natural number $k=k(m)$ such that the
base holes of $Y_i, Y_j, Y_k$ are at the distance $\le k$ in $Y$
and such that $k(m)>k(m-1)$. Consider the circles $y_i, y_j, y_k$
obtained by cutting $Y_i, Y_j, Y_k$ by planes parallel to the base
hole at distance $k$ from $O$, connect these circles with $O$ by
geodesics. Glue in an Euclidean disc $D_n$ to the circles $y_i,
y_j, y_k$ and connecting geodesics such that the boundary of
$D_n$ is glued, locally isometrically, according to the relation $r_m$. We supply
the resulting space $Y_C$ with the natural geodesic metric
$\dist$.

We keep the above notation for balls $B_n=\overline{B(O,n)}$, and
metric spaces $N_n$ and $\Gamma_n$ for this space $Y_C$.

The following properties of the space $(Y_C, \dist)$ are obvious.

\begin{lemma}\label{X} (1) The space $Y_C$ is geodesic and proper.

(2) For every $d>0$ there exists a number $r>0$ such that every
ball of radius $d$ in $Y_C$, whose center is outside $B(O,r)$, is
contractible.

(3) The fundamental group of $Y_C$ is isomorphic to $C$.
\end{lemma}

\begin{lemma} \label{UX}
The ultraball $\lio{B_n}_e$ of $Y_C$ with center $O$ is simply
connected if $\dist(e_n,O)$ is unbounded $\omega$-almost surely,
otherwise it is isometric to $Y_C$.
\end{lemma}

\proof Indeed, if a point $e=(e_n)$ from $X^\omega$ is such that
$\dist(e_n,O)$ is bounded $\omega$-almost surely then the
corresponding ultraball is isometric to $Y_C$ by Remark
\ref{uball}. Suppose that $$\lim_\omega \dist(e_n,O)=\infty.$$ Let
$U$ be the corresponding ultraball. Then every closed ball
$\overline{B_U(e,r)}$ in $U$ is the $\omega$-limit of
$B_{Y_C}(e_n,r)\cap B_n$. By Lemma \ref{X}, the balls
$B_{Y_C}(e_n,r)$ are contractible $\omega$-almost surely.
Therefore $\overline{B_U(e,r)}$ is contractible. Since every loop
in $U$ is contained in one of the balls $\overline{B_U(e,r)}$, $U$
is simply connected.\endproof

\subsection{Construction of the group}\label{cgr}

Let $A$ be an alphabet and $\free_A$ a free group generated by
$A$. For every $w\in \free_A$ we denote by $|w|$ the length of the
word $w$.

\begin{definition}[property $C^*(\lambda)$]
A set $\mathcal{W}$ of reduced words in $\free_A$, that is closed
under cyclic permutations and taking inverses, is said to satisfy
{\it property }$C^*(\lambda )$ if the following hold.
\begin{itemize}
\item[(1)] If $u$ is a subword in a word $w\in \mathcal{W}$ so
that $|u|\geq \lambda |w|$ then $u$ occurs only once in $w$;
\item[(2)] If $u$ is a subword in two distinct words $w_1,w_2\in
\mathcal{W}$ then $|u|\leq \lambda \min (|w_1|,|w_2|)$.
\end{itemize}
\end{definition}

 We need the following result from \cite{EO}.

\begin{proposition}\label{words}\cite{EO}
Let $A=\{ a,b \}$. For every $\lambda >0$ there exists a set
$\mathcal{W}$ of reduced words in $\free_A$, closed with respect
to cyclic permutations and taking inverses, satisfying the
following properties:
\begin{itemize}
\item[(1)] $\mathcal{W}$ satisfies $C^*(\lambda )$;

\item[(2)] for every $n\in \N$, the set $\{ w\in \mathcal{W} \mid
|w|\geq n \}$ satisfies $C^*(\lambda_n )$ with
$\lm_{n\to\infty}\lambda_n =0$;

\item[(3)] $\lm_{n\to \infty } {\rm{card}} \{ w\in \mathcal{W} \mid |w|=n
\}=\infty$.
\end{itemize}
\end{proposition}

\me

\Notat \quad Let us fix $\lambda =\frac{1}{500}$, and a set of
words $\mathcal{W}$ provided by Proposition \ref{words}.

\me

Let $\kappa (n)={\rm{card}} \{ w\in \mathcal{W} \mid |w|=n \}\, .$
We have that $\lm_{n\to \infty } \kappa(n) = \infty $.

\me

Fix a number $\zeta\in (0,1)$. For every $n\in \N$, let $\Gamma_n$
be a finite metric graph with edges of length at least $\zeta^n$
and at most $\zeta^{[n/2]}$ and diameter at most $10n$ for $n$
large enough. We endow $\Gamma_n$ with the length metric
$\dist_n$.  Let $N_n$ be the set of vertices of $\Gamma_n$ and let
$O_n$ be a fixed vertex in $N_n$. Let $E_n$ be the number of edges
of $\Gamma_n$.

\begin{definition}[fast increasing sequences]\label{seq}
An increasing sequence $(d_n)$ of positive numbers is called {\it
fast increasing with respect to the sequence of graphs
}$(\Gamma_n)$ if it satisfies the following:
\begin{itemize}
\item[(1)] for every $i\geq [\zeta^n d_n]$, $\kappa \left( i
\right)\geq E_n$; \item[(2)] $\lm_{n\to \infty }\frac{\zeta^n
d_n}{d_{n-1}}=\infty$; \item[(3)] $\lm_{n\to
\infty}\frac{E_n}{\zeta^n d_n}=0$.
\end{itemize}
\end{definition}

Fast increasing sequences of numbers clearly exist.

\me


Let us fix a fast increasing sequence $d=(d_n)$ with respect to
the sequence of graphs $(\Gamma_n)$.

To every edge $e=(x,y)$ in $\Gamma_n$ we attach a word $w_n(e)$ in
${\mathcal{W}}$ of length $\left[ d_n|e|\right]$ such that
\begin{itemize}
\item[(1)] $w_n(e\iv)=w_n(e)^{-1}$; \item[(2)] $w_n(e)\neq
w_n(e')$ if $e\neq e'$.
\end{itemize}

We can choose these words because for every edge $e=(x,y)$ in
$\Gamma_n$, we have $\left[ d_n \dist(x,y) \right] \geq [\zeta^n
d_n]$ and because we have enough words in $\mathcal{W}$ of any
given length (part (1) of Definition \ref{seq}).

\begin{definition}[the presentation of the group $G$]
\label{defwp} We define the set of relations $R_n$ as follows: for
every loop $p=e_1e_2...e_s$ in $\Gamma_n$ we include in $R_n$ the
free reduction of the word
$$
w_n(p)=w_n(e_1)w_n(e_2)\dots w_n(e_s).
$$
Let $R=\bigcup_{n\in \N } R_n$ and let $G=\langle a,b \mid R
\rangle$.
\end{definition}

\Notat \quad We denote by $\Cay(G)$ the left invariant Cayley graph
of $G$ with respect to the presentation $G=\la a,b \mid R \ra $,
that is the vertices are elements of $G$ and the (oriented) edges
are $(g,gx)$ for every $x\in \{a,b,a\iv,b\iv\}$. The edge $(g,gx)$
in $\Cay(G)$ is usually labeled by $x$, so $\Cay(G)$ can be viewed
as a labeled graph. Every path in $\Cay(G)$ is labeled by a word in
$a$ and $b$. The length of a path $p$ in $\Cay(G)$ is denoted by
$|p|$. The distance function in $\Cay(G)$ is denoted by $\dist$, it
coincides with the word metric on $G$.

\medskip

\Notat \quad For every word $w$ in the free group $F_{\{ a,b\} }$
we denote by $g_w$ the element in $G$ represented by $w$.

\medskip

As in \cite{EO} and \cite{O1}, we introduce the following types of
words.

\begin{definitions}[words of rank $n$] Every freely reduced product
\begin{equation}\label{prod1}
w=w_n(e_1)w_n(e_2) \dots w_n(e_m), \end{equation} where
$e_1,...,e_m$ are edges in $\Gamma_n$ is called a \textit{word of
rank} $n$. The words $w_n(e_i)$ will be called the {\em blocks} of
$w$.

Every freely reduced product
$$
w_n(p)=w_n(e_1)w_n(e_2) \dots w_n(e_m),
$$ where $p=e_1e_2...e_m$ is a path in
$\Gamma_n$, is called a \textit{net word of rank} $n$.
\end{definitions}

\begin{remark}\label{rk67}
 The words $w_n(e)$ have length at least
$[\zeta^nd_n]\geq [d_{n-1}]\geq \frac{d_1}{\zeta^{n-1}}-1\geq n$
for $n$ large enough. This and the small cancellation assumptions
from Proposition \ref{words} imply that at most $2\lambda_n$ of
the length of the block $w_n(e)$ can cancel in the product
(\ref{prod1}) provided none of its neighbor factors is
$w_n(e\iv)$. In particular, if a path $p$ in $\Gamma_n$ has no
backtracking, at most $2\lambda_n$ of the length of any factor
$w_n(e)$ cancels in the word $w_n(p)$.
\end{remark}

\Notat \quad For every path $p$ in $\Gamma_n$ starting at $O_n$
let $\bar p$ be the path in $\Cay(G)$ labeled by $w_n(p)$ starting
at $1$. We denote by $\Re_n \subset \Cay(G)$ the union of all these paths $\bar
p$. It is easy to see that $\Re_n$ consists of all prefixes of all
net words $w_n(p)$, where $p$ is a path in $\Gamma_n$ starting at
$O_n$.

\begin{definition}[cells of rank $n$]\label{crn} By definition of the set of
relations $R$, the boundary label of every cell in a \vk diagram
$\Delta$ over $R$ is a net word. Therefore a cell in $\Delta$ is
called a {\em cell of rank } $n$ if its boundary label is a net
word of rank $n$.
\end{definition}

\begin{definition}[minimal diagrams] A \vk diagram over $R$
is called {\em minimal} if it contains the minimal number of cells
among all \vk diagrams over $R$ with the same boundary label, and
the sum of perimeters of the cells is minimal among all diagrams
with the same number of cells and the same boundary label.
\end{definition}

\Notat \quad The boundary of any \vk diagram (cell) $\Delta$ is
denoted by $\partial\Delta$.

\begin{lemma}\label{vkd}
(1) Every minimal \vk diagram $\Delta$ over $R$ satisfies the
small can\-cel\-lation pro\-per\-ty $C'(1/10)$~(that is, the
length of any path contained in the boundaries of any two distinct
cells in $\Delta$ cannot be bigger than $1/10$ of the length of
the boundary of any of these cells).

(2) Every cell $\pi$ in a minimal \vk diagram $\Delta$ over $R$
satisfies $|\partial\pi |\leq 2|\partial\Delta|$.
\end{lemma}

\proof (1) is Lemma 4.2 in \cite{EO}.

(2) We prove the statement by induction on the number $n$ of cells
in $\Delta$. If $n=1$ then the statement is obviously true.
Suppose it is true for some $n$. We consider a minimal van Kampen
diagram $\Delta$ with $n+1$ cells. By Greendlinger's lemma
\cite{LS} and Part (1) there exists a cell $\pi$ and a common path
$p$ of $\partial\pi$ and $\partial\Delta$ whose length is bigger
than $\frac7{10}|\partial\pi |$. It follows that
$|\partial\pi|\leq 2|\partial\Delta |$. Removing $p$ and the
interior of $\pi$, we obtain a minimal diagram $\Delta'$ with
boundary length smaller than $|\partial\Delta|$ and with fewer
cells than $\Delta$. It remains to apply the induction assumption
to $\Delta'$. \endproof

\Notat \quad We shall denote the graphical equality of words by
$\equiv$.

\begin{lemma}\label{paths}
Let $u\equiv u_1u_2u_3$ be a  word of rank $n$ and $u'\equiv
u_1'u_2u_3'$ be a word of rank $m$, $n\ge m$. Suppose $|u_2|$ is
at least $5\lambda$ times the maximal length of a block in $u'$.
Then $m=n$. In addition, if $u=w_n(p)$ and $u'=w_n(q)$ are net
words then the paths $p$ and $q$ in $\Gamma_n$ have a common edge
$e$: $p=p_1ep_2$, $q=q_1eq_2$, and $u_1$ (resp. $u_1'$) is a
prefix of $w_n(p_1e)$ (resp. $w_n(q_1e)$), $u_3$ (resp. $u_3'$) is
a suffix of $w_n(ep_2)$ (resp. $w_n(eq_2)$).
\end{lemma}

\proof Indeed, the conditions of the lemma imply that one of the
blocks of $u$ that either contains $u_2$ or is contained in $u_2$
has in common with one of the blocks of $u'$ at least $\lambda$ of
its length. The small cancellation condition $C^*(\lambda)$
implies that the blocks coincide, so $m=n$.  The rest of the
statement follows immediately from the definition of net words and
Remark \ref{rk67}.
\endproof

\begin{lemma}\label{geod2}
Let $u$ and $v$ be two words in $\{a,b\}$ that are equal in $G$.
Suppose that $u$ is a (net) word of rank $n$ and $v$ is a shortest
word that is equal to $u$ in $G$. Then $v$ is also a (net) word of
rank $n$. In addition, if $u$ is a net word, $u=w_n(p)$, then
$v=w_n(q)$ for some simple path $q$ in $\Gamma_n$ having the same
initial and terminal vertices as $p$.
\end{lemma}

\proof Consider a \vk diagram $\Delta$ over $R$ with boundary
$\partial\Delta=st$ where $u$ labels $s$, $v\iv$ labels $t$.

By Greendlinger lemma, property $C'(1/10)$ implies that there
exists a cell $\pi$ in $\Delta$ such that $\partial\pi$ and
$\partial\Delta$ have a common subpath $r$ of length
$\frac{7}{10}|\partial\pi|$. Since $v$ is a shortest word that is
equal to $u$ in $G$, no more than $\frac12$ of $\partial\pi$ is a
subpath of $t$. Therefore $|r\cap s|\ge \frac15 |\partial\pi|$.
Notice that the label of $\partial\pi$ is the reduced form of a
product of at least two blocks. Therefore the label of $r\cap s$
contains at least $(1-4\lambda)/5$ of a block in $\partial\pi$.
Lemma \ref{paths} implies that $\pi$ is a cell of rank $n$. After
we remove the cell $\pi$ from $\Delta$ we obtain a diagram
$\Delta'$ corresponding to an equality $u'=v$ of the same type as
$u=v$, that is $u'$ is a word of rank $n$ representing the same
element in $G$ as $u$ and $v$, and if $u=w_n(p)$ then
$u'=w_n(p')$, where $p'$ is a path in $\Gamma_n$ with $p'_-=p_-,
p'_+=p_+$. Since $\Delta'$ has fewer cells than $\Delta$, it
remains to use induction on the number of cells in $\Delta$.
\endproof

\subsection{Tree-graded asymptotic cones}

Recall that we consider any sequence of metric graphs $\Gamma_n$,
$n\ge 1$, satisfying the properties listed before Definition
\ref{seq}, that the set of vertices of $\Gamma_n$ is denoted by
$N_n$, and that we fix basepoints $O_n$ in $N_n$.
 For every $x\in N_n$ let $p_x$ be a path from $O_n$ to $x$ in
$\Gamma_n$. We define
$$
\Phi_n\colon N_n \to \Re_n\, ,\, \Phi_n (x)=w_n(p_x) \hbox{ in }
G$$ (see notation before Definition \ref{crn}).

The value $\Phi_n (x)$ does not depend on the choice of the path
$p_x$, because $w_n(q)$ is equal to 1 in $G$ for every loop $q$ in
$\Gamma_n$ by the definition of the presentation of $G$. Hence
$\Phi_n$ is a map.

\begin{remark}\label{retubn}
Notice that every point in $\Re_n $ is at distance at most
$\zeta^{[n/2]}d_n(1+\lambda_n )$ from $\Phi_n(N_n)$.
\end{remark}

The sequence of maps $(\Phi_n)$ clearly defines a map
$$
(x_n)^\omega\mapsto \left( \Phi_n (x_n) \right)^\omega.
$$
from $\Pi N_n/\omega$ to $\Pi\Re_n/\omega$.

\begin{remark}\label{rre}
Let $a=\Phi_n (x)\, ,\, x\in N_n,$ and let $b\in G$ be such that $a$
and $b$ can be joined in $\Cay(G)$ by a path labeled by $w_n(q)$,
where $q$ is a path in $\Gamma_n$ with $q_-=x$ and $q_+=y$. Then
$b=\Phi_n(y)\in \Phi_n(N_n)$.
\end{remark}

\begin{lemma}\label{ren} Let $e=(e_n)^\omega \in \Pi N_n/\omega$,
$e'=(\Phi_n(e_n))^\omega$. The map $\Phi_\omega\colon \lio{N_n,
\dist_n}_e \to \lio{\Re_n, \dist/d_n }_{e'}$ such that
$$\Phi_\omega \left( \lio{x_n} \right)=\lio{\Phi_n(x_n)}$$ is a
surjective isometry.
\end{lemma}

\proof For every $x,y\in N_n$, let $p=e_1e_2...e_s$ be a shortest
path from $x$ to $y$ in $\Gamma_n$. Then $\Phi_n(x)$ and
$\Phi_n(y)$ are joined in $\Cay(G)$ by a path labeled by $w_n(p)$.
It follows that
$$\dist(\Phi_n(x), \Phi_n(y) )\leq \sum_{i=1}^s
|w_n(e_i)|\leq d_n\sum_{i=1}^{s} |e_i|=d_n\dist_n(x,y).$$

By Lemma \ref{geod2}, for every $x,y\in N_n$ there exists a
geodesic joining $\Phi_n(x)$ to $\Phi_n(y)$ labeled by a net word
$w_n(q)$ of rank $n$. If $q=e_1e_2...e_t$ then
$$
w_n(q)=w_n(e_1)\dots w_n(e_t).
$$ Therefore
$$
\begin{array}{l}
\dist(\Phi_n(x), \Phi_n(y)) = |w_n(q)|\geq
\sum_{i=1}^{t}\left( 1-2\lambda_n \right)|w_n(e_i)|\\
 \hspace{1.25in} \geq \left( 1-2\lambda_n \right) \sum_{i=1}^{t}(d_n|e_i| -1 )\geq
 \left(1-2\lambda_n \right) ( d_n \dist_n (x,y) - t)\\
\hskip 1.25 in\geq \left( 1-2\lambda_n \right) ( d_n \dist_n (x,y)
- E_n).
\end{array}
$$

Thus for every $x,y\in N_n$:

\begin{equation}\label{plonj1}
\left( 1-2\lambda_n \right) ( d_n \dist_n (x,y) - E_n)\leq
\dist(\Phi_n(x) , \Phi_n (y) )\leq d_n \dist_n (x,y).
\end{equation}

According to (\ref{plonj1}), for every $\lio{x_n}, \lio{y_n}\in
\lio{N_n, \dist_n}_e$ we have that
\begin{equation}\label{plonj2}
\lm_\omega \dist_n (x_n,y_n) - \lm_\omega \frac{E_n}{d_n}\leq
\lm_\omega \frac{\dist(\Phi_n(x_n) , \Phi_n (y_n) )}{d_n}\leq
\lm_\omega \dist_n (x_n,y_n)\, .
\end{equation}

Since $(d_n)_{n\in \N}$ is a fast increasing sequence we have that
$\lm_\omega \frac{E_n}{d_n}=0$. This implies that $\Phi_\omega$ is
well defined and that it is an isometry.

Remark \ref{retubn} implies the surjectivity of the
map~$\Phi_\omega$.\endproof

\me

\Notat \quad We denote by $e$ the element $(1)^\omega\in
G^\omega$.

\begin{proposition}\label{contrg}
Let $(\Gamma_n)_{n\in \N}$ be a sequence of metric graphs
satisfying the properties listed before Definition \ref{seq}, let
$(d_n)_{n\in \N}$ be a fast increasing sequence with respect to
$(\Gamma_n)_{n\in \N}$ and let $G=\la a,b \mid R\ra $ be the group
constructed as above. For every ultrafilter $\omega$ the
asymptotic cone $\co{G;e,d}$ is tree-graded with respect to the
set of pieces:
\begin{equation}\label{p}
\mathcal{P}=\left\{ \lio{g_n \Re_n } \mid  (g_n)^\omega \in G
^\omega\mbox{ such that } \lm_\omega \frac{\dist(e, g_n \Re_n )
}{d_n}< \infty \right\},
\end{equation}
in particular different elements $(g_n)^\omega$ correspond to
different pieces from $\pp$.
\end{proposition}

\proof {\bf Property $(T_1)$.} Suppose that $\lio{g_n \Re_n } \cap
\lio{g_n' \Re_n }$ contains at least two distinct points, where
$(g_n)^\omega , (g_n')^\omega \in G^\omega$. We may suppose that
$(g_n')^\omega =(1)^\omega$. Let
$$\lio{a_n},\lio{b_n} \in \lio{g_n \Re_n } \cap \lio{\Re_n }\, ,\,
\lio{a_n}\neq \lio{b_n}.$$ The inclusion $\lio{a_n},\lio{b_n} \in
\lio{\Re_n }$ implies that $$\lio{a_n} = \lio{ \Phi_n (x_n)}\, ,\,
\lio{b_n} =\lio{ \Phi_n (y_n)}.$$ where $x_n,y_n \in N_n \, ,\,
\lio{x_n}\neq \lio{y_n} $. The inclusion $\lio{a_n},\lio{b_n} \in
\lio{g_n \Re_n }$ implies that $\lio{a_n} =\lio{ g_n \Phi_n
(x_n')}\, ,\, \lio{b_n} =\lio{ g_n \Phi_n (y_n')}$, where
$x_n',y_n' \in N_n \, ,\, \lio{x_n'}\neq \lio{y_n'} $.

By Lemma \ref{geod2}, for every $n\ge 1$, there exists a geodesic
$\pgot_1^{(n)}$ in $\Cay(G)$ joining $\Phi_n(x_n)$ with
$\Phi_n(y_n)$ labeled by a net word $w_n(p_1^{(n)})$, where
$p_1^{(n)}$ is a simple path from $x_n$ to $y_n$ in $\Gamma_n$. It
follows that $\pgot_1^{(n)}\subset \Re_n$. Similarly, there exists
a geodesic $\pgot_2^{(n)}$ joining $g_n\Phi_n(x_n')$ to
$g_n\Phi_n(y_n')$ contained in $g_n\Re_n$. The label of this
geodesic is a net word $w_n(p_2^{(n)})$. Let $\q_n$ be a geodesic
joining $\Phi_n(x_n)$ to $g_n\Phi_n(x_n')$ and $\q'_n$ a geodesic
joining $\Phi_n(y_n)$ to $g_n\Phi_n(y_n')$ in $\Cay(G)$. Both
$\q_n$ and $\q'_n$ have length $o(d_n)$. The geodesics
$\pgot_1^{(n)}$ and $\pgot_2^{(n)}$ on the other hand have length
$O(d_n)$. We consider the geodesic quadrilateral composed of
$\pgot_1^{(n)}, \q_n, \pgot_2^{(n)}, \q'_n $ and a minimal \vk
diagram $\Delta_n$ whose boundary label coincides with the label
of this quadrangle. Then $\partial\Delta_n$ is a product of four
segments which we shall denote $s_n, t_n, s_n', t'_n $ (the labels
of these paths coincide with the labels of the paths
$\pgot_1^{(n)}, \q_n, \pgot_2^{(n)}, \q'_n$ respectively).

There exists a unique (covering) map $\gamma$ from $\Delta$ to
$\Cay(G)$ that maps the initial vertex of $s_n$ to $1$ and
preserves the labels of the edges. The map $\gamma$ maps $s_n$ to
$\pgot_1^{(n)}\subseteq \Re_n$ and $s_n'$ to
$\pgot_2^{(n)}\subseteq g\Re_n$.

Let $\Delta_n^1$ be the maximal (connected) sub-diagram of
$\Delta_n$ that contains $s_n$ and whose $\gamma$-image is
contained in $\Re_n$. Likewise, let $\Delta_n^2$ be the maximal
sub-diagram of $\Delta_n$ that contains $s'_n$ and whose
$\gamma$-image is contained in $g\Re_n$. The complement $\Delta_n
\setminus (\Delta_n^1 \cup \Delta_n^2)$ has several connected
components.

\begin{figure}[!ht]
\centering
\unitlength .5mm 
\linethickness{0.4pt}
\ifx\plotpoint\undefined\newsavebox{\plotpoint}\fi 


\caption{The diagram $\Delta_n$.} \label{fig6}
\end{figure}

Suppose that the complement contains cells, and let $\Theta_n$ be
one of the non-trivial components of the complement. The boundary
of $\Theta_n$ is contained in $\partial\Delta_n^1 \cup t_n \cup
\partial \Delta_n^2 \cup t'_n$. By Greendlinger's lemma, there exists a cell $\pi$ in
$\Theta_n$ such that $\partial\pi \cap \partial\Theta_n$ contains
a path $u_n$ of length at least $\frac{7}{10} |\partial\pi|$.
Suppose that $u_n$ has more than $15\lambda$ of its length in
common with $\partial\Delta_n^1$. Then the labels of $\partial\pi$
and $\partial\Delta_n^1$ contain a common subword of length at
least $5\lambda$ of the length of a block participating in the
label of $\partial\pi$. By Lemma \ref{paths}, $\pi$ has rank $n$
and the $\gamma$-image of $\Delta_n^1\cup\pi$ is in $\Re_n$, a
contradiction with the maximality of $\Delta_n^1$. Hence $|u_n\cap
\partial\Delta_n^1|\le 15\lambda|u_n|$. Similar argument applies
to $\Delta_n^2$.

Therefore $|u_n\cap \left(
\partial\Delta_n^1\cup\partial\Delta_n^2\right)|\leq
30\lambda|u_n|$. It follows that $u_n$ has more than
$\frac{6}{10}|\partial\pi|$ in common with $t_n\cup t'_n$. Since
$\gamma(t_n)$ and $\gamma(t_n')$ are both geodesics, $u_n$ must
intersect both of them. We have $|u_n|\le 30\lambda|u_n| +
|t_n|+|t'_n|$, hence $|u_n| =o(d_n)$. Therefore
$$
\dist(\Phi_n(x_n),\Phi_n(y_n))\leq |u_n|+|t_n|+|t'_n|=o(d_n),
$$
a contradiction.

\medskip

{\bf Property $(T_2)$.} According to Proposition \ref{approxul},
it suffices to study sequences of geodesic $k$-gons $P_n$ in
$\Cay(G)$ with all lengths of edges of order $d_n$, $k$ fixed and
$\lio{P_n}$ a simple geodesic triangle. We need to show that
$\lio{P_n}$ is contained in one piece.

We fix such a sequence $(P_n)_{n\in \N }$ of $k$-gons in
$\Cay(G)$. Let $\mathcal{V}_n$ be the set of vertices of $P_n$. We
consider minimal \vk diagrams $\Delta^{(n)}$ and covering maps
$\gamma_n\colon\Delta^{(n)}\to\Cay(G)$ such that
$\gamma_n(\partial\Delta^{(n)})$ is $P_n$. We can consider the
boundary of $\Delta^{(n)}$ also as a $k$-gon whose vertices and
sides correspond to the vertices and sides of $P_n$.

\medskip

{\bf (a) Properties of the diagrams $\Delta^{(n)}$.}

\me

By Lemma \ref{vkd}, each cell from $\Delta^{(n)}$ has boundary
length $\le O(d_n)$. On the other hand, the cells of rank $k\geq
n+1$ have boundary of length at least $[\zeta^{n+1}d_{n+1}]$.
Property (2) of the fast increasing sequence $(d_n)$ implies that
for $n$ large enough all cells from the diagram $\Delta^{(n)}$ are
of rank $k\leq n$.

Suppose that $\omega$-almost surely there exists a cell $\pi$ of
rank $m\leq n-1$ in $\Delta^{(n)}$ the boundary of which
intersects two edges $[x,y],\, [z,t]$ without common endpoint.
Recall that the diameter of a cell of rank $m$ is at most
$10md_m\leq 10(n-1)d_{n-1}$. Then there exist two points in
$\gamma_n[x,y]$ and in $\gamma_n[z,t]$ respectively, which are at
distance at most $10(n-1)d_{n-1}$ of each other. In the
$\omega$-limit of $P_n$ we obtain that two edges without common
endpoint intersect in a point. This contradicts the fact that
$\lio{P_n}$ is a simple loop. We conclude that $\omega$-almost
surely all cells whose boundaries intersect two edges without
common endpoint are of rank $n$.

Suppose that the boundary of one of the cells $\pi$ of rank $m$ in
$\Delta^{(n)}$ is not a simple path. Then by applying the
Greendlinger lemma to any hole formed by $\partial\pi$, we get a
cell $\pi'$ whose boundary has a common subpath $u$ with
$\partial\pi$ such that $|u|\ge \frac{7}{10}|\partial\pi'|$. Then
there exists a block $w$ in $\partial\pi'$ such that $|w\cap
\partial\pi | \ge \frac{7}{20} |w|$. We apply Lemma \ref{paths} to $\partial \pi$
 and $\partial \pi'$ and we obtain that the ranks of $\pi$ and $\pi'$ coincide
 and that the boundary label of the union $\pi\cup\pi'$ is a net word of
rank $m$ corresponding to a loop in $\Gamma_m$. Hence the union of
the cells $\pi$ and $\pi'$ can be replaced by one cell
corresponding to a relation from $R$, a contradiction with the
minimality of $\Delta^{(n)}$. Hence the boundary of each cell in
$\Delta^{(n)}$ is a simple path.

Suppose that the boundaries of two cells $\pi_1,\, \pi_2,$ in
$\Delta^{(n)}$, of rank $m_1$ and $m_2$ respectively, intersect in
several connected components. We apply the Greendlinger lemma to a
hole formed by $\partial\pi_1 \cup \partial\pi_2$ and we get a
cell $\pi'$ whose boundary has a common subpath, of length at
least $ \frac{7}{10}|\partial\pi'|$, with $\partial\pi_1 \cup
\partial\pi_2$. Therefore $\partial \pi'$ has a common subpath
with one $\partial\pi_i, \, i\in \{ 1,2 \} \, ,$ of length at
least $\frac{7}{20}|\partial\pi'|$. Lemma \ref{paths} implies that
the ranks of $\pi_i$ and $\pi'$ coincide and that the boundary
label of $\pi_i\cup\pi'$ is a net word of rank $m_i$ corresponding
to a loop in $\Gamma_{m_i}$. Hence $\pi_i \cup \pi'$ can be
replaced by one cell, a contradiction with the minimality of
$\Delta^{(n)}$. We conclude that the intersection of the
boundaries of two cells, if non-empty, is connected.

Suppose that the boundary of a cell $\pi$ in $\Delta^{(n)}$ of
rank $m$ intersects one side $[x,y]$ of $\partial \Delta^{(n)}$ in
several connected components. We consider a hole formed by
$\partial\pi \cup [x,y]$ and we apply the Greendlinger lemma to
it. We obtain a cell $\pi'$ whose boundary has a common subpath
$u$ with $\partial\pi \cup [x,y]$, such that $|u|\ge
\frac{7}{10}|\partial\pi'|$. Since $\gamma_n [x,y]$ is a geodesic,
$u$ cannot have more than $\frac{5}{7}|u|$ in common with $[x,y]$.
Hence $|u\cap \partial\pi |\ge \frac{1}{5}|\partial\pi'|$, which
implies that there exists a block $w$ in $\partial\pi'$ such that
$|w\cap \partial\pi |\ge \frac{1}{10}|w|$. We apply Lemma
\ref{paths} to $\pi$ and $\pi'$ and as previously we obtain a
contradiction of the minimality of $\Delta^{(n)}$. Consequently,
the intersection of the boundary of a cell in $\Delta^{(n)}$ with
a side of $\partial \Delta^{(n)}$, if non-empty, is connected.

\medskip

{\bf (b) Existence of a cell $\pi_n$ of rank $n$ in $\Delta^{(n)}$
such that $\dist (P_n, \gamma_n (\partial \pi_n )) =o(d_n)$.}

\me

Take any vertex $v=v_n$ of the $k$-gon $\partial\Delta^{(n)}$. Let
$[x,v], [v,y]$ be the two consecutive sides of the $k$-gon
$\partial\dn$. Let $x_n'\in [x,v]$ be such that $\gamma_n(x_n')$
is the last point on $[\gamma_n(v),\gamma_n(x)]$ (counting from
$\gamma_n(v)$) for which there exists a point $z$ on
$[\gamma_n(v),\gamma_n(y)]$ with $\dist(\gamma_n(x_n'),z)$ not
exceeding $\zeta^{n/2}d_n$. Since $\zeta^{n/2}d_n=o(d_n)$,
$\lio{x_n'}=\lio{\gamma_n v}$ (recall that the triangle
$\lio{P_n}$ is simple). Therefore $\dist(x_n',\gamma_n v)=o(d_n)$.

Similarly let $y_n'\in[y,v]$ be such that $\gamma_n(y_n')$ is the
last point on $[\gamma_n(v),\gamma_n(y)]$ for which there exists a
point $z$ on $[\gamma_n(v),\gamma_n(x)]$ with
$\dist(\gamma_n(y_n'),z)\le \zeta^{n/2}d_n$. Then
$\dist(y_n',\gamma_n v)=o(d_n)$.

Consider the set $\Pi_v$ of cells $\pi $ in $\Delta^{(n)}$ whose
boundaries have common points with both $[x,v]$ and $[v,y]$. The
boundary of $\pi$ naturally splits into four parts: a sub-arc of
$[x,v]$, a sub-arc of $[v,y]$, and two arcs $c(\pi), c'(\pi)$
which connect points on $[x,v]$ with points on $[v,y]$ and such
that $c(\pi)$ and $c'(\pi)$ do not have any common points with
$[x,v]\cup[v,y]$ others than their respective endpoints. We assume
that $c'(\pi)$ is closer to $v$ than $c(\pi)$.

The cells from $\Pi_v$ are ordered in a natural way by their
distance from $v$. Take the cell $\pi\in\Pi_v$ which is the
farthest from $v$ among all cells in $\Pi_v$ satisfying
$$\dist(\gamma_n(c(\pi)_-), \gamma_n(c(\pi)_+))\le [\zeta^{n/2}d_n].$$
Let us cut off the corner of $\dn$ bounded by the triangle
$\Theta_v=c(\pi)\cup [c(\pi)_-,v]\cup [v,c(\pi)_+]$. Notice that
by the definition of $x_n', y_n'$, we have  $c(\pi)_-\in
[x_n',v]$, $c(\pi)_+\in [v,y_n']$. Therefore the lengths of the
sides of $\Theta_v$ are $o(d_n)$. Also notice that $\omega$-almost
surely $\Theta_v$ contains all cells of rank $\le n-1$ from
$\Pi_v$. That follows from the fact that the diameter of
$\Re_{k},\, k\leq n-1,$ does not exceed $10(n-1)d_{n-1}$, hence
for $n$ large enough it does not exceed $[\zeta^{n/2}d_n]$ by
property (2) of the definition of a fast increasing sequence.

Let us do this operation for every vertex $v$ of the $k$-gon
$\dn$. As a result, we get a minimal diagram $\dn_1$ such that
$\gamma_n(\dn_1)$ is a $2k$-gon $P_n'$ with $k$ sides which are
sub-arcs of the sides of $P_n$ (we shall call them {\em long
sides}) and $k$ sides which are curves of type $c(\pi)$ whose
lengths are $o(d_n)$ ({\em short sides}). Some of the short sides
may have length $0$. The $\omega$-limit $\lio{P_n'}$ coincides
with $\lio{P_n}$. We shall consider $\partial\dn_1$ as a $2k$-gon
with long and short sides corresponding to the sides of $P_n'$.

Notice that by construction $\dn_1$ does not have cells of rank
$\le n-1$ which have common points with two long sides of the
$2k$-gon $\partial\dn_1$.

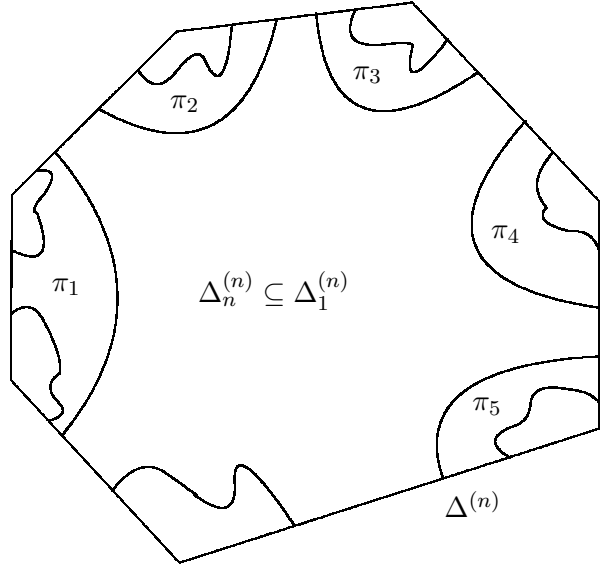
\begin{figure}[!ht]
\centering
\unitlength .7mm 
\linethickness{0.4pt}
\ifx\plotpoint\undefined\newsavebox{\plotpoint}\fi 
\begin{picture}(113.25,114)(0,0)
\multiput(33.5,7.5)(-.033719705,.0366174921){949}{\line(0,1){.0366174921}}
\multiput(1.5,42.25)(.03125,4.40625){8}{\line(0,1){4.40625}}
\multiput(1.75,77.5)(.0340043526,.0337323177){919}{\line(1,0){.0340043526}}
\multiput(33,108.5)(.27286585,.03353659){164}{\line(1,0){.27286585}}
\multiput(77.75,114)(.0337132004,-.0358499525){1053}{\line(0,-1){.0358499525}}
\put(113.25,76.25){\line(0,-1){43.25}}
\multiput(113.25,33)(-.105489418,-.0337301587){756}{\line(-1,0){.105489418}}
\qbezier(18.25,93.75)(44.38,79.25)(52,110.75)
\qbezier(59.5,112)(62.13,81.63)(90.25,100.75)
\qbezier(99.25,91.5)(73.25,61.75)(113.25,56)
\qbezier(113.25,46.75)(76.63,44.75)(83.5,23.75)
\qbezier(10,85.5)(32.88,58.38)(11.25,31.75)
\qbezier(25.75,101)(29.13,96)(33,101)
\qbezier(33,101)(36.88,105.88)(38.25,103.25)
\qbezier(38.25,103.25)(42.13,94.5)(43.5,109.75)
\qbezier(67,112.75)(63.38,105.13)(72.25,106)
\qbezier(72.25,106)(78.5,106.75)(77.75,105.5)
\qbezier(77.75,105.5)(75.38,95.25)(84.5,107)
\qbezier(6.5,82)(11.38,82.75)(7.75,77.5)
\qbezier(7.75,77.5)(4.88,74.25)(6.5,74)
\qbezier(6.5,74)(10.5,62.5)(1.5,67)
\qbezier(1.5,54.75)(5.63,58.38)(8.25,52.5)
\qbezier(8.25,52.5)(11.5,46)(10.75,43.5)
\qbezier(10.75,43.5)(7.38,41.25)(10.5,38)
\qbezier(10.5,38)(12.88,35.75)(8.75,34.5)
\qbezier(20.75,21.25)(26.13,28.13)(31,24.5)
\qbezier(31,24.5)(38.88,19)(36.25,20.5)
\qbezier(36.25,20.5)(44.38,14.5)(44,22.5)
\qbezier(44,22.5)(43.38,32.38)(55.25,14.75)
\qbezier(96.25,27.5)(90.38,31)(97,33.5)
\qbezier(97,33.5)(99.88,34.25)(100.25,38)
\qbezier(100.25,38)(100.75,41.75)(107.25,40.5)
\qbezier(107.25,40.5)(111.25,40.25)(113.25,38)
\qbezier(104.5,85.5)(99.13,79.75)(103.25,75)
\qbezier(103.25,75)(100.5,72.25)(107.75,70.5)
\qbezier(107.75,70.5)(112,70.88)(113.25,66.75)
\put(51.5,59.25){\makebox(0,0)[cc]{$\Delta_n^{(n)}\subseteq
\Delta_1^{(n)}$}} \put(12,60){\makebox(0,0)[cc]{$\pi_1$}}
\put(34.5,94.75){\makebox(0,0)[cc]{$\pi_2$}}
\put(69.25,100.25){\makebox(0,0)[cc]{$\pi_3$}}
\put(95.5,70){\makebox(0,0)[cc]{$\pi_4$}}
\put(92,37.25){\makebox(0,0)[cc]{$\pi_5$}}
\put(89.25,18.75){\makebox(0,0)[cc]{$\Delta^{(n)}$}}
\end{picture}
\caption{Diagram $\Delta^{(n)}$.} \label{fig7}
\end{figure}

Let $\pi_1$, $\pi_2$,...,$\pi_m$ be all Greendlinger
$\frac{6}{10}$-cells in $\dn_1$, i.e. for every $i=1,...,m$, the
intersection $\partial\pi_i\cap\partial\dn_1$ contains a subpath
$u_i$ of length at least $\frac6{10}|\partial\pi_i|$. Let $r_i$ be
the rank of the cell $\pi_i$, $i=1,...,m$. The path $u_i$ cannot
have more than $\frac56$ of its length in common with a long side
of the $2k$-gon $\partial\dn_1$ because the $\gamma_n$-images of
these sides are geodesics. By Lemma \ref{paths}, $u_i$ cannot have
a subpath of length bigger than $5\lambda$ times the length of a
block of rank $r_i$ in common with a short side of
$\partial\dn_1$. Since short sides and long sides in
$\partial\dn_1$ alternate $\omega$-almost surely, $u_i$ must have
points in common with two long sides of $\partial\dn_1$. Therefore
the number $m$ is at most $k$ and the rank $r_i$ is $n$ for every
$i=1,...,m$ ($\omega$-almost surely).

Let us cut off all cells $\pi_1,...,\pi_m$ from the diagram
$\dn_1$. The resulting diagram $\dn_2$ has a form of a polygon
where each side is either a part of a long side of $\dn_1$ (we
call it again {\em long}) or a part of $\partial\pi_i$ (we call it
{\em special}) or a part of a short side of $\dn_1$ (we call it
{\em short}). Notice that by the definition of $\dn_1$, the length
of any special side of $\dn_2$ cannot be smaller than
$[\zeta^{n/2}d_n]$ $\omega$-almost surely.

Suppose that the diagram $\dn_2$ contains cells $\omega$-almost
surely. Consider a Greendlinger $\frac{7}{10}$-cell $\pi$ of rank
$m$ in $\dn_2$ and the corresponding path
$u\subset\partial\pi\cap\partial\dn_2$. This path cannot have more
than $\frac57$ of its length in common with a long side of
$\dn_2$, more than $5\lambda$ times the length of a block of
$\partial \pi $ in common with a special or short side. Therefore
$u$ cannot contain a whole special side of $\dn_2$. Hence $u$ has
a subpath $u'$ of length at least
$(\frac{7}{10}-10\lambda)|\partial\pi|$ that intersects only long
and short sides of $\dn_2$. Hence $\pi$ is a Greendlinger
$\frac6{10}$-cell in $\dn_1$. This contradicts the fact that all
such cells were removed when we constructed $\dn_2$.

Thus $\dn_2$ contains no cells $\omega$-almost surely. In
particular, all cells in $\dn_1$ are of rank $n$ and all of them
are Greendlinger $\frac6{10}$-cells. For each cell $\pi_i$,
$i=1,...,m$, consider the decomposition $\partial\pi_i=u_iu_i'$.
Any two arcs $u_i', u_j'$ ($i\ne j$), have at most one maximal
sub-arc in common. The length of this sub-arc is at most
$5\lambda$ times the length of a maximal block of rank $n$ (by
Lemma \ref{paths} and the minimality of $\dn$). Hence
($\omega$-almost surely) the length of any arc $u_i'$ is at most
$5k\lambda[\zeta^{n/2}d_n]$. Therefore
$\lm_\omega\frac{|u_i'|}{d_n}=0$. Since $\lio{P_n'}$ is a simple
triangle, we can conclude that $\omega$-almost surely for all but
one $i\in \{ 1,...,m \}$ the length of $\partial\pi_i$ is
$o(d_n)$. Indeed otherwise we would have two points on $P_n'$ at
distance $O(d_n)$ along the boundary of $P_n'$ but at distance
$o(d_n)$ in $\Cay(G)$. The $\omega$-limits of these two points
would give us a self-intersection point of $\lio{P_n'}$.

Let us call this exceptional $i$ by $i_n$. Then $\lio{P_n'}$
coincides with $\lio{\gamma_n(\partial \pi_{i_n}})$. Since
$\gamma_n(\pi_{i_n})$ is contained in $g_n\Re_n$ for some $g_n$,
$\lio{P_n'}$ is contained in one piece $\lio{g_n\Re_n}$.
\endproof

\begin{proposition}[description of the set of pieces]\label{piec}
Consider the following two collections of metric spaces:
\begin{equation}\label{l1}\left\{\lio{g_n\Re_n}_e \mid
(g_n)^\omega\in G^\omega, \lm_\omega \frac{\dist(e, g_n \Re_n )
}{d_n}< \infty\right\}\end{equation} and
\begin{equation}\label{l2}\{\lio{N_n,\dist_n}_x \mid x\in \Pi
N_n/\omega\}\, .\end{equation}

We consider each $\lio{N_n,\dist_n}_x$ as a space with basepoint
$\lio{x_n}$ and each $\lio{g_n\Re_n}_e$ as a space with basepoint
$\lio{y_n}$, where $\lio{y_n}$ is the projection of $\lio{e}$ onto
$\lio{g_n \Re_n}$.

Then every space in one of these collections is isometric, as a
metric space with basepoint, to a space in the other collection.
Moreover every space in the second collection is isometric to
continuously many spaces in the first collection.
\end{proposition}

\proof Let $t_n=g_n\iv y_n$, $n\ge 1$. Let $y=(y_n)^\omega$ and
$t=(t_n)^\omega$. Then $\lio{g_n\Re_n}_e$ is isometric to
$\lio{g_n\Re_n}_y$ which, in turn, is isometric to
$\lio{\Re_n}_t$. Notice that $t_n\in \Re_n$, $\omega$-almost
surely. Remark \ref{retubn} implies that there exists a $u_n\in
\Phi_n(N_n)$ such that $\lim_\omega \frac{\dist(u_n,
t_n)}{d_n}=0$. Let $u=(u_n)^\omega$. For every $n\ge 1$, let
$x_n\in N_n$ be such that $u_n=\Phi_n(x_n)$, $x=(x_n)^\omega$.
Then by Lemma \ref{ren}, $\lio{g_n\Re_n}_e$ is isometric to
$\lio{N_n}_x$.

The fact that every limit set $\lio{N_n,\dist_n}_x$ is isometric
to a set $\lio{\Re_n, \dist/d_n }_{g}$ follows from Lemma
\ref{ren}. We write $g$ as $\lio{g_n^{-1}}$ for some $g_n^{-1}\in
\Phi_n(N_n)$. The set $\lio{g_n \Re_n, \dist/d_n }_{e}$ contains
$\lio{1}$ and with respect to this basepoint it is isometric to
$\lio{N_n,\dist_n}_x$.

We consider an arbitrary element $(\gamma_n)^\omega $ in $G^\omega
_e$ such that $\lim_\omega \frac{\dist (1,\gamma_n)}{d_n}=0$. The
set $\lio{\gamma_n g_n\Re_n }_{e}$ is distinct from the set
$\lio{g_n\Re_n }_{e}$, as the argument in Proposition \ref{contrg}
shows. On the other hand, the metric space $\lio{\gamma_n g_n\Re_n
}_{e}$ with basepoint $\lio{\gamma_n}=\lio{1}$ is isometric to the
metric space $\lio{g_n\Re_n }_{e}$ with basepoint $\lio{1}$, hence
to $\lio{N_n,\dist_n}_x$ with basepoint $\lio{x_n}$. We complete
the proof by noting that there are continuously many elements
$(\gamma_n)^\omega $ with $\lim_\omega \frac{\dist
(1,\gamma_n)}{d_n}=0$.
\endproof

\subsection{Free products appearing as fundamental groups of asymptotic
cones}\label{partc}

The following lemma is obvious.

\begin{lemma}
The collection of sets $\left\{ 2^{k}\N +2^{k-1} \mid k\in \N
\right\}$ is a partition of $\N$.
\end{lemma}



\me

\Notat \quad We denote the set $2^{k}\N +2^{k-1}$ by $\N_k$, for
every $k\in \N $. We denote by $k(n)$ the element $2^{k}n
+2^{k-1}$ of $\N_k$.

\me

Let $(M_k\, ,\, \dist_k)_{k\in \N }$ be a sequence of proper
geodesic locally uniformly contractible spaces, let $O_k$ be a
point in $M_k$ and let $\zeta$ be a real number in $(0,1)$. Fix
$k\in \N $. We apply Remark \ref{cresc} to the sequence of sets
$\left( B_n^{(k)} \right)_{n\in \N \cup \{0 \}}$, where
$B_0^{(k)}=\{ O_k\}$ and $B_n^{(k)}=\overline{B(O_k,n)},\, n\in \N
$, and to the sequence of numbers $\left(\zeta^{n} \right)_{n\in
\N}$. We obtain an increasing sequence
\begin{equation}\label{nets}
\{ O_k\} \subset N_1^{(k)}\subset N_2^{(k)} \subset  \cdots
\subset N_n^{(k)} \subset \cdots \, ,
\end{equation} such that $N_n^{(k)}$ is a $\zeta^{n}$-snet in $(B_n^{(k)} , \dist_k )$.
We consider the sequence of graphs $\Gamma_{\zeta^{[n/2]}}\left(
N_n^{(k)} \right)$ endowed with the length metric $\dist_n^{(k)}$.
We denote $\Gamma_{\zeta^{[n/2]}}\left( N_n^{(k)} \right)$ by
$\Gamma_n^{(k)}$.

\begin{remark}
Note that the diameter of $(N_n^{(k)} , \dist_k )$ is at most
$2n$, so by (\ref{eq2'}) the diameter of $(\Gamma_n^{(k)} ,
\dist_n^{(k)} )$ is at most $10n$, for $n$ large enough. Hence the
graphs $\Gamma_n^{(k)}$ satisfy the conditions listed before
Definition \ref{seq}.
\end{remark}

Now consider the sequence $(\Gamma_n\, ,\, \dist_n\, ,\, O_n)$ of
finite metric graphs endowed with length metrics and with
distinguished basepoints defined as follows: $(\Gamma_n\, ,\,
\dist_n\, ,\, O_n) \equiv (\Gamma_n^{(k)} \, ,\, \dist_n^{(k)}\,
,\, O_k )$ when $n\in \N_k$. We consider a sequence $(d_n)$ of
positive numbers which is fast increasing with respect to the
sequence of graphs $(\Gamma_n)$. We construct a group $G=\la a,b
\mid R\ra $ as in Section \ref{cgr}, associated to the sequences
$(\Gamma_n )$ and $(d_n)$.

For every $k\in \N$ let $\mu_k$ be an ultrafilter with the
property that $\mu_k(\N_k)=1$.

\begin{proposition}\label{contrgM}
The asymptotic cone $\mathrm{Con}^{\mu_k }(G;e,d)$ is tree-graded
with respect to a set of pieces $\mathcal{P}_k$ that are isometric
to ultraballs of $M_k$ with center $O_k$. Ultraballs with
different observation points correspond to different pieces from
$\pp_k$.
\end{proposition}

\proof By Proposition \ref{contrg}, $\mathrm{Con}^{\mu_k }(G;e,d)$
is tree-graded with respect to
\begin{equation}\label{l3}
\mathcal{P}_k=\left\{ \lm^{\mu_k }\left( g_n \Re_n \right) \mid
(g_n)^{ \mu_k} \in G^{\mu_k}\mbox{ such that } \lm_{\mu_k}
\frac{\dist(e, g_n \Re_n ) }{d_n}< \infty \right\}.
\end{equation}

By Proposition \ref{piec}, the collection of representatives up to
isometry of the set of pieces (\ref{l3}) coincides with the
collection of representatives up to isometry of the set of
ultralimits $\lm^{\mu_k }\left(N_n, \dist_n \right)_{x}$, $x\in
\Pi N_n/\mu_k$. The hypothesis that $\mu_k (\N_k)=1$ and the
definition of the sequence of graphs $(\Gamma_n)$ implies that
$\lm^{\mu_k }\left(N_n, \dist_n \right)_{x}=\lm^{\mu_k
}\left(N_n^{(k)}, \dist_n^{(k)} \right)_{x^{(k)}}$ for some
$x^{(k)} \in \Pi N_n^{(k)}/\mu_k$. It remains to apply Lemma
\ref{lim}.\endproof

\begin{cor}\label{cor56}
Suppose that the space $M_k$ is compact and locally uniformly
contractible. Then the asymptotic cone $\mathrm{Con}^{\mu_k
}(G;e,d)$ is tree-graded with respect to pieces isometric to
$M_k$, and the fundamental group of this asymptotic cone is the
free product of continuously many copies of $\pi_1(M_k)$.
\end{cor}

\proof It is a consequence of Proposition \ref{contrgM} and
Proposition \ref{pi1}.
\endproof

\begin{cor}\label{cth1}
There exists a $2$-generated group $\Gamma$ such that for every
finitely presented group $G$, the free product of continuously
many copies of $G$ is the fundamental group of an asymptotic cone
of $\Gamma$.
\end{cor}

\begin{theorem} \label{thcount} For every countable group $C$, there exists a
finitely generated group $G$ and an asymptotic cone $T$ of $G$
such that $\pi_1(T)$ is isomorphic to an uncountable free power of
$C$. Moreover, $T$ is tree-graded and each piece in it is
isometric either to a fixed proper metric space $Y_C$ with
$\pi_1(Y_C)=C$ or to a simply connected ultraball of $Y_C$.
\end{theorem}

\proof Let $C$ be a countable group. By Lemma \ref{X}, $C$ is the
fundamental group of a geodesic, proper, and locally uniformly
contractible space $Y_C$. Moreover, by Lemma \ref{UX}, there exists
a point $O$ in $Y_C$ such that every ultraball of $Y_C$ with center
$O$ either is isometric to $Y_C$ or is simply connected. It is easy
to see that the cardinality of the set of different ultraballs of
$Y_C$ with center $O$, that are isometric to $Y_C$, is continuum.
Consider the $2$-generated group $G=G(Y_C)$ obtained by applying the
above construction to $M_k=Y_C$ and $O_k=O$, $k\ge 1$. Then by
Proposition \ref{contrgM} there exists an asymptotic cone of $G$
that is tree-graded with respect to a set of pieces $\pp$ such that
the collection of representatives up to isometry of the pieces in
$\pp$ coincides with the collection of representatives up to
isometry of the set of ultraballs of $Y_C$ with center $O$. By
Proposition \ref{pi1}, the fundamental group of that asymptotic cone
is isomorphic to the free power of $C$ of cardinality continuum.
\endproof

\subsection{Groups with continuously
 many non-homeomorphic asymptotic cones}\label{nonh}

We use the construction in Section \ref{cgr} to obtain a
$2$-generated recursively presented group which has continuously
many non-$\pi_1$-equivalent (and thus non-homeomorphic) asymptotic
cones. Let us enumerate the set of non-empty finite subsets of
$\N$ starting with $\{1\}$ and $\{1,2\}$, then listing all subsets
of $\{1,2,3\}$ containing $3$, all subsets of $\{1,2,3,4\}$
containing $4$, etc. Let $F_k$, $k\in \N$, be the $k$-th set in
the sequence of subsets.

For every $n\ge 1$ let $\ttt^n$ be the $n$-dimensional torus
$\R^n/\Z^n$ with its natural geodesic metric and a basepoint
$O=(0,0,...,0)$.

For every $k\ge 1$ consider the bouquet of tori
$\calb_k=\bigvee_{n\in F_k} (\ttt^n,O)$. This is a compact locally
uniformly contractible geodesic metric space with a metric
$\dist_k$ induced by the canonical metrics on the tori and with
the basepoint $O_k=O$.

We repeat the construction of a group $G=\la a,b \mid R\ra$ in
Section \ref{partc} for the sequence of proper geodesic spaces
with basepoints $(\calb_k,  \dist_k, O_k)_{k\in \N }$.

Since all $\calb_k$ are bouquets of tori, we can choose the snets
$N^{(k)}_n$ coming from the same regular tilings of the tori of
different dimensions, and from their regular sub-divisions. There
is a recursive way to enumerate the snets $N^{(k)}_{k(n)}$. For an
appropriate choice of the set of words $\mathcal{W}$ in
Proposition \ref{words}, we obtain a recursively presented group
$G$. The group has the following property.

\begin{proposition}\label{spheres}
The asymptotic cone $\mathrm{Con}^{\mu_k }(G;e,d)$ is tree-graded
with respect to a set of pieces $\widetilde{\mathcal{P}}_k$ such
that every piece is isometric to one of the tori $\ttt^n$, $n\in
F_k$.
\end{proposition}

\proof Proposition \ref{contrgM} implies that the asymptotic cone
$\mathrm{Con}^{\mu_k }(G;e,d)$ is tree-graded with respect to a
set of pieces $\mathcal{P}_k$ such that all pieces are isometric
to $\calb_k$. It remains to use Lemma \ref{bouq}.\endproof

\me

\Notat \quad We denote $\mathrm{Con}^{\mu_k }(G;e,d)$ by $\calc_k$
and $\lm^{\mu_k }(e)$ by $e_k$.

\me

Let $I$ be an arbitrary infinite subset of $\N$, $I=\{i_1,\,
i_2,\dots ,\, i_n,\, \dots \}$. We consider the increasing
sequence of finite sets
$$
F_{k_1} \subset F_{k_2} \subset \cdots \subset F_{k_n}\subset
\cdots
$$ defined by $F_{k_n}=\{i_1,\,
i_2,\dots ,\, i_n \}$. Correspondingly we consider the sequence of
asymptotic cones $(\calc_{k_n})_{n\in \N }$. We consider an
ultrafilter $\omega $. The ultralimit
$\lio{\calc_{k_n}}_{(e_{k_n})_{n\in \N }}$ is also an asymptotic
cone of $G$, according to Corollary \ref{ulc}. We denote it by
$\calc_\omega (I)$.

\begin{lemma} \label{c6}
 Let $(\ttt^{k_i})$ be a sequence of tori
  $\ttt^{k_i}=\R^{k_i}/\Z^{k_i}$ with canonical flat metrics.
  Suppose that $\lio{k_i}=\infty$ for some ultrafilter $\omega$. Let
  $\ttt=\lio{\ttt^{k_i}}_e$ for some $e$. Then $\ttt$ contains
  isometric $\pi_1$-embedded copies of all tori $\ttt^n$, $n\ge 1$.
  \end{lemma}

  \proof Since tori are homogeneous spaces, we can assume that $e$
  is the sequence of points $(0,0,...)$. For every $n\ge 1$ the torus
  $\ttt^n$ isometrically embeds into $\ttt^{k_i}$ for
  $\omega$-almost all $i$ by the map $\phi_i\colon
  (x_1,...,x_n)\mapsto (x_1,...,x_n,0,0,...)$. Consequently $\ttt^n$
  isometrically embeds into $\ttt$ by $\phi_\omega \colon
  \bar{x}\mapsto\lio{\phi_i(\bar{x} )}$. Let $\cf$ be a
  non-0-homotopic loop in $\ttt^n$. Suppose that $\phi_\omega (\cf )$ is
$0$-homotopic in $\ttt$. Then there
  exists a continuous map $\psi \colon {\mathbb{D}}^2 \to \ttt$ with
  $\psi(\partial{\mathbb{D}}^2)=\phi_\omega (\cf )$. For every small
positive $\varepsilon $, there exists a triangulation of
${\mathbb{D}}^2$ such that if $e$ is an edge in the triangulation,
the images by $\psi$ of the endpoints of $e$ are at distance at
most $\varepsilon$. Let $\calv_\varepsilon$ be the set of vertices
of such a triangulation. The restricted map $\psi_\varepsilon \psi|_{\calv_\varepsilon}$ is an $\omega$-limit of maps $ \psi_i
\colon \calv_\varepsilon \to \ttt^{k_i}$. For every $i$ and for
every edge $e$ in the considered triangulation of ${\mathbb{D}}^2$
we join with a geodesic in $T^{k_i}$ the images by $\psi_i$ of the
endpoints of $e$. The length of this geodesic is $\omega$-almost
surely less than $2\varepsilon $. To each triangle of the
triangulation thus  corresponds a geodesic triangle in $T^{k_i}$
of perimeter smaller than
  $6\varepsilon$, $\omega$-almost surely. For $\varepsilon$ small enough
all these geodesic triangles are $0$-homotopic in some $T^{k_i}$.
But then $\cf$ is $0$-homotopic in $T^{k_i}$, a contradiction.
\endproof

\begin{lemma}\label{cI}
The asymptotic cone $\calc_\omega (I)$ is tree-graded with respect
to a set of pieces $\widetilde{\mathcal{P}}_{\omega}(I)$ such
that:
\begin{itemize}
  \item[(1)] All pieces are either isometric to one of the tori $\ttt^i$, $i\in
  I$, or they have the property that for every $n\in \N$ they contain an
isometric $\pi_1$-embedded copy of $\ttt^n$.
  \item[(2)] The fundamental group of every piece is Abelian.
\end{itemize}
\end{lemma}

\proof Proposition \ref{spheres} implies that for every $n\in
\N,\, \calc_{k_n}$ is tree-graded with respect to a set of pieces
$\widetilde{\mathcal{P}}_{k_n}$ such that every piece is isometric
to one of the tori $\{ \ttt^{i_1}, \ttt^{i_2},\, \dots ,\,
\ttt^{i_n} \}$. Theorem \ref{ultg} implies that $\calc_\omega
(I)=\lio{\calc_{k_n}}_{(e_{k_n})_{n\in \N }}$ is tree-graded with
respect to the set of pieces
\begin{equation}\label{pom}
\widetilde{\mathcal{P}}_\omega (I) = \left\{ \lio{M_n}\mid M_n \in
\widetilde{\pp}_{k_n}\, ,\,  \dist (e_{k_n},M_n) \mbox{ bounded
uniformly in }n\right\}\, .
\end{equation}

Let $\lio{M_n}$ be one of these pieces. Since $M_n \in
\widetilde{\pp}_{k_n}$, it follows that $M_n$ is isometric to one
of the tori $\{ \ttt^{i_1},\, \ttt^{i_2},\, \dots ,\, \ttt^{i_n}
\}$. Let $i(M_n)$ be the dimension of the torus $M_n$ and let
$\dist_n$ be the geodesic metric on $M_n$.

\me

\textbf{(1)} We have two possibilities.

\medskip

\noindent \textbf{I.} $\lio{i(M_n)}=\infty$. In this case we can
imply Lemma \ref{c6} and conclude that $\lio{M_n}$ contains
isometric and $\pi_1$-injective copies of tori $\ttt^N$ for every
$N$.

\medskip

\noindent \textbf{II.} $\lio{i(M_n)}<\infty$. It follows that
there exists a finite $m$ such that $i(M_n)\in \{ i_1,i_2,\dots ,
i_m\}$ $\omega$-almost surely. Remark \ref{udisj} implies that
there exists $j\in \{ 1,2,\dots , m\}$ such that $i(M_n)=i_j$
$\omega$-almost surely. Hence $\omega$-almost surely $M_n$ is
isometric to $\ttt^{i_j}$ and $\lio{M_n}$ is isometric to
$\ttt^{i_j}$.

\me

\textbf{(2)} Every torus $\ttt^n$ is a topological group, so it
admits a continuous binary operation and a continuous unary
operation satisfying the standard group axioms. It is not difficult to see that
$\omega$-limits of tori also are topological groups. Now the
statement follows from the fact that the fundamental group of
every topological group is Abelian \cite{Hat}.
\endproof

\begin{theorem}\label{theorem676}
The two-generated recursively presented group $G$ has continuously
many non-$\pi_1$-equivalent (and in particular non-homeomorphic)
asymptotic cones.
\end{theorem}

\proof Indeed, by Lemma \ref{cI} and Proposition \ref{pi1} the
fundamental group of $\calc_\omega (I)$ is a free product of
$\Z^{i}$, $i\in I$, and infinite dimensional Abelian groups. By
Kurosh's theorem \cite{LS}, if $j\not\in I$ then  $\Z^j$ cannot be a
free factor of that fundamental group. Hence the asymptotic cones
$\calc_\omega(I)$ for different subsets $I$ of $\N$ have different
fundamental groups.
\endproof

\begin{remark}\label{remark767} Each of the continuously many
asymptotic cones from Theorem \ref{theorem676} is a restrictive
asymptotic cone in the sense of Section \ref{adef}. Indeed by
Remark \ref{rem676}, each of the cones $\Con^{\mu_k}(G;e,d)$ is
isometric to a restrictive asymptotic cone
$\Con^{\nu_k}(G;e,(n))$. The map $\phi$ defined in Section
\ref{adef} just before Remark \ref{rem676} is in this case
injective. The images of the sets $\N_k$ under this map are
pairwise disjoint and $\nu_k(\phi(\N_k))=1$. It remains to use
Proposition \ref{prop676}.
\end{remark}

\section{Asymptotically tree-graded groups are relatively hyperbolic}
\label{secequiv}

Let $G$ be a finitely generated group and let $\{ H_1,\dots , H_m
\}$ be a collection of subgroups of $G$. Let $S$ be a finite
generating set of $G$ closed with respect to taking inverses.

\me


We denote by $\mh$ the set $\bigsqcup_{i=1}^m (H_i \setminus \{ e
\})$. We note that $\cgs$ is a subgraph of $\cgsh$, with the same
set of vertices but a smaller set of edges. We have that
$\dsh(u,v)\leq \ds (u,v)$, for every two vertices $u,v$.

For every continuous path $\pgot$ in a metric space $X$  we endow
the image of $\pgot$ with a pseudo-order $\le_\pgot$ (possibly not
anti-symmetric, but transitive and reflexive relation) induced by
the order on the interval of definition of $\pgot$. For every two
points $x,y$ we denote by $\pgot[x,y]$ the subpath of $\pgot$
composed of points $z$ such that $x\le_\pgot z\le_\pgot y$.

\me

\begin{definition}
Let $\pgot$ be a path in $\cgsh$. An $\mh$-\textit{component} of
$\pgot$ is a maximal sub-path of $\pgot$ contained in a left coset
$gH_i,\, i\in \{1,2, \dots ,m \},\, g\in G$ (i.e. this is a maximal
subpath with all labels edges belonging to $H_i$ for some $i$).

The path $\pgot$ is said to be \textit{without backtracking} if it
does not have two distinct $\mh$-components in the same left coset
$gH_i$.
\end{definition}

There are two notions of relative hyperbolicity. The weak relative
hyperbolicity has been introduced by B. Farb in \cite{Fa}. We use
a slightly different but equivalent definition. The proof of the
equivalence can be found in \cite{Osin}.

\begin{definition}
The group $G$ is \textit{weakly hyperbolic relative to} $\{
H_1,\dots , H_m \}$ if and only if the graph $\cgsh$ is
hyperbolic.
\end{definition}

The strong relative hyperbolicity has several equivalent
definitions provided by several authors. The definition that we
consider here uses the following property.

\begin{definition}
The pair $(G\, ,\, \{ H_1,\dots , H_m \})$ satisfies the
\textit{Bounded Coset Penetration (BCP) property} if for every
$\lambda \geq 1$ there exists $a=a(\lambda )$ such that the
following holds. Let $\pgot$ and $\q$ be two
$\lambda$-bi-Lipschitz paths without backtracking in $\cgsh$ such
that $\pgot_-=\q_-$ and $\ds (\pgot_+,\q_+)\leq 1$.
\begin{itemize}
  \item[(1)] Suppose that $s$ is an $\mh$-component of $\pgot$ such that $\ds (s_-,s_+)\geq
  a$. Then $\q$ has an $\mh$-component contained in the same left
  coset as $s$;
  \item[(2)] Suppose that $s$ and $t$ are two $\mh$-components of $\pgot$ and $\q$, respectively,
   contained in the same left
  coset. Then $\ds (s_-,t_-)\leq a$ and $\ds (s_+,t_+)\leq a$.
\end{itemize}
\end{definition}

\begin{definition}
The group $G$ is \textit{(strongly) hyperbolic relative to} $\{
H_1,\dots , H_m \}$ if it is weakly hyperbolic relative to $\{
H_1,\dots , H_m \}$ and if $(G\, ,\, \{ H_1,\dots , H_m \})$
satisfies the BCP property.
\end{definition}

We are going to prove the following theorem.

\begin{theorem}\label{dir}
A finitely generated group $G$ is asymptotically tree-graded with
respect to subgroups $\{ H_1,\dots , H_m \}$ if and only if $G$ is
(strongly) hyperbolic relative to $\{ H_1,\dots , H_m \}$ and each
$H_i$ is finitely generated.
\end{theorem}

This section is devoted to the proof of the ``only if" statement.
Note that the fact that each $H_i$ is finitely generated has been
proved before (Proposition \ref{propfg}).

The ``if" statement is proved in the Appendix.

\subsection{Weak relative hyperbolicity}

The most difficult part of Theorem \ref{dir} is the following
statement.

\begin{theorem}\label{weak}
If $G$ is asymptotically tree-graded with respect to $\{ H_1,\dots
, H_m \}$ then $G$ is weakly hyperbolic relative to $\{ H_1,\dots
, H_m \}$.
\end{theorem}

The main tool is a characterization of hyperbolicity due to
Bowditch \cite[Section 3]{Bow3*}. For the sake of completeness we
recall the results of Bowditch here.

\subsubsection{A characterization of hyperbolicity}\label{carhip}

Let $\cg$ be a connected graph, with vertex set $\vv$ and distance
function $\dist$, such that every edge has length $1$.

We assume that to every pair $u,v\in \vv$ we have associated a
subset $\La{uv}$ of $\vv$. Assume that each $\La{uv}$ is endowed
with a relation $\luv$ such that the following properties are
satisfied.
\begin{itemize}
  \item[($l_1$)] $\luv$ is reflexive and transitive;
  \item[($l_2$)] for every $x,y\in \La{uv}$ either $x\luv y$ or $y\luv x$;
  \item[($l_3$)] for every $u,v\in \vv$ we have  $\La{uv}=\La{vu}$ and
$\luv=\geq_{vu}$.
\end{itemize}

We note that the relations $\luv$ may not be anti-symmetric.

\me

\Notat \quad For $x,y\in \La{uv}$ with $x\luv y$, we write
$$
\La{uv}[x,y]=\La{uv}[y,x]=\{ z\in \La{uv} \mid x \luv z\luv y\}\,
.
$$

We also assume that we have a function $\phi : \vv \times \vv
\times \vv \to \vv$ with the following properties.
\begin{itemize}
  \item[($c_1$)](symmetry) $\phi \circ \sigma =\phi$ for every
$3$-permutation $\sigma$;
  \item[($c_2$)] $\phi (u,u,v)=u$ for all $u,v\in \vv$;
  \item[($c_3$)] $\phi (u,v,w) \in \La{uv} \cap \La{vw} \cap
  \La{uw}$.
\end{itemize}

Suppose moreover that there exists a constant $K>0$ such that the
following conditions are satisfied.
\begin{itemize}
  \item[(I)] For every $u,v,w \in \vv$, the Hausdorff distance
  between the sets $\La{uv} [u, \phi (u,v,w) ]$ and $\La{uw} [u, \phi (u,v,w)
  ]$ is at most $K$;
  \item[(II)] If $p,q\in \vv$ are such that $\dist (p,q)\leq 1$
  then $\diam\, \La{uv} [\phi (u,v,p)\, ,\, \phi (u,v,q)]$ is at most $K$;
  \item[(III)] If $w\in \La{uv}$ then $\diam \, \La{uv} [w\, ,\, \phi (u,v,w)]$
   is at most $K$.
\end{itemize}

 We call $(\La{uv}\, ,\, \luv )$ a \textit{line
from $u$ to $v$}. We call $\phi (u,v,w)$ the \textit{center of
$u,v,w$}.

\begin{proposition}[\cite{Bow3*}, Proposition 3.1]\label{hypb}
If the graph $\cg$ admits a system of lines and centers satisfying
the conditions above then $\cg$ is hyperbolic with the
hyperbolicity constant depending only on $K$. Moreover, for every
$u,v\in \vv$, the line $\La{uv}$ is at uniformly bounded Hausdorff
distance from any geodesic joining $u$ to $v$, where the previous
bound depends only on $K$.
\end{proposition}

\subsubsection{Generalizations of already proven results and new results}

\begin{lemma}
Let $\q : [0,t]\to X$ be an $(L,C)$-quasi-geodesic. Let $x$ be a
point in its image and let $a,b$ be its endpoints. Then
\begin{equation}\label{qgeod}
\dist (a,b)\geq \frac{1}{L_1}[ \dist (a,x)+\dist (x,b)]-C_1 \, ,
\end{equation} where $L_1=L^2$ and $C_1=C \left(
\frac{2}{L}+1\right)$.
\end{lemma}

\proof Let $s\in [0,t]$ be such that $\q(s)=x$. We have that
$\dist (a,b)\geq \frac{1}{L}t-C$. On the other hand $s\geq
\frac{1}{L}\dist (a,x)-C$ and $t-s\geq \frac{1}{L}\dist (x,b)-C$
imply that $t\geq \frac{1}{L}[ \dist (a,x)+\dist (x,b)]-2C$.
\endproof

Let $(X , \dist )$ be a metric space asymptotically tree-graded
with respect to a collection of subsets $\aaa$. Given $L\geq 1$
and $C\geq 0$ we denote by $M(L,C)$ the constant given by
$(\alpha_2')$ for $\theta=\frac13$.

\begin{definition}[parameterized saturations]\label{psatt}
Given $\q$ an $(L,C)$-quasi-geodesic and $\mu \geq 0$, we define
the $\mu$-\textit{saturation} $\Sat^\mu(\q)$ as the union of $\q$
and all $A\in\aaa$ with $\nn_\mu(A)\cap \q\ne\emptyset$.
\end{definition}

Notice that if a metric space $X$ is asymptotically tree-graded
with respect to a collection $\aaa=\{A_i\mid i\in I\}$ then $X$ is
also asymptotically tree-graded with respect to
$\nn_\mu(\aaa)=\{\nn_\mu(A_i)\mid i\in I\}$ for every number
$\mu>0$. This immediately follows from the definition of
asymptotically tree-graded spaces. One can also easily see that
properties $(\alpha_1), (\alpha_2), (\alpha_3)$ are preserved.
Hence the following two lemmas follow from Lemmas \ref{lsat},
\ref{limrsat} and \ref{gat}.

\begin{lemma}[uniform quasi-convexity of parameterized saturations]\label{plqqcun}
For every $L\geq 1$, $C\geq 0$ and $\mu \geq M(L,C)$, and for
every $\lambda \geq 1$, $\kappa \geq 0$, there exists $\tau =\tau
(L,C,\mu, \lambda , \kappa ) $ such that for every $R\ge1$, for
every $(L,C)$-quasi-geodesic $\q$, the saturation $\Sat^\mu(\q)$
has the property that every $(\lambda , \kappa )$-quasi-geodesic
$\cf$ joining two points in its $R$-tubular neighborhood is
entirely contained in its $\tau R$-tubular neighborhood.
\end{lemma}

\begin{lemma}[parameterized saturations of polygonal lines]\label{plimrsat}
The statements in Lemmas \ref{limrsat} and \ref{gat} remain true
if we replace the saturations by $\mu$-saturations, for every
$\mu>0$.
\end{lemma}

\begin{lemma}\label{polig}
Let $\q=\q_1 \cup \q_2 \cup \dots \cup \q_n$ be such that
\begin{itemize}
  \item[(1)] $\q_i$ is an $(L,C)$-almost-geodesic in $X$ for $i=1,2,\dots ,n$;
  \item[(2)] $\q_{i}\cap \q_{i+1}=\{ x_i\}$ is an
   endpoint of $\q_{i}$ and of $\q_{i+1}$ for $i=1,\dots n-1$;
  \item[(3)] $x_{i-1}$ and $x_{i}$ are the two
   endpoints of $\q_{i}$ for $i=2,\dots n-1$;
  \item[(4)] each $\q_i$ satisfies one of the following two
  properties:
\begin{itemize}
  \item[(i)] the endpoints of $\q_i$ are in a set $A_i\in \aaa$ or
  \item[(ii)] $\q_i$ has length
at most $\ell$, where $\ell$ is a fixed constant;
\end{itemize}
  \item[(5)] $A_i\neq A_j$ if $i\neq j$.
\end{itemize}

Then there exists $L_n\geq L\, ,\, C_n\geq C$ depending on $n,\,
\ell$ and $(L,C)$, such that $\q$ is an
$(L_n,C_n)$-almost-geodesic.
\end{lemma}

\proof Clearly $\q$ is an $L$-Lipschitz map. We prove by induction
on $n$ that $\dist (\q (t),\q (s))\geq \frac{1}{L_n} |t-s| - C_n$
for some $L_n\geq L$ and $C_n\geq C$.

The statement is true for $n=1$. Assume it is true for some $n$.
Let $\q=\q_1 \cup \q_2 \cup \dots \cup \q_n \cup \q_{n+1}$ be as
in the statement of the Lemma. Let $\q'=\q_1 \cup \q_2 \cup \dots
\cup \q_n$ which, by the induction hypothesis, is an
$(L_n,C_n)$-almost-geodesic.

Suppose that $\q_{n+1}$ satisfies (4), (ii). Then $\q$ is an
$(L_n\, ,\, 2 (\ell +C_n ) )$-almost-geodesic.

Suppose that $\q_{n+1}$ satisfies (4), (i). Let $A=A_{n+1}$. We
take $M_n=M(L_n,C_n)$. Let $y$ be the farthest point from $x_{n}$
in the intersection $\overline{\nn}_{M_n}(A) \cap \q'$. Consider
$\q_y$ a sub-almost-geodesic of $\q'$ of endpoints $y$ and $x_n$.
By Lemma \ref{qqc}, $\q_y$ is contained in the $\tau_n
M_n$-tubular neighborhood of $A$. On the other hand, $\q_y=\q_i'
\cup \q_{i+1} \cup \dots \cup \q_n$, where $\q_i'$ is a
sub-almost-geodesic of $\q_i$. Again Lemma \ref{qqc} implies that
every $\q_j$ satisfying (4), (i), is contained in $\nn_\tau (A_i)$
for some uniform constant $\tau$. Therefore, every such $\q_j$
composing $\q_y$ has endpoints at distance at most the diameter of
$\nn_\tau (A_i) \cap \nn_{\tau_n M_n} (A)$, hence at most $D_n$,
for some $D_n =D_n (\tau_n M_n )$. It follows that the distance
$\dist (y, x_{n})$ is at most $n( \ell +D_n)$. Lemma \ref{lemaC}
implies that if the endpoints of $\q_{n+1}$ are at distance at
least $D'=D'(L_n,C_n, D_n)$, then $\q$ is an $(L_n+1\, ,\,
2D')$-almost-geodesic.

If the endpoints of $\q_{n+1}$ are at distance at most $D'$ then
the length of $\q_{n+1}$ is at most $LD'+C$ and $\q$ is an $(L_n\,
,\, 2(LD'+C +C_n) )$-almost-geodesic.\endproof

\begin{lemma}\label{triangle1}
For every $L\geq 1$, $C\geq 0$, $M\geq M(L,C)$ and $\delta>0$
there exists $D_0>0$ such that the following holds. Let $A\in
\aaa$ and let $\q_i\colon [0,\ell_i]\to X$, $i=1,2$, be two
$(L,C)$-quasi-geodesics with one common endpoint $b$ and the other
two respective endpoints $a_i\in \nn_M (A)$, such that the
diameter of $\q_i\cap \overline{\nn}_M (A)$ does not exceed
$\delta$ for $i=1,2$. Then one of the two situations occurs:
\begin{itemize}
  \item[(a)] either $\dist(a_1,a_2)\leq D_0$ or
  \item[(b)] $b\in \nn_M (A)$ and  $\ell_i \leq L \delta +C$.
\end{itemize}
\end{lemma}

\proof Let $\dist(a_1,a_2)=D$. We show that if $D$ is large enough
then we are in situation (b). Remark \ref{buragos} implies that we
may suppose that $\q_i$ are $(L+C,C)$-almost geodesics.

According to Lemma \ref{lemaC}, there exists $D'$ such that if
$D\geq D'$ then $\q_1\sqcup [a_1,a_2 ]$ is an $(L+C+1, 2D'
)$-quasi-geodesic. Suppose that $D\geq D'$.

Suppose that $b$ is not contained in $\nn_M (A)$. Let $t \in
[0,\ell_2]$ be such that $\q_2 (t)\in \overline{\nn}_M (A)$ and
$\q_2|_{[0,t ]}$ does not intersect $\nn_M (A)$. The sub-arc
$\q_2|_{[t , \ell_2 ]}$ has endpoints at distance at most $\delta
$, hence it has length at most $L\delta +C$. It follows that
$\q_1\sqcup [a_1,a_2 ] \sqcup \q_2|_{[t , \ell_2 ]}$ is an $(L+C+1
\, ,\, C_1 )$-quasi-geodesic, where $C_1=C_1 (D',\delta )$. Lemma
\ref{lqqcun} implies that $\q_1\sqcup [a_1,a_2 ] \sqcup \q_2|_{[t,
\ell_2 ]}$ is contained in the $\tau $-tubular neighborhood of
$\Sat (\q_2|_{[0,t ]})$, where $\tau = \tau (L,C,D',\delta )$.
This implies that $\nn_M(A)\cap \nn_\tau \left( \Sat (\q_2|_{[0,t
]})\right)$ has diameter at least $D$. By Lemma \ref{diam}, for
$D$ large enough we must have that $A\subset \Sat (\q_2|_{[0,t
]})$. This contradicts the choice of $t$.

It follows that $b$ is contained in $\nn_M (A)$, which implies
that $\ell_i \leq L\dist (a_i,b) +C \leq  L\delta +C$.\endproof

\begin{cor}\label{12}
For every $L\geq 1$, $C\geq 0$, $M\geq M(L,C)$ and $\delta>0$
there exists $D_1>0$ such that the following holds. Let $A\in
\aaa$ and let $\q_i\colon [0,\ell_i]\to X$, $i=1,2$, be two
$(L,C)$-quasi-geodesics with one common endpoint $b$ and the other
two respective endpoints $a_i\in \nn_M (A)$, such that the
diameter of $\q_i\cap \overline{\nn}_M (A)$ does not exceed
$\delta$. Then $\dist(a_1,a_2)\leq D_1$.
\end{cor}

\begin{figure}[!ht]
\centering
\unitlength .5mm 
\linethickness{2pt}
\ifx\plotpoint\undefined\newsavebox{\plotpoint}\fi 

\caption{Corollary \ref{12} and Lemma \ref{restr}.} \label{fig8}
\end{figure}

\begin{lemma}\label{restr}
For every $L\geq 1$, $C\geq 0$ and $M\geq M(L,C)$ there exists
$\dg =\dg(L,C,M)>0$ such that the following holds. Let $A\in \aaa$
and let $\q \colon [0,\ell]\to X$ be an $(L,C)$-almost-geodesic
with endpoints $x$ and $y\in \nn_M (A)$. There exists a sub-arc
$\q'$ of $\q$ with one endpoint $x$ and the second endpoint in
$\nn_M (A)$ such that the diameter of $\q'\cap \overline{\nn}_M
(A)$ is at most $\dg$.
\end{lemma}

\proof If $x\in \nn_M (A)$ then we take $\q' = \{ x \}$. Suppose
that $x\not \in \nn_M (A)$. Let $t=\inf \{t'\in [0,\ell] \mid
t'\in \q^{-1}(\nn_M (A)) \}$. Then $\q (t) \in \overline{\nn}_M
(A)$. Let $s_i\in [0,t]$ be such that $\q (s_i) \in
\overline{\nn}_M (A),\, i=1,2$. If $|s_1-s_2| \geq 3L(M+1)$ then
property $(\alpha_2')$ implies that $\q([s_1,s_2])\cap \nn_M (A)
\neq \emptyset $. This contradicts the choice of $t$. Therefore
$|s_1-s_2| \leq 3L(M+1)$. We deduce that $\q ([0,t]) \cap
\overline{\nn}_M (A) $ has diameter at most $3L^2(M+1)$.

The definition of $t$ implies that there exists $t_1 >t$ with
$t_1-t\leq \frac{1}{L}$ and $\q (t_1) \in \nn_M (A)$. We take
$\q'=\q|_{[0,t_1]}$. The diameter of $\q' \cap \overline{\nn}_M
(A)$ is at most $\dg = 3L^2(M+1)+1$.\endproof

\subsubsection{Hyperbolicity of $\cgsh$}

Let $G$ be a finitely generated group that is asymptotically
tree-graded with respect to the finite collection of subgroups $\{
H_1,\dots , H_m \}$. This means that $\cgs$ is asymptotically
tree-graded with respect to the collection of subsets $\aaa = \{
gH_i \mid g\in G \, ,\, i=1,2,\dots , m \}$. We prove that $\cgsh
$ is hyperbolic, using Proposition \ref{hypb}. The following
result is central in the argument.

\begin{proposition}\label{int3}
Let $L\geq 1,\, C\geq 0$, let $\mu \geq M(L,C)$ and let
$\q_1,\q_2,\q_3$ be three $(L,C)$-almost-geodesics composing a
triangle in $\cgs$. We consider the set
$$
\mathcal{C}_\kappa ^\mu (\q_1,\q_2,\q_3)=\nn_\kappa
(\su{\q_1})\cap \nn_\kappa (\su{\q_2})\cap \nn_\kappa
(\su{\q_3})\, .
$$
\begin{itemize}
  \item[(1)] There exists $\kappa_0=\kappa_0(L,C,\mu )$ such that for every
  $\kappa \geq \kappa_0$ the set $\mathcal{C}_\kappa ^\mu
  (\q_1,\q_2,\q_3)$ intersects each of the almost-geodesics
  $\q_1,\q_2,\q_3$. In particular it is non-empty.
  \item[(2)] For every $\kappa \geq \kappa_0$ there exists $D_\kappa$ such that
the set $\mathcal{C}_\kappa ^\mu (\q_1,\q_2,\q_3)$ has diameter at
most $D_\kappa$ in $\cgsh$.
\end{itemize}
\end{proposition}

\me

\noindent\textit{Proof of (1).} Let $\{i,j,k \}=\{1,2,3 \}$.
According to Lemma \ref{plimrsat}, the result in Lemma
\ref{lqqcun} is true if we replace $\Sat (\q )$ by $\su{\q_i}\cup
\su{\q_j}$. In particular there exists $\tau =\tau (L,C, \mu ) $
such that $\q_k \subset \nn_\tau (\su{\q_i}) \cup \nn_\tau
(\su{\q_j})$. The traces on $\q_k$ of the two sets $\nn_\tau
(\su{\q_i})$ and $\nn_\tau (\su{\q_j})$ compose a cover of two
open sets, none of them empty. Since $\q_k$ is an almost geodesic,
it is connected, hence $\q_k \cap \nn_\tau (\su{\q_i})$ and $\q_k
\cap \nn_\tau (\su{\q_j})$ intersect. The intersection is in
$\mathcal{C}_\kappa ^\mu (\q_1,\q_2,\q_3)$ for every $\kappa \geq
\tau$.\hspace*{\fill }$\Box$

\me

We need several intermediate results before proving (2). In the
sequel we work with the data given in the statement of the
Proposition \ref{int3}, without mentioning it anymore.

\begin{lemma}\label{s1}
There exist positive constants $\alpha, \beta$ depending only on
$L,C,\mu$ and $\kappa $ such that every point $x\in
\mathcal{C}_\kappa ^\mu (\q_1,\q_2,\q_3)$ is in one of the two
situations:
\begin{itemize}
 \item[(i)] the ball $B(x, \alpha ) $ intersects each of the three almost-geodesics
  $\q_1,\q_2,\q_3$;
 \item[(ii)] $x\in \nn_\kappa (A)$ and $\nn_\beta (A)$
  intersects each of the three almost-geodesics $\q_1,\q_2,\q_3$.
\end{itemize}
\end{lemma}

\proof Let $x$ be an arbitrary point in $\mathcal{C}_\kappa ^\mu
(\q_1,\q_2,\q_3)$. The inclusion $x\in \nn_\kappa (\su{\q_i})$,
$i\in \{ 1,2,3\}$, implies that there are two possibilities:
\begin{itemize}
  \item[($I_i$)] $ x\in \nn_\kappa (\q_i)$ or
  \item[($II_i$)] $ x\in \nn_\kappa (A)$, where $A\in \aaa$, $\nn_\mu (A) \cap \q_i \neq
  \emptyset$.
\end{itemize}

If we are in case ($I$) for the three edges then this means that
(i) is satisfied with $\beta = \kappa $.

Suppose that only one edge is in case ($II$). Suppose it is
$\q_3$. Then $ x\in \nn_\kappa (\q_1) \cap \nn_\kappa (\q_2)$ and
there exists $A\in \aaa$ with $\nn_\mu (A) \cap \q_3 \neq
  \emptyset$ such that $ x\in \nn_\kappa (A)$. It follows that
  $\nn_\beta (A)$ intersects the three edges for
  $\beta = \max (\mu , 2 \kappa)$, so (ii) is satisfied.

Suppose that two edges are in case ($II$), for instance $\q_2$ and
$\q_3$. Consequently, $ x\in \nn_\kappa (\q_1)$ and $ x\in
\nn_\kappa (A_2)\cap \nn_\kappa (A_3)$, with $\nn_\mu (A_i) \cap
\q_i \neq \emptyset$, where $i=2,3$. If $A_2=A_3=A$ then
$\nn_\beta (A)$ intersects the three edges for
  $\beta = \max (\mu , 2 \kappa)$, so (ii) is satisfied. If $A_2\neq
  A_3$ then, according to Lemma \ref{plimrsat} (more precisely to
  Lemma \ref{gat} which also holds for $\mu$-saturations) we have that
  $x\in \nn_{\varkappa } (\q_2 \cup \q_3 )$,
   where $\varkappa =\varkappa (\mu , \kappa )$. Suppose that $x\in \nn_{\varkappa }
  (\q_2)$. Then $\nn_\beta (A_3)$ intersects the three edges for
  $\beta = \max (\mu , 2 \kappa , \kappa + \varkappa )$, so (ii) is satisfied.

Suppose that the three edges are in case ($II$). It follows that $
x\in \nn_\kappa (A_1)\cap \nn_\kappa (A_2)\cap \nn_\kappa (A_3)$,
with $\nn_\mu (A_i) \cap \q_i \neq \emptyset$, where $i=1,2,3$.

If the cardinality of the set $\{A_1,A_2,A_3 \}$ is 1 then we are in
situation (ii) with $\beta =\mu $. Suppose the cardinality of the
set is 2. Suppose that $A_1=A_2\neq A_3$. Lemma \ref{gat} for
$\mu$-saturations implies that $x\in \nn_{\varkappa } (\q_2 \cup
\q_3 ) \cap \nn_{\varkappa } (\q_1 \cup \q_3 )$. If $x\in
\nn_{\varkappa } (\q_3)$ then $A=A_1=A_2$ has the property that
$\nn_\beta (A)$ intersects the three edges for $\beta = \max (\mu ,
\kappa + \varkappa )$, and we are in case (ii). Otherwise $x\in
\nn_{\varkappa } (\q_1 ) \cap \nn_{\varkappa } (\q_2 )$, hence
$\nn_\beta (A_3)$ intersects the three edges for $\beta = \max (\mu
, \kappa + \varkappa )$.

Assume that the cardinality of the set $\{A_1,A_2,A_3 \}$ is 3. Then
$x\in \nn_{\varkappa } (\q_1 \cup \q_2 ) \cap \nn_{\varkappa } (\q_2
\cup \q_3 ) \cap \nn_{\varkappa } (\q_1 \cup \q_3 )$. It follows
that $x$ is in the $\varkappa$-tubular neighborhood of at least two
edges. Suppose these edges are $\q_1$ and $\q_2$. Then $\nn_\beta
(A_3)$ intersects the three edges for $\beta = \max (\mu , \kappa +
\varkappa )$.\endproof

\begin{lemma}\label{s3}
For every $r >0$ there exists $\varrho =\varrho (r , L,C)$ such
that the following holds. Let $A\neq B$ be such that $A,B \in
\aaa$, and both $\nn_r(A)$ and $\nn_r (B)$ intersect each of the
three almost-geodesic edges of the triangle. Then there exists $x$
such that $B(x,\varrho )$ intersects each of the edges of the
triangle.
\end{lemma}

\proof Let $y\in \nn_r (A)$ and $z\in \nn_r (B)$. Lemma
\ref{restr} implies that up to taking a subsegment of $[y,z]$, we
may suppose that the diameters of $[y,z]\cap \nn_r (A)$ and of
$[y,z]\cap \nn_r (B)$ are at most $\dg$, where $\dg = \dg (r)$. We
apply Lemma \ref{gat} for $r$-saturations and for each $\q_i$,
$i\in \{1,2,3\}$, and we obtain that both $B(y, \varrho )$ and
$B(z, \varrho )$ intersect $\q_i$, where $\varrho = \varrho (r
)$.\endproof

\begin{lemma}\label{s4}
There exists $R=R(L,C)$ such that for every triangle with
$(L,C)$-almost-geodesic edges, one of the following two situations
holds.
\begin{itemize}
  \item[(C)] There exists $x$ such that $B(x,R)$ intersects each
  of the three edges of the triangle;
  \item[(P)] There exists a unique $A\in \aaa$ such that $\nn_R (A)$ intersects each
  of the three edges of the triangle.
\end{itemize}
\end{lemma}

\proof Let $\q_1,\q_2,\q_3$ be the three edges. For $\mu =M(L,C) $
and $\kappa_0=\kappa_0(L,C)$ we have that $\mathcal{C}_\kappa^\mu
(\q_1,\q_2,\q_3)$ is nonempty. It remains to apply Lemmas \ref{s1}
and \ref{s3}.\endproof

\me

\noindent \textit{Notation}: We denote the vertices of the
triangle by $O_1,O_2,O_3$, such that $\q_i$ is opposite to $O_i$.

\me

\begin{lemma}\label{s2}
For every $r>0$ there exists $D=D(r, L,C)$ such that the following
holds. Let $x$ be such that $B(x,r)$ intersects the three edges of
the triangle.
\begin{itemize}
  \item[(a)] If $y$ is such that $B(y,r)$ intersects the three edges
  then $\dsh (x,y)\leq D$.
  \item[(b)] If $A\in \aaa$ is such that $\nn_r (A)$ intersects the three edges
  then $\dsh (x,A)\leq D$.
\end{itemize}
\end{lemma}

\proof Let $x_i$ be nearest points to $x$ in $\q_i,\, i=1,2,3$.

\sm

\textbf{(a)} We denote $\dsh (x,y)$ by $D$. Let $y_i$ be nearest points to
$y$ in $\q_i,\, i=1,2,3$. Then $\dsh (x_i,y_j)\geq D-2r$ for every
$i,j\in \{ 1,2,3 \}$. Suppose that $D>2r$. Without loss of
generality we may assume that $y_1\in \q_1[x_1,O_3]$. We have $\ds
(x_1,x_2)\leq 2r$, hence $\q_1[x_1,O_3]\subset \nn_{2\tau r}
\left( \Sat \left( \q_2[x_2,O_3]\right) \right)$, where $\tau
=\tau (L,C)$. In particular $y_1$ is contained either in
$\nn_{2\tau r} ( \q_2[x_2,O_3])$ or in $\nn_{2\tau r} (B)$ for
$B\in \aaa$ such that $\nn_M (B)$ intersects $\q_2[x_2,O_3]$.

\me

\textbf{Case (a)I.} Suppose that $y_2\in \q_2[x_2,O_1]$.

\textbf{Case (a)I.1.} Suppose that $y_1\in \nn_{2\tau
r}(\q_2[x_2,O_3])$. Then there exists $u\in \q_2[x_2,O_3]$ such
that $\ds(y_1,u)\leq 2\tau r $. It follows that $\ds(u, x_2)\geq
\dsh(u, x_2) \geq D-2r-2\tau r$. Inequality (\ref{qgeod}) implies
that
$$
\ds(u, y_2)\geq \frac{1}{L_1}[ \dist (u,x_2)+\dist
(x_2,y_2)]-C_1\geq \frac{1}{L_1}(2D-4r-2\tau r)-C_1 \, .
$$

On the other hand $\dist (u,y_2)\leq 2\tau r + 2r$. Hence $D\leq
2r+\tau r +L_1 (r+\tau r +C_1/2)$.

\me

\textbf{Case (a)I.2.} Assume that $y_1\in \nn_{2\tau r} (B)$,
where $B\in \aaa$ is such that $\nn_M (B)$ intersects
$\q_2[x_2,O_3]$. Let $w_2$ be a point in $\nn_M (B)\cap
\q_2[x_2,O_3]$.

Suppose that $\q_2[x_2,y_2]\cap \overline{\nn}_{2\tau r}(B)\neq
\emptyset $. Let $z_2$ be a point in the previous intersection.
Then $\q_2[w_2,z_2]$ has its endpoints in $\nn_{\chi }(B)$, with
$\chi =\max (M, 2\tau r +1)$. Consequently $\q_2[w_2,z_2]\subset
\nn_{\tau \chi }(B)$. In particular $x_2$ is contained in
$\nn_{\tau \chi}(B)$ and $\dsh (y_1,x_2)\leq \tau (2r+\chi) +1$,
hence $D\leq \tau (2r+\chi)+2r+1$.

Suppose that $\q_2[x_2,y_2]\cap \overline{\nn}_{2\tau
r}(B)=\emptyset $. We have that $x_2$ is in $\q_2[w_2,y_2]$. Also,
$\q_2[w_2,y_2]$ has its endpoints in $\nn_{\chi }(B)$, with $\chi
=\max (M, 2r (\tau +1))$. Consequently $\q_2[w_2,y_2]\subset
\nn_{\tau \chi }(B)$. In particular $x_2$ is contained in
$\nn_{\tau \chi}(B)$ and $\dsh (y_1,x_2)\leq \tau (2r+\chi) +1$,
hence $D\leq \tau (2r+\chi)+2r+1$.

\me

\textbf{Case (a)II.} Suppose that $y_2\in \q_2[x_2,O_3]$. If
$y_3\in \q_3[x_3,O_1]$ then we repeat the previous argument with
$y_1$ replaced by $y_3$. If $y_3\in \q_3[x_3,O_2]$ then we repeat
the previous argument with $(y_1,y_2)$ replaced by $(y_3,y_1)$.

\me

\textbf{(b)} We denote $\dsh (x,A)$ by $D$. We note that for every point
$y$ in $\nn_r(A) \cap (\q_1 \cup \q_2 \cup \q_3)$ we have that
$\ds (x_i,y)\geq \dsh (x_i,y)\geq D-2r$ for $i=1,2,3$. We choose
$y_i \in \nn_r(A) \cap \q_i,\, i=1,2,3$. Suppose $y_1\in
\q_1[x_1,O_3]$. Like in case (a), we have that $y_1$ is contained
either in $\nn_{2\tau r}(\q_2[x_2,O_3])$ or in $\nn_{2\tau r}(B)$
for some $B\in \aaa$ such that $\nn_M (B)$ intersects
$\q_2[x_2,O_3]$.

\me

\textbf{Case (b)I.} Suppose that $y_2\in \q_2[x_2,O_1]$.

\textbf{Case (b)I.1.} Assume that $y_1\in \nn_{2\tau r} (
\q_2[x_2,O_3])$. Then there exists $u\in \q_2[x_2,O_3]$ such that
$\ds(y_1,u)\leq 2\tau r $. It follows that $u\in \nn_{r(1+2\tau)}
(A)$ which together with $y_2\in \nn_r(A)$ implies that
$\q_2[u,y_2]\in \nn_{\tau r(1+2\tau)} (A)$. In particular $x_2\in
\nn_{\tau r(1+2\tau)} (A)$, therefore $D\leq r + \tau r(1+2\tau)$.

\me

\textbf{Case (b)I.2} Suppose $y_1\in \nn_{2\tau r} (B)$, with
$B\in \aaa$ such that $\nn_M (B)$ intersects $\q_2[x_2,O_3]$. Let
$w_2$ be a point in $\nn_M (B)\cap \q_2[x_2,O_3]$.

Suppose that $\q_2[x_2,y_2]\cap \overline{\nn}_{2\tau r}(B)\neq
\emptyset $. As in the proof of part (a), Case I.2, we obtain that
$\dsh (y_1,x_2)\leq \tau (2r+\chi)+1$, whence $D\leq \tau
(2r+\chi)+2r+1$.

Suppose that $\q_2[x_2,y_2]\cap \overline{\nn}_{2\tau
r}(B)=\emptyset $. Then $x_2$ is in $\q_2[w_2,y_2]$. On the other
hand, $\q_2[w_2,y_2]$ has its endpoints in the $M$-tubular
neighborhood of $\Sat^{2\tau r}([y_1,y_2])$. It follows that
$\q_2[w_2,y_2]$, in particular $x_2$, is in the $tM$-tubular
neighborhood of $\Sat^{2\tau r}([y_1,y_2])$. In $\cgsh$,
$\Sat^{2\tau r}([y_1,y_2])$ is contained in the $(2\tau
r+1)$-tubular neighborhood of $[y_1,y_2]$. Since in $\cgs$ we have
that $[y_1,y_2]\subset \nn_{\tau r} (A)$, we deduce that in
$\cgsh$, $x_2$ is in the $(tM+3\tau r+1)$-tubular neighborhood of
$A$. Hence $D\leq tM+(3\tau +1)r+1$.

\me

\textbf{Case (b)II.} Suppose that $y_2\in \q_2[x_2,O_3]$. Then we
can use the same argument as in Case II of part (a).\endproof

\me

\noindent\textit{Proof of Proposition \ref{int3}, (2).} By Lemma
\ref{s4} we are either in case (C) or in case (P).

\me

\textbf{Case (C)}. Let $y\in \mathcal{C}_\kappa ^\mu
(\q_1,\q_2,\q_3)$. According to Lemma \ref{s1} we have either (i)
or (ii). Suppose that (i) is satisfied. Then, by Lemma \ref{s2},
(a), $\dsh (x,y)\leq D$, where $D=D(\alpha , R,L,C )$.

Suppose that (ii) is satisfied, that is $y\in \nn_\kappa (B)$ and
$\nn_\beta (B)$ intersects each of the three almost-geodesics
$\q_1,\q_2,\q_3$. Lemma \ref{s2}, (b), implies that $\dsh
(x,B)\leq D$, where $D=D(\beta , R,L,C)$. Therefore $\dsh
(x,y)\leq D + \kappa +1$.

\me

\textbf{Case (P)}. Let $y\in \mathcal{C}_\kappa ^\mu
(\q_1,\q_2,\q_3)$. Suppose that $y$ satisfies (i). Lemma \ref{s2},
(b) implies that $\dsh (y,A)\leq D$, with $D=D(\alpha , R, L,C)$.

If $y$ satisfies (ii) of Lemma \ref{s1}, then the unicity stated
in (P) implies that $y\in \nn_\kappa (A)$, hence that $\dsh
(y,A)\leq \kappa$.

We may conclude that in all cases the diameter of the set
$\mathcal{C}_\kappa ^\mu (\q_1,\q_2,\q_3)$ in the metric $\dsh$ is
uniformly bounded.\hspace*{\fill }$\Box$

\me

We now define a system of lines and centers in $\cgsh$ such that
the properties in Section \ref{carhip} are satisfied.

First of all, for every pair of vertices $u,v$ in $\cgsh$ we
choose and fix a geodesic $[u,v]$ in $\cgs $ joining the two
points. Let $M_0=M(1,0)$ and let $\kappa_0$ be the constant given
by Proposition \ref{int3} for $\mu=M_0$. We may suppose that
$\kappa_0\geq M_0$. For every pair of vertices $u,v$ in $\cgsh$,
we define $\La{uv}$ as $\nn_{\kappa_0} (\Sat ([u,v]))$. The
relation on it is defined as follows: to every $x\in
\nn_{\kappa_0} (\Sat ([u,v]))$ we associate one nearest point
(projection) $x'\in [u,v]$ and we put $x\luv y$ if $x'$ is between
$u$ and $y'$. Properties ($l_1$), ($l_2$), ($l_3$) are obviously
satisfied.

We define the function $\phi $ by choosing, for every three
vertices $u,v,w$ in $\cgs$ a point $C_{uvw}$ in
$\mathcal{C}^{M_0}_{\kappa_0}([u,v],[u,w],[v,w])$ and defining
$\phi (u,v,w)=\phi \circ \sigma (u,v,w)=C_{uvw}$ for every
$3$-permutation $\sigma$. We choose $C_{uuv}=u$.

Properties ($c_1$), ($c_2$), ($c_3$) are satisfied. Before
proceeding further, we prove some intermediate results.

\begin{lemma}\label{pic}
For every $\alpha >0$ there exists $\lambda =\lambda(\alpha )$
such that the following holds. Let $[u,v]$ be a geodesic and let
$A\in \aaa$ be such that $\nn_\alpha (A) \cap [u,v]\neq
\emptyset$. Let $x$ be a point in $\nn_{\alpha} (A)$ and let
$x'\in [u,v]$ be a projection of $x$. Then $x' \in \nn_\lambda
(A)$.
\end{lemma}
\begin{figure}[!ht]
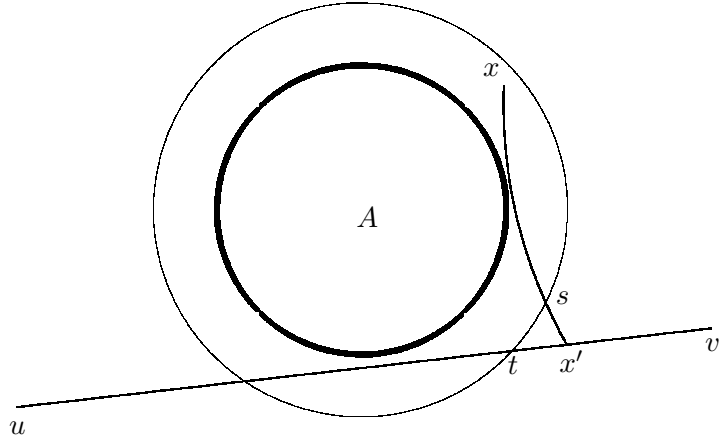

\centering
\unitlength .6mm 
\linethickness{.05pt}
\ifx\plotpoint\undefined\newsavebox{\plotpoint}\fi 

\caption{Projection of a point onto the saturation.} \label{fig9}
\end{figure}

\proof Suppose that $x'\not \in \nn_\alpha (A)$. Lemma \ref{restr}
implies that there exist $t\in [u,v]\cap \nn_\alpha (A)$ and $s\in
[x',x]\cap \nn_\alpha (A)$ such that the sets $[x',t] \cap
\overline{\nn}_\alpha (A)$ and $[x',s] \cap \overline{\nn}_\alpha
(A)$ have diameters at most $\dg$, where $\dg =\dg (\alpha )$.
Corollary \ref{12} implies that $\dist (s,t)\leq D_1$. On the
other hand, since $\dist (x,x')\leq \dist (x, t)$, it follows that
$\dist (s,x')\leq \dist (s, t)\leq D_1$. We conclude that $\ds
(x',A)\leq D_1+\alpha $.\endproof

\begin{cor}\label{c1pic}
Let $x$ be a point in $\nn_{\kappa} (\Sat^\mu ([u,v]))$ and let
$x'\in [u,v]$ be a projection of $x$. Then $\dsh (x,x')\leq \chi
$, where $\chi =\chi (\kappa , \mu )$.
\end{cor}

\proof Since $x\in \nn_{\kappa}(\Sat^\mu([u,v]))$ it follows that
either $x\in \nn_{\kappa} ([u,v])$ or $x\in \nn_{\kappa} (A)$,
where $\nn_\mu (A) \cap [u,v]\neq \emptyset$. In the first case it
follows that $\dsh (x,x') \leq \kappa$, while in the second case
we may apply Lemma \ref{pic}.\endproof

\begin{cor}\label{c2pic}
Let $u,v$ be a pair of vertices in $\cgsh$ and let $x,y \in
\La{uv}$ and $x',y'$ their chosen respective projections on
$[u,v]$. Then, in $\cgsh$, $\La{uv} [x,y]= \La{uv} [y,x]$ is at
Hausdorff distance $\chi$ of $[x',y']\subset [u,v]$, where
$\chi=\chi(G)$.
\end{cor}

Before proving properties (I), (II), (III), we make some remarks
and introduce some notations.

\begin{remarks} \label{rem999}
1) For every quasi-geodesic $\q$ in $\cgs$, we have that $\Sat^\mu
(\q )$ is in the $(\mu+1)$-tubular neighborhood of $\q$ in
$\cgsh$.

\me

2) Lemma \ref{s1} implies that there exist two constants $\eta$
and $c$ such that for every three geodesics $[u,v], [v,w], [u,w]$
in $\cgs$ every point $x \in \calc_{\kappa_0}^{M_0}([u,v], [v,w],
[u,w])$ satisfies one of the following two properties:
\begin{itemize}
 \item[(i)] the ball $B(x, \eta) $ intersects each of the three geodesics
  $[u,v], [v,w], [u,w]$;
 \item[(ii)] $x\in \nn_{\kappa_0} (A)$ and $\nn_c (A)$
  intersects each of the three geodesics $[u,v], [v,w], [u,w]$.
\end{itemize}

We note that the constants $\eta$ and $c$ depend on $M_0$ and
$\kappa_0$, so they depend only on $G$. We may suppose without
loss of generality that $c\geq M_0$.

\me

3) Lemma \ref{pic} implies that there exists $\xi$ such that if
$[u,v]$ is a geodesic, $A\in \aaa$ is such that $\nn_c (A)$
intersects $[u,v]$ and $x$ is a point in $\nn_{\kappa_0} (A)$,
then any projection of $x$ on $[u,v]$ is in $\nn_\xi (A)$. The
constant $\xi$ depends on $\max (c, \kappa_0)$, so it depends only
on $G$. Without loss of generality we may suppose that $\xi \geq
M_0$.

\me

4) In the sequel we denote the constant $\dg (1,0,c)$ provided by
Lemma \ref{restr} simply by $\dg$.
\end{remarks}

\me

\noindent\textit{Proof of properties} (I), (II), (III).

\quad \textbf{(I).} Let $x=\phi (u,v,w)$ and let $x_1$ and $x_2$
be the chosen projections of $x$ on $[u,v]$ and on $[u,w]$,
respectively. According to Corollary \ref{c2pic}, it suffices to
prove that $[u, x_1 ]$ and $[u,x_2]$ are at uniformly bounded
Hausdorff distance in $\cgsh$. The point $x=\phi (u,v,w)$
satisfies either (i) or (ii) from Remark \ref{rem999}, part 2.

Suppose $x$ is in case (ii). Then $x\in \nn_{\kappa_0}(A)$ such
that $\nn_{c}(A)$ intersects the three geodesic edges. Lemma
\ref{pic} implies that $x_1,x_2 \in \nn_\xi(A)$. The geodesic
$[u,x_1]$ has its endpoints in $\nn_\xi (\Sat^\xi [u,x_2])$. Lemma
\ref{plqqcun} implies that $[u,x_1]$ is entirely contained in
$\nn_{\tau \xi }(\Sat^\xi [u,x_2])$. It follows that $[u,x_1]$ is
in the $[(\tau +1)\xi +1]$-tubular neighborhood of $[u,x_2]$ in
$\cgsh$. A similar argument done for $[u,x_2]$ allows to conclude
that (I) is satisfied.

Suppose $x$ is in case (i). Then $\ds (x,x_i)\leq \eta$ for
$i=1,2$. Hence $\ds (x_1,x_2)\leq 2\eta $ and $[u,x_i]$ has its
endpoints in $\nn_{2\eta } (\Sat [u,x_j])$, for $\{ i,j \}=\{ 1,2
\}$. We repeat the previous argument.

\me

\textbf{(II)} The fact that $\dsh (p,q)\leq 1$ means that either
$\ds (p,q)\leq 1$ or $p,q\in A_0$, where $A_0\in \aaa $. Let
$x=\phi (u,v,p)$ and $y=\phi (u,v,q)$. We have to show that
$\La{uv}[x,y]$ has uniformly bounded diameter in $\cgsh$. Let
$x_0$ and $y_0$ be the respective projections of $x$ and $y$ on
$[u,v]$. Corollary \ref{c2pic} implies that it suffices to prove
that $[x_0,y_0]$ has uniformly bounded diameter in $\cgsh$, where
by $[x_0,y_0]$ we denote the sub-arc of $[u,v]$ of endpoints
$x_0,y_0$.

\me

Suppose that both $x$ and $y$ are in case (i). We have that $x_0
\in \nn_{2\eta }[u,p]\cap \nn_{2\eta }[v,p]$ and that $y_0 \in
\nn_{2\eta }[u,q]\cap \nn_{2\eta }[v,q]$. Since $[u,p]\subset
\nn_{\tau }(\Sat [u,q])$ and $[v,p]\subset \nn_{\tau }(\Sat
[v,q])$, we conclude that $x_0,y_0 \in \calc_{2\eta +\tau
}^{M_0}([u,q],[v,q],[u,v])$, hence that $[x_0,y_0]\subset
\calc_{\tau (2\eta +\tau )}^{M_0}([u,q],[v,q],[u,v])$. We complete
the proof by applying Proposition \ref{int3}.

\me

Suppose $x$ is in case (i) and $y$ is in case (ii). The case when
$x$ is in case (ii) and $y$ is in case (i) is discussed similarly.
As above we have that $x_0 \in \calc_{2\eta +\tau
}^{M_0}([u,q],[v,q],[u,v])$. We have that $y\in \nn_{\kappa_0}
(A)$ such that $\nn_{c} (A)$ intersects $[u,q],[v,q],[u,v]$. Lemma
\ref{pic} implies that $y_0\in \nn_\xi (A)$. Then $y_0\in
\calc_{\xi }^{c}([u,q],[v,q],[u,v])$. As previously we obtain that
$[x_0,y_0]\subset \calc_{\tau' r }^{s}([u,q],[v,q],[u,v])$, where
$r=\max (2\eta +\tau \, ,\,  \xi)$, $s=\max (M_0, c)$ and $\tau'=\tau' (s)$. Proposition \ref{int3} allows to complete the
argument.

\me

Suppose that both $x$ and $y$ are in case (ii). Then $x\in
\nn_{\kappa_0} (A)$ such that $\nn_{c} (A)$ intersects $[p,u],
[p,v],[u,v]$. Let $p_1\in [u,p]\cap \nn_{c} (A)$ and $p_2\in
[v,p]\cap \nn_{c} (A)$ be such that $[p,p_i] \cap
\overline{\nn}_{c} (A)$ has diameter at most $\dg$, $i=1,2$.
Likewise we consider $u_1\in [u,v]\cap \nn_{c} (A)$ and $u_2\in
[u,p]\cap \nn_{c} (A)$ so that $[u,u_i] \cap \overline{\nn}_{c}
(A)$ has diameter at most $\dg$, and $v_1\in [p,v]\cap \nn_{c}
(A)$ and $v_2\in [u,v]\cap \nn_{c} (A)$ so that $[v,v_i] \cap
\overline{\nn}_{c} (A)$ has diameter at most $\dg$. Corollary
\ref{12} implies that $\ds (p_1,p_2), \ds (u_1,u_2)$ and $\ds
(v_1,v_2)$ are at most $\zeta $, where $\zeta =\zeta (G)$.

We have that either $A\subset \Sat [u,q ]$ or $\nn_c(A) \cap
\nn_{\tau }(\Sat [u,q])$ has diameter at most $\gamma$, where
$\gamma = \gamma (G) $. The latter case implies, together with the
inclusion $[u,p]\subset \nn_{\tau }(\Sat [u,q])$, that $\dist
(p_1,u_2)\leq \gamma$. Thus, we have that either $A\subset \Sat
[u,q ]$ or $\dist (p_1,u_2)\leq \gamma$. Likewise, we obtain that
either $A\subset \Sat [v,q ]$ or $\dist (p_2,v_1)\leq \gamma$.

Suppose that $\dist (p_1,u_2)\leq \gamma$. Then $\dist
(p_1,u_1)\leq \gamma+\zeta$, hence $B(p_1, \gamma+\zeta )$
intersects  $[p,u], [p,v],[u,v]$. We can argue similarly to the
case above when $x$ is in case (i) and $y$ is in case (ii), with
$x$ replaced by $p_1$ and $\eta$ by $\gamma+\zeta $. We obtain
that if $p_1'$ is the chosen projection of $p_1$ on $[u,v]$ then
$[p_1',y_0]$ has the diameter bounded in $\cgsh$ by a constant
depending on $G$. Since $[x_0,y_0]\subset [x_0,p_1']\cup
[p_1',y_0]$, it remains to prove that $[x_0,p_1']$ has bounded
diameter in $\cgsh$. Lemma \ref{pic} provides for $\alpha = \max
(\kappa_0 , c )$ a constant $\tilde{\lambda} $. We have that $x_0$
and $p_1'$ are in $\nn_{\tilde{\lambda}} (A)$, hence that
$[x_0,p_1']\subset \nn_{\tau \tilde{\lambda} }(A)$. We conclude
that the diameter of $[x_0,p_1']$ in $\cgsh$ is at most $2\tau
\tilde{\lambda} +1$. A similar argument works if $\dist
(p_2,v_1)\leq \gamma$.

Now suppose that $A\subset \Sat [u,q ] \cap \Sat [v,q ]$. Lemma
\ref{pic} implies that $x_0 \in \nn_\xi (A)$. Since $y$ is also in
case (ii), we have that $y\in \nn_{\kappa_0} (B)$ such that
$\nn_{c} (B)$ intersects the three geodesic edges $[q,u],
[q,v],[u,v]$ and that $y_0 \in \nn_\xi (B)$. We have that $A\cup B
\subset \Sat^c [u,q ] \cap \Sat^c [v,q ]\cap \Sat^c [u,v]$. Lemma
\ref{plqqcun} implies that $[x_0,y_0]\subset \calc_{ \tau \xi }^c
([q,u], [q,v],[u,v])$ and Proposition \ref{int3} allows to finish
the argument.

\me

\textbf{(III)} Let $u,v,w$ be three vertices such that $w\in
\nn_{\kappa_0}(\Sat [u,v])$. Let $x=\phi (u,v,w)$. Let $w_0$ and
$x_0$ be the projections of $w$ and $x$ respectively on $[u,v]$.
We bound the diameter of $[x_0,w_0]$ in $\cgsh$.

We have $x,w \in \calc_{\kappa_0}^{M_0} ([u,v],[u,w],[v,w])$.
Suppose both $x$ and $w$ are in case (i). Then $x_0,w_0 \in
\calc_{\kappa_0+\eta }^{M_0} ([u,v],[u,w],[v,w])$, consequently
$[x_0,w_0] \subset  \calc_{\tau (\kappa_0+\eta )}^{M_0}
([u,v],[u,w],[v,w])$ and we apply Proposition \ref{int3} to obtain
the conclusion.

Suppose that $x$ is in case (i) and $w$ in case (ii). The case
when $x$ is in case (ii) and $w$ in case (i) is similar. The ball
$B(x, \eta )$ intersects the three edges and $w\in \nn_{\kappa_0}
(A)$ such that $\nn_{c} (A)$ intersects the three edges. Lemma
\ref{pic} implies that $w_0\in \nn_\xi (A) \subset \calc_\xi^c
([u,v], [u,w], [v,w])$. The point $x_0$ is in $\calc_{\eta
+\kappa_0}^{M_0} ([u,v], [u,w], [v,w])$. It follows that
$[x_0,w_0]\subset \calc_{\tau' r }^{s}([u,v],[u,w],[v,w])$, where
$r=\max (\eta +\kappa_0\, ,\, \xi)$, $s=\max (M_0, c)$ and $\tau'
= \tau' (s)$. We apply Proposition \ref{int3}.

Suppose that $x$ and $w$ are both in case (ii). We have that $x\in
\nn_{\kappa_0} (A)$ and $w\in \nn_{\kappa_0} (B)$ such that both
$\nn_{c} (A)$ and $\nn_{c} (B)$ intersect the three edges. We also
have that $x_0\in \nn_\xi (A)$ and $w_0\in \nn_\xi (B)$, hence
$[x_0,w_0]\subset \calc_{\tau \xi }^{c}([u,v],[u,w],[v,w])$. We
end the proof by applying Proposition \ref{int3}.\hspace*{\fill
}$\Box$

\me

Proposition \ref{hypb} implies that $\cgsh$ is hyperbolic. Moreover
we have that $\La{uv}$ is at bounded Hausdorff distance from every
geodesic connecting $u$ and $v$ in $\cgsh$. Since in the previous
argument the choice of the geodesics $[u,v]$ in $\cgs$ was
arbitrary, we have the following.

\begin{proposition}\label{udist}
Every geodesic in $\cgs$ joining two points $u$ and $v$ is at
bounded Hausdorff distance in $\cgsh$ from any geodesic joining
$u$ and $v$ in $\cgsh$.
\end{proposition}

\subsection{The BCP Property}

Given two vertices $u,v$ in $\cgsh$, we denote by $[u,v]$ a
geodesic joining them in $\cgs$ and by $\g_{uv}$ a geodesic
joining them in $\cgsh$.

\begin{definition}\label{lift}
For a path $\pgot$ in $\cgsh$, we denote by $\tilde{\pgot}$ a path
in $\cgs$ obtained by replacing every $\mh$-component $s$ in
$\pgot$ by a geodesic in $\cgs$ connecting $s_-$ and $s_+$. We
call $\tilde{\pgot}$ a \textit{lift} of $\pgot$.
\end{definition}

We now prove the following.

\begin{proposition}\label{pbcp}
If $G$ is asymptotically tree-graded with respect to $\{ H_1,\dots
, H_m \}$ and $G$ is weakly hyperbolic relative to $\{ H_1,\dots ,
H_m \}$ then the pair $(G,\{ H_1,\dots , H_m \})$ satisfies the
BCP-property.
\end{proposition}

\proof  Let $\lambda \geq 1$. Let $\pgot$ and $\q$ be two
$\lambda$-bi-Lipschitz paths without backtracking in $\cgsh$ such
that $\pgot_-=\q_-$ and $\ds (\pgot_+,\q_+)\leq 1$.

\me

\textbf{(1)} Let $s$ be an $\mh$-component of $\pgot$ contained in
a left coset $A\in \aaa$, and let $\ds (s_-,s_+)=D$. We show that
if $D$ is large enough then $\q$ has an $\mh$-component contained
in $A$.

\me

\noindent\textit{Notations}: In this section $M$ denotes
$M(\lambda , 0)$, the constant given by $(\alpha_2')$ for
$\theta=\frac13$ and $(L,C)=(\lambda ,0)$.

The graph $\cgsh$ is hyperbolic. Therefore for the given $\lambda$
there exists $\varkappa = \varkappa (\lambda )$ such that two
$\lambda$-bi-Lipschitz paths $\pgot$ and $\q$ in $\cgsh$ with
$\dsh (\pgot_-,\q_-)\leq 1$ and $\dsh (\pgot_+,\q_+)\leq 1$ are at
Hausdorff distance at most $\varkappa$.

\me

\textbf{Step I}. We show that for $D\geq D_0(G)$, some lift
$\tilde{\q}$ of $\q$ intersects $\nn_{M'}(A)$, where $M'=M'(G)$.

We choose $u$ on the arc $\pgot[\pgot_-,s_-]$ such that either the
length of $\pgot[u,s_-]$ is $2 \lambda (\varkappa + 1)$ or, if the
length of $\pgot[\pgot_-,s_-]$ is less than $2\lambda (\varkappa
+1)$, $u=\pgot_-$. Likewise we choose $v$ on the arc
$\pgot[s_+,\pgot_+]$ such that either the length of $\pgot[s_+,v]$
is $2\lambda (\varkappa + 1)$ or $v=p_+$. We have that $\dsh(u,
s_-),\, \dsh (s_+, v)\in [2(\varkappa + 1)\, ,\, 2\lambda^2
(\varkappa + 1)]$, in the first cases.

There exist $w$ and $z$ on $\q$ such that $\dsh (u,w)\leq
\varkappa$ and $\dsh (v,z)\leq \varkappa$. We consider $\g_{uw}$
and $\g_{vz}$ geodesics in $\cgsh$.

Let $u'$ be the farthest from $u$ point on $\g_{uw}$ which is
contained in the same left coset $B\in \aaa$ as an $\mh$-component
$\sigma $ of $\pgot[u,v]$. Suppose that $\sigma \cap
\pgot[s_-,v]\neq \emptyset $. We have that
$$
\dsh (u,u')\geq \dsh (u,\sigma_+)-1\geq \frac{1}{\lambda }
\mathrm{length}\, (\pgot[u,\sigma_+])-1\geq \frac{1}{\lambda }
\mathrm{length}\, (\pgot[u,s_-])-1\geq 2\varkappa +1 \, .
$$

This contradicts the inequality $\dsh (u,u')\leq \varkappa$.
Therefore $\sigma $ is contained in $\pgot[u,s_-]\setminus \{ s_-
\}$. We choose $v'$ the farthest from $v$ point on $\g_{vz}$
contained in the same left coset as a component $\sigma' $ of
$\pgot[u,v]$. In a similar way we prove that $\sigma' $ is
contained in $\pgot[s_+,v]\setminus \{ s_+ \}$. It is possible
that $u'=u,\, \sigma =\{ u \}$ and/or $ v'=v$, $\sigma' =\{ v \}$.
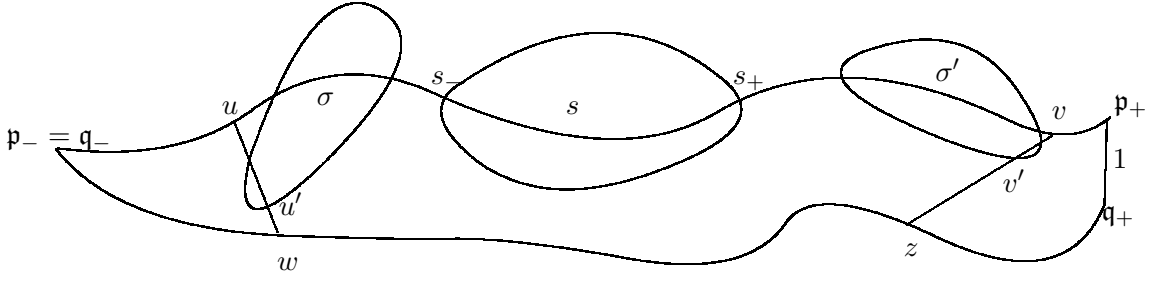
\begin{figure}[!ht]
\centering
\unitlength .75mm 
\linethickness{0.4pt}
\ifx\plotpoint\undefined\newsavebox{\plotpoint}\fi 
\begin{picture}(195.75,48.63)(0,0)
\qbezier(5.5,25.5)(26.5,22.75)(39.5,32)
\qbezier(39.5,32)(54.88,43.5)(72.75,35)
\qbezier(72.75,35)(104.38,21)(122.5,32)
\qbezier(122.5,32)(144.75,45)(175,30)
\qbezier(175,30)(185.5,25.5)(192,31)
\qbezier(5.25,25.25)(16.13,11.38)(45.5,10)
\qbezier(45.5,10)(68.63,9.25)(76.25,9.5)
\qbezier(76.25,9.5)(88,9.38)(107.75,5.75)
\qbezier(107.75,5.75)(127.88,2.5)(134.5,12.25)
\qbezier(134.5,12.25)(139.88,19.88)(160.75,10)
\qbezier(160.75,10)(184.13,-1.25)(191,15.5)
\multiput(191,15.5)(.033333,1.016667){15}{\line(0,1){1.016667}}
\qbezier(39.75,24.5)(51.75,58.63)(63.75,49.25)
\qbezier(63.75,49.25)(72.38,42.5)(52.5,22.75)
\qbezier(52.5,22.75)(34.88,5.75)(39.75,24.75)
\qbezier(77.25,37.25)(102.88,54.38)(123,38)
\qbezier(123,38)(133.75,29.63)(112.5,21.75)
\qbezier(112.5,21.75)(92.38,14)(79.75,23.25)
\qbezier(79.75,23.25)(68.5,32.25)(77.25,37.25)
\qbezier(151,43.5)(166.88,48.13)(175.25,36.25)
\qbezier(175.25,36.25)(188.5,17.38)(162.75,27)
\qbezier(162.75,27)(133.13,38.75)(151,43.5)
\put(5.75,27){\makebox(0,0)[cc]{$\pgot_-=\q_-$}}
\put(36,32.75){\makebox(0,0)[cc]{$u$}}
\put(46.25,5){\makebox(0,0)[cc]{$w$}}
\multiput(36.75,30.25)(.03369565,-.08586957){230}{\line(0,-1){.08586957}}
\put(74.5,37.5){\makebox(0,0)[cc]{$s_-$}}
\put(128,37.5){\makebox(0,0)[cc]{$s_+$}}
\put(96.75,32.75){\makebox(0,0)[cc]{$s$}}
\put(53,34.25){\makebox(0,0)[cc]{$\sigma$}}
\put(163.25,39.5){\makebox(0,0)[cc]{$\sigma'$}}
\put(183,31.75){\makebox(0,0)[cc]{$v$}}
\put(156.75,7.5){\makebox(0,0)[cc]{$z$}}
\multiput(156,12)(.055139186,.03372591){467}{\line(1,0){.055139186}}
\put(193.75,23.75){\makebox(0,0)[cc]{$1$}}
\put(195.75,33){\makebox(0,0)[cc]{$\pgot_+$}}
\put(193.5,13.5){\makebox(0,0)[cc]{$\q_+$}}
\put(46.75,15.75){\makebox(0,0)[cc]{$u'$}}
\put(175,19.75){\makebox(0,0)[cc]{$v'$}}
\end{picture}

\caption{Proof of (1) in BCP Property.} \label{fig10}
\end{figure}

We consider the path in $\cgsh$ defined as $\rrr=\g_{wu'}\sqcup
\g_{u'\sigma_+}\sqcup \pgot[\sigma_+, \sigma_-']\sqcup
\g_{\sigma_-'v'}\sqcup \g_{v'z}$, where $\g_{wu'}$ and $\g_{v'z}$
are sub-geodesics of $\g_{wu}$ and $\g_{vz}$, respectively, and
$\g_{u'\sigma_+}$ and $\g_{\sigma_-'v'}$ are composed of one edge.
The length of $\rrr$ is at most $N= \lambda (4\varkappa
+5)+2\varkappa $. It contains the component $s$. We show that it
has no backtracking. By construction and the fact that geodesics
do not have backtracking (\cite[Lemma 2.23]{Osin}), we have that
the sub-arcs $\rrr[w,v']$ and $\rrr[u',z]$ do not have
backtracking. Suppose that $\g_{wu'}$ and $\g_{v'z}$ have
$\mh$-components in the same left coset. It follows that there
exists $x\in \g_{wu'}$ and $y\in \g_{v'z}$ with $\dsh (x,y)\leq
1$. Then $\dsh (u,v)\leq 2\varkappa +1$. By construction either
$\mathrm{length}\, \pgot[u,v]\geq 2\lambda (\varkappa +1 )+1$ or
$u=\pgot_-$ and $v=\pgot_+$. In the latter case, the geodesic
$\g_{wu}$ is trivial, $\g_{vz}$ is an edge $e$ in $\cgs$, and
$\rrr=\pgot \cup e$ has no backtracking. In the former case we
have that $\dsh (u,v) > 2\varkappa +2$, which contradicts the
previous inequality.

We conclude that $\rrr$ is without backtracking. A lift
$\tilde{\rrr}$ of it is composed of $n$ consecutive sub-paths,
\begin{equation}\label{rrr}
\tilde{\rrr}=\tilde{\rrr}_1\cup \dots \cup \tilde{\rrr}_n \, ,
\end{equation} with
$n\leq N$, such that each $\tilde{\rrr}_i$ is either
\begin{itemize}
  \item[(R$_1$)] a $\lambda$-bi-Lipschitz arc both in $\cgs$ and in $\cgsh$
  of length at most $N$ or
  \item[(R$_2$)] a geodesic in $\cgs $ with
endpoints in some left coset $A_i\in \aaa$.
\end{itemize}

Since $\rrr$ is without backtracking, we have that $A_i\neq A_j$
when $i\neq j$. Lemma \ref{polig} implies that $\tilde{\rrr}$ is
an $(L_N,C_N)$-almost geodesic.

On the other hand, $\dsh (w,z)\leq \mathrm{length}\, \rrr \leq N$.
Hence the length of $\q[w,z]$ is at most $N_1$, where $N_1=\lambda
N$. As above, a lift $\tilde{\q}[w,z]$ decomposes into $m$
consecutive sub-paths,
\begin{equation}\label{tq}
\tilde{\q}[w,z]=\tilde{\q}_1\cup \dots \cup \tilde{\q}_m \, ,
\end{equation} with
$m\leq N_1$, such that each $\tilde{\q}_i$ is either
\begin{itemize}
  \item[(Q$_1$)] a $\lambda$-bi-Lipschitz arc both in $\cgs$ and $\cgsh$, of length at most $N_1$ or
  \item[(Q$_2$)] a geodesic in $\cgs $ with
endpoints in some left coset $B_i\in \aaa$.
\end{itemize}

Since $\q$ is without backtracking, we have that $B_i\neq B_j$
when $i\neq j$. Lemma \ref{polig} implies that $\tilde{\q}[w,z]$
is an $(L_{N_1},C_{N_1})$-almost geodesic. We denote $L'=\max
(L_N,L_{N_1})$ and $C'=\max (C_N,C_{N_1})$. We denote
$M'=M(L',C')$. Lemma \ref{lqqcun} implies that in $\cgs$ the path
$\tilde{\rrr}$ is contained in the $\tau'$-tubular neighborhood of
$\Sat (\tilde{\q}[w,z])=\Sat^{M'} (\tilde{\q}[w,z])$, where $\tau'
=\tau'(L',C')$. In particular the component $s$ is contained in
$\nn_{\tau'} (\Sat (\tilde{\q}[w,z]))$, hence the set $\nn_{\tau'}
(\Sat (\tilde{\q}[w,z]))\cap A$ has diameter at least $D$. Lemma
\ref{diam} implies that for $D\geq D_0(L',C',\tau ')$ we must have
that $\nn_{M'}(A) \cap \tilde{\q}[w,z] \neq \emptyset$.

\me

\textbf{Step II}. We show that there exist two points $w_1$ and
$z_1$ on $\tilde{\q}[w,z]$ such that $\ds (w_1,s_-)\leq D_1$ and
$\ds (z_1,s_+)\leq D_1$, where $D_1=D_1(G)$. We do this by means
of Corollary \ref{12}.

Lemma \ref{restr} implies that there exist $w_1,z_1\in
\tilde{\q}[w,z]\cap \nn_{M'}(A)$ such that $\tilde{\q}[w,w_1]$ and
$\tilde{\q}[z_1,z]$ intersect $\overline{\nn}_{M'}(A)$ in two sets
of diameter at most $\dg_1$, where $\dg_1 = \dg_1(L',C',M')$.

We show that $\tilde{\rrr}[w,s_-]$ and $\tilde{\rrr}[s_+,z]$
intersect $\overline{\nn}_{M'}(A)$ in two sets of bounded
diameter. We prove it only for $\tilde{\rrr}[w,s_-]$, the same
argument works for $\tilde{\rrr}[s_+,z]$. Let $x\in
\tilde{\rrr}[w,s_-]\cap \overline{\nn}_{M'}(A)$ and let $\ds(x,
s_-)= \delta$. According to the decomposition (\ref{rrr}), we have
that $\tilde{\rrr}[x,s_-]= \tilde{\rrr}_i'\cup
\tilde{\rrr}_{i+1}\cup \dots \cup \tilde{\rrr}_j$, where $i\leq
j,\, i,j \in \{1,2,\dots ,n \}$ and $\tilde{\rrr}_i'$ is
eventually a restriction of $\tilde{\rrr}_i$ such that $x$ is an
endpoint of it. If all the components are of type (R$_1$), then
$\tilde{\rrr}[x,s_-]$ has length at most $N$ and $\delta \leq N$.
Suppose that at least one component is of type (R$_2$). We have at
most $N$ such components. Then at least one component
$\tilde{\rrr}_{k}$ of type (R$_2$) has the distance between its
endpoints at least $\frac{\delta -N}{N}$. On the other hand since
$x,s_-\in \nn_{M'+1}(A)$ and $\tilde{\rrr}[x,s_-]$ is an
$(L',C')$-almost-geodesic, it follows that
$\tilde{\rrr}[x,s_-]\subset \nn_{\tau' (M'+1)}(A)$. In particular
$\tilde{\rrr}_{k}$ is contained in the same tubular neighborhood,
therefore the diameter of $A_k\cap \nn_{\tau'(M'+1)}(A)$ is at
least $\frac{\delta -N}{N}$. There exists $\delta_0=\delta_0
(L',C', N)$ such that if $\delta \geq \delta_0$ then $A_k=A$. This
contradicts the fact that $\rrr$ is without backtracking. We
conclude that $\delta \leq \delta_0$.

We apply Corollary \ref{12} to $\tilde{\q}[w,w_1]$ and to
$\tilde{\rrr}[w,s_-]$ and we obtain that $\ds(w_1,s_-)\leq D_1$,
where $D_1=D_1(L',C',\delta_0 )$. With a similar argument we
obtain that $\ds(z_1,s_+)\leq D_1$.

\me

\textbf{Step III}. We show that $\q$ has a component in $A$.

We have that $\ds (w_1,z_1)\geq D-2D_1$ and that
$\tilde{\q}[w_1,z_1]\subset \nn_{\tau' D_1 }(A)$.
 The decomposition (\ref{tq}) implies that $\tilde{\q}[w_1,z_1]
 = \tilde{\q}_k'\cup
\tilde{\q}_{k+1}\cup \dots \cup \tilde{\q}_{l-1}\cup
\tilde{\q}_{l}'$, where $k\leq l,\, k,l \in \{1,2,\dots ,N_1 \}$
and $\tilde{\q}_k'$, $\tilde{\q}_{l}'$ are eventually restrictions
of $\tilde{\q}_k$, $\tilde{\q}_{l}$, respectively, with endpoints
$w_1$ and $z_1$. If $D-2D_1>N_1$ it follows that
$\tilde{\q}[w_1,z_1]$ has at least a component of type (Q$_2$).
Since it has at most $N_1$ such components, we may moreover say
that $\tilde{\q}[w_1,z_1]$ has at least a component
$\tilde{\q}_{i}$ with endpoints at distance at least
$\frac{D-2D_1-N_1}{N_1}$. Consequently the diameter of $B_i\cap
\nn_{\tau' D_1 } (A)$ is at least $\frac{D-2D_1-N_1}{N_1}$. For
$D$ large enough we obtain that $B_i=A$. We conclude that $\q$ has
a component in $A$.

\me

\textbf{(2)} Suppose that $s$ and $t$ are $\mh$-components of
$\pgot$ and $\q$, respectively, contained in a left coset $A\in
\aaa$. We show that $\ds (s_-,t_-)$ and $\ds (s_+,t_+)$ are
bounded by a constant depending on $G$.

We take $u\in \pgot[\pgot_-, s_-]$ either such that the length of
$\pgot[u,s_-]$ is $2\lambda (\varkappa +1)$ or, if the length of
$\pgot[\pgot_-, s_-]$ is less than $2\lambda (\varkappa +1)$,
$u=\pgot_-$. Likewise we take $v\in \pgot[s_+,\pgot_+]$ either
such that the length of $\pgot[s_+,v]$ is $2\lambda (\varkappa
+1)$ or, if the length of $\pgot[s_+,\pgot_+]$ is less than
$2\lambda (\varkappa +1)$, $v=\pgot_+$.

Since $\dsh (s_-,t_-)\leq 1$ and $\cgsh$ is hyperbolic, there
exists $w\in \q[\q_-\, ,\, t_-]$ such that $\dsh(u,w)\leq
\varkappa $. Similarly, $\dsh (s_+,t_+)\leq 1$ implies the
existence of $z\in \q[t_+,\q_+]$ such that $\dsh (v,z)\leq
\varkappa $. We consider two geodesics $\g_{uw}$ and $\g_{vz}$. As
in Step 1 of the proof of (1), we show that the path $\g_{wu}\cup
\pgot[u,v]\cup \g_{vz}$ can be modified to give a path $\rrr$ with
endpoints $w$ and $z$ and of length at most $N$, without
backtracking, containing $s$, such that any of its lifts,
$\tilde{\rrr}$, decomposes as in (\ref{rrr}) and it is an
$(L',C')$-almost-geodesic. Again as in Step I of the proof of (1),
we show that the length of $\q[w,z]$ is at most $N_1$ and that any
lift $\tilde{\q}[w,z]$ decomposes as in (\ref{tq}) and it is an
$(L',C')$-almost-geodesic.

With an argument as in Step II of the proof of (1), we show that
$\tilde{\rrr}[w,s_-]$ and $\tilde{\rrr}[s_+,z]$ intersect
$\nn_{M'} (A)$ in sets of diameter at most $\delta_0$. The same
argument can be used to show that $\tilde{\q}[w,t_-]$ and
$\tilde{\q}[t_+,z]$ intersect $\nn_{M'}(A)$ in sets of diameter at
most $\delta_0'=\delta_0'(L',C',N_1)$. Corollary \ref{12} implies
that $\ds(s_-,t_-)$ and $\ds (s_+,t_+)$ are at most $D_1$, where
$D_1=D_1 (L',C', \delta_0 , \delta_0')$.\endproof

\subsection{The Morse Lemma}

Proposition \ref{udist} can be strengthened to the following
statement.

\begin{proposition}\label{prudist}
Let $\q :[0,\ell] \to \cgs$ be an $(L,C)$-quasi-geodesic and let
$\pgot$ be a geodesic in $\cgsh$ joining the endpoints of $\q$. In
$\cgs$ the quasi-geodesic segment $\q$ is contained in the
$T$-tubular neighborhood of the $M$-saturation of the lift
$\tilde{\pgot }$ of $\pgot$. Conversely, the lift $\tilde{\pgot }$
is contained in the $T$-tubular neighborhood of the $M$-saturation
of $\q$. The constants $T$ and $M$ depend on $L,C$ and $S$.
\end{proposition}

\proof According to Proposition \ref{udist}, the Hausdorff distance
from $\q$ to $\pgot$ in $\cgsh$ is at most $\varkappa$. We divide
$\pgot $ into arcs of length $3(\varkappa +2)$, with the exception
of the two arcs at the endpoints, which can be shorter. Let $s$ be
one of these arcs. Consider $u$ on the sub-arc of $\pgot$ between
$\pgot_-$ and $s_-$ such that either $\dsh (u,s_-)=\varkappa +2$ or
$u=\pgot_-$. Likewise let $v$ be a point on the sub-arc of $\pgot$
between $s_+$ and $\pgot_+$ such that either $\dsh (s_+,v)=\varkappa
+2$ or $v=\pgot_+$. Let $w$ and $z$ be two points on $\q$ at
distance at most $\varkappa$ from $u$ and $v$ respectively, in
$\cgsh$.

We repeat the argument from the proof of Proposition \ref{pbcp},
Step I. Consider $\g_{uw}$ and $\g_{vz}$ geodesics in $\cgsh$.
Consider $u'$ the farthest from $u$ point on $\g_{uw}$ contained
in the same left coset as an $\hh$-component $\sigma$ of $\pgot$.
Likewise let $v'$ be the farthest from $v$ point on $\g_{vz}$
contained in the same left coset as an $\hh$-component $\sigma'$
of $\pgot$. Then $\sigma$ does not intersect $s$, otherwise the
distance from $u$ to $s$ in $\cgsh$ would be at most $\varkappa
+1$. Similarly, $\sigma'$ does not intersect $s$.

Consider the path in $\cgsh$ defined as $\rrr=\g_{wu'}\sqcup
\g_{u'\sigma_+}\sqcup \pgot[\sigma_+, \sigma_-']\sqcup
\g_{\sigma_-'v'}\sqcup \g_{v'z}$, where $\g_{wu'}$ and $\g_{v'z}$
are sub-geodesics of $\g_{wu}$ and $\g_{vz}$, respectively, and
$\g_{u'\sigma_+}$ and $\g_{\sigma_-'v'}$ are composed of one edge.
It has no backtracking and its lift $\tilde{\rrr}$ is an
$(L',C')$-quasi-geodesic, where $L'$ and $C'$ depend on the length
of $\rrr$, hence on $\varkappa$. It has the same endpoints as a
sub-quasi-geodesic of $\q$ of endpoints $w$ and $z$, hence it is
contained in the $T$-tubular neighborhood of the $M$-saturation of
it, where $M=M(L,C)$ and $T=T(L,C,\varkappa )$. In particular this
is true for the lift of $s$. Since $s$ is arbitrary, we have
obtained that the lift $\tilde{\pgot }$ is contained in the
$T$-tubular neighborhood of the $M$-saturation of $\q$.

We now consider $\q$ endowed with the order from $[0,\ell]$. We
consider the path $\widehat{\q}$ in $\cgsh$ obtained by deleting the
part of $\q$ between the first and the last point in $\q$ contained
in the same left coset, replacing it with an edge, and performing
this successively for every coset intersecting $\q$ in more than one
point. Then, for a constant $D$ to be chosen later, we divide
$\widehat{\q}$ into arcs $t$ such that $t_+$ is the first vertex on
$\widehat{\q}$ (in the order inherited from $\q$) which is at
distance $D$ of $t_-$. We start constructing these arcs from $\q_-$
and we end in $\q_+$ by an arc which possibly has endpoints at
distance smaller than $D$. Consider $t$ one of these arcs. Let $u$
be a point on $\q$ between $\q_-$ and $t_-$ with the property that
it is at distance $\varkappa +2$ of $t$. If no such point exists,
take $u=\q_-$. Similarly, take $v$ a point on $\q$ between $t_+$ and
$\q_+$ with the property that it is at distance $\varkappa +2$ of
$t$, or $v=\q_+$. There exist $w$ and $z$ respectively on $\pgot$ at
distance at most $\varkappa$ from $u$ and $v$. Then $\dsh (w,z)\leq
2\varkappa +2(\varkappa +2)+6D$. It follows that the lift
$\tilde{\pgot}_{wz} $ of the sub-geodesic $\pgot_{wz}$ of $\pgot$ of
endpoints $w$ and $z$ is an $(L'',C'')$-quasi-geodesic, where $L''$
and $C''$ depend on $\varkappa $ and $D$.

As above we choose $u'\in \g_{uw}$ and $\sigma$ an $\hh$-component
of $\widehat{\q}$ in the same left class as $u'$. The choice of
$u$ implies that $\sigma $ does not intersect $t$, otherwise $u$
would be at distance at most $\varkappa +1$ of $t$. Likewise we
choose $v'$ and $\sigma'$, and we construct the path
$\rrr'=\g_{wu'}\sqcup \g_{u'\sigma_+}\sqcup \widehat{\q}[\sigma_+,
\sigma_-']\sqcup \g_{\sigma_-'v'}\sqcup \g_{v'z}$ in $\cgsh$ of
bounded length, with $\widehat{\q}[\sigma_+, \sigma_-']$
containing $t$. As in the proof of Proposition \ref{pbcp}, Step I,
the sub-arcs $\rrr'[w,v']$ and $\rrr' [u',z]$ do not have
backtracking. Suppose that $\g_{wu'}$ and $\g_{v'z}$ have
$\hh$-components in the same left coset. Let $w'$ and $z'$ be the
nearest points to $u'$ and respectively $v'$ contained in the same
left coset. Lemma \ref{polig} implies that
${\frak{l}}=\g_{\sigma_+u'}\sqcup \g_{u'w'}\sqcup \g_{w'z'} \sqcup
\g_{z'v'}\sqcup \g_{v'\sigma_-'}$, which has length at most
$2\varkappa +3$, lifts to an $(L_1, C_1)$-quasi-geodesic, where
$(L_1, C_1)$ depends of $\varkappa$. It follows that the sub-arc
of $\q$ between $\sigma_+$ and $\sigma_-'$ is contained in the
$\tau $-neighborhood of the $M'$-saturation of
$\widetilde{{\frak{l}}}$, where $M'=M'(\varkappa )$ and $\tau
=\tau (\varkappa ,L,C)$. It follows that the diameter of
$\widehat{\q}[\sigma_+, \sigma_-']$ is at most $2\tau +2
+2M'+\mathrm{length}\, {\frak{l}} $. Hence $D\leq 2(\tau
+1+M'+\varkappa)+3$. Thus, if we take $D=2(\tau
+1+M'+\varkappa)+4$, we get a contradiction.

We conclude that $\rrr'$ has no backtracking, hence it lifts to a
quasi-geodesic, by Lemma \ref{polig}. We make a slight change when
lifting it to a path $\tilde{\rrr}'$, in that the sub-arcs in
$\widehat{\q}$ are lifted to the corresponding sub-arcs of $\q$.
We obtain a quasi-geodesic $\tilde{\rrr}'$ with the same endpoints
as $\tilde{\pgot}_{wz} $, hence contained in the $T$-tubular
neighborhood of the $M$-saturation of it, where $M=M(L,C,\varkappa
)$ and $T=T(L,C,\varkappa )$. In particular this applies to the
lift of $t$. Since $t$ was arbitrary, this allows to obtain the
desired statement for $\q$ and $\tilde{\pgot}$.\endproof

Proposition \ref{prudist} together with Proposition \ref{udist} and
Lemmas \ref{lqqcun}, \ref{limrsat} and \ref{gat} imply Theorem
\ref{morse}.

\subsection{Undistorted subgroups and outer automorphisms of relatively hyperbolic
groups}\label{sundis}

\begin{theorem}\label{thundis} Let
$G=\la S\ra$ be a finitely generated group that is hyperbolic
relative to subgroups $H_1,...,H_n$. Let $G_1=\la S_1\ra$ be an
undistorted finitely generated subgroup of $G$. Then $G_1$ is
relatively hyperbolic with respect to subgroups $H_1',...,H_m'$,
where each $H_i'$ is one of the intersections $G_1\cap gH_jg\iv$,
$g\in G$.
\end{theorem}

\proof Since $G_1$ is undistorted, there exists a constant $D\ge
1$ such that for every element $g\in G_1$, $|g|_{S_1} \le D|g|_S$.
Here by $|g|_S$ and $|g|_{S_1}$ we denote the length of $g$ in $G$
and $G_1$ respectively. We can assume that $S_1\subseteq S$ so
that the graph $\Cay(G_1,S_1)$ is inside $\Cay(G,S)$. Then every
geodesic in $\Cay(G_1,S_1)$ is a $(D,0)$-quasi-geodesic of
$\Cay(G,S)$.

\textbf{Step I.} Let us prove that for every coset $gH_i$ and
every constant $C>0$ there exists $C'=C'(C,g,i)>0$ such that
$G_1\cap \nn_C(gH_i)\subseteq \nn_{C'}(G_1\cap gH_ig\iv)$. By
contradiction, let $(x_j)_{j\in \N}$ be a sequence of elements in
$G_1$ such that $x_j=gh_jp_j\in G_1$, $h_j\in H_i$, $|p_j|_S<C$,
and $\dist(x_j,G_1\cap gH_ig\iv)\ge j$ for every $j$. Without loss
of generality we can assume that $p_j=p$ is constant. Then
$x_jx_1\iv\in G_1\cap gH_ig\iv$. Hence $\dist(x_j, G_1\cap
gH_ig\iv)\le |x_1|_S$, a contradiction.

\me

\textbf{Step II.} Let $R>0$ and let $gH_i$ be such that $\nn_R
(gH_i)\cap G_1\neq \emptyset$.

We prove that for every $K>0$ there exists $K'=K'(K ,R)$ such that
$$
G_1\cap \nn_{K} (gH_i) \subset \nn_{K'}(G_1\cap g_1 \gamma H_i
\gamma^{-1})
$$ for some $g_1\in G_1$
and some $\gamma\in G$ with $|\gamma|_S\le R$.

Fix $K>0$ and define $K'$ as the maximum of numbers $C'(K,\gamma ,
i)$ defined in Step I taken over all $i\in \{ 1,2,\dots ,n\}$ and
all $\gamma\in G$ with $|\gamma|_S\le R$.

Let $g\in G$ be such that $G_1\cap \nn_{R} (gH_i)\neq \emptyset$.
Let $g_1$ be an element of the intersection. Then $g_1\iv\nn_R
(gH_i)= \nn_R (g_1\iv gH_i)$ contains $1$, hence $g_1\iv gH_i
=\gamma H_i$ where $|\gamma|_S\le R$.

Step I and the choice of $K'$ imply that
$$
G_1 \cap \nn_{K} (\gamma H_i) \subset \nn_{K'}(G_1\cap \gamma H_i
\gamma^{-1})\, .
$$

Multiplying this inclusion by $g_1$ on the left, we obtain
$$
G_1 \cap \nn_K (gH_i)\subset \nn_{K'}(G_1\cap g_1 \gamma H_i
\gamma^{-1}).
$$

\me

\textbf{Step III.} Let $R=M(D, 0,\frac{1}{3} )$ be the constant
given by the property $(\alpha_2')$ satisfied by the left cosets
$\{ g H_i \mid g\in G,\, i=1,2,\dots ,n\}$ in $\Cay (G,S)$.

For every $i\in\{1,...,n\}$ consider the following equivalence
relation on the ball $B(1,R)$ in $G$: $$\gamma\sim_i\gamma' \hbox{
iff } G_1\gamma H_i=G_1\gamma'H_i.$$ For each pair
$(\gamma,\gamma')$ of $\sim_i$-equivalent elements in $B(1,R)$ we
choose one $g_1\in G_1$ such that $\gamma\in g_1\gamma'H_i$. Let
$\tilde C$ be the maximal length of these elements $g_1$.

Let ${\cal M}$ be the collection of all nontrivial subgroups of
$G_1$ in the set
$$
\{G_1\cap
 \gamma H_i \gamma^{-1}  \mid i\in \{ 1,2,\dots ,n\}\, ,\,
 |\gamma|_S\le R\}.
 $$

By Step II, this collection of subgroups has the property that for
every $K>0$ there exists $K'=K'(K,R)$ such that for every $g\in G$
with $\nn_R(gH_i)\cap G_1 \neq \emptyset$, we have

\begin{equation}\label{eqstar} G_1\cap \nn_K (gH_i) \subset \nn_{K'} (g_1 H)
\end{equation}
for some $g_1\in G_1$ and $H\in \mathcal{M}$.

We say that two non-trivial subgroups $G_1\cap \gamma
H_i\gamma\iv$ and $G_1\cap \beta H_i\beta\iv$ from $\cal M$ are
equivalent if $\gamma\sim_i\beta$.

Let $H_1',...,H_m'$ be the set of representatives of equivalent
classes in $\cal M$. If $\cal M$ is empty, we set $m=1$,
$H_1'=\{1\}$.

Notice that for every $H\in \cal M$ there exists $j\in
\{1,...,m\}$ such that $H$ is at Hausdorff distance at most
$\tilde C$ from a left coset $gH_j'$ from $G_1$. Indeed, $H=\gamma
H_i\gamma\iv\cap G_1$. Let $H_j'=\beta H_i\beta\iv\cap G_1$ be
equivalent to $H$. Then $\gamma=g\beta h$ for some $g\in G_1$,
$h\in H_i$, where $|g|\le \tilde C$. Then $$H=g\beta hH_ih\iv
\beta\iv g\iv\cap G_1 = g H_j' g\iv \, ,$$ from which we deduce
that $H$ is at Hausdorff distance at most $\tilde C$ from $gH_j'$.

Hence (\ref{eqstar}) remains true if we replace ${\cal M}$ by the
smaller set $\{ H_1',...,H_m'\}$ and $K'$ by $K'+\tilde C$.

We shall prove that $G_1$ is relatively hyperbolic with respect to
$\{ H_1',\dots ,H_m' \}$ by checking properties $(\alpha_1)$,
$(\alpha_2^{\frac{1}{6D}})$, $(\alpha_3)$ from Theorem \ref{tgi}
and Remark \ref{tgi3} for the collection of left cosets $\{ g_1
H_j' \mid g_1\in G_1, j=1,2,\dots ,m\}$.

\me

{\em Property $(\alpha_1)$.} Consider $g_1 H_j' \neq g_1' H_k'$.
We have that $$\nn_\delta(g_1 H_j')\cap \nn_\delta (g_1'
H_k')\subset \nn_\delta (g_1\gamma H_{i_j} \gamma^{-1})\cap
\nn_\delta (g_1'\gamma' H_{i_k} (\gamma')^{-1})\subset \nn_{\delta
+R}(g_1\gamma H_{i_j})\cap \nn_{\delta +R} (g_1'\gamma' H_{i_k})\,
.$$

Suppose that $g_1\gamma H_{i_j}= g_1'\gamma' H_{i_k}$. Then
$(g_1\gamma )\iv g_1'\gamma'\in H_{i_j}$ hence $g_1\gamma H_{i_j}
=g_1'\gamma' H_{i_j}$. We deduce that $ H_{i_j}= H_{i_k}$.
Therefore $g_1\gamma=g_1'\gamma'h$ for some $h\in H_{i_j}$. Hence
$\gamma\sim_{i_j}\gamma'$, so $\gamma=\gamma'$. We deduce that
$g_1\gamma H_{i_j}\gamma\iv = g_1' \gamma H_{i_j}\gamma\iv$. So
$g_1H_j'=g_1'H_k'$, a contradiction.

Thus, $g_1\gamma H_{i_j}\neq g_1'\gamma' H_{i_k}$. Property
$(\alpha_1)$ satisfied by the left cosets $\{ g H_i \mid g\in G,\,
i=1,2,\dots ,n\}$ allows to complete the proof.

\me

{\em Property $(\alpha_2^{\frac{1}{6D}})$.} Let $\theta_1 \in
\left[ 0,\frac{1}{6D}\right)$. We may write $\theta_1 \frac{\theta}{D}$, with $\theta \in \left[ 0,\frac{1}{6}\right)$.
Let $\g\colon [0,\ell ]\to \Cay(G_1,S_1)$ be a geodesic of length
$\ell$ in $\Cay (G_1, S_1)$ with endpoints in $\nn_{\theta_1 \ell}
(g_1 H_j')\subset \nn_{\theta_1 \ell +R} (g_1 \gamma H_i)$, where
$|\gamma|_S\le R$ and $i\in \{ 1,2,\dots ,n\}$. Then $\g$ is a
$(D,0)$-quasi-geodesic in $\Cay(G,S)$.

Suppose that $\ell \leq 6DR$. Then $\g$ is contained in the
$(3DR+R)$-tubular neighborhood of $g_1 H_j'$ in $\Cay(G_1,S_1)$.

Suppose that $\ell > 6DR$. Then the endpoints of $\g$ are
contained in $\nn_{(\theta +\frac{1}{6})\frac{\ell }{D}} (g_1
\gamma H_i)\subset \nn_{\frac{1}{3}\frac{\ell }{D}} (g_1 \gamma
H_i)$ in $\Cay (G,S)$. Since the property $(\alpha_2')$ is
satisfied by the cosets of $H_i$ in $G$, it follows that $\g$
intersects $\nn_R (g_1 \gamma H_i )$. Hence $\g$ intersects
$G_1\cap \nn_{R} (g_1 \gamma H_i)=g_1 [G_1\cap \nn_{R} (\gamma
H_i)]\subset g_1 \nn_{R'} (H_j')$ where $R'=R'(R,R)$ is given by
Step II.

We conclude that $\g$ intersects $\nn_{M'} (g_1 H_j')$ in $\Cay
(G_1,S_1)$, for $M'=\sup (DR', 3DR+R)$.

\me

{\em Property $(\alpha_3)$.} We use the property (\ref{eqstar}) of
$\{ H_1',\dots , H_m' \}$ and the property $(\alpha_3')$ satisfied
by the cosets of groups $H_i$.

Fix an integer $k\ge2$. Let $P$ be a $(\vartheta, 2, 8D)$-fat
geodesic $k$-gon in $\Cay(G_1,S_1)$ for some $\vartheta$. Then $P$
has $(D,0)$-quasi-geodesic sides in $\Cay(G,S)$ and it is $\left(
\frac{\vartheta}{D}, 2D, 8D \right)$-fat. Consequently, for
$\vartheta $ large enough, by property $(\alpha_3')$ satisfied by
the left cosets $\{ gH_i \mid g\in G,\, i=1,\dots ,n\}$, the
$k$-gon $P$ is contained in a tubular neighborhood
$\nn_\varkappa(gH_i)$ in $\Cay(G,S)$ for some $\varkappa>0$.

Suppose that all edges of $P$ have lengths at most $3D\varkappa$
in $\Cay(G_1,S_1)$. Then $P$ has diameter at most $3kD\varkappa $
in the same Cayley graph.

Suppose that one edge $\g$ of $P$ has length at least $3D\varkappa
$. This, the fact that $P\subset \nn_\varkappa (gH_i)$
 and property $(\alpha_2')$ satisfied by the left
cosets $\{ gH_i\}$ implies that $\g$ intersects $\nn_R (gH_i)$,
therefore $\nn_R (gH_i) \cap G_1 \neq \emptyset$.

Then by (\ref{eqstar}) there exists
$\varkappa'=\varkappa'(\varkappa ,R)$ such that
     $$
G_1\cap \nn_\varkappa (gH_i) \subset \nn_{\varkappa'} (g_1 H_j')
     $$ for some $g_1\in G_1$ and $j\in \{ 1,2,\dots ,m
     \}$.
We conclude that in this case $P\subset \nn_{\varkappa'} (g_1
H_j')$.

Thus we can take $\xi$ needed in $(\alpha_3)$ to be the maximum of
$3kD\varkappa$ and $\varkappa'$. \endproof

\begin{remarks}\label{rlast}
(1) If in Theorem \ref{thundis} the subgroup $G_1$ is unconstricted
then $G_1$ is inside a conjugate of one of the subgroups $H_i$.

(2) If the subgroup $G_1$ intersects with all conjugates of the
subgroups $H_1,...,H_n$ by hyperbolic subgroups then $G_1$ is
hyperbolic.

\end{remarks}

\proof (1) Indeed, Corollary \ref{cutp1} implies that $G_1$ is
contained in the $K$-tubular neighborhood of a left coset $gH_i$,
where $K$ depends only on the non-distortion constants. For every
$g_1\in G_1$, $G_1=g_1 G_1$ is contained in the $K$-tubular
neighborhoods of $g_1gH_i$ and of $gH_i$. Since $G_1$ is infinite,
property $(\alpha_1)$ implies that $g_1gH_i=gH_i$. We conclude
that $G_1$ is contained in $gH_ig\iv$.

(2) By Theorem \ref{thundis} $G_1$ is relatively hyperbolic with
respect to hyperbolic subgroups, so every asymptotic cone of $G_1$
is tree-graded with respect to $\R$-trees, whence it is an
$\R$-tree itself. Therefore $G_1$ is hyperbolic
\cite{Gr2}.\endproof

\begin{cor}\label{829} Let $G$ be a finitely generated group that
 is relatively hyperbolic with
respect to  subgroups $H_1,...,H_m$. Suppose that $H_1$ is
unconstricted and that each $H_i\, ,\, i\in \{2,\dots ,m\}$, is
infinite and either unconstricted or does not contain a copy of
$H_1$. Let $\mathrm{Fix}(H_1)$ be the subgroup of the automorphism
group of $G$ consisting of the automorphisms that fix $H_1$ as a
set. Then:
\begin{itemize}
\item[(1)] $\mathrm{Inn}(G)\mathrm{Fix}(H_1)$ has index at most $m!$ in
$\mathrm{Aut}(G)$ (in particular, if $m=1$, these two subgroups
coincide).

\item[(2)] $\mathrm{Inn}(G)\cap \mathrm{Fix}(H_1)
= \mathrm{Inn}_{H_1}(G)$, where $\mathrm{Inn}_{H_1}(G)$ is by
definition $\{i_h\in \mathrm{Inn}(G) \mid h\in
  H_1\}\, .$
\item[(3)] There exists a natural homomorphism from a subgroup of index
at most $m!$ in $\mathrm{Out}(G)$ to $\mathrm{Out}(H_1)$ given by
$\phi\mapsto i_{g_\phi}\phi|_{H_1}$, where $g_\phi$ is an element
of $G$ such that $i_{g_\phi}\phi\in \mathrm{Fix}(H_1)$, and
$|_{H_1}$ denotes the restriction of the automorphism to $H_1$.
\end{itemize}
\end{cor}

\proof (1) Indeed, every automorphism $\phi$ of $G$ is a
quasi-isometry of the Cayley graph of $G$. Hence $\phi(H_1)$ is an
undistorted subgroup of $G$ that is isomorphic to $H_1$. By Remark
\ref{rlast}, (1), we have that $\phi(H_1)\subset g H_jg\iv$ for some
$g\in G$ and $j\in \{1,2,\dots ,m\}$. In particular $i_g\iv \phi
(H_1)\subset H_j$. By hypothesis $H_j$ is unconstricted. If we
denote by $\psi$ the automorphism $i_g\iv \phi$, we have that $\psi
\iv (H_j)\subset \gamma H_k \gamma\iv$, for some $\gamma\in G$ and
$k\in \{1,2,\dots ,m\}$. Consequently $H_1\subset \gamma H_k
\gamma\iv$. We deduce from the fact that $H_1$ is infinite and from
property $(\alpha_1)$ that $H_1= \gamma H_k \gamma\iv$ and $i_g\iv
\phi (H_1)= H_j$. In particular every automorphism of $G$ induces a
permutation of the set
$$\{H_i\mid H_i \hbox{ is isomorphic to } H_1\}.$$ Therefore we have an action of $\mathrm{Aut}(G)$ on a
subset of $\{H_1,...,H_m\}$. Let $S$ be the kernel of this action.
Then $|\mathrm{Aut}(G):S|\le m!$. The composition of any $\phi\in
S$ with an inner automorphism $i_{g}\iv$ induced by $g\iv$ is in
$\mathrm{Fix}(H_1)$. Therefore $S$ is contained in
$\mathrm{Inn}(G)\mathrm{Fix}(H_1)$.

(2) Let $i_g$ be an element in $\mathrm{Inn}(G)\cap \mathrm{Fix}(H_1)$. Then $g$ normalizes $H_1$, hence by
\cite{Osin}, $g\in H_1$.

(3) This immediately follows from (1) and (2).
\endproof

\section{Appendix. Relatively hyperbolic groups are asymptotically
tree-graded. By Denis Osin and Mark Sapir}

Here we prove the ``if" statement in Theorem \ref{dir}.

\begin{theorem}\label{main}
Let $G$ be a group generated by a finite set $S$, that is
relatively hyperbolic with respect to finitely generated subgroups
$H_1,...,H_m$. Then $G$ is asymptotically tree-graded with respect
to these subgroups.
\end{theorem}

Throughout the rest of this section we assume that $G$,
$H_1,...,H_m$, an ultrafilter $\omega$, and a sequence of numbers
$d=(d_i)$ are fixed, $G$ is generated by a finite set $S$ and is
hyperbolic relative to $H_1,...,H_m$. We denote the asymptotic
cone $\co{G;e,d}$ by $C$.

If $(g_i), (h_i)$ are sequences of numbers, we shall write
$g_i\le_\omega h_i$ instead of ``$g_i\le h_i$ $\omega$-almost
surely". The signs $=_\omega$, $\in_\omega$ have similar meanings.

As before, $\hh=(\bigcup_{i=1}^m H_i) \setminus \{ e \}$. For
every $i=1,...,m$, in every coset of $H_i$ ($i=1,...,m$) we choose
a smallest length representative. The set of these representatives
is denoted by $T_i$. Let $\ttt_i$ be the set $\{(g_j)^\omega\mid
\lio{|g_i|_S}<\infty\}$. For each $\gamma=(g_j)^\omega\in \ttt_i$
we denote by $M_\gamma$ the $\omega$-limit $\lio{g_jH_i}_e$. We
need to show that $C$ is tree-graded with respect to all
$\pp=\{M_\gamma\mid \gamma\in\ttt_i, i=1,...,,m\}$.

We use the notation $\dx $ and $\dxh$ for combinatorial metrics on
$\cgs$ and $\cgsh$. When speaking about geodesics in $\cgsh$ we
always assume them to be geodesic with respect to $\dxh$.

The lemma below can be found in \cite[Theorem 3.23]{Osin}.

\begin{lemma} \label{rips}
There exists a constant $\nu >0$ such that the following condition
holds. Let $\Delta = pqr$ be a geodesic triangle in $\Cay(G, S\cup
\mathcal H)$ whose sides are geodesic (with respect to the metric
$\dist _{S\cup \mathcal H}$). Then for any vertex $v$ on $p$,
there exists a vertex $u$ on the union $q\cup r$ such that
$$\dx(u,v)\le \nu .$$
\end{lemma}

\begin{lemma}\label{log}
Let $p$ and $q$ be paths in $\G $ such that $p_-=q_-$, $p_+=q_+$,
and $q$ is geodesic. Then for any vertex $v\in q$, there exists a
vertex $u\in p$ such that $$\dx (u,v)\le (1+\nu ) \log\limits_2
|p|.$$
\end{lemma}

\begin{proof}
Let $f:\mathbb N\to \mathbb N$ be the smallest function such that
the following condition holds. Let $p$ and $q$ be paths in $\G $
such that $p_-=q_-$, $p_+=q_+$, $q$ is geodesic, and $|p|\le n$.
Then for any vertex $v\in q$, there exists a vertex $u\in p$ such
that $\dx (u,v)\le f(n).$ Clearly $f(n)$ is finite for each value
of the argument. By dividing $p$ into two parts and applying Lemma
\ref{rips}, we obtain $f(m+n)\le \max \{f(m), \; f(n)\} +\nu $. In
particular, $f(2n)\le f(n)+\nu $ and $f(n+1)\le f(n)+\nu $.

Suppose that $$ n=\varepsilon _0 +2\varepsilon _1 +\ldots +
2^k\varepsilon _k, $$ where $\varepsilon _i\in \{ 0,1\} $ and
$\varepsilon _k=1$. Then
$$
\begin{array}{rl}
f(n)= & f\big( \varepsilon _0 +2( \varepsilon _1 + \ldots +2
(\varepsilon _{k-1} +2)\ldots )\big) \le \\ & \\
 & \underbrace{\nu+\nu +\ldots +\nu }\limits_{2k\; {\rm times}} +f(1) \le
2\nu  \log\limits_2 n.
\end{array}
$$
\end{proof}

The next lemma can be found in \cite[Lemma 3.1]{Osin}.

\begin{lemma}\label{length-of-comp}
There is a constant $\alpha $ such that for any cycle $q$ in $\G
$, and any set of isolated $\hh$-components of $p_1, \ldots , p_k$
of $q$, we have
$$
\sum\limits_{i=1}^k \dx ((p_i)_-, (p_i)_+)\le \alpha |q|.
$$
\end{lemma}

The following lemma which holds for any (not necessarily
relatively hyperbolic) finitely generated group $G$ and any
subgroup $H\le G$.

\begin{lemma}\label{int}
For any $i=1,...,m$, $\theta , \sigma\in \ttt_i$, if $\theta\ne
\sigma$ then the intersection $M_\theta \cap M_\sigma $ consists
of at most $1$ point.
\end{lemma}

\begin{proof}
Suppose that $x,y\in M_\theta \cap M_\sigma $. Suppose that
$\theta =(f_j)^\omega $, $\sigma =(g_j)^\omega $. Then $$x=\lio{
f_ja_j},\ y=\lio{f_js_j}$$ for some $a_j, s_j\in H_i$ and
$$x=\lio{g_jb_j}, \ y=\lio{g_jt_j} $$ for some $b_j, t_j\in
H_i$. Since the sequences $(f_ja_j)^\omega$ and $(g_jb_j)^\omega$
are equivalent, we have $f_ja_j=g_jb_ju_j$, where $|u_j|_S=_\omega
o(d_j)$. Similarly $f_js_j=g_jt_jv_j$, where $|v_j|_S=_\omega
o(d_j)$. From these equalities we have
$$
a_j^{-1}s_j=u_j^{-1}b_j^{-1}t_jv_i.
$$

Let $U_j$, $V_j$ be shortest words over $S$ representing $u_j$ and
$v_j$ respectively. Let also $h_j=a_j^{-1}s_i$ and
$k_j=b_j^{-1}t_j$. Then there exists a quadrangle
$$q^j=p_1^jp_2^jp_3^jp_4^j$$ in $\cgs$ such that $\phi
(p_1^j)\equiv U_i$, $\phi (p_3^j)\equiv V_j^{-1}$ and $p_2^j$,
$p_4^j$ are edges of $\cgs$ labelled $h_j$ and $k^{-1}_j$
respectively. Note that the cycle $q^j$ contains only two
components, namely, $p_2^j$ and $p_4^j$, as the labels of $p_1^j$
and $p_3^j$ are words over $S$. Let $A\subseteq \N$ be the set of
all $j$ such that the components $p_2^j$ and $p_4^j$ are
connected. There are two cases to consider.

{\it Case 1.} $\omega(A)=1$. Note that $\phi(p_1^j)$ represents an
element of $H_i$ in $G$ for any $j\in A$, i.e., $s_j\in_\omega
H_i$. It follows that $\theta =\sigma$.

{\it Case 2.} $\omega (A)=0$. Note that $p_2^j$ is an isolated
component of $q^j$ for any $j\in \N\setminus A$. Since $\omega
(\N\setminus A)=1$, applying Lemma \ref{length-of-comp}, we obtain
$$
|h_j|_S= \dx((p_2)_-, (p_2)_+)\le \alpha |q^j|\le _\omega \alpha
(2+2o(d_j))=o(d_j).
$$
This yields
$$
\dist (x,y)=\lio{\frac{1}{d_j} \dx (f_ja_j, f_js_j)}
=\lio{\frac{1}{d_j} |h_j|_S} =0,
$$
i.e., $x=y$.
\end{proof}

The following lemma does use the relative hyperbolicity of $G$.

\begin{lemma}\label{int1} For every $i\ne i'$ and every $\theta\in\ttt_i$,
$\sigma\in\ttt_{i'}$, the intersection $M_\theta\cap M_\sigma$
consists of at most one point.
\end{lemma}

\proof Indeed, repeating the argument from the proof of Lemma
\ref{int}, we immediately get contradiction with the BCP property.
\endproof

Lemmas \ref{int}, \ref{int1} show that the asymptotic cone $C$
satisfies the property $(T_1)$ with respect to the set $\pp$. Now
we are going to prove $(T_2)$.

\begin{lemma}\label{40}
Let $\g$ be a simple loop in $C$. Suppose that $\g =\lio{\g_j}$
for certain loops $\g_j$ in $\cgs$, $|\g_j|\le Cd_j$ for some
constant $C$. Then there exists $i=1,...,m$ and $\theta \in
\ttt_i$ such that $\g$ belongs to $M_\theta$.
\end{lemma}

\begin{proof}
Let $a\ne b$ be two arbitrary points of $\g$, $$a=\lio{a_j},
\;\;\;\;\; b=\lio{b_j}, $$ where $a_j$, $b_j$ are vertices on
$\g_j$. For every $j$, we consider a geodesic path $\q_j$ in
$\cgsh$ such that $(\q_j)_-=a_j$, $(\q_j)_+=b_j$.

According to Lemma \ref{log}, for every vertex $v\in \q_j$, there
exist vertices $x_j=x_j(v)\in \g_j[a_j, b_j]$ and $y_j=y_j(v)\in
\g_j[b_j, a_j]$ (here $\g_j[a_j, b_j]$ and $\g_j[b_j, a_j]$ are
segments of $\g_j=\g_j[a_j,b_j]\g_j[b_j,a_j]$) such that
\begin{equation}\label{50}
\dx(v, x_j)\le 2\nu \log\limits_2 |\g_j[a_j, b_j]|< 2\nu
\log\limits_2 (Cd_j)=o(d_j)
\end{equation}
and similarly
\begin{equation}\label{60}
\dx (v, y_j)\le 2\nu \log\limits_2 |\g_j[b_j, a_j]|< 2\nu
\log\limits_2 (Cd_j)=o(d_j).
\end{equation}

Summing (\ref{50}) and (\ref{60}), we obtain $$ \dx (x_j,y_j)\le
\dx (x_j, v) +\dx (v, y_j)=o(d_j).$$ Thus for any $j$, there are
only two possibilities: either $\lio{x_j}=\lio{y_j}=a$ or $\lio{
x_j}=\lio{y_j}=b$, otherwise the loop $\g$ is not simple.

For every $j$, we take two vertices $v_j, w_j\in \q_j$ such that
$$
\lio{x_j(v_j)}=\lio{y_j(v_j)}=a,
$$
$$
\lio{x_j(w_j)}=\lio{y_j(w_j)}=b,
$$
and $\dxh(v_j,w_j)=1$. Since $\lio{x_j(v_j)}=a$, we have $\dx
(x_j(v_j), a)=_\omega o(d_j)$. Hence $$\dx(v_j, a_j)\le \dx(v_j,
x_j(v_j))+\dx (x_j(v_j), a_j)=_\omega o(d_j).$$ Similarly,
$$\dx (w_j, b)=_\omega o(d_j).$$ This means that
\begin{equation}\label{avwb}
\lio{a_j}=\lio{v_j}, \;\;\; {\rm and} \;\;\; \lio{b_j}=\lio{w_j}.
\end{equation}
For every $i=1,...,m$, set $A_i=\{ j\in \N\; |\; v_j^{-1}w_j\in
H_i\} $. Let us consider two cases.

{\it Case 1.} $\omega(A_i)=1$ for some $i$. Set
$\theta=(t_i(v_j))^\omega \in\ttt_i$ where $t_i(v_j)$ is the
representative of the coset $v_jH_i$ chosen in the definition of
$\ttt_i$. Then $a,b\in_\omega M_\theta$. Indeed, this is obvious
for $a$ since $\lio{a_j}=\lio{v_j} \in M_\theta $. Further, since
$v_j^{-1}w_j\in _\omega H_i$, we have $t(w_j)=_\omega t(v_j)$.
Hence $\lio{b_j}=\lio{w_j}\in M_\theta $.

{\it Case 2.} $\omega (A_i)=0$ for every $i$. Recall that
$v_j^{-1}w_j\in S\cup \hh$. Thus we have $v_j^{-1}w_j\in _\omega
S$. This implies $|v_j^{-1}w_j|_S=_\omega 1$ and
$\lio{v_j}=\lio{w_j}$. Taking into account (\ref{avwb}), we obtain
$\lio{a_j}=\lio{b_j}$, i.e., $a=b$.

Since $a$ and $b$ were arbitrary points of $\beta$, the lemma is
proved.
\end{proof}

Now property $(T_2)$ immediately follows from Proposition
\ref{approxul} and Lemma \ref{40}.

\me

{\bf Acknowledgement.} The authors are grateful to Anna Erschler and
Denis Osin for many helpful conversations. We are also grateful to
Chris Hruska for a discussion that lead us to a new version of
Theorem \ref{morse}.

\me

\addtocontents{toc}{\contentsline {section}{\numberline { }
References \hbox {}}{\pageref{bibbb}}}


\begin{thebibliography}{ECHPT}
\label{bibbb}

\bibitem[AN]{AN} S. A. Adeleke, P. M. Neumann. Relations related to
  betweenness: their structure and automorphisms, Mem. Amer. Math. Soc. 131
  (1998), no. 623.



\bibitem[BGS]{BGS} W. Ballmann, M. Gromov and V. Schroeder. {\em Manifolds of
Non-positive Curvature}, Springer, 1999.

\bibitem[BHV]{BHV} B. Bekka, P. de la Harpe, A. Valette. {\em Kazhdan's
Property (T)}. Preprint 2002.

\bibitem[Bou]{Bou} N. Bourbaki. {\em Topologie g\'en\'erale}, quatri\`eme
\'edition, Hermann, Paris, 1965.

\bibitem[Bow$_1$]{Bow1*} B. Bowditch. Relatively hyperbolic
groups. Preprint Southampton, 1997.


\bibitem[Bow$_2$]{Bow2*} B. Bowditch. Treelike structures arising from continua and convergence
groups. Memoirs Amer. Math. Soc. Volume 662 (1999).


\bibitem[Bow$_3$]{Bow3*} B. Bowditch. Intersection numbers and the hyperbolicity of the curve
complex. Preprint. 2003.


\bibitem[Bow$_4$]{Bow4*} B. Bowditch, \textit{private communications}.

\bibitem[Bo\.{z}]{Boz} M. Bo\.{z}ejko. Uniformly amenable discrete
groups, Math. Ann. 251 (1980), 1-6.

\bibitem[Bri]{Bri} M. Bridson. Asymptotic cones and polynomial
isoperimetric inequalities. Topology 38 (1999), no. 3, 543--554.

\bibitem[BrH]{BH} M. R. Bridson, A. Haefliger. {\em Metric Spaces of Non-positive
Curvature}, Sprin\-ger, 1999.

\bibitem[Bo]{Buragos} Dmitri Burago, Yuri Burago, Sergei Ivanov.
{\em A course in metric geometry,} Graduate Studies in
Mathematics, 33. American Mathematical Society, Providence, RI,
2001.

\bibitem[Bu]{Bu} J. Burillo.
{{\it Dimension and fundamental groups of asymptotic cones}}, PhD
thesis, University of Utah, june 1996.



\bibitem[Chis]{Chis} Ian Chiswell. {\em
Introduction to $\Lambda$-trees.} World Scientific Publishing Co.,
Inc., River Edge, NJ, 2001.

\bibitem[Del]{Del}  T. Delzant. Sous-groupes distingu\'es
et quotients des groupes hyperboliques, Duke Math. J. 83:3 (1996),
661–682.

\bibitem[Dah$_1$]{Dah} F. Dahmani. Combination of convergence groups. Geometry and Topology 7 (2003), 933-963.

\bibitem[Dah$_2$]{Dah2} F. Dahmani. {\em Les groupes relativement hyperboliques et leurs bords}, PhD Thesis, University of Strasbourg.

\bibitem[dH]{dlH} P. de la Harpe. Moyennabilit\'e de quelques
 groupes topologiques de dimension infinie. C. R. Acad. Sci. Paris s\'er. I,
 277 (1973), 1037-1040.

\bibitem[DP$_1$]{EP} A. Dioubina, I. Pol\-te\-ro\-vich.
Explicit constructions of universal R-trees and asymptotic
geometry of hyperbolic spaces. preprint, math.DG/9904133.




\bibitem[Dr$_2$]{Dr2} C. Dru\c{t}u. Quasi-isometric classification of non-uniform lattices in semisimple
groups of higher rank. Geom. Funct. Anal. {\bf 10} (2000), no. 2,
327-388.


\bibitem[Dr$_3$]{Dr1} C. Dru\c{t}u. C\^ones asymptotiques et invariants de
quasi-isom\'etrie pour des espaces m\'etriques hyperboliques. Ann.
Inst. Fourier {\bf 51} (2001), no. 1, 81-97.


\bibitem[Dr$_4$]{Dr3} C. Dru\c{t}u. Quasi-isometry invariants and
asymptotic cones. Int. J. of Algebra and Computation 12 (2002),
no. 1 and 2, 99-135.

\bibitem[EO]{EO} A. Erschler, D. Osin.
Fundamental groups of asymptotic cones. preprint
arXiv:math.GR/0404111 5 april 2004.


\bibitem[Fa]{Fa} B. Farb. Relatively hyperbolic groups. Geom. Funct. Analysis {\bf{8}}
(1998), 810-840.

\bibitem[Gr$_1$]{Gr1} M. Gromov. Groups of polynomial growth
and expanding maps.  Publ. Math. IHES {\bf{53}} (1981), 53-73.

\bibitem[Gr$_2$]{Gr1p} M. Gromov. Hyperbolic groups.
Essays in group theory, 75-263, Math. Sci. Res. Inst. Publ., 8,
Springer, New York, 1987.

\bibitem[Gr$_3$]{Gr2} M. Gromov. Asymptotic Invariants of Infinite
Groups. Geometric Group Theory(vol. 2), G. A. Niblo, M. A. Roller
(eds), Proc. of the Symposium held in Sussex, LMS Lecture Notes
Series 181, Cambridge University Press 1991.


\bibitem[GLP]{GLP} M. Gromov, J. Lafontaine, P. Pansu. {\em Structures
m\'et\-riques pour les vari\'et\'es riemanniennes}. Cedic/Fernand
Nathan, Paris (1981).


\bibitem[Hat]{Hat} Allen Hatcher. {\em Algebraic topology}. Cambridge University Press,
Cambridge, 2002.

\bibitem[KaL$_1$]{KaL1} M. Kapovich, B. Leeb. On asymptotic cones and
quasi-isometry classes of fundamental groups of nonpositively
curved manifolds. Geom. Funct. Analysis {\bf{3}} (1995), vol.5,
582-603.

\bibitem[KaL$_2$]{KaL2} M. Kapovich, B. Leeb. Quasi-isometries
preserve the geometric decomposition of Haken manifolds. Invent.
Math. {\bf 128} (1997), no.2, 393-416.

\bibitem[KaL$_3$]{KaL3} M. Kapovich, B. Leeb. $3$-manifold groups and nonpositive curvature. Geom. Funct. Anal.
8 (1998), no. 5, 841-852.

\bibitem[Kel]{Kel} G. Keller. Amenable groups and varieties of
groups, Illinois J. Math. 16 (1972), 257-268.

\bibitem[KlL]{KlL} B. Kleiner, B. Leeb. Rigidity of quasi-isometries
for symmetric spaces and Euclidean buildings. Publ. Math. IHES
{\bf{86}} (1997), 115-197.


\bibitem[KSTT]{KSTT} L. Kramer, S. Shelah, K. Tent, S.
Thomas. Asymptotic cones of finitely presented groups. preprint,
arXive, math.GT/0306420.

\bibitem[KT]{KT} L. Kramer, K. Tent. Asymptotic cones and
ultrapowers of Lie groups. preprint, arXive, math.GT/0311101.

\bibitem[LS]{LS} R. Lyndon and P. Schupp.
\newblock {\it Combinatorial group theory}.
\newblock Springer-Verlag, 1977.

\bibitem[Ols$_1$]{O1} A.Yu. Olshanskii. ${\rm SQ}$-universality of hyperbolic groups. (Russian)
Mat. Sb. 186 (1995), no. 8, 119--132; translation in Sb. Math. 186
(1995), no. 8, 1199--1211.


\bibitem[Ols$_2$]{Ols}
A.Yu.Olshanskii. Distortion functions for subgroups. In: Group
Theory Down Under (ed. J.Cossey, C.F. Miller, W.D. Neumann, M.
Shapiro), de Gruyter, 1999, 281--291.

\bibitem[Os]{Osin} D.V. Osin. Relatively hyperbolic groups. Preprint, 2003.

\bibitem[Pa]{Pa} P. Pansu. Croissance des boules et des
g\'e\-o\-d\'e\-si\-ques ferm\'ees  dans les nilvari\'et\'es.
Ergod. Th. Dynam. Syst. 3 (1983), 415-445.

\bibitem[Pp]{Pp} P. Papasoglu.
On the asymptotic cone of groups satisfying a quadratic
isoperimetric inequality. J. Differential Geom. 44 (1996), no. 4,
789-806.

\bibitem[PW]{PW} P. Papasoglu, K. Whyte.
Quasi-isometries between groups
with infinitely many ends. Comment. Math. Helv. 77 (2002), no. 1,
133--144.

\bibitem[Ri]{Ri} T.R. Riley. Higher connectedness of asymptotic cones. Topology
42 (2003), no. 6, 1289--1352.

\bibitem[RiSe]{RipsSela} E.Rips, Z. Sela. Canonical
representatives and equations in hyperbolic groups. Invent. Math.
120 (1995), no. 3, 489--512.

\bibitem[SBR]{SBR} M. V. Sapir, J. C. Birget, E. Rips.
\newblock Isoperimetric and isodiametric functions of groups,
Annals of Mathematics, 157, 2 (2002), 345-466.

\bibitem[Sch]{Sch} R. E. Schwartz. The quasi-isometry
classification of rank one lattices. Inst. Hautes Etudes Sci.
Publ. Math. no. 82 (1995), 133-168 (1996).



\bibitem[Sh]{She} S. Shelah. Classification theory and the number
of non-isomorphic models. Studies in Logic and the Foundations of
Mathematics, 92. North-Holland Publishing Co., Amsterdam-New York,
1978.

\bibitem[TV]{TV} S. Thomas, B. Velickovic. Asymptotic cones of finitely
generated groups. Bull. Lon- don Math. Soc. 32 (2000), no. 2,
203-208.

\bibitem[Tr]{Trofimov} V.I. Trofimov.
Some asymptotic characteristics of groups. (Russian) Mat. Zametki
46 (1989), no. 6, 85--93, 128; translation in Math. Notes 46
(1989), no. 5-6, 945--951 (1990).


\bibitem[VDW]{VDW} L. van den Dries, A. J. Wilkie. On Gromov's
theorem concerning groups of polynomial growth and elementary
logic. J. of Algebra {\bf{89}} (1984), 349-374.

\bibitem[Wys]{Wys} J. Wysocz\'{a}nski. On uniformly amenable groups, Proc. Amer. Math. Soc. 102 (1988),
no. 4, 933-938.


\bibitem[Ya]{Ya} A. Yaman.
A topological characterisation of relatively hyperbolic groups.
Journal fur die reine und angewandte Mathematik (Crelle's
Journal), to appear.




\end{thebibliography}
\end{document}